\documentclass[twoside, 10pt]{article}

\usepackage[hmargin= 0.8in, height=9.3in]{geometry}
\usepackage{amssymb,amsthm, amsfonts,amsmath,mathrsfs, amsxtra, url, hyperref}

\usepackage{fancyhdr}
\fancypagestyle{plain}{
 \fancyhf{}

}

\renewcommand{\theequation}{\thesection.\arabic{equation}}
 \numberwithin{equation}{section}

\newtheorem {thm}{Theorem}[section]
\newtheorem {prop}{Proposition}[section]
\newtheorem {lemm}{Lemma}[section]
\newtheorem {deff}{Definition}[section]
\newtheorem {cor}{Corollary}[section]
\newtheorem {rem}{Remark}[section]

\def\ba{\begin{array}}
\def\ea{\end{array}}
\def\bea{\begin{eqnarray}}
\def\eea{\end{eqnarray}}
\def\beas{\begin{eqnarray*}}
\def\eeas{\end{eqnarray*}}
\def\bi{\begin{itemize}}
\def\ei{\end{itemize}}
\def\bc{\begin{cases}}
\def\ec{\end{cases}}

\def\a{\alpha}
\def\g{\gamma}
\def\d{\delta}
\def\e{\varepsilon}
\def\z{\zeta}
\def\k{\kappa}
\def\l{\lambda}

\def\si{\sigma}
\def\vs{\varsigma}
\def\t{\tau}

\def\o{\omega}

\def\vf{\varphi}

\def\D{\Delta}
\def\G{\Gamma}
\def\L{\Lambda}
\def\O{\Omega}
\def\F{\Phi}
\def\P{\Psi}

\def\Th{\Theta}
\def\U{\Upsilon}

\def\bF{{\bf F}}
\def\bG{{\bf G}}

\def\bU{{\bf U}}
\def\bV{{\bf V}}

\def\cA{{\cal A}}

\def\cC{{\cal C}}
\def\cD{{\cal D}}
\def\cE{{\cal E}}
\def\cF{{\cal F}}
\def\cG{{\cal G}}

\def\cI{{\cal I}}
\def\cJ{{\cal J}}
\def\cK{{\cal K}}
\def\cL{{\cal L}}
\def\cM{{\cal M}}
\def\cN{{\cal N}}
\def\cO{{\cal O}}
\def\cP{{\cal P}}
\def\cQ{{\cal Q}}

\def\cU{{\cal U}}
\def\cV{{\cal V}}

\def\cX{{\cal X}}
\def\cY{{\cal Y}}
\def\cZ{{\cal Z}}

\def\hC{\mathbb{C}}

\def\hE{\mathbb{E}}

\def\hG{\mathbb{G}}
\def\hH{\mathbb{H}}

\def\hK{\mathbb{K}}
\def\hL{\mathbb{L}}
\def\hM{\mathbb{M}}
\def\hN{\mathbb{N}}

\def\hQ{\mathbb{Q}}
\def\hR{\mathbb{R}}
\def\hS{\mathbb{S}}

\def\hU{\mathbb{U}}
\def\hV{\mathbb{V}}

\def\sB{\mathscr{B}}

\def\sK{\mathscr{K}}
\def\sL{\mathscr{L}}

\def\sN{\mathscr{N}}
\def\sO{\mathscr{O}}
\def\sP{\mathscr{P}}

\def\sT{\mathscr{T}}

\def\sV{\mathscr{V}}

\def\sX{\mathscr{X}}
\def\sY{\mathscr{Y}}
\def\sZ{\mathscr{Z}}

 \def\fA{\mathfrak{A}}
\def\fB{\mathfrak{B}}

\def\fE{\mathfrak{E}}
\def\fF{\mathfrak{F}}

\def\fI{\mathfrak{I}}

\def\fL{\mathfrak{L}}

\def\fN{\mathfrak{N}}
\def\fO{\mathfrak{O}}
\def\fP{\mathfrak{P}}

\def\fX{\mathfrak{X}}

\def\ff{\mathfrak{f}}

\def\fz{\mathfrak{z}}
\def\fx{\mathfrak{x}}
\def\fl{\mathfrak{l}}
\def\fh{\mathfrak{h}}
\def\fw{\mathfrak{w}}

\def\fk{\mathfrak{k}}
\def\fn{\mathfrak{n}}

\def\({\textnormal{(}}
\def\){\textnormal{)}}
\def\[{[\neg[}
\def\]{]\neg]}
\def\lan{\langle}
\def\ran{\rangle}

\def\no{\noindent}

\def\ss{\smallskip}
\def\ms{\medskip}

\def\q{\quad}
\def\qq{\qquad}

\def\neg{\negthinspace}
\def\dneg{\neg \neg}
\def\tneg{\neg \neg \neg}

\def\ol{\overline}
\def\ul{\underline}
\def\ua{\mathop{\uparrow}}
\def\da{\mathop{\downarrow}}

\def\lra{\mathop{\longrightarrow }}

\def\wt{\widetilde}
\def\wh{\widehat}

\def\pas{{\hbox{$P-$a.s.}}}

\def\hb{\hbox}
\def\dis{\displaystyle}
\def\cd{\cdot}
\def\cds{\cdots}

\def\fa{\,\forall \,}
\def\pa{\partial}
\def\es{\emptyset}

\def\dfnn{\stackrel{\triangle}{=}}
\def\b1{{\bf 1}}
\def\qed{\hfill $\Box$ \medskip}

\def\liminf{\mathop{\ul{\rm lim}}}
\def\limsup{\mathop{\ol{\rm lim}}}

\newcommand{\lsup}[1]{ \underset{#1}{\limsup}}
\newcommand{\linf}[1]{ \underset{#1}{\liminf}}
\newcommand{\lmt}[1]{ \underset{#1}{\lim}}
\newcommand{\lmtu}[1]{ \underset{#1}{\lim} \neg \ua \,}
\newcommand{\lmtd}[1]{ \underset{#1}{\lim} \neg \da \,}

\begin{document}

\title{\bf On Zero-Sum Stochastic Differential Games
}

\author{
  Erhan Bayraktar\thanks{ \noindent Department of
  Mathematics, University of Michigan, Ann Arbor, MI 48109; email:
{\tt erhan@umich.edu}.}  \thanks{E. Bayraktar is supported in part by the National Science Foundation under applied mathematics research grants and a Career grant, DMS-0906257, DMS-1118673, and DMS-0955463, respectively, and in part by the Susan M. Smith Professorship.} $\,\,$,
$~~$Song Yao\thanks{
\noindent Department of
  Mathematics, University of Pittsburgh, Pittsburgh, PA 15260; email: {\tt songyao@pitt.edu}. } }
\date{ }

\maketitle

  \begin{abstract}
We generalize the results of Fleming and Souganidis \cite{Fleming_1989} on zero-sum stochastic differential games to the case when the controls are unbounded. We do this by proving a dynamic programming principle using a covering argument instead of relying on a discrete approximation (which is used along with a comparison principle in  \cite{Fleming_1989}). Also, in contrast with \cite{Fleming_1989}, we define our pay-off through a doubly reflected backward stochastic differential equation. The value function (in the degenerate case of a single controller) is closely related to the second order doubly reflected BSDEs.

\ss  \noindent   {\bf Keywords: }\: Zero-sum stochastic differential games, Elliott-Kalton strategies, dynamic programming principle, stability under pasting, doubly reflected backward stochastic differential equations, viscosity solutions, obstacle problem for fully non-linear PDEs, shifted processes, shifted SDEs, second-order  doubly reflected backward stochastic differential equations.

\end{abstract}

    \tableofcontents

\section{Introduction}

    \setcounter{equation}{0}

\label{sec:intro}

In this paper  we use doubly reflected backward stochastic differential equations (DRBSDEs) to generate
payoffs for  a zero-sum  stochastic differential game introduced by the seminal work of Fleming and Souganidis \cite{Fleming_1989}. In our setting, the two players compete by choosing 
  square-integrable controls.
  \if{0}
 \cite{Fleming_1989} initiated a formal study of  zero-sum SDGs between two players. 
  They showed that the lower and  upper values of such games satisfy the dynamic programming principle
  and  solve   the associated Bellman-Isaacs equations in the viscosity sense.
  \fi

DRBSDEs were first analyzed by Cvitani\'c and Karatzas \cite{Cvitanic_Karatzas_1996}, who showed that the solution of a DRBSDE
 is the value of a certain Dynkin game, a zero-sum stochastic game of optimal stopping.  Then Hamad\`ene et al. \cite{Hamadene_Lepeltier_1995a, Hamadene_Lepeltier_1995b, Hamadene_Lepeltier_Wu_1999, Hamadene_Lepeltier_2000, EH-03}  added controls into DRBSDEs to
 study mixed  control and stopping  games and saddle point problems, only when the drift is controlled. Recent advances for Dynkin games and controller and stopper games were made by Karatzas et al. \cite{MR1872738, MR2435857, OS_CRM}, Bayraktar et al. \cite{OSNE1, OSNE2, BH11, BS12}, by \cite{Peng_XuXm_2010}. 
 On the other hand, when there are two competing controllers who can also control the diffusion coefficient, there are a lot of technicalities involved as it is demonstrated by \cite{Fleming_1989}. In particular, \emph{Elliott-Kalton strategies} needs to be used for the controller with lower priority. 
 Recently Buckdahn and Li \cite{Buckdahn_Li_1,Buckdahn_Li_2,Buckdahn_Li_3} and Hamad\`ene et al. \cite{Hamadene_Rotenstein_Zlinescu_2009} made some significant advances to these types of games. However, as in \cite{Fleming_1989}, they assumed that the control spaces are compact. Also, the analysis in these papers 
 is different than \cite{Fleming_1989} and ours in that they work with a uniform  canonical space $\O = \big\{\o  \in  \hC([0,T]; \hR^d ) \neg : \o(0)  =  0 \big\} $
 regardless of the starting time of the game. 

One encounters tremendous technical difficulties when the compactness assumption of the control spaces is removed,
 since the approximation tool of  \cite{Fleming_1989} (also see Fleming and Hern\'andez Hern\'andez \cite{Fleming_2011}) is not applicable any longer. There are some exceptions to this rule: Square-integrable controls was considered by Browne \cite{Browne_2000} for a
 specific zero-sum investment game between two small investors whose controls are
  in form of their portfolios. The PDEs in this special case have smooth solutions, therefore the problem can be solved by relying on a verification theorem instead of the dynamic programming principle. In a more general setting, Chapter 6 of Krylov \cite{Krylov_CDP} considered square-integrable controls. However, the analysis was done only for cooperative games (i.e. the so called $\sup$ $\sup$ case). It is also worth mentioning that inspired by the ``tug-of-war" (a discrete-time random turn game, see e.g.
 \cite{PSSW_2009} and \cite{LLPU_1994}),  Atar and Budhiraja  \cite{Atar_Budhiraja_2010} studied
 a zero-sum  stochastic differential game with $\hU=\hV=\{x \in \hR^n: |x|=1 \} \times  [0,\infty) $   played until the state process exits a given domain. The authors showed that the value of such a game is the unique viscosity solution
to the inhomogenous infinity Laplace equation. As in Chapter 6 of \cite{Krylov_CDP}, they depend on an approximating the game with unbounded controls with a sequence of games with bounded controls. They prove a dynamic programming principle for the latter case and prove the equicontinuity of the approximating sequence to conclude that the value function is a viscosity solution to the infinity Laplace equation. Instead of relying on approximation argument, we directly prove a dynamic programming principle for the game with square-integrable controls.

 \ss  In this paper,  the controls of respective players take values in two separable metric spaces  $\hU$ and $\hV$.   
  We  follow  the probabilistic setting of \cite{Fleming_1989} and rely on the existence of the regular conditional probability distributions. When  the game starts from time $t \in [0,T]$,
  we consider
   the canonical space  $\O^t   \dfnn   \big\{\o     \in     \hC([t,T]; \hR^d ) \neg : \o(t)    =    0 \big\}$,
 whose coordinator  process $ B^t $   is a Brownian motion under the  Wiener measure $P^t_0$. 
 Denote by $\cU^t$ (resp.\;$\cV^t$) the set of all $\hU -$valued (resp.\;$\hV -$valued) square-integrable  process.  If   player I
  chooses a  $\mu \in \cU^t$ and  player  II  selects a $\nu \in \cV^t$ as controls,
 the state  process    $X^{t,x,\mu,\nu}$ starting from $x$ 
 will evolve according to the following  SDE:
    \bea
  X_s= x + \int_t^s b(r, X_r,\mu_r,\nu_r ) \, dr + \int_t^s \si(r, X_r,\mu_r,\nu_r) \, dB^t_r,  \q s \in [t, T] ,  \label{FSDE}
    \eea
where the drift $b$ and the diffusion $\si$ are Lipschitz in $x$ and have linear growth in $(\mu,\nu)$.
 The payoff player I will receive from  player II is determined by the first component  of the solution
  $ \big(Y^{t,x,\mu,\nu},Z^{t,x,\mu,\nu},\ul{K}^{t,x,\mu,\nu}, \ol{K}^{t,x,\mu,\nu} \big) $ to the following DRBSDE:
   \bea   \label{DRBSDE00}
    \left\{\ba{l}
 \dis  Y_s= h \big( X^{t,x,\mu,\nu}_T \big)  \neg + \neg  \int_s^T  \neg  f (r, X^{t,x,\mu,\nu}_r,  Y_r, Z_r, \mu_r, \nu_r)  \, dr
 \neg+ \neg  \ul{K}_{\,T} \neg-\neg \ul{K}_{\,s}
  \neg- \neg \big( \ol{K}_T \neg - \neg \ol{K}_s \big)
   \neg-\neg \int_s^T \neg Z_r d B^t_r  , \q    s \in [t,T] ,   \vspace{1mm}      \\
   \dis \ul{l}(s,X^{t,x,\mu,\nu}_s) \le Y_s  \le  \ol{l}(s,X^{t,x,\mu,\nu}_s)       , \q    s \in [t,T] , \vspace{1mm} \q \\
  \dis       \int_t^T \neg  \big(  Y_s   -  \ul{l}(s,X^{t,x,\mu,\nu}_s)  \big)     d \ul{K}_{\,s}
  = \int_t^T \neg  \big(  \ol{l}(s,X^{t,x,\mu,\nu}_s) - Y_s  \big)     d \ol{K}_s   = 0  ,
     \ea \right.
    \eea
  with two separate obstacle functions $\ul{l}  < \ol{l} $
  satisfying $ \ul{l}(T,\cd) \le h(\cd) \le \ol{l} (T,\cd)$. When $ \ul{l}, h , \ol{l} ,f $ are all $2/q-$H\"older continuous
  in $x $ for some $q \in (1,2]$,    $ Y^{t,x,\mu,\nu}$ is $q-$integrable by El Asri et al. \cite{EHW_2011}.
  As we see from \eqref{FSDE} and \eqref{DRBSDE00} that the controls $\mu$ and $\nu$
  influence the game in two aspects: either affect  \eqref{DRBSDE00} via the state process $X^{t,x,\mu,\nu}$
  or appear directly in the generator $f$ of   \eqref{DRBSDE00} as parameters.

  \ss
  We use Elliott-Kalton strategies as in \cite{Fleming_1989}. In this game, one player (e.g. Player I) has priority
   and chooses a control and its opponent (e.g. Player II) will react by selecting  a corresponding strategy.
  We specify Player II's strategy by
   a  measurable mapping  $\beta :  [t,T] \neg \times \neg  \O^t  \neg  \times \neg  \hU \to \hV $ if the game starts from
   time $t$. This additional specification is  to accommodate a particular measurability issue, see Remark \ref{r_strategy}.
   Under a linear growth condition on the $\mu-$variable,
 $\beta$     induces a mapping $\beta \lan \cd \ran : \cU^t \to \cV^t$  by
 $
  \big( \beta \lan \mu \ran \big)_r (\o)  \dfnn  \beta \big(r, \o, \mu_r(\o)\big)  $, $  \fa \mu \neg  \in  \neg  \cU^t,
  ~  (r,\o)  \neg \in  \neg  [t,T]  \neg  \times  \neg  \O^t $,
   which is exactly an Elliott-Kalton strategy.
  Then  $ w_1 (t,x) \dfnn \underset{\beta \in \fB^t  }{\inf} \; \underset{\mu \in \cU^t }{\sup}
  \; Y^{t,x,\mu, \beta \lan \mu \ran }_t $
  represents  Player I's  priority  value of the  game starting   from time $t  $ and state $x$,
  where $\fB^t$ collects  all admissible strategies for Player II.
 Player II's priority value $w_2 (t,x)$  is defined similarly.

  \ss  Although value functions $w_1 (t,x)$, $w_2(t,x)$ are still $2/q-$H\"older continuous in $x$, they are no longer $1/q-$H\"older
 continuous in $t$ as in  the case of compact control spaces. 
 Hence we are not able to use the approach of \cite{Fleming_1989} to show the dynamic programming principle for $w_1$ and $w_2$;
 see Remark  \ref{r_strategy}.  Instead, we use the continuity of   $Y^{t,x,\mu, \nu  }$ in controls $(\mu,\nu)$,
 properties of shifted processes (especially shifted SDEs) as well as
 stability under pasting of controls/strategies (as listed below)  to prove a dynamic  programming principle, say for $w_1$:
  \bea   \label{eq:xxa111}
  w_1(t,x)  =      \underset{\beta \in \fB^t}{\inf} \; \underset{\mu \in \cU^t}{\sup}    \;
     Y^{t,x,\mu, \beta \lan \mu \ran}_t \Big(\t_{\mu,\beta},
    w_1\big(\t_{\mu,\beta},   X^{t,x,\mu, \beta \lan \mu \ran}_{\t_{\mu,\beta}} \big) \Big)
    \eea
 for any family $\{\t_{\mu,\beta} \neg : \mu \neg \in \neg \cU^t, \beta \neg \in  \neg \fB^t\}$ of $\hQ  -$valued stopping times. The crucial ingredients of the proof of dynamic programming principle \eqref{eq:xxa111} are:

  \ss  \no i) When $b$, $\si$ are $\l-$H\"older continuous in $\mu$ (or $\nu$) and $f $ is $2\l-$H\"older continuous in $\mu $ (or $\nu$)
   for some $\l \in (0,1)$,  applying an a priori estimate \eqref{eq:n211} for DRBSDEs,
  we obtain a continuous dependence result for $X^{t,x,\mu,\nu}$ and $Y^{t,x,\mu,\nu}$ on control $\mu$ (or $\nu$), see
 Lemma \ref{lem_DX_varpi} and Lemma \ref{lem_estimate_Y}. This dependence  together with
 two nice topological properties of the canonical spaces $\O^t$, namely separability and Lemma \ref{lem_ocset_Omega}  are crucial in the covering argument
 which is used to construct $\e$ optimal strategies starting at any stopping time.

 \ss \no ii) Let $0\le t \le s \le T$.   For any random variable $\xi$ on $\O^t$ we define a shifted random variable
 $\xi^{s, \o} (\wt{\o}) \dfnn  \xi (\o \otimes \wt{\o})$, $ \fa  \wt{\o} \in  \O^s$ as a projection of $\xi$ onto $\O^s$ along
 a give path $\o $ of $ \O^t$, where $\o \otimes_s \wt{\o}$ is the concatenation of paths $\o$ and $\wt{\o}$ at time $s$; see
 \eqref{def_concatenation}. Its discrete-time finite-state counterpart is
 a restriction of a binomial/trinomial tree  of asset prices to one of its branches. Similarly, one can introduce
 shifted processes and shifted random fields, in particular, shifted controls and shifted strategies.
 Soner et al. \cite{STZ_2011b, STZ_2011c} as well as  our generalization in Subsection \ref{subsection:shift_measurability}
 \& \ref{subsection:shift_integrability}  show that these shifted random objects
 almost surely inherit measurability and integrability.

   In Proposition  \ref{prop_FSDE_shift}, we extend a result of \cite{Fleming_1989} on shifted forward SDEs:
    For $P^t_0-$a.s. $\o \in \O^t$,
  the process obtained by shifting $X^{t,x,\mu,\nu}$ solves \eqref{FSDE}  with  parameters $\big(s, X^{t,x,\mu,\nu}_s (\o), \mu^{s,\o},\nu^{s,\o} \big) $  on the probability space $\big(\O^s, \cF^s_T,P^s_0 \big)$.
   Similarly, the process obtained by shifting 
 $ \big(Y^{t,x,\mu,\nu},Z^{t,x,\mu,\nu},\ul{K}^{t,x,\mu,\nu}, \ol{K}^{t,x,\mu,\nu} \big) $
 solves \eqref{DRBSDE00} with the parameters $\big(s, X^{t,x,\mu,\nu}_s (\o), \mu^{s,\o},\nu^{s,\o} \big) $  on
 $\big(\O^s, \cF^s_T,P^s_0 \big)$; see Proposition \ref{prop_DRBSDE_shift}. These  two propositions
  are also crucial in demonstrating   \eqref{eq:xxa111}.

\no iii) In constructing the $\e-$optimal strategies above, we use
pasting of controls and strategies. Our sets of controls and strategies are closed under pasting; see Proposition \ref{prop_paste_control} \& \ref{prop_paste_strategy}.
  In the latter proof we show that an additional path-continuity of the strategies that we use to prove
  the dynamic programming principle is also closed under pasting.

   \ss     Next, using  the dynamic programming principle,  the continuity of   $Y^{t,x,\mu, \nu  }$ in controls $(\mu,\nu)$ as well as
  the separability of  $\hU$, $\hV$   we deduce that the value functions $w_1$ and $w_2$ are (discontinuous)
  viscosity solutions of  the corresponding obstacle problem of  fully non-linear PDEs, see Theorem \ref{thm_viscosity}.

 \ss  When $\hV$ becomes a singleton,  the zero-sum stochastic differential
 game degenerates into a classical stochastic control problem including only one player. In particular,
 when $\hU =\{$all symmetric $d-$dimensional
 matrices\}, $b(t,x,u ) =b(t,x)$ and $\si(t,x,u)= u $,
  the  value function $w$ of the optimization problem coincides with
 that  of the second-order DRBSDEs and is related to the one in Nutz \cite{Nutz_2011}
  via a probability transformation of strong form \eqref{eq:a155}.
  Motivated by applications in financial mathematics and probabilistic numerical methods,
 Cheridito et al. \cite{CSTV_2007} introduced second-order BSDEs. Later, Soner et al. \cite{STZ_2011c} refined this notion and
 Soner et al. \cite{STZ_2011d} related it to  $G-$expectations of Peng \cite{Peng_G_2007a, Peng_G_2007b}.
 Quite recently Matoussi et al. \cite{Matoussi_Possamai_Zhou_2012} analyzed the second order reflected BSDEs.

 The rest of the paper is organized as follows:
 After listing the notations used, we present two basic properties of DRBSDEs in Section \ref{sec:intro}.
 In Section \ref{sec:zs_drgame}, we set up the zero-sum stochastic differential games based on DRBSDEs and
 present a dynamic programming principle in Theorem \ref{thm_DPP}, for priority values of both players defined via Elliott-Kalton strategies.
 Using the dynamic programming principle, we show in Section \ref{sec:PDE} that  the   priority values
   are (discontinuous) viscosity solutions of the corresponding obstacle problem of  fully non-linear   PDEs; see Theorem~\ref{thm_viscosity}.
 In Section \ref{sec:shift_prob}, we explore the properties of shifted processes (including measurability/integrability),
 shifted SDEs and pasting of controls/strategies.  The contents of this section are all technical necessities in proving our main results, Theorems \ref{thm_DPP}
 and \ref{thm_viscosity}.  In Section \ref{sec:one_control}, we will discuss the classical stochastic control problem as a degenerate case and   connect it  to second order doubly reflected BSDEs. The proofs of our results are deferred to Section \ref{sec:Proofs}.

\subsection{Notation and Preliminaries} \label{subsec:preliminary}


      We let
              $\hE$ denote a generic Euclidian space and let   $\hM$ be a generic metric space  with metric $\rho_{\overset{}{\hM}}$
   and denote by $\sB(\hM)$  the Borel $\si-$field on $\hM$. For any $x \in \hM $     and $\d >0$,
           $O_\d(x) \dfnn \{x' \in \hM : \rho_{\overset{}{\hM}}(x,x')  < \d \}$ denotes
         the open   ball    centered at $x   $     with radius $\d  $ and its closure is
              $\ol{O}_\d(x) \dfnn \{x' \in \hM : \rho_{\overset{}{\hM}}(x,x')  \le  \d \}$.
                For any function $\phi: \hM \to \hR$, we define
    \beas
      \linf{x' \to x} \phi(x') \dfnn \lmtu{n \to \infty} \underset{x' \in O_{1/n}(x)}{\inf} \phi(x')
     \q \hb{ and }  \q  \lsup {x' \to x} \phi(x') \dfnn \lmtd{n \to \infty} \underset{x' \in O_{1/n}(x)}{\sup} \phi(x')  ,
     \q \fa  x \in \hM .
    \eeas

  \ss   Fix $d    \in \neg  \hN$.  For any  $0 \le t \le T < \infty $,
  we  set     $\hQ_{t,T} \dfnn  \big(  [t,T) \cap \hQ \big) \cup \{T\} $, $\hQ$ being the rational numbers,
   and  let $\O^{t,T} \neg \dfnn \neg \big\{\o  \neg  \in  \neg  \hC([t,T]; \hR^d) \neg : \o(t)  \neg = \neg  0 \big\}$
 be the   canonical space  over the period $[t,T]$,  which is   equipped with the uniform norm
 $\|\o\|_{t,T}  \dfnn \underset{s \in [t,T]}{\sup} |\o(s)| $.
 We let   $O_\d(\o) \dfnn \{\o' \in \O^{t,T} : \| \o'- \o \|_{t,T} < \d \}$ denote  the open ball centered
 at $ \o \in \O^{t,T} $  with radius $\d >0$,    
   and  let   $\sB(\O^{t,T})$  be   the correspondingly Borel $\si-$field of $\O^{t,T}$. 
   We denote by  $B^{t,T} $    the   canonical process   on $\O^{t,T}$, and by $P^{t,T}_0$ the   Wiener measure
   on $\big(\O^{t,T},  \sB(\O^{t,T})\big)$ under which $ B^{t,T} $    is a  $d-$dimensional Brownian motion.
     Let      $\bF^{t,T} =  \Big\{ \cF^{t,T}_s \dfnn \si \big(B^{t,T}_r; r \in [t,s]\big) \Big\}_{s \in [t,T]}$
   be the   filtration generated by $ B^{t,T} $ and let $\cC^{t,T} $   collect    all {\it cylinder} sets in $\cF^{t,T}_T$, i.e.
      $
   \cC^{t,T}  \dfnn \left\{  \underset{i=1}{\overset{m}{\cap}}  \big(  B^{t,T}_{t_i} \big)^{-1}     ( \cE_i ): m \in \hN, \, t < t_1 < \cds < t_m \le T , \, \{ \cE_i\}^m_{i=1} \subset \sB(\hR^d)  \right\}$.
         It is well-known that
   \bea    \label{eq:xxc023}
   \sB(\O^{t,T})  =\si(\cC^{t,T}  )  = \si \Big\{ \big(B^{t,T}_r\big)^{-1} (\cE ) : r \in [t,T], \cE \in \sB(\hR^d) \Big\} = \cF^{t,T}_T .
   \eea
   For any $\bF^{t,T}-$stopping time $\t$, we define two stochastic intervals
   $\[t,\t\[ \; \dfnn   \big\{ (r,\o) \in [t,T] \times \O^t: r  <  \t(\o) \big\} $,
     $\[\t,T\] \dfnn \big\{ (r,\o) \in [t,T] \times \O^t: r \ge \t(\o) \big\} $
   and set  $ \[\t,T\]_A \dfnn \big\{ (r,\o) \in [t,T] \times A: r \ge \t(\o) \big\}$  for any
   $A \in \cF^{t,T}_\t$.

 \ss  The following two results are basic, see \cite{SDGVU} for proofs.

  \begin{lemm}  \label{lem_countable_generate1}
 Let $0 \neg \le \neg t \neg \le \neg T \neg < \neg \infty$, for any $s \neg \in \neg [t,T]  $,
 the $\si-$field $  \cF^{t,T}_s  $ is countably generated by
  \beas
   \cC^{t,T}_s \dfnn \Big\{  \underset{i=1}{\overset{m}{\cap}} \big( B^{t,T}_{t_i} \big)^{-1} \big( O_{\l_i} (x_i) \big) :   \,    m \in \hN, \,  t_i \in \hQ
   \hb{ with } t \le t_1 < \cds < t_m \le s ,\, x_i \in \hQ^d , \, \l_i \in \hQ_+ \Big\}  .
   \eeas
   \end{lemm}

\begin{lemm}  \label{lem_shift_inverse} 
 Let $0 \le t \le s \le S \le T < \infty$.
 The truncation  mapping $\Pi^{T,S}_{t,s}: \O^{t,T} \to \O^{s, S}$ defined by
 \beas
  \big(\Pi^{T,S}_{t,s}(\o)\big)(r) \dfnn \o (r) - \o(s) , \q \fa \o \in \O^{t,T}, ~  \fa s \in [s, S]
 \eeas
 is continuous (under uniform norms) and $\cF^{t,T}_r \big/ \cF^{s,S}_r-$measurable for any $r \in [s,S]$.
Moreover, we have
 \beas
    P^{t,T}_0 \Big( \big(\Pi^{T,S}_{t,s} \big)^{ -1}(A) \Big) = P^{s,S}_0(A), \q  \fa   A \in \cF^{s,S}_S.
\eeas
\end{lemm}

 From now on, we fix  a  time horizon $T  \in (0,\infty)  $ and
   shall drop it from  the above notations, i.e.
 \big($  \O^{t,T}$, $\|~\|_{t,T}$, $B^{t,T}$,  $\bF^{t,T}$, $P^{t,T}_0$,   $\cC^{t,T}_s$\big) $\lra$
  \big($\O^t$, $\|~\|_t$, $B^t$,  $\bF^t$, $P^t_0$,   $\cC^t_s$\big).
  The expectation under $P^t_0$ will be denoted by   $E_t $. When $S=T$ we simply denote $\Pi^{T,T}_{t,s}$
   by $\Pi_{t,s}$ in Lemma \ref{lem_shift_inverse}.

   Given $t \in [0,T]$, we  let $\cP^t$  denote  the set of all probability measures on $\big(\O^t,  \sB(\O^t)  \big) = \big(\O^t,  \cF^t_T\big)$
   by \eqref{eq:xxc023}.    For any   $P \in \cP^t$,    we set
    $      \sN^P \dfnn \big\{ \cN \subset \O^t: \cN \subset A \hb{ for some  } A \in \cF^t_T \hb{ with } P(A ) =0   \}    $
        as the collection of all  $P-$null sets.  The $P-$augmentation $ \bF^P$ of $\bF^t$ consists of
     $ \cF^P_s   \dfnn \si \big(  \cF^t_s  \cup  \sN^P  \big)  $,   $ s \in [t,T]$.
      \Big(In particular, we will write $\ol{\bF}^t  =  \big\{ \ol{\cF}^t_{\neg s} \big\}_{s \in [t,T]} $
      for $\bF^{P^t_0} =  \Big\{ \cF^{P^t_0}_s  \Big\}_{s \in [t,T]}$.\Big)
         The completion of  $\big(\O^t,  \cF^t_T, P\big)$ is the probability space
        $\Big(\O^t,  \cF^P_T, \ol{P}\Big)$ with $\ol{P} \Big|_{\cF^t_T} \neg  = P $.
       For convenience, we will simply write $P$ for   $\ol{P}$.   

    \ss     Similar to Lemma 2.4 of \cite{STZ_2011a}, we have the following result:

  \begin{lemm}  \label{lem_F_version}
    Let \,$t \in [0,T]$ and  $P \in \cP^t$.

\ss  \no    1)  For any $\xi   \neg  \in   \neg    \hL^1 \big( \cF^P_T, P\big)$ and $s   \neg  \in \neg  [t,T]$,
  $  
    E_{P} \neg \big[ \xi \big| \cF^P_s \big]    \neg  =  \neg   E_P \neg \big[ \xi \big|  \cF^t_s  \big] $,     \pas ~
      Consequently,   a martingale \(resp.\;local martingale or semi-martingale\)  with respect to $(\bF^t, P)$
 is also a martingale \(resp.\;local martingale or semi-martingale\) with respect to $\left(\bF^P, P\right)$.

 \ss  \no    2)    For any $\hE-$valued, $\bF^P-$adapted  continuous process $\{X_s\}_{s \in [t,T]}$, there exists a  unique
   \big(in sense of $P-$evanescence\big) $\hE-$valued, $\bF^t-$adapted  continuous  process $\{\wt{X}_s\}_{s \in [t,T]}$  such that $ P \big(\wt{X}_s= X_s, ~\fa s \in [t,T]\big) =1$.
   For any $\hE-$valued, $\bF^P-$progressively measurable  process $\{X_s\}_{s \in [t,T]}$, there exists a  unique
 \big(in $ds \times d P-$a.s. sense\big) $\hE-$valued, $\bF^t-$ progressively measurable  process $\{\wt{X}_s\}_{s \in [t,T]}$
  such that $ \wt{X}_s  (\o)= X_s (\o) $ for $ds \times d P-$a.s. $(s,\o) \in [t,T] \times \O^t$.
        In both cases, we call  $\wt{X}$  the  $\bF^t-$version of $X$.

\end{lemm}

   For any $p \in [1,\infty)$,   $t \in [0,T]$ and $P \in \cP^t$,     we introduce some   spaces of functions:

  \ss \no 1)  For any sub$-\si-$field  $\cF$  of $\cF^P_T$,  let
     $\hL^p(\cF,\hE, P) $ be  the space of all  $\hE-$valued,
$  \cF-$measurable random variables $\xi$ such that $\|\xi\|_{\hL^p(\cF,P)}
   \dfnn  \Big\{E_P\big[|\xi|^p\big]\Big\}^{1/p  } < \infty$.

\ss \no 2)  For any  filtration $\fF=\{\fF_s\}_{s \in [t,T]}$ on $\big(\O^t, \cF^P_T\big)$,
  $\sP  (\fF )$  will denote  the   $\fF-$progressively  measurable $\si-$field of $ [t,T] \times \O^t$.
     Let  $\hC^0_\fF   ([t,T], \hE, P)$ be the space of all $\hE-$valued, $\fF-$adapted processes
  $\{X_s\}_{s \in [t,T]}$ with \pas ~ continuous paths.  We define the following subspaces of $\hC^0_\fF   ([t,T], \hE, P)$:


  \ss \no ~   $\bullet$ $\hC^p_\fF  ( [t,T],\hE,P) \dfnn    \bigg\{  X \in \hC^0_\fF  ( [t,T],\hE,P) :\,  \|X\|_{\hC^p_\fF( [t,T],P)}
  \dfnn \bigg\{ E_P \Big[\,\underset{s \in [t,T]}{\sup}|X_s|^p \Big]  \bigg\}^{ 1 / p }<\infty \bigg\}$;

  \ss \no ~  $\bullet$ $ \hC^{\pm, p}_\fF  ( [t,T] ,P) \neg \dfnn  \neg    \Big\{  X  \neg \in  \neg  \hC^0_\fF ( [t,T], \hR, P)
   \neg  : X^\pm \dfnn ( \pm  X ) \vee 0     \in     \hC^p_\fF  ( [t,T], P) \Big\}$;

    \ss \no ~  $\bullet$ $   \sV_\fF( [t,T],P) \dfnn \big\{  X \in \hC^0_\fF ( [t,T], \hR, P) :\,  X$ has
    $P-$a.s. finite variation\big\};

  \ss \no ~   $\bullet$  $ \hK_\fF( [t,T],P)  \dfnn  \big\{  X \in \hC^0_\fF  ( [t,T], \hR, P) : \,  X_t=0$
  and $X$ has  \pas ~ increasing paths\big\};

   \ss \no ~   $\bullet$  $ \hK^p_\fF( [t,T],P) \dfnn \big\{  X \in \hK_\fF( [t,T],P) :\, E_P \big[ X^p_T \big] < \infty$\big\}.

\ss \no 3)     Let  $ \hH^{p,loc}_\fF   ([t,T], \hE, P) $  be the space of all $\hE-$valued,
 $\fF-$progressively   measurable processes $\{X_s\}_{s \in [t,T]}$ with
   $ \int_t^T \neg |X_s|^p   \,  ds <\infty $, $P^t_0-$a.s.
    And for any $\wh{p} \in [1,\infty)$, we let    $ \hH^{p,\wh{p}}_\fF   ([t,T], \hE, P) $  denote the space of all $\hE-$valued,
 $\fF-$progressively   measurable processes $\{X_s\}_{s \in [t,T]}$ with
   $
   \|X\|_{\hH^{p,\wh{p}}_\fF([t,T], \hE, P)} \dfnn \left\{  E_P \left[  \big( \hb{$\int_t^T \neg |X_s|^p   \,  ds$}
  \big)^{ \wh{p} /p } \right] \right\}^{ 1/\wh{p} }<\infty $.

 \ss  Also, we set  $ \hG^q_\fF( [t,T],P) \dfnn \hC^q_\fF( [t,T], \hR,P) \times  \hH^{2,q}_\fF([t,T], \hR^d,P)
        \times \hK^q_\fF( [t,T],P)  \times \hK^q_\fF( [t,T],P) $.

 \ss   If $\hE=\hR$ (resp. $P=P^t_0$), we will drop it from the above notations. Moreover, we use the convention $ \inf \es \dfnn \infty$.

 \subsection{Doubly Reflected Backward Stochastic Differential Equations}

  Let \,$t \in [0,T]$.  A  $t-$parameter set $\big(\xi,\ff,\ul{L}, \ol{L} \big)$ consists of
  a random variable   $\xi \in \hL^0 \big(\ol{\cF}^{\, t}_{\neg T} \big)$,
  a function $\ff: [t,T] \times \O^t \times \hR \times \hR^d \to \hR $, and
  two processes $ \ul{L}, \ol{L}  \in \hC^0_{\ol{\bF}^t} \big([t,T]\big)$ such that
  $\ff$ is $\sP \big(  \ol{\bF}^t  \big)   \otimes \sB(\hR)  \otimes \sB(\hR^d)/\sB(\hR)-$measurable
  and that $ \ul{L}_T \le  \xi \le \ol{L}_T $, $P^t_0-$a.s. In particular, $\big(\xi,\ff,\ul{L}, \ol{L} \big)$
 is called a $(t,q)-$parameter set if $\xi \in \hL^q \big(\ol{\cF}^{\, t}_{\neg T} \big)$, $ \ul{L} \in \hC^{+,q}_{\ol{\bF}^t} \big([t,T]\big)$  and $ \ol{L} \in \hC^{-,q}_{\ol{\bF}^t} \big([t,T]\big)$.

 \begin{deff}
 Given $t \in [0,T]$ and a $t-$parameter set $\big(\xi,\ff,\ul{L}, \ol{L} \big)$, a quadruplet $(Y,Z,\ul{K} , \ol{K})
 \in \hC^0_{\ol{\bF}^t}([t,T]) \times \hH^{2,loc}_{\ol{\bF}^t}([t,T],\hR^d) \times \hK_{\ol{\bF}^t}([t,T]) \times \hK_{\ol{\bF}^t}([t,T])$ is called a solution of the
 doubly reflected backward stochastic differential equation  on the probability space $(\O^t, \ol{\cF}^t_{\neg T}, P^t_0)$
    with terminal condition $\xi$, generator $\ff$, lower obstacle $\ul{L}$ and upper obstacle $\ol{L}$
    \big(DRBSDE$ \big( P^t_0,\xi,\ff,\ul{L}, \ol{L} \big)$ for short\big) if it holds $P^t_0-$a.s. that
     \bea   \label{DRBSDE01}
    \left\{\ba{l}
 \dis   Y_s= \xi \neg + \neg  \int_s^T  \neg  \ff  (r,   Y_r, Z_r)  \, dr  \neg+ \neg  \ul{K}_{\,T} \neg-\neg \ul{K}_{\,s}
  \neg- \neg \big( \ol{K}_T \neg - \neg \ol{K}_s \big)
   \neg-\neg \int_s^T \neg Z_r d B^t_r  , \q    s \in [t,T] ,   \vspace{1mm} \q \\
   \dis \ul{L}_s \le Y_s  \le  \ol{L}_s       , \q    s \in [t,T]   \q \hb{and}  \q  \int_t^T \neg  \big(  Y_s   -  \ul{L}_s  \big)     d \ul{K}_{\,s}
  = \int_t^T \neg  \big(  \ol{L}_s - Y_s  \big)     d \ol{K}_s   = 0  .
     \ea \right.
    \eea
\end{deff}

 The last two equalities in \eqref{DRBSDE01} are known as the {\it flat-off} conditions corresponding to $\ul{L}$ and $\ol{L}$
 respectively, under which the two increasing processes $\ul{K}$, $\ol{K}$ keep process $Y$  between $\ul{L}$ and $\ol{L}$ at the minimal effort: i.e.,  only when $Y$ tends to drop below   $\ul{L}$ (resp.\;rise above $\ol{L}$), $\ul{K}$ (resp.\;$\ol{K}$)
   generates an upward (resp.\;downward) momentum.

\ss    We first have the following comparison result and a priori estimate for DRBSDEs, which generalize those in
 \cite{Buckdahn_Li_3} and \cite{Hamadene_Hassani_2005}.

  \begin{prop} \label{prop_comp_DRBSDE}
 Given $t \in [0,T]$ and two $(t,q)-$parameter sets  $ \big(\xi_1,\ff_1,\ul{L}^1, \ol{L}^1 \big), \big(\xi_2,\ff_2,\ul{L}^2, \ol{L}^2 \big)$ with  $P^t_0(\xi_1 \le \xi_2) = P^t_0 \big(\ul{L}^1_s \le \ul{L}^2_s, \, \ol{L}^1_s \le \ol{L}^2_s, \fa s \in [t,T] \big) =1 $,  let  $ \big(Y^i ,Z^i ,\ul{K}^i , \ol{K}^i  \big) \in  \hC^q_{\ol{\bF}^t}([t,T]) \times
 \hH^{2,q}_{\ol{\bF}^t}([t,T],\hR^d) \times \hK_{\ol{\bF}^t}([t,T]) \times \hK_{\ol{\bF}^t}([t,T]) $,  $i=1,2$ be a solution
   of DRBSDE$ \big( P^t_0,\xi_i,\ff_i,\ul{L}^i, \ol{L}^i \big)$.
    For either $i=1$ or $i=2$,      if \, $\ff_i$   is  Lipschitz continuous in $(y,z)$:
  i.e. for some $\g > 0 $, it holds for $ds \times dP^t_0-$a.s.    $(s, \o) \in [t,T] \times \O^t $ that
   \bea   \label{ff_Lip}
 \big| \ff_i(s, \o,y ,z ) - \ff_i(s,\o,y', z') \big| \le \g \big( |y-y'| + |z-z'| \big), \q \fa y,y' \in \hR,~ \fa z,z' \in \hR^d ,
   \eea
  and if $  \ff_1 (s,  Y^{3-i}_s ,Z^{3-i}_s ) \le  \ff_2 (s, Y^{3-i}_s ,Z^{3-i}_s )  $,  $ds \times dP^t_0-$a.s.,
     then  it holds $P^t_0-$a.s. that   $ Y^1_s \le Y^2_s$  for any $s \in [t,T]$.

\end{prop}

 \begin{prop} \label{prop_apriori_DRBSDE}
Given $t \in [0,T]$ and two $(t,q)-$parameter sets  $ \big(\xi_1,\ff_1,\ul{L}^1, \ol{L}^1 \big), \big(\xi_2,\ff_2,\ul{L}^2, \ol{L}^2 \big)$ with  $P^t_0 \big(\ul{L}^1_s = \ul{L}^2_s, \, \ol{L}^1_s = \ol{L}^2_s, \fa s \in [t,T] \big) =1$,
 let  $ \big(Y^i ,Z^i ,\ul{K}^i , \ol{K}^i  \big) \in  \hC^q_{\ol{\bF}^t}([t,T]) \times
 \hH^{2,q}_{\ol{\bF}^t}([t,T],\hR^d) \times \hK_{\ol{\bF}^t}([t,T]) \times \hK_{\ol{\bF}^t}([t,T]) $,  $i=1,2$ be a solution
   of DRBSDE$ \big( P^t_0,\xi_i,\ff_i,\ul{L}^i, \ol{L}^i \big)$.
   If  \, $\ff_1$ satisfies \eqref{ff_Lip}, then for any   $\varpi \neg \in \neg (1,q]$
\bea      \label{eq:n211}
     E_t \bigg[ \underset{s \in [t,T]}{\sup} \big|    Y^1_s \neg - \neg  Y^2_s   \big|^\varpi  \bigg]
       \neg \le   \neg  C(T,\varpi, \g )   \Bigg\{ E_t \big[  \big|   \xi_1   \neg -  \neg  \xi_2    \big|^\varpi \big]
           \neg +  \neg  E_t \bigg[  \bigg( \neg \int_t^T   \neg   \ff_1 (r, Y^2_r ,Z^2_r )
 \neg - \neg  \ff_2 (r, Y^2_r,Z^2_r)   dr    \bigg)^{\neg \varpi}  \bigg] \Bigg\}   .
  \eea

\end{prop}

Given   a $(t,q)-$parameter set  $ \big(\xi,\ff,\ul{L}, \ol{L} \big)$ such that
 $\ff$     is  Lipschitz continuous in $(y,z)$.
 If $    E_t   \Big[  \Big(\int_t^T  \neg \big| \ff    (s, 0,0 )   \big| \, ds \Big)^q \Big] < \infty $
 and if $P^t_0( \ul{L}_s   < \ol{L}_s, \fa s \in [t,T]) =1$, then we know from
     Theorem 4.1  of  \cite{EHW_2011}   that the DRBSDE$ \big( P^t_0,\xi,\ff,\ul{L}, \ol{L} \big)$
        admits a unique solution $ \big(Y ,Z ,\ul{K} , \ol{K}  \big) \in \hG^q_{\ol{\bF}^t } \big([t,T]\big) $.

 \section{Stochastic Differential Games with Square-Integrable Controls} \label{sec:zs_drgame}

    Let  $(\hU, \rho_{\overset{}{\hU}})$ and $(\hV, \rho_{\overset{}{\hV}})$ be two  separable metric spaces,
   whose  Borel$-\si-$fields will be denoted by $\sB(\hU)$  and $\sB(\hV)$ respectively.
   For some $u_0 \in \hU$ and $v_0 \in \hV$, we define
  \beas
   [u]_{\overset{}{\hU}} \dfnn \rho_{\overset{}{\hU}}(u,u_0), ~\fa u \in \hU \q \hb{and} \q  [v]_{\overset{}{\hV}} \dfnn \rho_{\overset{}{\hV}}(v,v_0) , ~\fa v \in \hV .
  \eeas

    Fix    a non-empty  $\hU_0 \subset  \hU$ and  a non-empty    $\hV_0 \subset  \hV$.
   We shall   study a  zero-sum stochastic differential game between two players,
   player I and player II,  who  choose square-integrable $\hU-$valued controls
    and   $\hV-$valued controls respectively  to compete:

  \begin{deff}
 Given $t \in [0,T]$, an admissible control process $\mu=\{\mu_r\}_{r \in [t,T]}$  for player I   over   period $[t,T]$
 is a $\hU-$valued, $\bF^t-$progressively measurable process
 such that $\mu_r \in \hU_0$,  $dr \times dP^t_0 -$a.s. and that
 $ E_t  \int_t^T [ \mu_r ]^2_{\overset{}{\hU}} \,dr < \infty$.
 Admissible control processes    for player  II  over   period $[t,T]$ are defined similarly.
 The set of all admissible controls for player I \(resp.\;II\) over period $[t,T]$ is denoted by $\cU^t$ \(resp.\;$\cV^t$\).

   \end{deff}

 Our zero-sum stochastic differential game is formulated via a (decoupled) SDE$-$DRBSDE  system
   with the following parameters: Fix $k \in \hN$, $\g >0$ and $q \in (1,2]$.

  \ss  \no 1)   Let  $b: [0,T] \times  \hR^k  \times \hU \times \hV  \to  \hR^k  $
  be a $\sB([0,T]) \otimes \sB(\hR^k) \otimes \sB(\hU) \otimes \sB(\hV)/\sB(\hR^k)-$measurable
  function and let $\si: [0,T] \times  \hR^k \times \hU \times \hV  \to  \hR^{k \times d} $
  be a $\sB([0,T]) \otimes \sB(\hR^k) \otimes \sB(\hU) \otimes \sB(\hV)/\sB(\hR^{k \times d})-$measurable function  such  that
     for any $(t,u,v)  \neg \in \neg  [0,T]  \neg \times \neg  \hU  \neg \times \neg  \hV$ and any $ x, x'    \in    \hR^k$
   \bea
   && \hspace{2.05cm} |b(t,0,u,v)  |     +    |\si(t,0,u,v) |
       \neg \le  \neg  \g \big(1 + [u]_{\overset{}{\hU}} + [v]_{\overset{}{\hV}}  \big)   \label{b_si_linear_growth}   \\
    \hb{and}  && |b(t,x,u,v) \neg - \neg b(t,x',u,v)|    +    |\si(t,x,u,v) \neg - \neg \si(t,x',u,v) |
       \neg \le  \neg  \g |x-x'|      ;  \qq     \label{b_si_Lip}
   \eea

      \no 2) Let   $\ul{l}, \, \ol{l}: [0,T]  \neg  \times   \hR^k  \neg \to \neg  \hR$ be two continuous functions such that
            $ \ul{l}(t,x)  \neg < \neg  \ol{l}(t,x)$, $  \fa (t,x)  \neg \in \neg  [0,T]  \neg \times \neg  \hR^k$ and that
    \bea
       |\ul{l}(t,x ) \neg - \neg \ul{l}(t,x' )|   \vee   |\ol{l}(t,x ) \neg - \neg \ol{l}(t,x' ) |
       \neg \le  \neg  \g |x-x'|^{2/q} ,  \q   \fa t \in [0,T] ,~ \fa  x, x'    \in    \hR^k;       \label{2l_growth}
   \eea

     \no 3) Let $h \neg : \hR^k  \neg \to \neg  \hR $ be a $ 2/q-$H\"older continuous function with coefficient $\g$ such that
     $\ul{l}(T,x)  \neg \le \neg h(x)  \neg \le \neg  \ol{l}(T,x)$, $  \fa  x  \neg \in \neg  \hR^k$;

  \ss \no 4)  Let   $f: [0,T] \times \hR^k \times \hR \times \hR^d \times \hU \times \hV \to \hR$  be
  $\sB([0,T]) \otimes \sB(\hR^k) \otimes \sB(\hR) \otimes \sB(\hR^d)  \otimes \sB(\hU) \otimes \sB(\hV)/\sB(\hR)-$measurable
  function 
%
%
 such that 
   for any  $(t,  u,v   ) \in [0,T]   \times \hU \times \hV $  and any $(x,y,z),  (x',y',z') \in   \hR^k \times \hR \times \hR^d$
     \bea
   |f(t,0,0,0,u,v)  |   & \dneg  \dneg  \le & \dneg \dneg  \g \Big(1 + [u]^{2/q}_{\overset{}{\hU}}
 + [v]^{2/q}_{\overset{}{\hV}}  \Big)   \label{f_linear_growth} \\
 \hb{and} \q \big|f(t,x,y,z,u,v)-f(t,x', y',z',u,v)\big| & \dneg \dneg\le & \dneg \dneg   \g \big(|x-x'|^{2/q} + |y-y'|+ |z-z'|\big)   .  \qq  \label{f_Lip}
 \eea

        For any $\l \ge 0$,    we let  $c_\l$ denote a generic constant,  depending on $\l$, $T$, $q$,  $\g $,
        $ \ul{l}_*  \neg  \dfnn  \neg  \underset{s \in [t,T]}{\sup}   \big| \ul{l} (s, 0 ) \big|$
      and $ \ol{l}_*  \neg \dfnn \neg  \underset{s \in [t,T]}{\sup}   \big| \ol{l} (s, 0 ) \big|$,
      whose form may vary from line to line. (In particular, $c_0$ stands for a generic constant depending
      on $ T $, $q$, $ \g $, $ \ul{l}_*$   and  $\ol{l}_*$.)

  \ms  Fix $t \in [0,T]$. Assume that when  player I (resp.\;II) select admissible control
  $\mu \in \cU^t$ (resp.\;$\nu \in \cV^t$),
    the corresponding state process  starting from time $t$ at point $x \in \hR^k$ 
    is  driven  by  the  SDE \eqref{FSDE} on the probability space $\big( \O^t, \ol{\cF}^t_{\neg T}  , P^t_0 \big)$.
   Clearly,     the measurability of $b$\, and $\si$ implies that for any $x' \in \hR^k$,
 both $\big\{b(s, x',\mu_s,\nu_s )\big\}_{s \in [t,T]} $
   and  $ \big\{\si(s, x',\mu_s,\nu_s )\big\}_{s \in [t,T]} $ are   $\bF^t-$progressively measurable processes.
 Also, we see from \eqref{b_si_Lip} that  for any $(s,\o) \in [t,T] \times \O^t$, both $b(s, \cd,\mu_s(\o),\nu_s(\o) )  $
 and $\si(s, \cd,\mu_s(\o),\nu_s(\o) )  $  are Lipschitz with coefficient $ \g $.  Since
   \bea   \label{eq:p311}
    E_t   \int_t^T  \neg \big(| b(s, 0,\mu_s,\nu_s )  |^2 + | \si (s, 0,\mu_s,\nu_s )  |^2 \big)
  \, ds \le c_0 + c_0 E_t   \int_t^T  \big(   [ \mu_s ]^2_{\overset{}{\hU}} + [ \nu_s ]^2_{\overset{}{\hV}} \big) ds < \infty
   \eea
 by \eqref{b_si_linear_growth}, it is well-known (see e.g. Theorem 2.5.7 of \cite{Krylov_CDP}) that
      \eqref{FSDE} admits a  unique  solution  $\left\{ X^{t,x,\mu,\nu}_s \right\}_{s \in [t,T]}   \\
    \in \hC^2_{\ol{\bF}^t }([t,T], \hR^k)$.

\ss  Applying Theorem 2.5.9 of \cite{Krylov_CDP}
  with $(\xi_s , \wt{\xi}_s, \wt{b}_s(0),\wt{\si}_s(0) \big)  \equiv  (x, 0, 0,0 )$ (thus $\wt{x}_s \equiv  0 $ therein)
   and using \eqref{b_si_linear_growth} yield
  \bea
     E_t \left[  \underset{s \in [t,T]}{\sup}  \big| X^{t,x,\mu,\nu}_s  \big|^2 \right]
      & \tneg \le& \tneg   c_0  \bigg( 1 + |x|^2 + E_t \int_t^T \big( [ \mu_s ]^2_{\overset{}{\hU}} + [ \nu_s ]^2_{\overset{}{\hV}}  \big) \, ds  \bigg) < \infty .                 \label{eq:esti_X_1}
 \eea

\subsection{Continuous Dependence Results}

  \begin{lemm}  \label{lem_DX_varpi}
 Let $\varpi  \in [1,2]$,  $t \in [0,T]$ and $(x, \mu,\nu) \in \hR^k \times \cU^t \times \cV^t$.
   \bea
\hb{\(1\) For any  $s \in [t,T]$,  } \,  E_t  \left[  \underset{r \in [  t,s]}{\sup}  \big| X^{t,x,\mu,\nu}_r \neg - \neg x  \big|^2 \right]       \le    c_0 (1 \neg + \neg |x|^2 )  (s \neg -\neg t )
 \neg + \neg  c_0 \, E_t  \neg  \int_t^s  \dneg \big( [\mu_r ]^2_{\overset{}{\hU}}   \neg + \neg   [\nu_r ]^2_{\overset{}{\hV}} \big) dr  . \hspace{2.5cm}    \label{eq:esti_X_2}   \\
\hb{\(2\) For any  $x' \in \hR^k$,  } \,  E_t  \left[  \underset{s \in [  t,T]}{\sup}  \big| X^{t,x,\mu,\nu}_s \neg - \neg X^{t,x',\mu,\nu}_s  \big|^{\varpi } \right]       \le    c_{\varpi }  |x-x'  |^{\varpi }  . \hspace{6.1cm} \label{eq:esti_X_3}
 \eea

 \no \(3\)  If $b $ and $ \si $ are  $\l-$H\"older continuous   in $u$ for some $\l \neg \in  \neg (0,1]$, i.e.,
 for any
    $   \big(\ol{t},\ol{x}   , \ol{v} \big)  \neg \in \neg  [0,T]  \neg \times \neg   \hR^k    \neg \times \neg  \hV$
    and  $\ol{u}_1, \ol{u}_2  \neg \in \neg  \hU$
     \bea
          \big|b(\ol{t}, \ol{x} ,\ol{u}_1  , \ol{v} ) \neg - \neg b(\ol{t}, \ol{x} ,\ol{u}_2  , \ol{v} )\big|
            +   \big|\si(\ol{t}, \ol{x} ,\ol{u}_1  , \ol{v} ) \neg - \neg \si(\ol{t}, \ol{x} ,\ol{u}_2  , \ol{v} ) \big|
   \le     \g   \,   \rho^\l_{\overset{}{\hU}} \big( \ol{u}_1, \ol{u}_2 \big)   ,   ~          \label{b_si_Lip_u}
   \eea
   then  for any  $    \mu'  \in    \cU^t$
  \bea
       E_t \bigg[ \underset{s \in [t,T]}{\sup} \Big|    X^{t,x,\mu,\nu}_s
    \neg -  \neg  X^{t,x,\mu',\nu}_s \Big|^{\varpi } \, \bigg]     \le
 c_{\varpi} E_t \Bigg[ \bigg( \int_t^T   \rho^{2\l}_{\overset{}{\hU}}(\mu_r,\mu'_r) \,  dr \bigg)^{\neg   \varpi/2 } \, \Bigg] .        \label{eq:s015}
\eea

   Similarly, if  $b $ and $ \si $ are  $\l-$H\"older continuous   in $u$ for some $\l \in (0,1]$, i.e.,  for any
    $   \big(\ol{t},\ol{x}   , \ol{u} \big) \in [0,T] \times  \hR^k   \times \hU$
    and  $\ol{v}_1, \ol{v}_2 \in \hV$
   \bea
          \big|b(\ol{t}, \ol{x} ,\ol{u}  , \ol{v}_1 ) \neg - \neg b(\ol{t}, \ol{x} ,\ol{u}  , \ol{v}_2 )\big|
            +   \big|\si(\ol{t}, \ol{x} ,\ol{u}  , \ol{v}_1 ) \neg - \neg \si(\ol{t}, \ol{x} ,\ol{u}  , \ol{v}_2 ) \big|
   \le     \g   \,   \rho^\l_{\overset{}{\hV}} \big( \ol{v}_1, \ol{v}_2 \big)   ,   ~          \label{b_si_Lip_v}
   \eea
    then for any  $    \nu\,'  \in    \cV^t$
\bea    \label{eq:s015b}
       E_t \bigg[ \underset{s \in [t,T]}{\sup} \Big|    X^{t,x,\mu,\nu}_s
    \neg -  \neg  X^{t,x,\mu,\nu\,'}_s \Big|^{\varpi } \, \bigg]     \le
 c_{\varpi} E_t \Bigg[ \bigg( \int_t^T   \rho^{2\l}_{\overset{}{\hV}}(\nu_r,\nu\,'_{\neg r} ) \,  dr \bigg)^{\neg \varpi/2 } \, \Bigg] .
\eea
\end{lemm}

   By  Lemma \ref{lem_F_version} (2), $ X^{t,x,\mu,\nu} $ admits a unique $\bF^t-$version $\wt{X}^{t,x,\mu,\nu}$,
  which clearly belongs to $  \hC^2_{\bF^t }([t,T], \hR^k)$ and  also  satisfies \eqref{eq:esti_X_1}, \eqref{eq:esti_X_2}, \eqref{eq:esti_X_3},
  \eqref{eq:s015} and \eqref{eq:s015b}.



 \ss If   $\big( \wt{\mu}, \wt{\nu} \big)$ is another pair of $\cU^t \neg \times \neg  \cV^t$ such that
  $ \big(  \mu,  \nu \big)  \neg = \neg  \big( \wt{\mu}, \wt{\nu} \big)  $ $dr  \neg \times \neg  d P^t_0-$a.s.
  on $\[t,\t\[$ for some $\bF^t-$stopping time $\t$, then both
  $\big\{X^{t,x,\mu,\nu}_{\t \land s}\big\}_{s \in [t,T]}$ and $\big\{X^{t,x,\wt{\mu},\wt{\nu}}_{\t \land s}\big\}_{s \in [t,T]}$
  satisfy the same SDE:
    \bea \label{eq:p605}
  X_s= x + \int_t^s   b_\t (r, X_r   ) \, dr + \int_t^s   \si_\t (r, X_r ) \, dB^t_r,  \q s \in [t, T] ,
    \eea
  with $b_\t(r,\o, x) \dfnn \b1_{\{r < \t(\o)\}} b\big(r, x,\mu_r(\o),\nu_r(\o) \big)$
  and $\si_\t(r,\o, x) \dfnn \b1_{\{r < \t(\o)\}} \si\big(r, x,\mu_r(\o),\nu_r(\o) \big)$,
  $ \fa (r,\o,x) \in [t,T] \times \O^t \times \hR^k$. Clearly, for any $ x \in \hR^k$  both $b_\t(\cd,\cd,x)$ and $\si_\t(\cd,\cd,x)$ are $\bF^t-$progressively measurable processes, and for any $(r,\o) \in [t,T] \times \O^t $ both
  $b_\t(r,\o,\cd)$ and $\si_\t(r,\o,\cd)$ are Lipschitz continuous with coefficient $\g$.   Thus \eqref{eq:p605} has a unique solution
  in $\hC^2_{\ol{\bF}^t}\big([t,T],\hR^k\big)$, i.e.
   \bea  \label{eq:p611}
  P^t_0 \big( X^{t,x,\mu,\nu}_{\t \land s} = X^{t,x,\wt{\mu},\wt{\nu}}_{\t \land s}, ~ \fa  s \in [t,T ] \big) = 1 .
   \eea

   Let $\Th$ stand for the  quadruplet $(t,x,\mu,\nu)$.
  By the continuity of   $\ul{l}$ and $\ol{l}$,    $\ul{L}^\Th_s \dfnn \ul{l} \big(s,\wt{X}^\Th_s \big)$  and
  $\ol{L}^{\, \raisebox{-0.5ex}{\scriptsize $\Th$}}_s \dfnn \ol{l} \big(s,\wt{X}^\Th_s \big)$, $s \in [t,T]$ are  two  real$-$valued, $ \bF^t-$adapted  continuous processes such that
    $ \ul{L}^\Th_s   < \ol{L}^{\, \raisebox{-0.5ex}{\scriptsize $\Th$}}_s  $, $\fa s \in [t,T]$.

 \ss  Given  an  $\bF^t-$stopping time  $\t$,   the measurability of $\big( f ,  \wt{X}^\Th , \mu, \nu \big)$
 and \eqref{f_Lip} imply that
   \beas
     f^\Th_\t   (s,\o,y,z ) \dfnn \b1_{\{s < \t(\o)\}}  f  \Big(s,  \wt{X}^\Th_s(\o) , \, y, z, \mu_s (\o),\nu_s (\o) \Big) ,
    \q   \fa (s,\o,y,z  ) \in   [t,T] \times \O^t \times \hR \times \hR^d
    \eeas
    is a    $\sP \big(  \bF^t  \big)   \otimes \sB(\hR)  \otimes \sB(\hR^d)/\sB(\hR)-$measurable function
      that is Lipschitz continuous in $(y,z)$ with coefficient $\g$.
  And  one can deduce from    \eqref{2l_growth}, \eqref{f_linear_growth}, \eqref{f_Lip}, H\"older's inequality  and \eqref{eq:esti_X_1} that
 \bea
      E_t \neg  \bigg[   \underset{s \in [t,T]}{\sup}    \big| \ul{L}^\Th_{\t \land s}\big|^q
   \neg +   \neg  \underset{s \in [t,T]}{\sup} \big| \ol{L}^{\, \raisebox{-0.5ex}{\scriptsize $\Th$}}_{\t \land s}\big|^q
  \neg+  \neg    \Big(\int_t^T  \neg \big| f^\Th_\t   (s, 0,0 )   \big|  ds \Big)^q\bigg]
   \neg  \le  \neg  c_0    \neg + \neg  c_0   E_t \bigg[   \underset{s \in [t,T]}{\sup}   \big| \wt{X}^\Th_s \big|^2
   \neg + \neg \int_t^T \dneg \big( [ \mu_s ]^2_{\overset{}{\hU}} \neg+\neg [ \nu_s ]^2_{\overset{}{\hV}}  \big)  ds \bigg]
    \neg < \neg  \infty .   \q
  \label{eq:s031}
  \eea
  Thus, for any $\cF^t_\t-$measurable random variable $\xi$ with $ \ul{L}^\Th_\t \le  \xi \le \ol{L}^\Th_\t $, $P^t_0-$a.s., it follows that  $  E_t \big[  | \xi  |^q  \big]     < \infty $, i.e.  $\xi \in \hL^q \big( \cF^t_\t \big)$.
   Then Theorem 4.1  of  \cite{EHW_2011} shows that
    the DRBSDE$\big(P^t_0,\xi,  f^\Th_\t,\ul{L}^\Th_{\t \land \cd},\ol{L}^{\,\raisebox{-0.5ex}{\scriptsize $\Th$}}_{\t \land \cd}\big)$
    admits a unique solution $ \big(Y^\Th (\t, \xi),Z^\Th (\t, \xi),\ul{K}^{\Th}(\t, \xi),  \ol{K}^{\, \raisebox{-0.5ex}{\scriptsize $\Th$}} (\t, \xi) \big) \neg \in \neg \hG^q_{\ol{\bF}^t } \big([t,T]\big) $. Clearly,
      its   $\bF^t-$version  $ \Big(\wt{Y}^\Th(\t, \xi) ,  \wt{Z}^\Th(\t, \xi),  \wt{\ul{K}}^{\,\raisebox{-0.5ex}{\scriptsize $\Th$}}(\t, \xi), \\
    \wt{\ol{K}}^{\, \raisebox{-0.8ex}{\scriptsize $\Th$}} (\t, \xi) \Big) $ by Lemma \ref{lem_F_version} (2) belongs to $  \hG^q_{\bF^t} \big([t,T]\big)$.   As $ \cF^t_t = \{\es, \O^t\} $,   $\wt{Y}^\Th_t (\t, \xi)$   is  a constant.

 \ss      Given another $\bF^t-$stopping time  $\z$ such that $\z \le \t  $, $P^t_0-$a.s.,  one can easily show that
 $  \Big\{ \Big( \wt{Y}^\Th_{\z \land s}   (\t,\xi), \b1_{\{s < \z \}}\wt{Z}^\Th_s   (\t,\xi) ,  \\   \wt{\ul{K}}^{\,\raisebox{-0.5ex}{\scriptsize $\Th$}}_{\z \land s}   (\t,\xi),
  \wt{\ol{K}}^{\, \raisebox{-0.8ex}{\scriptsize $\Th$}}_{\z \land s}
    (\t,\xi) \Big) \Big\}_{s \in [t,T]} \in \hG^q_{\bF^t } \big([t,T]\big)$
 solves the DRBSDE$\Big( P^t_0, \wt{Y}^\Th_\z  (\t,\xi),f^\Th_\z ,
 \ul{L}^\Th_{\z \land \cd}, \ol{L}^{\,\raisebox{-0.5ex}{\scriptsize $\Th$}}_{\z \land \cd}   \Big)$. To wit, we have
 \bea
\Big(\wt{Y}^\Th_s \big(\z,  \wt{Y}^\Th_\z  (\t,\xi)\big) ,  \wt{Z}^\Th_s \big(\z,  \wt{Y}^\Th_\z  (\t,\xi)\big),  \wt{\ul{K}}^{\,\raisebox{-0.5ex}{\scriptsize $\Th$}}_s \big(\z,  \wt{Y}^\Th_\z  (\t,\xi)\big),
 \wt{\ol{K}}^{\, \raisebox{-0.8ex}{\scriptsize $\Th$}}_s \big(\z,  \wt{Y}^\Th_\z  (\t,\xi)\big) \Big)
 \nonumber    \\
  =  \Big( \wt{Y}^\Th_{\z \land s}   (\t,\xi), \b1_{\{s < \z \}}\wt{Z}^\Th_s   (\t,\xi) ,  \wt{\ul{K}}^{\,\raisebox{-0.5ex}{\scriptsize $\Th$}}_{\z \land s}   (\t,\xi),
  \wt{\ol{K}}^{\, \raisebox{-0.8ex}{\scriptsize $\Th$}}_{\z \land s}   (\t,\xi) \Big), \q s \in [t,T].   \label{eq:p677}
 \eea

   The continuity of functions $ h$  implies that
   $ h\big( \wt{X}^\Th_T \big)$ is a real$-$valued, $  \cF^t_T -$measurable random variables
   such that   $  \ul{L}^\Th_T = \ul{l} \big(T,\wt{X}^\Th_T \big) \le h \big( \wt{X}^\Th_T \big) \le
      \ol{l} \big(T,\wt{X}^\Th_T \big) = \ol{L}^\Th_T $.
    Hence, we can use \eqref{DRBSDE01} to obtain  that
   \beas
   \ul{l}(t,x)
   = \ul{l}\big(t, \wt{X}^\Th_t \big) 
   \le  Y^\Th_t   \big(T,h \big(\wt{X}^\Th_T \big)  \big)
   \le \ol{l}\big(t, \wt{X}^\Th_t \big)
    =  \ol{l}(t,x) , \q P^t_0-a.s.,
 \eeas
  which leads to  that
       \bea   \label{eq:j701}
 \ul{l}(t,x) \le \wt{Y}^\Th_t   \big(T,h \big(\wt{X}^\Th_T \big)  \big) \le \ol{l}(t,x) .
   \eea

 Inspired by Proposition 6.1 (ii) of \cite{Buckdahn_Li_3}, we have  the following a priori estimates about the dependence of $Y^{t,x,\mu,\nu}_s  \Big(T, h\big(\wt{X}^{t,x,\mu,\nu}_T \big) \Big) $
on initial state $x$ and on  controls $(\mu,\nu)$.

\begin{lemm} \label{lem_estimate_Y}
 Let $\varpi \in (1,q]$,   $t \in [0,T]$ and $(x, \mu,\nu) \in \hR^k \times \cU^t \times \cV^t$.
  \bea
    \hb{(1) For any  } x' \in \hR^k, ~
   \, E_t \bigg[ \, \underset{s \in [t,T]}{\sup} \Big|  Y^{t,x,\mu,\nu}_s  \big(T, h\big(\wt{X}^{t,x,\mu,\nu}_T \big) \big)
     - Y^{t,x',\mu,\nu}_s  \big(T, h\big(\wt{X}^{t,x',\mu,\nu}_T \big) \big) \Big|^\varpi \, \bigg]
    \le c_\varpi | x \neg - \neg x' |^{\frac{2 \varpi}{q}}   . \hspace{1cm} \label{eq:s025}
     \eea
 (2)   Let $\ul{l},\ol{l}$ and $h$   satisfy
   \bea
       |\ul{l}(t,x ) \neg - \neg \ul{l}(t,x' )|    \vee    |\ol{l}(t,x ) \neg - \neg \ol{l}(t,x' ) |
       \vee |h( x ) \neg - \neg h( x' ) |
       \neg \le  \neg  \g \, \psi \big( |x-x'| \big) ,  \q   \fa t \in [0,T] ,~ \fa  x, x'    \in    \hR^k          \label{2l_growth2}
   \eea
   for an increasing $\hC^2  (\hR_+  ) $ function $\psi $ such that for some $0<R_1<1< R_2$
   \beas
         \psi(a)    = \frac12 a^2   \,   \hb{ if } a \in [0, R_1] , \;~
        \psi(a) \le a^{2/q}     \,    \hb{ if } a \in (R_1, R_2) ,  \hb{\; and \;}
    \psi(a) = a^{2/q}    \,  \hb{ if } a > R_2 .
         \eeas
      \(Clearly,  $\psi(a) \le a^{2/q}$ for any $a \ge 0$. So   \eqref{2l_growth2} implies \eqref{2l_growth}
     and the $2/q-$H\"older continuity of $h$.\)

  \ss Let $\l \in (0,1]$.   If $b, \si $ are $\dis \l - $H\"older continuous in $u$ \big(see \eqref{b_si_Lip_u}\big)
 and  if $f$ is $\dis 2\l - $H\"older continuous in $u$, i.e.,  for any
    $   \big(\ol{t},\ol{x} ,\ol{y} , \ol{z} , \ol{v} \big) \in [0,T] \times  \hR^k \times \hR \times \hR^d    \times \hV$
    and  $\ol{u}_1, \ol{u}_2 \in \hU$
   \bea
            \big|f(\ol{t}, \ol{x} ,\ol{y} , \ol{z} ,\ol{u}_1  , \ol{v} ) \neg -   \dneg   f(\ol{t}, \ol{x} ,\ol{y} , \ol{z} ,\ol{u}_2  , \ol{v})\big|
       \neg \le  \neg  \g   \,     \rho^{ 2  \l }_{\overset{}{\hU}}(\ol{u}_1, \ol{u}_2)    ,   ~          \label{f_Lip_u}
   \eea
 then  for any   $  \mu'    \in    \cU^t$
  \bea
  &&  E_t \bigg[ \underset{s \in [t,T]}{\sup} \Big|   Y^{t,x,\mu,\nu}_s  \big(T, h\big(\wt{X}^{t,x,\mu,\nu}_T \big) \big)
     - Y^{t,x,\mu',\nu}_s  \big(T, h\big(\wt{X}^{t,x,\mu',\nu}_T \big) \big) \Big|^{\varpi} \, \bigg]  \nonumber  \\
  && \hspace{3cm}  \le   c_\varpi  \k^\varpi_\psi   \Bigg\{   E_t \bigg[ \Big( \int_t^T \rho^{2\l}_{\overset{}{\hU}}(\mu'_s , \mu_s)    ds  \Big)^{\frac{\varpi}{2}}\bigg]     + E_t \bigg[ \Big( \int_t^T \rho^{2\l}_{\overset{}{\hU}}(\mu'_s , \mu_s)    ds  \Big)^{\varpi}\bigg] \Bigg\}   ,       \label{eq_estimate_Y_u}
  \eea
   where $\k_\psi \dfnn  \Big( 2   + R^{-1}_1 \neg \underset{a \in [R_1,R_2]}{\sup } 1 \vee \psi'(a)
   + \neg \underset{a \in [R_1,R_2]}{\sup } \big| \psi''(a) \big| \Big) R^{2-\frac{2}{q}}_2 $.

 Similarly, if $b, \si $ are $\dis \l - $H\"older continuous in $v$ \big(see \eqref{b_si_Lip_v}\big) and
 if $f$ is additionally $\dis 2\l - $H\"older continuous in $v$, i.e.,  for any
    $   \big(\ol{t},\ol{x} ,\ol{y} , \ol{z} , \ol{u} \big) \in [0,T] \times  \hR^k \times \hR \times \hR^d    \times \hU$
    and  $\ol{v}_1, \ol{v}_2 \in \hV$
     \bea          \label{f_Lip_v}
  \big|f(\ol{t}, \ol{x} ,\ol{y} , \ol{z} ,\ol{u}  , \ol{v}_1 ) \neg -   \dneg   f(\ol{t}, \ol{x} ,\ol{y} , \ol{z} ,\ol{u}  , \ol{v}_2 )\big| \neg \le  \neg  \g   \,     \rho^{ 2  \l }_{\overset{}{\hV}}  ( \ol{v}_1,\ol{v}_2  )    ,   ~
     \eea
 then  for any   $  \nu \, '    \in    \cV^t$
  \bea
   && E_t \bigg[ \underset{s \in [t,T]}{\sup} \Big|   Y^{t,x,\mu,\nu}_s  \big(T, h\big(\wt{X}^{t,x,\mu,\nu}_T \big) \big)
     - Y^{t,x,\mu,\nu \, '}_s  \big(T, h\big(\wt{X}^{t,x,\mu,\nu \, '}_T \big) \big) \Big|^{\varpi} \, \bigg] \nonumber  \\
  && \hspace{3cm}  \le   c_\varpi  \k^\varpi_\psi   \Bigg\{   E_t \bigg[ \Big( \int_t^T \rho^{2 \l}_{\overset{}{\hV}}(\nu\,'_s , \nu_s)    ds  \Big)^{\frac{\varpi}{2}}\bigg]     + E_t \bigg[ \Big( \int_t^T \rho^{2 \l}_{\overset{}{\hV}}(\nu\,'_s , \nu_s)    ds  \Big)^{ \varpi }\bigg] \Bigg\} .   \label{eq_estimate_Y_v}
  \eea

\end{lemm}

\subsection{Definition of the value functions and the Dynamic Programming Principle}

\ss  Now, we are ready to introduce   values of the zero-sum stochastic differential games via the following
  notion of admissible strategies.

   \begin{deff}
  Given $t \in [0,T]$, an admissible strategy $\a$  for player I over period $[t,T]$
  is a $\hU -$valued function $\a$ on $[t,T] \neg \times \neg  \O^t  \neg \times \neg  \hV$
  that is $\sP \big(  \bF^t  \big)  \neg \otimes \neg  \sB(\hV)    \big/    \sB(\hU) -$measurable and satisfies:
   \(i\) $ \a(r,  \hV_0) \neg \subset \neg  \hU_0  $, $dr  \neg \times \neg  d P^t_0-$a.s.
   \(ii\)  For a $\k   >   0$ and a  non-negative  measurable process $\Psi$ on $(\O^t,\cF^t_T) $
   with $ E_t \neg \int_t^T \neg   \Psi^2_r  \,  dr   <   \infty$,
   it holds  $dr  \neg \times \neg  d P^t_0-$a.s. that
   \bea    \label{eq:r503}
    \big[ \a(r,\o,v)  \big]_{\overset{}{\hU}} \le \Psi_r (\o) + \k [v]_{\overset{}{\hV}} \, , \q \fa v \in \hV .
   \eea
  Admissible strategies $\beta: [t,T] \times \O^t \times \hU \to \hV$   for player II  over period $[t,T]$ are defined similarly.
  The set of all  admissible strategies for player I \(resp.\;II\)  on $[t,T]$ is denoted by $\cA^t$ \(resp.\;$\fB^t$\).
 \end{deff}

 Given $t \in [0,T]$, an admissible strategy $\a \in \cA^t$  induces a mapping $\a\lan \cd \ran : \cV^t \to \cU^t$ by
  \beas
  \big( \a\lan \nu \ran \big)_r (\o) \dfnn \a \big(r, \o, \nu_r(\o)\big), \q  \fa \nu \in \cV^t, ~  (r,\o) \in [t,T] \times \O^t .
  \eeas
  To see this, let $ \nu \in \cV^t $. Clearly, $ \a\lan \nu \ran $ is a $\hU-$valued, $\bF^t-$progressively measurable process.
   Since $ \big\{ (r,\o) \in [t,T] \times \O^t:    \nu_r (\o) \in   \hV_0 \big\}
    \cap \big\{(r,\o) \in [t,T] \times \O^t: \a(r,\o,  \hV_0) \subset   \hU_0 \big\}
    \subset \big\{ (r,\o) \in [t,T] \times \O^t: \big( \a\lan \nu \ran \big)_r (\o) \in   \hU_0 \big\}$,
    it holds $dr \times dP^t_0-$a.s. that $\a\lan \nu \ran   \in   \hU_0$.
    On the other hand, one can deduce that
   \beas
   E_t \neg \int_t^T \dneg \big[ \big( \a\lan \nu \ran \big)_r \big]^2_{\overset{}{\hU}} \, dr
   \le 2 E_t \neg \int_t^T \neg   \Psi^2_r  \,  dr
   +   2 \k^2 E_t \neg \int_t^T  [  \nu_r  ]^2_{\overset{}{\hV}} \, dr  < \infty   .
   \eeas
    Thus, $\a\lan \nu \ran \in \cU^t$.
    If $\nu^1  \in \cV^t $ is equal to $\nu^2 \in \cV^t$, $d r \times d P^t_0-$a.s.  on $\[t,\t\[   $
    for any   $\bF^t-$stopping time $\t$,
   then $\a \lan \nu^1 \ran =\a  \lan \nu^2 \ran $,   $d r \times d P^t_0-$a.s.  on $\[t,\t\[  $. So $\a\lan \cd \ran$
   is exactly an Elliott$-$Kalton strategy  considered in e.g. \cite{Fleming_1989}.
   Similarly, any $\beta \in \fB^t$  gives rise to  a mapping $\beta \lan \cd \ran : \cU^t \to \cV^t$.

 \begin{deff} \label{def_tilde_strategy}
 Given $t \in [0,T]$,  an $  \cA^t -$strategy $\a$   is said to be of $\wh{\cA}^t$  if for any $\e  \neg > \neg  0$,
 there exist  a $\d  >0$ and  a closed subset $F$ of $\O^t$ 
 with    $   P^t_0\big( F \big)  \neg > \neg  1-\e$ such that for any   $\o ,\o'   \in    F$ with $ \|\o-\o' \|_t < \d $
  \bea   \label{eq:r742}
 \underset{r \in [t,T]}{\sup} \, \underset{v \in  \hV}{\sup} \; \rho_{\overset{}{\hU}} \big( \a(r,\o,v) , \a(r,\o',v) \big) < \e .
   \eea
We define $\wh{\fB}^t \subset \fB^t$ similarly.

 \end{deff}

     For any $(t,x) \in  [0,T] \times  \hR^k $,    
   we define
  \beas
\q  w_1 (t,x) \dfnn \underset{\beta \in \fB^t }{\inf} \; \underset{\mu \in \cU^t}{\sup} \; \wt{Y}^{t,x,\mu, \beta \lan \mu \ran }_t   \Big(T,h \Big(\wt{X}^{t,x,\mu, \beta \lan \mu \ran}_T \Big)  \Big) ~\; \hb{and } ~\;   \wh{w}_1 (t,x) \dfnn \underset{\beta \in \wh{\fB}^t }{\inf} \; \underset{\mu \in \cU^t}{\sup} \; \wt{Y}^{t,x,\mu, \beta \lan \mu \ran }_t   \Big(T,h \Big(\wt{X}^{t,x,\mu, \beta \lan \mu \ran}_T \Big)  \Big)
  \eeas
as    player I's  {\it priority value} and {\it intrinsic priority  value} of the zero-sum stochastic differential  game that starts   from time $t  $ and state $x$. Correspondingly, we define
\beas
\q  w_2 (t,x )  \dfnn    \underset{\a \in \cA^t}{\sup} \, \underset{\nu \in \cV^t}{\inf} \; \wt{Y}^{t,x,\a \lan \nu \ran,\nu}_t  \Big(T, h \Big(\wt{X}^{t,x,\a \lan \nu \ran,\nu}_T  \Big) \Big) ~\; \hb{and } ~\;   \wh{w}_2 (t,x )  \dfnn    \underset{\a \in \wh{\cA}^t}{\sup} \, \underset{\nu \in \cV^t}{\inf} \; \wt{Y}^{t,x,\a \lan \nu \ran,\nu}_t  \Big(T, h \Big(\wt{X}^{t,x,\a \lan \nu \ran,\nu}_T  \Big) \Big)
  \eeas
  as    player  II's {\it priority value} and {\it intrinsic priority   value} of the zero-sum stochastic differential  game that starts   from time $t  $ and state $x$.   By  \eqref{eq:j701}, one has
  \bea   \label{eq:n411}
  \ul{l}(t,x) \le w_1 (t,x) \le \wh{w}_1 (t,x)  \le \ol{l}(t,x)
  \q \hb{and} \q \ul{l}(t,x) \le \wh{w}_2 (t,x)  \le  w_2 (t,x) \le  \ol{l}(t,x) .
  \eea

  The two obstacle functions $\ul{l}$, $\ol{l}$ as well as  the DRBSDE structure  prevent the value functions
 from taking $\pm \infty$ values. 
The values $w_1(t,x )$ and $\wh{w}_1(t,x )$, otherwise,  might blow up unless
  we   impose additional integrability conditions on $\cU^t$ and $\fB^t$ analogous to e.g.
 Assumption 5.7 of \cite{STZ_2011c}.


\begin{rem}
  Given $t \in [0,T]$, we can regard  $\mu \in \cU^t$ as  a member of $ \cA^t $ since
    \beas
     \a^\mu (r,\o,v) \dfnn  \mu_r(\o), \q
  \fa (r,\o,v) \in [t,T]   \times    \O^t    \times    \hV
  \eeas
    is clearly  a
  $\sP \big(  \bF^t  \big)    \otimes \sB(\hV)    \to    \sB(\hU) -$measurable function
  such that $\a^\mu (r,\hV_0) = \mu_r \in \hU_0  $, $dr \times d P^t_0-$a.s. and  that
  \eqref{eq:r503} holds for $\Psi^\a  =[ \mu  ]_{\overset{}{\hU}}$ and any $\k_\a>0$. Similarly, $ \cV^t $  can be embedded into
  $ \fB^t$.  Then it follows   that
  \beas
 \q  w_1(t,x) \le \underset{\nu \in \cV^t}{\inf} \, \underset{\mu \in \cU^t}{\sup} \; \wt{Y}^{t,x,\mu, \nu}_t \Big(T,h\big(\wt{X}^{t,x,\mu, \nu}_T \big)  \Big) \q \hb{and} \q
  w_2(t,x ) \ge \underset{\mu \in \cU^t}{\sup}\, \underset{\nu \in \cV^t}{\inf} \; \wt{Y}^{t,x,\mu, \nu}_t \Big(T,h\big(\wt{X}^{t,x,\mu, \nu}_T \big)  \Big)   .
   \eeas
 However, the fact that $ \underset{\mu \in \cU^t}{\sup}\, \underset{\nu \in \cV^t}{\inf} \; \wt{Y}^{t,x,\mu, \nu}_t \Big(T,h\big(\wt{X}^{t,x,\mu, \nu}_T \big)  \Big) \le \underset{\nu \in \cV^t}{\inf} \, \underset{\mu \in \cU^t}{\sup} \; \wt{Y}^{t,x,\mu, \nu}_t \Big(T,h\big(\wt{X}^{t,x,\mu, \nu}_T \big)  \Big) $
 does not necessarily imply that $ w_2(t,x )   \le w_1(t,x)$.
\end{rem}


  Let \,$t \in [0,T]$  and  $x_1,x_2 \in \hR^k$.    For any $\beta \in \fB^t$ and $\mu \in \cU^t  $,
   \eqref{eq:s025} shows that
  \beas   
    E_t \bigg[ \underset{s \in [t,T]}{\sup} \Big| \, \wt{Y}^{t,x_1,\mu,\beta \lan \mu \ran}_s
     \big(T, h\big(\wt{X}^{t,x_1,\mu,\beta \lan \mu \ran}_T \big) \big)
     - \wt{Y}^{t,x_2,\mu,\beta \lan \mu \ran}_s  \big(T, h\big(\wt{X}^{t,x_2,\mu,\beta \lan \mu \ran}_T \big) \big) \Big|^q
     \, \bigg]  \le c_0 | x_1 \neg - \neg x_2 |^2  .
  \eeas
  It then follows that
   \bea
 && \hspace{-2cm}  \wt{Y}^{t,x_2,\mu,\beta \lan \mu \ran}_t  \big(T, h\big(\wt{X}^{t,x_2,\mu,\beta \lan \mu \ran}_T \big) \big)
  \neg  - \neg  c_0 | x_1 \neg - \neg x_2 |^{2/q}
    \neg  \le \neg  \wt{Y}^{t,x_1,\mu,\beta \lan \mu \ran}_t  \big(T, h\big(\wt{X}^{t,x_1,\mu,\beta \lan \mu \ran}_T \big) \big) \nonumber \\
    &&    \neg  \le \neg  \wt{Y}^{t,x_2,\mu,\beta \lan \mu \ran}_t  \big(T, h\big(\wt{X}^{t,x_2,\mu,\beta \lan \mu \ran}_T \big) \big)
     \neg + \neg  c_0 | x_1 \neg - \neg x_2 |^{2/q} .   \label{eq:k115}
  \eea
   Taking supremum over $\mu \in \cU^{t }$ and then taking infimum over $\beta \in \fB^{t }$ yield that
  \beas
   w_1(t,x_2)    -    c_0 | x_1 \neg - \neg x_2 |^{2/q}  \le    w_1(t,x_1)
            \le    w_1(t,x_2)   +    c_0 | x_1 \neg - \neg x_2 |^{2/q}    .
  \eeas
  Thus $ \big|   w_1(t,x_1) \neg - \neg  w_1(t,x_2) \big| \le c_0 | x_1 \neg - \neg x_2 |^{2/q} $,
  and one can deduce the similar inequalities   for     $\wh{w}_1 $, $w_2 $  and $\wh{w}_2 $:

 \begin{prop}  \label{prop_w_conti}
 For any $t \in [0,T]$  and  $x_1,x_2 \in \hR^k$, we have
\beas
  \q  \big| w_1(t,x_1) \neg - \neg  w_1(t,x_2) \big| +\big| \wh{w}_1(t,x_1) \neg - \neg  \wh{w}_1(t,x_2) \big|
  +\big| w_2(t,x_1) \neg - \neg  w_2(t,x_2) \big| +\big| \wh{w}_2(t,x_1) \neg - \neg  \wh{w}_2(t,x_2) \big|
   \le c_0 | x_1 \neg - \neg x_2 |^{2/q} .
\eeas

 \end{prop}


   However, these value functions 
   are generally not $1/q-$H\"older continuous in $t$ unless the 
control spaces are compact. 


 \begin{rem} \label{r_strategy}

 When trying to directly prove the dynamic programming principle, 
 \cite{Fleming_1989} encountered a measurability issue; see  page 299 therein.
 To overcome this technical difficulty, they  first proved that the value functions
 are unique viscosity solutions to the associated Bellman-Isaacs equations
 by a   time-discretization approach \(assuming that the limiting Isaacs equation has a comparison principle\), which relies on the following regularity of the approximating  values $v_\pi$
   \beas
   |v_\pi(t,x)- v_\pi(t',x')| \le C \big( |t-t'|^{1/2} + |x-x'|\big) \q \fa (t,x), (t',x') \in [0,T] \times \hR^k
   \eeas
 with a uniform coefficient $C>0$  for all partitions $\pi$  of $[0,T]$.
 Since our value functions are  not $1/2-$H\"older continuous in $t$ given $q=2$, 
 this method does not work in general under our assumptions.
 Instead, we specify Elliott$-$Kalton strategies as  measurable random fields from one control space to another in order to
 avoid similar measurability issues when pasting strategies \(see Proposition \ref{prop_paste_strategy}\).
 This is a crucial ingredient in the proof of the supersolution (resp. subsolution) side of the dynamic programming principle \(Theorem \ref{thm_DPP}\) for $w_1$(resp. $w_2$).
 \end{rem}

 Given   $i=1,2$, 
 since    $w_i(t,\cd)$ is continuous   for any $t \in [0,T]$,
 one can deduce that for any  $ \bF^t-$stopping time $\t $ with countably many values $\{t_n\}_{n \in \hN} \subset [t,T]$,
 and any $\hR^k-$valued, $ \cF^t_\t -$measurable random variable $\xi$
  \bea \label{eq:n333}
    w_i(\t, \xi) =\sum_{n \in \hN} \b1_{\{ \t = t_n \}} w_i (t_n,\xi)  \hb{  is     $\cF^t_\t -$measurable. }
  \eea
 Similarly, $\wh{w}_i(\t, \xi) $   is  also   $\cF^t_\t -$measurable.


   \ss   Then we have the following dynamic programming principle for value functions.

 \begin{thm}  \label{thm_DPP}
 Let    $(t,x) \neg \in \neg [0,T]   \times   \hR^k $. For any family $\{\t_{\mu,\beta} \neg : \mu \neg \in \neg \cU^t, \beta \neg \in  \neg \fB^t\}$ of $\hQ_{t,T}  -$valued, $\bF^t-$stopping times
  \bea  \label{eq:DPP0}
  w_1(t,x)& \le &    \underset{\beta \in \fB^t}{\inf} \; \underset{\mu \in \cU^t}{\sup}    \;
    \wt{Y}^{t,x,\mu, \beta \lan \mu \ran}_t \Big(\t_{\mu,\beta},
    w_1\big(\t_{\mu,\beta},  \wt{X}^{t,x,\mu, \beta \lan \mu \ran}_{\t_{\mu,\beta}} \big) \Big) \\
  \hb{and}   \q    \wh{w}_1(t,x)& \le &    \underset{\beta \in \wh{\fB}^t}{\inf} \; \underset{\mu \in \cU^t}{\sup}    \;
    \wt{Y}^{t,x,\mu, \beta \lan \mu \ran}_t \Big(\t_{\mu,\beta},
    \wh{w}_1 \big(\t_{\mu,\beta},  \wt{X}^{t,x,\mu, \beta \lan \mu \ran}_{\t_{\mu,\beta}} \big) \Big) ;    \label{eq:DPP1}
    \eea
 the reverse inequality (of \eqref{eq:DPP1}) holds if
  \beas
 ( \bV_\l )  \qq  \bc \hb{$\ul{l},\ol{l}$ and $h$   satisfy \eqref{2l_growth2};  $b, \si$ are
 $\l-$H\"older continuous  in $v$ \big(see \eqref{b_si_Lip_v}\big); and} \\
  \hb{$f$ is $2\l-$H\"older continuous  in $v$
  \big(see \eqref{f_Lip_v}\big) for some }\l \in (0,1).
  \ec
  \eeas

 \ss On the other hand, for    any family $\{\t_{\nu,\a} \neg: \nu \neg\in \neg \cV^t, \a \neg\in \neg\cA^t \}$
  of $\hQ_{t,T} -$valued, $\bF^t-$stopping times
    \bea \label{eq:DPP2}
     w_2 (t,x) & \ge &    \underset{\a \in \cA^t }{\sup} \, \underset{\nu \in \cV^t }{\inf} \;
      \wt{Y}^{t,x,\a \lan \nu \ran,\nu}_t \Big(\t_{\nu,\a},  w_2 \big(\t_{\nu,\a},  \wt{X}^{t,x,\a \lan \nu \ran,\nu}_{\t_{\nu,\a}} \big) \Big) \\
   \hb{and}   \q     \wh{w}_2 (t,x) & \ge &    \underset{\a \in \wh{\cA}^t }{\sup} \, \underset{\nu \in \cV^t }{\inf} \;
      \wt{Y}^{t,x,\a \lan \nu \ran,\nu}_t \Big(\t_{\nu,\a},  \wh{w}_2 \big(\t_{\nu,\a},  \wt{X}^{t,x,\a \lan \nu \ran,\nu}_{\t_{\nu,\a}} \big) \Big);  \label{eq:DPP3}
 \eea
    the reverse inequality of \eqref{eq:DPP3} holds if
     \beas
 ( \bU_\l )  \qq  \bc \hb{$\ul{l},\ol{l}$ and $h$   satisfy \eqref{2l_growth2};  $b, \si$ are
 $\l-$H\"older continuous  in $u$ \big(see \eqref{b_si_Lip_u}\big); and} \\
  \hb{$f$ is $2 \l-$H\"older continuous  in $u$
  \big(see \eqref{f_Lip_u}\big) for some }\l \in (0,1).
  \ec
  \eeas

 \end{thm}

 Note that each $\wt{Y}^{t,x,\mu, \beta \lan \mu \ran}_t \Big(\t_{\mu,\beta},
    w_1\big(\t_{\mu,\beta},  \wt{X}^{t,x,\mu, \beta \lan \mu \ran}_{\t_{\mu,\beta}} \big) \Big)$ in \eqref{eq:DPP0}
 is well-posed since $ w_1 \Big(\t_{\mu,\beta},  \wt{X}^{t,x,\mu, \beta \lan \mu \ran}_{\t_{\mu,\beta}} \Big)$ is
$\cF^t_{\t_{\mu,\beta}}-$ measurable by \eqref{eq:n333} and since
 \beas
 \ul{L}^{t,x,\mu, \beta \lan \mu \ran}_{\t_{\mu,\beta}} = \ul{l} \Big( \t_{\mu,\beta},  \wt{X}^{t,x,\mu, \beta \lan \mu \ran}_{\t_{\mu,\beta}}\Big)  \le  w_1\Big(\t_{\mu,\beta},  \wt{X}^{t,x,\mu, \beta \lan \mu \ran}_{\t_{\mu,\beta}} \Big)
 \le \ol{l} \Big( \t_{\mu,\beta},  \wt{X}^{t,x,\mu, \beta \lan \mu \ran}_{\t_{\mu,\beta}}\Big)
 =  \ol{L}^{t,x,\mu, \beta \lan \mu \ran}_{\t_{\mu,\beta}}
 \eeas
 by \eqref{eq:n411}. The proof of Theorem \ref{thm_DPP} (see Subsection \ref{subsection:DPP}) relies on
  \eqref{eq_estimate_Y_u}, \eqref{eq_estimate_Y_v},  properties
 of shifted processes (especially shifted SDEs) as well as stability under pasting of controls/strategies,
  the latter two  of which will be  discussed  in Section \ref{sec:shift_prob}.

\section{An Obstacle Problem for Fully non-linear PDEs}

\label{sec:PDE}

\ms In this section,  we show that  the (intrinsic) priority values
   are (discontinuous) viscosity solutions of the following obstacle problem of a  PDE with a
fully non-linear Hamiltonian $H$:
  \bea \label{eq:PDE}
 \hspace{-0.3cm}  \min \neg \Big\{ \neg (w  \neg -  \neg \ul{l})(t,x),   \max \neg \Big\{ \dneg - \dneg \frac{\pa }{\pa t} w(t,x)
     \neg  - \neg  H \big(t,x, w(t,x), D_x w(t,x), D^2_x w(t,x)\big),
    (w \neg - \neg \ol{l})(t,x)  \Big\} \neg  \Big\} \neg = \neg  0 ,    \,    \fa (t,x)
 \neg \in  \neg  (0,T)  \neg \times  \neg  \hR^k .    ~ \;       
\eea

 \if{0}
 \bea \label{eq:PDE}
   \min \Big\{ (w-\ul{l})(t,x), ~ \max \Big\{ - \neg \frac{\pa }{\pa t} w(t,x)
    - H \big(t,x, w(t,x), D_x w(t,x), D^2_x w(t,x)\big),
    (w-\ol{l})(t,x)  \Big\} \Big\}= 0 , ~ \fa (t,x) \in (0,T) \times \hR^k  .       
\eea
 \fi

\begin{deff}  \label{def:viscosity_solution}
 Let   $H: [0,T] \times   \hR^k \times \hR  \times   \hR^k
  \times \hS_k \to [-\infty, \infty]$ be an  upper  \(resp.\;lower\)  semicontinuous functions with  $\hS_k$ denoting
   the set of all  $\hR^{k \times k}-$valued symmetric matrices.
  An upper \(resp.\;lower\) semicontinuous function $w :  [0, T] \times \hR^k \to \hR$    is called  a viscosity subsolution \(resp. supersolution\) of \eqref{eq:PDE} if   $w(T, x) \le $ \(resp. $\ge$\) $ h (x)$,   $\fa x \in   \hR^k$,
  and if    for any $(t_0,x_0, \vf ) \in (0,T) \times \hR^k \times \hC^{1,2}\big([0,T] \times \hR^k\big)$
  such that $w(t_0,x_0) = \vf (t_0,x_0)$ and that
  $w - \vf$ attains a strict local  maximum \(resp.\;strict local  minimum\) at $(t_0,x_0)$, we have
 \beas
     \min  \neg  \Big\{ ( \vf \neg - \neg \ul{l})(t_0,x_0),   \max  \neg  \Big\{  \neg  - \dneg \frac{\pa }{\pa t} \vf (t_0,x_0)
     \neg - \neg  H \big(t_0,x_0, \vf (t_0,x_0), D_x \vf (t_0,x_0), D^2_x \vf (t_0,x_0)\big),   ( \vf  \neg -  \neg  \ol{l})(t_0,x_0)  \Big\} \Big\}  \neg \le \neg   (\hb{resp.}  \ge)  \;  0.
 \eeas

  \end{deff}

 Although the function $H$ in Definition \ref{def:viscosity_solution} may take $\pm \infty$ values,
 the left-hand-side of the inequality above  is between $(w-\ol{l})(t_0,x_0)$ and $(w-\ul{l})(t_0,x_0)$ and thus finite.

  \ss    For any $(t,x,y,z,\G, u,v ) \neg \in \neg  [0,T]  \neg \times \neg  \hR^k  \neg \times \neg  \hR
   \neg \times  \neg  \hR^d  \neg \times \neg \hS_k \neg \times \neg  \hU_0  \neg \times \neg  \hV_0   $,
   we set 
      \beas
    H  (t,x,y,z,\G,u,v) \dfnn \frac12 trace\big(\si \si^T(t,x,u,v) \,\G\big)+ z \cd b(t,x,u,v)
    +  f\big(t,x, y, z \cd \si(t,x,u,v), u,v \big)
    \eeas
 and  consider   the following Hamiltonian functions:
  \beas
    \ul{H}_1(\Xi)   & \dfnn &         \underset{u \in \hU_0}{\sup} ~\;
 \linf{ \Xi'  \to \Xi }  ~\;   \underset{v \in \hV_0}{\inf}  \;  H (\Xi',u,v) , \q
 \ol{H}_1(\Xi)  \dfnn  \lmtd{n \to \infty} \,    \underset{u \in \hU_0}{\sup} \;
 \underset{v \in \sO^n_u}{\inf}  ~\, \lsup{\hU_0 \ni u' \to u}\; \underset{\Xi'     \in O_{1/n} (\Xi)}{\sup} \;  H (\Xi',u',v) ,  \\
 \hb{and}  \q  \ol{H}_2(\Xi)  & \dfnn &     \underset{v \in \hV_0}{\inf}  ~\;
 \lsup{ \Xi'  \to \Xi }   ~\;   \underset{u \in \hU_0}{\sup} \;  H (\Xi',u,v)  ,  \q
\ul{H}_2(\Xi)   \dfnn      \lmtu{n \to \infty} \,    \underset{v \in \hV_0}{\inf} \, \underset{u \in \sO^n_v}{\sup} ~\,
 \linf{\hU_0 \ni u' \to u}\; \underset{ \Xi'     \in O_{1/n} (\Xi)}{\inf} \;  H (\Xi',u',v),
    \eeas
 where $\Xi = (t,x,y,z,\G)$, $ \sO^n_u \dfnn \big\{ v \in \hV_0: [v]_{\overset{}{\hV}}
  \le n + n [u]_{\overset{}{\hU}} \big\}$ and
    $ \sO^n_v \dfnn \big\{ u \in \hU_0: [u]_{\overset{}{\hU}} \le n + n [v]_{\overset{}{\hV}} \big\}$.

 \ss  For any $ (t,x) \in [0,T] \times \hR^k $,     Proposition \ref{prop_w_conti} implies that
 \beas
           w^*_1 (t,x )    \dfnn     \lsup{ t' \to t}  w_1  \big(t',x\big)
    =    \lsup{ (t',x') \to (t,x) }  w_1  (t',x')
   \q \hb{and}  \q  w^*_2 (t,x )    \dfnn     \linf{ t' \to t}  w_2 \, \big(t',x\big)
    =    \linf{ (t',x') \to (t,x) }  w_2   (t',x')   .
 \eeas
  In fact,     $w^*_1$ is the smallest upper semicontinuous function above $w_1$
 (also known as the upper  semicontinuous envelope of $w_1$), while $w^*_2$  is the largest lower semicontinuous function
 below  $w_2$    (also known as the lower  semicontinuous envelope of $w_2$).
   Similarly, for $i=1,2$,
\beas
 \ul{w}_i (t,x ) \dfnn \linf{ t' \to t}  \wh{w}_i \, \big(t',x\big)   \q \hb{and} \q
 \ol{w}_i (t,x ) \dfnn \lsup{ t' \to t}  \wh{w}_i \, \big(t',x\big)  , \q    \fa    (t,x) \in [0,T] \times \hR^k
\eeas
are    the lower and upper semicontinuous envelopes of $\wh{w}_i$  respectively.

   Given   $x \in \hR^k$, though $w_i(T,x) = \wh{w}_i(T,x)  = h(x) $, probably neither
 of $w^*_i (x)$, $\ul{w}_i (x)$, $\ol{w}_i (x)$ equals to $h(x)$ as
  the  value functions $w_i$, $\wh{w}_i$ may not be continuous in $t$.

\begin{thm} \label{thm_viscosity}
 ~

\no  1\) If $\hU_0$ \(resp.\;$\hV_0$\) is a countable union of closed subsets of $\hU$ \(resp.\;$\hV$\),
  then  $\ol{w}_1 $ and $w^*_1 $ \(resp.\;$ \ul{w}_2 $ and $w^*_2 $\)
   are  two  viscosity subsolutions \(resp.\;supersolutions\)
  of   \eqref{eq:PDE} with the fully nonlinear Hamiltonian $\ol{H}_1$ \(resp.\;$\ul{H}_2$\).

\no  2\) On the other hand,  if  \($\bV_\l$\) \big(resp.\;\($\bU_\l$\)\big) holds for some $\l \in (0,1)$,
    then  $\ul{w}_1 $ \(resp.\;$ \ol{w}_2 $\)  is a  viscosity supersolution
 \(resp.\;subsolution\)  of   \eqref{eq:PDE} with the fully nonlinear
 Hamiltonian $ \ul{H}_1 $ \(resp.\;$ \ol{H}_2 $\).

\end{thm}

 \section{Shifted Processes}

 \label{sec:shift_prob}

In this section, we fix   $0 \le t \le s   \le  T$ and explore  properties of shifted processes
 from $\O^t$ to $\O^s$, which are necessary  for Section \ref{sec:zs_drgame} and Section \ref{sec:PDE}.

 \subsection{Concatenation   of Sample Paths}

   We  concatenate  an $\o \in \O^t$
 and an $ \wt{\o} \in \O^s$ at time $s$ by:
 \bea   \label{def_concatenation}
 \big(\o \otimes_s  \wt{\o}\big)(r)  \dfnn   \o(r) \, \b1_{\{r \in [t,s)\}}   + \big(\o(s) + \wt{\o}(r) \big) \, \b1_{\{r \in [s,T]\}} , \q \fa  r \in [t,T] ,
 \eea
  which is   still  of $\O^t$.
 Clearly, this   concatenation    is an associative operation: i.e., for any $r \in [s, T]$ and $\wh{\o}   \in   \O^r$
  \beas  
     (\o \otimes_s  \wt{\o})\otimes_r  \wh{\o} =  \o \otimes_s ( \wt{\o} \otimes_r  \wh{\o} ) .
 \eeas

  Given  $\o \in \O^t$.     we set $\o \otimes_s \es  =\es $ and $  \o \otimes_s \wt{A} \dfnn
   \big\{ \o  \otimes_s \wt{\o}: \wt{\o} \in  \wt{A} \big\}$ for any non-empty     $\wt{A} \subset   \O^s$.
  The next result shows that   $A \in \cF^t_s$  consists of all branches $\o \otimes_s \O^s $ with $\o \in A$.

\begin{lemm}  \label{lem_element}
  Let  $ A \in \cF^t_s$.  
  If $\o \in A  $, then   $  \o \otimes_s \O^s   \subset A  $ \(i.e.  $A^{s,\o}=\O^s$\).
  Otherwise,  if $\o \notin A  $, then $   \o \otimes_s \O^s   \subset    A^c $
  \(i.e. $A^{s,\o}=  \es $\).
  \end{lemm}

  Also,  for  any   $A \subset \O^t$ we set $A^{s, \o} \dfnn   \{ \wt{\o} \in \O^s: \o \otimes_s \wt{\o} \in A  \} $
   as the  projection of $A$  on $\O^s $ along $\o$. In particular, $\es^{s,\o} = \es$.
   For any $A \subset \wt{A} \subset \O^t$ and any collection $\{A_i \}_{i \in \cI}$ of subsets of $\O^t$,
   One  can deduce that
    \bea
     \big(A^c\big)^{s,\o}
    & \tneg =& \tneg  \{ \wt{\o} \in \O^s: \o \otimes_s \wt{\o} \in A^c  \}=
  \O^s \big\backslash \{ \wt{\o} \in \O^s: \o \otimes_s \wt{\o} \in A  \}
  =  \O^s \big\backslash A^{s,\o}    = \big( A^{s,\o} \big)^c   , \hspace{1cm}  \label{lem_basic_complement} \\
         A^{s,\o}
  & \tneg =& \tneg \{ \wt{\o} \in \O^s: \o \otimes_s \wt{\o} \in A   \} \subset
  \big\{ \wt{\o} \in \O^s: \o \otimes_s \wt{\o} \in \wt{A}   \big\}
  =  \wt{A}^{s,\o}         , \label{lem_basic_subset} \\
        \hb{ and }  \;    \Big( \underset{i \in \cI}{\cup} A_i \Big)^{s,\o}
       &  \tneg =& \tneg  \Big\{ \wt{\o} \in \O^s: \o \otimes_s \wt{\o} \in \underset{i \in \cI}{\cup} A_i  \Big\}=
   \underset{i \in \cI}{\cup} \big\{ \wt{\o} \in \O^s: \o \otimes_s \wt{\o} \in  A_i  \big\}
 =  \underset{i \in \cI}{\cup} \;  A^{s,\o}_i      .      \label{lem_basic_union}
 \eea


 \begin{lemm} \label{lem_concatenation}
 Let $\o \in \O^t$.   For any   open \(resp.\;closed\) subset $A$ of $\O^t$, $A^{s,\o}$ is
  an open \(resp.\;closed\) subset   of $\O^s$.
 Moreover, given    $r \neg \in  \neg  [s,T]$.  we have  $A^{s, \o}  \neg  \in   \neg    \cF^s_r$
 for any $A  \neg \in \neg  \cF^t_r$ and        $\o \otimes_s \wt{A} \in \cF^t_r$ for any  $ \wt{A} \in \cF^s_r$.

 \end{lemm}

    For any $\cD \neg \subset  \neg  [t,T]    \times    \O^t$,
  we accordingly set $\cD^{s, \o}  \neg \dfnn  \neg  \big\{\big(r, \wt{\o}\big)  \neg \in \neg  [s, T]  \neg \times \neg  \O^s: \big(r, \o \otimes_s \wt{\o}\big)  \neg \in \neg  \cD \big\}$.
    Similar to   \eqref{lem_basic_complement}-\eqref{lem_basic_union},
   for any $\cD \subset \wt{\cD} \subset [t,T] \times \O^t$ and any collection $\{\cD_i \}_{i \in \cI}$ of subsets of $[t,T] \times \O^t$, one has
   \bea
   ~\;  \big( ([t,T] \times \O^t) \big\backslash \cD\big)^{s,\o}
     = ( [s,T] \times \O^s ) \big\backslash \cD^{s,\o} = \big( \cD^{s,\o} \big)^c   ,  ~\,
         \cD^{s,\o} \neg \subset  \wt{\cD}^{s,\o}  ~\,
  \hb{and}  ~\, \Big( \underset{i \in \cI}{\cup} \cD_i \Big)^{s,\o} \neg =  \underset{i \in \cI}{\cup}  \; \cD^{s,\o}_i  . \label{lem_basic2}
 \eea

  \subsection{Measurability of Shifted Processes}
  \label{subsection:shift_measurability}


   For any $\hM-$valued random variable $\xi$ on $\O^t$, we define a shifted random variable $\xi^{s,\o}$ on $\O^s$ by  $ \xi^{s,\o}(\wt{\o}) \dfnn \xi ( \o \otimes_s  \wt{\o} ) $,  $  \fa  \wt{\o} \in \O^s $.
 And for any $\hM-$valued process $X = \{X_r\}_{r \in [t,T]}$, its corresponding shifted process
 with respect to $s$ and $\o$  consists of   $ X^{s,\o}_r = (X_r)^{s,\o} $, $ \fa  r \in [s,T]  $.
  In light of Lemma \ref{lem_concatenation},   shifted random variables and shifted processes ``inherit" measurability
  in the following way:

 \begin{prop}  \label{prop_shift1}
  If  $\xi $ is   $\cF^t_r-$measurable for some $r \in [s,T]$,  then   $\xi^{s,\o} $ is $  \cF^s_r-$measurable.
   Moreover, for any $\hM-$valued,  $\bF^t-$adapted   process $  \{X_r \}_{r \in [t, T]}$,
the shifted  process $   \big\{X^{s,\o}_r \big\}_{  r \in [s,T]}$ is $\bF^s-$adapted.
 \end{prop}

\begin{prop} \label{prop_shift2}
 Given $T_0 \in [s,T] $, $ \cD^{s,\o} \neg \in  \neg  \sB\big([s,T_0]\big)     \otimes \cF^s_{T_0}$
 for any $ \cD  \neg \in \neg  \sB([t,T_0])    \otimes \cF^t_{T_0}$.
   Consequently, if $  \{X_r \}_{r \in [t, T]}$  is an   $\hM-$valued,   measurable process
   on $ \big(\O^t,\cF^t_T\big) $ \big(resp. an $\hM-$valued, $\bF^t-$progressively measurable process\big),
   then    the shifted  process $ \big\{X^{s,\o}_r \big\}_{  r \in [s,T]}$
  is a   measurable process on $ \big(\O^{s},\cF^s_T\big)$
  \big(resp. an $\bF^s-$progressively measurable process\big).
   Moreover, we have $ \cD^{s,\o} \in \sP\big(\bF^s\big)$ for any $\cD \in \sP\big(\bF^t\big)$.

\end{prop}


  \ss  For any $\cJ  \neg \subset  \neg  [t,T]  \neg \times  \neg  \O^t    \times     \hM$,
  we set $\cJ^{s, \o}  \neg \dfnn   \neg  \big\{\big(r, \wt{\o},x\big)  \neg \in  \neg  [s, T] \times \O^s \times \hM: \big(r, \o \otimes_s \wt{\o},x\big)  \neg \in  \neg  \cJ \big\} $.

 \begin{cor} \label{cor_shift2}

      For any $ \cJ    \in    \sP\big(\bF^t\big) \otimes \sB(\hM)$, $ \cJ^{s,\o}
    \in \sP\big(\bF^s\big) \otimes \sB(\hM)$.
   Let $\wt{\hM}$ be another generic metric space.  If  a function $f :  [t,T] \times \O^t \times \hM \rightarrow  \wt{\hM}  $
   is    $\sP\big(\bF^t\big) \otimes \sB(\hM)/\sB\big(\wt{\hM} \big)-$measurable, then
  the function $f^{s,\o}\big(r, \wt{\o},x\big) \dfnn f\big(r, \o \otimes_s \wt{\o},x\big)$,
   $\fa \big(r, \wt{\o},x\big) \in [s, T] \times \O^s \times \hM$ is $ \sP\big(\bF^s\big) \otimes \sB(\hM)/\sB\big(\wt{\hM} \big)-$measurable.
       \end{cor}


 \ss    When $s=\t(\o)$ for some $\bF^t-$stopping time $\t$, we shall simplify the above notations by:
 \beas
      \o \otimes_\t \wt{\o}   =    \o \otimes_{\t(\o)} \wt{\o}  ,\q  A^{\t,\o}  =   A^{\t(\o),\,\o} ,
       \q  \cD^{\t,\o}  =   \cD^{\t(\o),\,\o}, \q
   \xi^{\t,\o}    =     \xi^{\t(\o),\,\o}    \q \hb{and} \q    X^{\t,\o} = X^{\t(\o),\,\o}  .
 \eeas

   \if{0}

\begin{cor}  \label{cor_element}
  Let  $ A \in \cF^t_\t$ for some $\bF^t-$stopping time $\t$.  
  If $\o \in A  $, then   $  \o \otimes_\t \O^{\t(\o)}   \subset A  $ \(i.e.  $A^{\t,\o}=\O^{\t(\o)}$\).
  Otherwise,  if $\o \notin A  $, then $   \o  \otimes_\t \O^{\t(\o)}    \subset    A^c $  \(i.e. $A^{\t,\o}=  \es $\).
  \end{cor}

  \ss \no {\bf Proof:}  Fix $\o \in A$ and set $s = \t(\o)$. Since $\o \in A \cap \{\t \le s\} \in \cF^t_s$, Lemma
  \ref{lem_element} shows that
     $  \o \otimes_\t \O^{\t(\o)} =  \o \otimes_s \O^s  \subset  A \cap \{\t \le s\} \subset A $. \qed

   \fi

   The following lemma shows that given an $\bF^t-$stopping time $\t$, an $\cF^t_\t-$measurable random variable only depends on what happens before $\t$:

  \begin{lemm}   \label{lem_bundle}
 For any $\bF^t-$stopping time $\t$ and $\xi \in \cF^t_\t$,
    $\xi^{\t,\o} \equiv \xi (\o ) $. In particular, $\t\big(\o \otimes_\t \O^{\t(\o)}\big)=\t(\o)$.

     \end{lemm}

 \if{0}

 For another  $\bF^t-$stopping time $\z$, we define $    \z^\o_\t (\wt{\o}) \dfnn
   \t(\o) \vee \z \big( \o \otimes_\t  \wt{\o} \big)$, $\fa \wt{\o} \in \O^{\t(\o)}$. The following Lemma extends Lemma
   \ref{lem_concatenation} for stopping times.

 \begin{lemm}    \label{lem_concatenation_tau}
 Let $\t$ and $\z$ be two $\bF^t-$stopping times.  For any $\o \in \O^t$, $\z^\o_\t$ is an $\bF^{\t(\o)}-$stopping time.    For any $A \in \cF^t_{\t(\o) \vee \z}$,   $A^{\t, \o}  \in     \cF^{\t(\o)}_{\z^\o_\t}$;
        for any  $ \wt{A} \in \cF^{\t(\o)}_{\z^\o_\t}$,   $\o \otimes_\t \wt{A} \in \cF^t_{\t(\o) \vee \z}$.

     \end{lemm}

  \ss \no {\bf Proof:} For any $r \in \big[\t(\o),T\big] $, we can deduce that
 \beas
  \big\{  \z^\o_\t \le r\big\}
    = \big\{\wt{\o} \in \O^{\t(\o)}:   \z \big( \o \otimes_\t  \wt{\o} \big) \le r\big\}
       = \big\{\wt{\o} \in \O^{\t(\o)}:  \o \otimes_\t  \wt{\o} \in \{ \z  \le r \} \big\}= \{ \z  \le r \}^{\t,\o} .
 \eeas
 As $ \{ \z  \le r \} \in \cF^t_r$,   Lemma \ref{lem_concatenation} implies that $\big\{  \z^\o_\t \le r\big\}
 =   \{ \z  \le r \}^{\t,\o} \in \cF^{\t(\o)}_r$. Thus $ \z^\o_\t $ is an $\bF^{\t(\o)}-$stopping time.

 \ss       Let $A \in \cF^t_{\t(\o) \vee \z}$.  For any $r \in \big[\t(\o),T\big] $,
 since $ A  \cap\{ \t(\o) \vee \z \le r \} \in \cF^t_r $, Lemma \ref{lem_concatenation} implies  that
 \beas
  A^{\t, \o} \cap \big\{ \z^\o_\t   \le r \big\}
    &=& \big\{ \wt{\o} \in \O^{\t(\o)}: \o  \otimes_\t \wt{\o} \in A,   \z^\o_\t (\wt{\o} )   \le r  \big\}
  = \big\{ \wt{\o} \in \O^{\t(\o)}: \o  \otimes_\t \wt{\o} \in A \cap \{ \t(\o) \vee \z \le r \}   \big\} \\
 & = & \big( A \cap \{ \t(\o) \vee  \z \le r \}\big)^{\t,\o} \in  \cF^{\t(\o)}_r   .
 \eeas
 Hence,  $A^{\t, \o}  \in     \cF^{\t(\o)}_{\z^\o_\t}$.   On the other hand, let  $ \wt{A} \in \cF^{\t(\o)}_{\z^\o_\t}$.   For any $r \in \big[\t(\o),T\big] $,
 since $ \wt{A} \cap \{ \z^\o_\t   \le r \} \in \cF^{\t(\o)}_r$,  one can deduce from Lemma \ref{lem_concatenation}   that
     \beas
   \big( \o \otimes_\t \wt{A} \big) \cap \{\t(\o) \vee \z \le r \}
        & =& \big\{ \o  \otimes_\t \wt{\o}: \wt{\o} \in  \wt{A},\t(\o) \vee \z (\o  \otimes_\t \wt{\o}) \le r  \big\}
                    =    \big\{ \o  \otimes_\t \wt{\o}: \wt{\o} \in  \wt{A} \cap \{ \z^\o_\t   \le r \} \big\} \\
                   &  = &     \o  \otimes_\t  \big(  \wt{A} \cap \{ \z^\o_\t   \le r \} \big)\in \cF^t_r   .
       \eeas
        Thus  $\o \otimes_\t \wt{A} \in \cF^t_{\t(\o) \vee \z}$.\qed

   \fi

  \subsection{Integrability of Shifted Processes}
  \label{subsection:shift_integrability}

  In this subsection,  let  $\t$ be    an $\bF^t-$stopping time  with countably many values.
 Using the regular conditional probability distribution (see e.g. \cite{Stroock_Varadhan}),
 we show below that shifted random variables ``inherit"  integrability property.

\begin{prop} \label{prop_rcpd_L1}
  For any   $\xi \in \hL^1 \big( \cF^t_T \big)$, it holds for $P^t_0-$a.s.     $\o \in \O^t$
  that   $\xi^{\t,\o} \in \hL^1 \big( \cF^{\t(\o)}_T ,  P^{\t(\o)}_0   \big) $ and
  \bea   \label{eq:f475}
    E_{\t(\o)} \big[ \xi^{\t,\o} \big]= E_t \big[\xi\big| \cF^t_\t\big](\o) \in  \hR  ,
    \eea
    where 
    $ E_{\t(\o)}$ stands for $E_{P^{\t(\o)}_0}$.
 Consequently, for any $p \in [1,\infty)$ and  $\xi \in \hL^p \big( \cF^t_T \big)$, it holds for $P^t_0-$a.s.  $\o \in \O^t$  that $\xi^{\t,\o} \in \hL^p \big( \cF^{\t(\o) }_T, P^{\t(\o)}_0 \big) $.
\end{prop}

\begin{cor} \label{cor_rcpd11}
 For any   $P^t_0-$null set    $  \cN          $,
  it holds  for  $P^t_0-$a.s.  $\o  \in \O^t  $ that $\cN^{\t,\o}$ is a $P^{\t(\o)}_0-$null set.
 Consequently, for any two real-valued random variables $\xi_1$ and $ \xi_2$, if  $\xi_1 \le \xi_2$,   $P^t_0-$a.s.,
 then it holds for  $P^t_0-$a.s. ~  $\o  \in \O^t  $ that  $\xi^{\t,\o}_1 \le \xi^{\t,\o}_2$,   $ P^{\t(\o)}_0 -$a.s.
\end{cor}

  Next, let us extend   Proposition \ref{prop_rcpd_L1}  to   $\hE-$valued  measurable processes.

\begin{prop} \label{prop_rcpd3}
 Let $\{X_r\}_{r \in [t,T]}$ be an $\hE-$valued  measurable process on $\big(\O^t,\cF^t_T\big)$ such that
   $ E_t \Big[\big( \neg \int_\t^T \neg \big|X_r  \big|^p  dr\big)^{ \wh{p}/p}\Big]   <  \infty$ for some $p,\wh{p} \in [1,\infty)$.
   It holds for $P^t_0-$a.s. ~$\o \in \O^t$
      that $\{X^{\t,\o}_r\}_{r \in [\t(\o),T]}$ is a measurable process on $\big(\O^{\t(\o)},\cF^{\t(\o)}_T\big)$ with
   $E_{  \t(\o) }\Big[ \big( \neg \int_{\t(\o)}^T \neg \big| X^{\t,\o}_r \big|^p  dr \big)^{\wh{p}/p} \Big] < \infty$.

  \end{prop}

 \if{0}
 Next, let $X^{(i)}$ be the $i-$th component of $X$. Since the integral
 $   \int_t^T \b1_{\{r \ge \t(\o')\}} X^{(i)}_r(\o')  dr $ is also well-defined for any $\o' \in \cN^c$,
 we see that $  \xi_i \dfnn \b1_{ \cN^c} \int_t^T \b1_{\{r \ge \t \}}  X^{(i)}_r    dr
 \in \hL^1\big(\cF^t_T  \big)$. Similar to \eqref{eq:f479}, it holds any $\o \in \O^t$ and $\wt{\o} \in \O^{\t(\o)}$ that
 $ \xi^{\t,\o}_i (\wt{\o})
 = \b1_{\big( \cN^{\t,\o}\big)^c}  (\wt{\o}) \cd \int_{\t(\o)}^T   \big(X^{(i)}\big)^{\t,\o}_r  (\wt{\o}) \, dr $.
   Since $\cN^{\t,\o} \in \sN^{P^{\t(\o)}_0}$ for $P^t_0-$a.s.~  $\o \in \O^t$,  we
       can   deduce from    \eqref{eq:f475}  that for $P^t_0-$a.s. $\o \in \O^t$
\beas
  \hspace{1.8cm}    E_{ P^{\t(\o)}_0}\bigg[ \int_{\t(\o)}^T   \big(X^{(i)} \big)^{\t,\o}_r    dr \bigg]
 = E_{ P^{\t(\o)}_0} \big[    \xi^{\t,\o}_i \big]
=  E_{ P}\big[   \xi_i \big| \cF^t_\t\big](\o)
=  E_{ P}\bigg[  \int_\t^T     X^{(i)}_r    dr  \Big| \cF^t_\t\bigg](\o)  \in \hR .  \hspace{1.8cm} \hb{\qed}
\eeas
 \fi

 \begin{cor} \label{cor_rcpd3}
 Given $p,\wh{p} \in [1,\infty)$,  if  $\{X_r\}_{r \in [t,T]} \in  \hH^{p,\wh{p}}_{\bF^t}([t,T],\hE) $
 \big(resp.\;$ \hC^p_{\bF^t}([t,T], \hE)$\big), then it holds for $P^t_0-$a.s.~$\o \in \O^t$,
 $\{X^{\t,\o}_r\}_{r \in [\t(\o),T]} \in \hH^{p,\wh{p}}_{\bF^{\t(\o)}}\big([\t(\o),T],\hE,P^{\t(\o)}_0\big) $
 \Big(resp.\;$ \hC^p_{\bF^{\t(\o)}}\big([\t(\o),T], \hE , P^{\t(\o)}_0\big)$\Big).
         \end{cor}


  Similar to Corollary \ref{cor_rcpd11}, a shifted $dr \times dP^t_0-$null set still has zero product measure:

\begin{prop} \label{prop_rcpd5}
 For any $\cD  \in   \sB\big([t,T]  \big) \otimes \cF^t_T $ with $(dr \times dP^t_0) \big(\cD \cap \[\t,T\]\big) = 0 $,
  it holds  for    $P^t_0-$a.s. ~ $\o  \in \O^t  $ that $\cD^{\t,\o} \in \sB\big([\t(\o),T]\big) \otimes \cF^{\t(\o)}_T$
  with   $\big(dr \times d P^{\t(\o)}_0\big)\big(\cD^{\t,\o}\big)=0 $.

 \end{prop}

The following analyzes the admissibility of controls and strategies when they are shifted.

 \begin{prop} \label{prop_dwarf_strategy}
 \(1\) For any $\mu \in \cU^t$ \(resp.\;$\nu \in \cV^t$\),  it holds for $P^t_0-$a.s. $\o \in \O^t$ that  $\mu^{\t,\o} = \{\mu^{\t,\o}_r\}_{r \in [\t(\o),T]} \in \cU^{\t(\o)} $
 \big(resp.\;$\nu^{\t,\o} =  \{\nu^{\t,\o}_r\}_{r \in [\t(\o),T]} \in \cV^{\t(\o)}$\big).

 \ss \no \(2\)  For any    $\a \in  \cA^t$ \big(resp.\;$\a \in  \wh{\cA}^t$, $\beta \in  \fB^t$ and $\beta \in \wh{\fB}^t$\big),
    it holds for $P^t_0-$a.s. $\o \in \O^t$ that
   $\a^{s,\o}   \in  \cA^s $  \big(resp.\;$\a^{s,\o}   \in \wh{\cA}^s $, $\beta^{s,\o}   \in  \fB^s$ and $\beta^{s,\o}   \in \wh{\fB}^s$\big).
 \end{prop}

 \subsection{Shifted Stochastic Differential Equations}

 In this subsection, we still consider   an $\bF^t-$stopping time $\t$ with countably many values.

 \ss  Fix $x \in \hR^k$, $ \mu  \in \cU^t  $,   $ \nu  \in   \cV^t$ and  set $\Th=(t,x,\mu,\nu)$.
 For $P^t_0-$a.s. $\o \in \O^t$,  Proposition \ref{prop_dwarf_strategy} (1) shows  that
  $\big( \mu^{\t,\o} , \nu^{\t,\o} \big)  \in \\ \cU^{\t(\o)} \neg \times \neg  \cV^{\t(\o)}$, and thus
     we know from Section \ref{sec:zs_drgame}   that
 the following SDE  on the probability space $\Big(\O^{\t(\o)}, \ol{\cF}^{\t(\o)}_{\neg T}, P^{\t(\o)}_0 \Big)$:
    \bea  \label{eq:p201}
  X_s=   \wt{X}^{\Th}_{\t(\o)}  (\o) + \int_{\t(\o)}^s  b(r, X_r,\mu^{\t,\o}_r,\nu^{\t,\o}_r ) \, dr
   + \int_{\t(\o)}^s \si(r, X_r,\mu^{\t,\o}_r,\nu^{\t,\o}_r) \, dB^{\t(\o)}_r,  \q s \in [\t(\o), T]
    \eea
  admits a  unique  solution  $\Big\{ X^{ \Th^\o_\t }_s \Big\}_{s \in [\t(\o),T]} \neg  \in  \neg  \hC^2_{\ol{\bF}^{\t(\o)} } \big([\t(\o),T], \hR^k \big) $  with  $\Th^\o_\t  \neg \dfnn  \neg  \Big(\t(\o), \wt{X}^{\Th}_{\t(\o)}  (\o) ,\mu^{\t,\o},\nu^{\t,\o}\Big) $.  As shown below, 
   the $\bF^{\t(\o)}-$version of $\Big\{ X^{ \Th^\o_\t }_s \Big\}_{s \in [\t(\o),T]}$
   is exactly   the shifted process $ \big\{  \big( \wt{X}^\Th \big)^{\t,\o}_s \big\}_{s \in [\t(\o),T]} $.

 \ss

\begin{prop} \label{prop_FSDE_shift}
 It holds
   for $P^t_0-$a.s.~$\o \in \O^t$ that 
  $   \wt{X}^{\Th^\o_\t}_s =  \big(\wt{X}^\Th \big)^{\t,\o}_s   $,  
     $ \fa    s \in [\t(\o),T]$.

\end{prop}

 This result has appeared  in \cite{Fleming_1989} for case of compact control spaces (see the paragraph below (1.16) therein)
 and appeared in  Lemma 3.3 of \cite{Nutz_2011} where only one unbounded control is considered.
     The proof of Proposition \ref{prop_FSDE_shift}  depends on the following result
  about the convergence of shifted random variables in probability.

 \begin{lemm} \label{lem_shift_converge_proba}
  For any   $ \{\xi_i\}_{i \in \hN} \subset \hL^1(\cF^t_T)$ that converges to 0 in probability $P^t_0$,
 we can find a subsequence $ \big\{  \wh{\xi}_{\,i} \big\}_{i \in \hN}  $ of it such that for $P^t_0-$a.s. $\o \in \O^t$,
 $ \big\{ \wh{\xi}^{\,\t,\o}_{\,i} \big\}_{i \in \hN}  $   converges to 0 in probability $P^{\t(\o)}_0$.

 \end{lemm}

 For any $\cF^t_T-$measurable random variable $ \xi $ with $\ul{L}^\Th_T \le \xi \le \ol{L}^\Th_T $, $P^t_0-$a.s.,
   Proposition \ref{prop_shift1},   Corollary \ref{cor_rcpd11} and Proposition \ref{prop_FSDE_shift} imply that
    for $P^t_0-$a.s. $\o \in \O^t$, $ \xi^{\t,\o} \in \cF^{\t(\o)}_T $ and
     \beas
   \q  \ul{L}^{\Th^\o_\t}_T =  \ul{l} \Big(T , \wt{X}^{\Th^\o_\t}_T   \Big)
    =  \ul{l}  \Big(T , \big(\wt{X}^\Th_T \big)^{\t,\o}  \Big)  \le \big(\ul{L}^\Th_T\big)^{\t,\o} \le \xi^{\t,\o} \le  \big(\ol{L}^\Th_T\big)^{\t,\o}
       \le \ol{l}  \Big(T , \big(\wt{X}^\Th_T \big)^{\t,\o}  \Big)  = \ol{l} \Big(T , \wt{X}^{\Th^\o_\t}_T   \Big)
       = \ol{L}^{\Th^\o_\t}_T .
     \eeas
   Then Section    \ref{sec:zs_drgame} also shows that for $P^t_0-$a.s. $\o \in \O^t$,
    the  DRBSDE$\Big(P^{\t(\o)}_0, \xi^{\t,\o}, f^{\Th^\o_\t}_T, \ul{L}^{\Th^\o_\t}, \ol{L}^{\Th^\o_\t} \Big)$
    on the probability space $\Big(\O^{\t(\o)}, \ol{\cF}^{\t(\o)}_{\neg T}, P^{\t(\o)}_0 \Big)$
    admits a unique solution    $ \Big(Y^{\Th^\o_\t} \big(T, \xi^{\t,\o}\big),Z^{\Th^\o_\t} \big(T, \xi^{\t,\o}\big), \ul{K}^{\Th^\o_\t} \big(T, \xi^{\t,\o}\big), \ol{K}^{\, \raisebox{-0.5ex}{\scriptsize $\Th^\o_\t$}} \big(T, \xi^{\t,\o}\big) \Big) \in \hG^q_{\ol{\bF}^{\t(\o)} } \big([\t(\o),T]\big) $.
    Similar to Proposition \ref{prop_FSDE_shift}, the  $\bF^{\t(\o)}-$version of $Y^{\Th^\o_\t} \big(T, \xi^{\t,\o}\big)$
    coincides with the shifted process $ \Big\{  \big( \wt{Y}^\Th (T,\xi) \big)^{\t,\o}_s \Big\}_{s \in [\t(\o),T]} $.

  \begin{prop} \label{prop_DRBSDE_shift}
 It holds  for $P^t_0-$a.s.~$\o \in \O^t$ that 
  $  \wt{Y}^{\Th^\o_\t}_s \big(T, \xi^{\t,\o}\big) = \big( \wt{Y}^\Th (T,\xi) \big)^{\t,\o}_s  $,
     $ \fa    s \in [\t(\o),T]$. In particular, $ \wt{Y}^{\Th^\o_\t}_{\t(\o)} \big(T, \xi^{\t,\o}\big) = \big( \wt{Y}^\Th_\t (T,\xi) \big)  (\o) $ for $P^t_0-$a.s.~$\o \in \O^t$.

\end{prop}

 \ss  Proposition \ref{prop_DRBSDE_shift} can also be shown by  {\it Picard} iteration, see
 (4.15) of \cite{STZ_2011b} for a BSDE version.

  \subsection{Pasting of Controls and Strategies}

   We define $ \wh{\Pi}_{t,s} (r,\o) \dfnn \big(r, \Pi_{t,s}(\o)\big) $,
   $\fa (r,\o) \in [s,T] \times \O^t$.   Analogous to Lemma \ref{lem_shift_inverse}, one has


  \begin{lemm}  \label{lem_shift_inverse2}
     Let $r \in [s,T]$. For any
 $\cD \in \sB([s,r])\otimes \cF^s_r$, $ \wh{\Pi}^{\,-1}_{\,t,s}(\cD)
  \in \sB([s,r])\otimes \cF^t_r  $ and $    \big( dr \times dP^t_0 \big) \big( \wh{\Pi}^{-1}_{t,s} (\cD)\big)
    = \big( dr \times dP^s_0 \big) (\cD) $. Consequently, the mapping $ \wh{\Pi}_{t,s}:
 [s,T] \times \O^t \to [s,T] \times \O^s$ is $ \sP_s(\bF^t)/\sP(\bF^s)-$measurable, where
$\sP_s(\bF^t) \dfnn  \{\cD \in \sP(\bF^t): \cD \subset  [s,T] \times \O^t \} 
$ is a $\si-$field of $ [s,T] \times \O^t $.

   \end{lemm}

  \if{0}

   \begin{cor} \label{cor_concatenation_inverse}
    For any  $  P^s_0 - $null set $A$, \;$  \Pi^{-1}_{t,s} (A) $ is a $  P^t_0 - $null set;
 and for any  $dr \times dP^s_0 - $null set $\cD$, \;$ \wh{\Pi}^{-1}_{t,s} (\cD) $ is a $dr \times dP^t_0 - $null set.

   \end{cor}

  \fi


Now, we are ready to   discuss   pasting of controls and strategies.

  \begin{prop}  \label{prop_paste_control}
   Let $\mu \in \cU^t $ for some $t \in [0,T]$ and let $\t$ be an $\bF^t-$stopping time taking values in a countable subset $\{t_n\}_{n \in \hN}$ of \,$[t,T]$. Given $N \in \hN$,
  let  $\{ A^n_i \}^{\ell_n}_{i=1} \neg \subset \neg \cF^t_{t_n} \neg $ be  disjoint subsets of $\{\t = t_n\}$
  for $n =1,\cds \neg ,  N$ and set $A_0 \dfnn  \O^t \Big\backslash \Big( \underset{n=1}{\overset{N}{\cup}} \overset{\ell_n}{\underset{i=1}{\cup}} A^n_i \Big) $. Then
   for any    $\{ \mu^n_i \}^{\ell_n}_{i=1}   \subset \cU^{t_n} $, $n =1,\cds \neg ,  N$
     \bea       \label{paste_control}
  \wh{\mu}_r(\o)   \dfnn   \begin{cases}
    \big(\mu^n_i\big)_r\big(\Pi_{t,t_n} (\o)\big),\q & \hb{if $(r,\o ) \in
   \[\t,T\]_{A^n_i} = [t_n,T] \times  A^n_i$    for $ n=1 \cds \neg , N$ and
    $ i=1,\cds \neg , \ell_n $} ,   \\
     \mu_r(\o),  & \hb{if $(r,\o) \in \[t,\t\[ \, \cup \, \[\t,T\]_{ A_0 }$     }
    \end{cases}
    \eea
    defines a $\cU^t-$control such that  for any $(r,\o) \in      \[ \t , T \]$
     \beas
  \wh{\mu}^{\t,\o}_r = \begin{cases}
    \big(\mu^n_i\big)_r,\q & \hb{if $ \o   \in    A^n_i$    for $ n=1 \cds \neg , N$ and
    $ i=1,\cds \neg , \ell_n $} ,   \\
     \mu^{\t,\o}_r ,  & \hb{if $ \o   \in A_0 $.   }
    \end{cases}
    \eeas
      We can paste  $\{ \nu^n_i \}^{\ell_n}_{i=1}   \subset \cV^{t_n} $, $n =1,\cds,N$ to a $\,\nu \in \cV^t $ with respect to
      $\{ A^n_i \}^{\ell_n}_{i=1} $, $n =1,\cds, N$ in the same manner.
       \end{prop}

  \begin{prop}  \label{prop_paste_strategy}
   Let $\a \in \cA^t $ \big(resp.\;$\wh{\cA}^t$\big) for some $t \in [0,T]$ and let $\t$ be an $\bF^t-$stopping time taking values in a countable subset $\{t_n\}_{n \in \hN}$ of \,$ \hQ_{t,T} $. Given $N \in \hN$,
  let  $\{ A^n_i \}^{\ell_n}_{i=1} \neg \subset \neg \cF^t_{t_n} \neg $ be  disjoint subsets of $\{\t = t_n\}$
  for $n =1,\cds \neg ,  N$ and set $A_0 \dfnn  \O^t \Big\backslash \Big( \underset{n=1}{\overset{N}{\cup}} \overset{\ell_n}{\underset{i=1}{\cup}} A^n_i \Big) $. Then
   for any    $\{ \a^n_i \}^{\ell_n}_{i=1}   \subset \cA^{t_n} $ \big(resp.\;$\wh{\cA}^{\, t_n}$\big), $n =1,\cds \neg ,  N$
     \bea       \label{paste_strategy}
  \wh{\a} (r,\o,v)  \neg  \dfnn  \neg    \begin{cases}
     \a^n_i \big( r, \Pi_{t,t_n}(\o),v \big), \dneg  & \hb{if $(r,\o ) \neg \in \neg
   \[\t,T\]_{A^n_i}  \neg = \neg  [t_n,T]  \neg \times \neg   A^n_i$    for $ n  \neg = \neg  1 \cds \neg , N$ and
    $ i \neg = \neg 1,\cds \neg , \ell_n $} ,   \\
     \a(r,\o,v),  & \hb{if $(r,\o) \in \[t,\t\[ \, \cup \, \[\t,T\]_{ A_0 }$, }
    \end{cases}     ~    \fa v \in \hV
    \eea
      is an $\cA^t-$strategy \big(resp.\;$\wh{\cA}^t-$strategy\big)   such that
      given $\nu \in \cV^t$, it holds for $P^t_0-$a.s. $\o \in \O^t$ that for any $r \in [\t(\o),T]$
     \bea   \label{eq:r401}
  \big( \wh{\a}\lan \nu \ran \big)^{\t,\o}_r = \begin{cases}
     \big( \a^n_i \lan \nu^{\, t_n,\o} \ran \big)_r,   ~    & \hb{if $ \o   \in    A^n_i$    for $ n=1 \cds \neg , N$ and
    $ i=1,\cds \neg , \ell_n $} ,   \\
     \big(  \a \lan \nu \ran \big)^{\t,\o}_r  ,  & \hb{if $ \o   \in A_0 $.      }
    \end{cases}
    \eea
      We can paste  $\{ \beta^n_i \}^{\ell_n}_{i=1}   \subset \fB^{t_n} $ \big(resp.\;$\wh{\fB}^{\, t_n}$\big), $n =1,\cds,N$ to a $\beta \in \fB^t $ \big(resp.\;$\wh{\fB}^t$\big) with respect to
      $\{ A^n_i \}^{\ell_n}_{i=1} $, $n =1,\cds, N$ in the same manner.
       \end{prop}

 \section{ Optimization Problems with Square-Integrable Controls
}

 \label{sec:one_control}

 In this section, we will remove $v-$controls (or take $\hV = \hV_0 = \{v_0\}$) so that
 the zero-sum stochastic differential game discussed above degenerates as an one-control optimization problem for player I.

   \subsection{General Results}

  We will follow the setting of Section \ref{sec:zs_drgame} except that we  take away
     the $v-$controls from all notations and definitions. In particular,  
  $\cV^t$, $\fB^t$ and $\wh{\fB}^t$ disappear (or become singletons) while  $\cA^t$ is equivalent  to   $\cU^t $.
  Then for any $(t,x) \in  [0,T] \times  \hR^k $, $w_1(t,x)$, $\wh{w}_1(t,x)$,   $w_2(t,x)$  coincide as
  \bea  \label{def_w_fcn}
   w  (t,x) \dfnn   \underset{\mu \in \cU^t}{\sup} \; \wt{Y}^{t,x,\mu  }_t   \Big(T,h \big(\wt{X}^{t,x,\mu }_T \big)  \Big) .
  \eea

   \if{0}

   \fi

       As $(\bV_\l)$ trivially  holds   for any $\l \in (0,1)$, the one-control version of Theorem \ref{thm_DPP}  reads as:

 \begin{prop}   \label{prop_1control_DPP}
 Let    $(t,x) \neg \in \neg [0,T]   \times   \hR^k $. For any family $\{\t_{\mu } \neg : \mu \neg \in \neg \cU^t\}$
 of $\hQ_{t,T}  -$valued, $\bF^t-$stopping times
  \bea \label{eq:xxc011}
  w (t,x)& = &      \underset{\mu \in \cU^t}{\sup}    \;
    \wt{Y}^{t,x,\mu }_t \Big(\t_{\mu },
    w \big(\t_{\mu },  \wt{X}^{t,x,\mu }_{\t_{\mu }} \big) \Big)  .
    \eea
\end{prop}

  Moreover, for  $\Xi = (t,x,y,z,\G) \in [0,T] \times \hR^k \times \hR \times \hR^d \times \hS_k $,
     $\ul{H} (\Xi) $ and $\ol{H} (\Xi)$ simplify respectively as $\ul{H}(\Xi) \dfnn \underset{u \in \hU_0}{\sup}~ \linf{\Xi'     \to \Xi } \;  H (\Xi',u )$ and
 \beas
 \ol{H} (\Xi) &=& \lmtd{n \to \infty} \,    \underset{u \in \hU_0}{\sup}
   ~\, \lsup{\hU_0 \ni u' \to u}\; \underset{\Xi'     \in O_{1/n} (\Xi)}{\sup} \;  H (\Xi',u' )
  = \lmtd{n \to \infty} \,    \underset{u \in \hU_0}{\sup}
   ~\,   \underset{\Xi'     \in O_{1/n} (\Xi)}{\sup} \;  H (\Xi',u ) \\
 &=&   \lmtd{n \to \infty}
   ~  \underset{\Xi'     \in O_{1/n} (\Xi)}{\sup} ~\, \underset{u \in \hU_0}{\sup}\;  H (\Xi',u )
 = \lsup{\Xi'     \to \Xi } ~ \underset{u \in \hU_0}{\sup}\;  H (\Xi',u )   ,
 \eeas
  where we used  the fact that  $\underset{u \in \hU_0}{\sup}
   ~\, \lsup{\hU_0 \ni u' \to u} = \underset{u \in \hU_0}{\sup}    $ in the second equality.
  Then we have the following   one-control version of     Theorem \ref{thm_viscosity}:

\begin{prop}   \label{prop_1control_viscosity}

     The lower semicontinuous envelopes of $w$:  $ \ul{w}  (t,x ) \dfnn \linf{ t' \to t}  w \big(t',x\big) $,
 $(t,x) \in [0,T] \times \hR^k$   is a  viscosity supersolution
    of   \eqref{eq:PDE} with the fully nonlinear
 Hamiltonian $ \ul{H}  $\,. On the other hand, if $\hU_0$ is a countable union of closed subsets of $\hU$,
  then the upper semicontinuous envelopes of $w$: $\ol{w}  (t,x ) \dfnn \lsup{ t' \to t}  w \big(t',x\big)  $, $(t,x) \in [0,T] \times \hR^k$
  is a  viscosity subsolutions
  of   \eqref{eq:PDE} with the fully nonlinear Hamiltonian $\ol{H} $.

 \end{prop}

 \begin{rem} \label{rem:semiconti}
  Similar to  \cite{STZ_2011c},  
  we only need to assume the measurability \(resp.\; lower  semi-continuity\)
 of the terminal function $h$
 for the ``\,$\le$" \(resp.\;``\,$\ge$"\) inequality of \eqref{eq:xxc011} and thus
 for the viscosity subsolution \(resp.\;supersolution\) part of
 Proposition \ref{prop_1control_viscosity}.
 \end{rem}

   \subsection{Connection to the \emph{Second-Order} Doubly Reflected  BSDEs}   \label{subsection:DR2BSDE}

 Now  let us further take $k=d$,  $\hU = \hS_d $, $u_0 = 0$
 and $\hU_0 = \hS^{>0}_d \dfnn \{\G \in \hS_d: det( \G ) > 0  \}$.  

 \begin{lemm}  \label{lem_hSd_property}
     $\hS_d$ is a separable normed vector space on which the determinant $det(\cd)$ is continuous.
 \end{lemm}

  It follows   that $\hU_0 = \hS^{>0}_d $ consists of closed subsets of $ \hU   $:
  $F_n \dfnn \{ u \in \hU:  det(u) \ge 1/n \}$, $ n \in \hN$.
  We also  specify:
    \bea   \label{eq:xd133}
b(t,x,u ) =b(t,x) \q \hb{and} \q \si(t,x,u)= u  , \q \fa (t,x,u) \in [0,T] \times \hR^d \times \hU.
    \eea
for a function $b(t,x):  [0,T] \times \hR^d \to \hR^d$ that is $\sB([0,T]) \otimes \sB(\hR^d) / \sB(\hR^d)-$measurable
 and Lipschitz continuous in $x$ with coefficient $\g > 0$.
 Via the transformation \eqref{eq:a155}, we will show that the value function $w$ defined in \eqref{def_w_fcn} is
  the value function of second-order doubly reflected  BSDEs.


  \ss   Given \,$t \in [0,T]$,  we say that a  $P \in \cP^t$  is a  {\it semi-martingale measure}
  if   $B^t$ is a  continuous  semi-martingale with respect to $(\bF^t, P)$.  Let $\cQ^t$ be the collection
  of all  semi-martingale measures on $(\O^t, \cF^t_T )$.

    \begin{lemm}  \label{lem_pqv}
  Let \,$t \in [0,T]$.  For $i,j \in \{1,\cds,d\}$, there exists an $\hR \cup \{\infty\}-$valued,
  $\bF^t-$progressively measurable process  $\hat{a}^{t,i,j}$
  such that for any $P \in \cQ^t$, it holds \pas ~   that
   \bea\label{eq:a075}
     \hat{a}^{t,i,j}_s = \hat{a}^{t,j,i}_s  = \lsup{m \to \infty}  \,  m      \Big(      \lan B^{t,i},  B^{t,j}  \ran^P_s
  - \lan B^{t,i},  B^{t,j}  \ran^P_{(s-1\neg /m)^{\neg {}^+} }  \Big)
     ,\q  s  \in [t,T],
  \eea
  where $\lan B^{t,i}, B^{t,j} \ran^P$'s denote  the $P-$cross variance between the $i-$th and $j-$th components of $B^t $.
\end{lemm}

  \ms  Similar to \cite{STZ_2011a}, we  let  $\cQ^t_W$  collect
  all   $P \in \cQ^t$ such that \pas  
 \bea    \label{eq:xd011}
          \lan B^t  \ran^P_s   \hb{  is absolutely continuous  in $ s $
 \; and \; $\hat{a}^t_s  \in \hS^{>0}_d    
 $ for a.e. $s \in [t,T]$.}
 \eea
  In general, two different probabilities $P_1, P_2 $ of $  \cQ^t_W$    are mutually singular, see  Example 2.1 of \cite{STZ_2011a}.
    \begin{lemm} \label{lem_hat_a_half}
  For any \,$t \in [0,T]$,  there exist a unique $\hS^{>0}_d -$valued, $\bF^t-$progressively measurable process $\hat{q}^t$ such that
 for any $P \in \cQ^t_W$, it holds $P-$a.s. that $\big(\hat{q}^t_s\big)^2
= \hat{q}^t_s\cd \hat{q}^t_s = \hat{a}^t_s$ for a.e. $s \in [t,T]$.

 \end{lemm}

 \ss For any   $P \neg \in  \neg  \cQ^t_W$, we define
   $  
   \fI^P_s  \dfnn   \int^{P}_{[t,s]}          \big(\hat{q}^t_r\big)^{-1  }    d B^t_r  $, $    s \in [t,T ]$,
       which is   a continuous semi-martingale with respect to $\big(  \bF^P, P  \big) $.
  Since the first part of \eqref{eq:xd011} and    \eqref{eq:a075} imply  that \pas,
   \bea  \label{eq:xc073}
    \lan B^t  \ran^P_s   \neg = \neg  \int_t^s  \hat{a}^t_r dr  , \q  \fa    s  \neg \in \neg   [t,T] ,
   \eea
   one can deduce from Lemma \ref{lem_hat_a_half} that \pas  
      \beas   
       \big\lan  \fI^P  \big\ran^{P}_s
      = \int_t^s    \big(\hat{q}^t_r\big)^{-1} \neg  \cd  \neg  \big(\hat{q}^t_r\big)^{-1}    d \lan B^t \ran^{P}_r
      = \int_t^s   \big(\hat{q}^t_r\big)^{-1}  \neg \cd \neg  \big(\hat{q}^t_r\big)^{-1}
   \neg \cd    \hat{a}^t_{r} \, d r = (s-t) \, I_{d \times d}    ,    \q  \fa  s \in [t,T] .
     \eeas
 In light of  L\'evy's characterization, the martingale part $W^P $ of $ \fI^P $  is a  Brownian motion
 under $P$.  Let $\bG^P =\big\{ \cG^P_s \big\}_{s \in [t,T]}$ denote 
 the $P-$augmented filtration generated
 by the $P-$Brownian motion $W^P$, i.e.
 \beas
  \cG^P_s  \dfnn \si \Big( \si \big( W^P_r, r \in [t,s] \big)    \cup  \sN^P  \Big)  , \q  \fa  s \in [t,T] .
 \eeas

   \ss     Let $(t,x) \in [0,T] \times \hR^d$ and $\mu   \in   \cU^t$. According to our specification \eqref{eq:xd133},
 $\left\{ X^{t,x,\mu }_s \right\}_{s \in [t,T]} \in \hC^2_{\ol{\bF}^t }([t,T], \hR^k)$
 stands for the unique solution of
  the following  SDE on the probability space $\big( \O^t, \cF^t_T, P^t_0 \big)$:
   \bea
  X_s= x + \int_t^s b(r, X_r  ) \, dr + \int_t^s  \mu_r  \, dB^t_r,  \q s \in [t, T] .      \label{FSDE3}
    \eea
 Similar to \eqref{eq:esti_X_1}, we have
  \bea
     E_t \left[  \underset{s \in [t,T]}{\sup}  \big| X^{t,x,\mu }_s  \big|^2 \right]
      & \tneg \le& \tneg   c_0  \bigg( 1 + |x|^2 + E_t \int_t^T  | \mu_s |^2   \, ds  \bigg) < \infty .                 \label{eq:esti_X_11}
 \eea
 By Lemma \ref{lem_F_version} (2), $X^{t,x,\mu }$   admits   a unique $\bF^t-$version $\wt{X}^{t,x,\mu } \in \hC^2_{\bF^t }([t,T], \hR^k)$  which also satisfies \eqref{eq:esti_X_11}.
 As $\wt{X}^{t,x,\mu }$ has continuous paths except on some $\cN^{\,t,x}_\mu \in \cF^t_T$ with $P^t_0 \big(\cN^{\,t,x}_\mu\big)=0$,
  we can view
 \bea    \label{def_X_alpha}
  \cX^{t,x,\mu}  \dfnn \b1_{ \big( \cN^{\,t,x}_\mu \big)^c} \Big( \wt{X}^{ \raisebox{0.3ex}{\scriptsize $t,x,\mu$}} -x \Big)
 \eea
as a mapping from $ \O^t$ to $ \O^t$. We claim that $ \cX^{t,x,\mu}$ is actually a measurable mapping
 from $\big( \O^t,   \cF^t_T \big)$ to $\big( \O^t,   \cF^t_T \big)$:
    To see this, we pick up  an arbitrary pair
  $(s ,  \cE) \in [t,T] \times \sB(\hR^d)$.     The $\bF^t-$adaptness of  $ \wt{X}^{t,x,\mu} $    implies that
   \bea
   \big(\cX^{t,x,\mu}\big)^{-1}\Big( \big(B^t_s\big)^{-1}(\cE)\Big) &=& \big\{\o \in \O^t: \cX^{t,x,\mu} (\o) \in \big(B^t_s\big)^{-1}(\cE) \big\}
   =\big\{\o \in \O^t: \cX^{t,x,\mu}_s (\o)   \in \cE \big\}  \nonumber \\
   &=& \left\{
      \ba{ll}
       \cN^{\,t,x}_\mu \cup \Big( (\cN^{\,t,x}_\mu)^c  \cap \big\{\o \in \O^t: \wt{X}^{ \raisebox{0.3ex}{\scriptsize $t,x,\mu$}}_s  (\o)   \in \cE_x \big\}  \Big)
       \in    \cF^t_T \, ,   \q \; & \hb{if }0 \in \cE , \ms \\
       (\cN^{\,t,x}_\mu)^c \cap \big\{\o \in \O^t: \wt{X}^{ \raisebox{0.3ex}{\scriptsize $t,x,\mu$}}_s  (\o)   \in \cE_x \big\}
       \in  \cF^t_T \, ,  \q & \hb{if }0 \notin \cE  ,
       \ea
       \right.  \qq   \qq     \label{eq:a247}
   \eea
where $\cE_x = \{x+x': x' \in \cE\} \in \sB(\hR^d)$.
   Thus $\big(B^t_s\big)^{-1}(\cE)  \in \L^t \dfnn  \Big\{A \subset \O^t: \big(\cX^{t,x,\mu}\big)^{-1}(A) \in   \cF^t_T  \Big\}$.
   Clearly, $\L^t$ is
     a $\si-$field  of $\O^t$.
     It follows that
  \bea   \label{eq:a151}
    \cF^t_T  =\si \left(\big(B^t_s\big)^{-1}(\cE) ; \,  s \in [t,T] , ~ \cE   \in   \sB(\hR^d)  \right)  \subset \L^t   ,
 \eea
    proving the measurability of the mapping $ \cX^{t,x,\mu} $.
    Consequently, we can induce a  probability measure
       \bea   \label{eq:a155}
   P^{t,x,\mu} \dfnn P^t_0  \circ \big( \cX^{t,x,\mu} \big)^{-1}
   \eea
   on $\big(\O^t,  \cF^t_T\big)  $, i.e. $P^{t,x,\mu} \in \cP^t$.
   Similar to \cite{STZ_2011a}, we set  $  \cQ^{t,x}_S \dfnn \{ P^{t,x,\mu} \}_{ \mu \in \cU^t } $.

  \if{0}
 There exists a $P^t_0-$null set $\cN^{\,t,x}_\mu$ such that
 for each $\o \in (\cN^{\,t,x}_\mu)^c $
 \beas
  \wt{X}^{t,x,\mu }_s (\o) = X^{t,x,\mu }_s (\o) , \q \fa s \in [t,T]
 \q \hb{and} \q \hb{the path } s \to \wt{X}^{t,x,\mu }_s (\o)       \hb{ is continuous.}
 \eeas
  \fi

     \begin{lemm}  \label{lem_X_alpha}
  Given  $(t,x) \in [0,T] \times \hR^d$ and $\mu \in \cU^t$, let $\cX^{t,x,\mu}: \O^t \to \O^t$  be the mapping    defined in \eqref{def_X_alpha}. It holds  for any   $s \in [t,T]$ that
    $\big(\cX^{t,x,\mu}\big)^{-1} \Big(\cF^{P^{t,x,\mu}}_s    \Big) \subset  \ol{\cF}^t_{\neg s}  $\,.
     Moreover, we have
    \bea  \label{eq:d265}
    P^{t,x,\mu} = P^t_0 \circ \big(\cX^{t,x,\mu}\big)^{-1}   \; \hb{ on }    \;    \cF^{P^{t,x,\mu}}_T   .
    \eea
    \end{lemm}

         \begin{prop} \label{prop_PS-in-PW}
  For any $(t,x) \in [0,T] \times \hR^d$,   we have  $\cQ^{t,x}_S \subset \cQ^t_W$.
   \end{prop}

The following result about $\cQ^{t,x}_S$ is inspired by  Lemma 8.1 of \cite{STZ_2011a}.

        \begin{prop} \label{prop_filtration_coincide}
 Let $(t,x) \in [0,T] \times \hR^d$. For any $P    \in    \cQ^{t,x}_S $,  $ \bF^P$  coincides with  $ \bG^P$,
 the $P-$augmented filtration generated  by the $P-$Brownian motion $W^P$.

    \end{prop}

    Fix $(t,x)  \in [0,T] \times \hR^d$ and  set $B^{t,x} \dfnn x + B^t$. By the continuity of $\ul{l}$ and $\ol{l}$,
 $ \ul{\cL}^{t,x}_s \dfnn \ul{l} \big(s, B^{t,x}_s \big) $  and
 $    \ol{\cL}^{t,x}_s \dfnn  \ol{l} \big(s, B^{t,x}_s \big)$, $   s \in [t,T] $
 are two real-valued, $\bF^t-$adapted continuous processes satisfying $\ul{\cL}^{t,x}_s < \ol{\cL}^{t,x}_s $, $\fa s \in [t,T]$.
 The measurability of $\big(f, B^{t,x}\big)$, the measurability of $\hat{q}^t$
 by Lemma \eqref{lem_hat_a_half} as well as the Lipschitz continuity of $f$ in $(y,z)$ imply that
 \beas
 \hat{f}^{t,x} (s,\o,y,z) \dfnn   f(s,B^{t,x}_s (\o) , y,z, \hat{q}^t_s(\o) \big),
\q \fa (s,\o,y,z) \in [t,T] \times \O^t \times \hR \times \hR^d
 \eeas
 is a    $\sP \big(  \bF^t  \big)   \otimes \sB(\hR)  \otimes \sB(\hR^d)/\sB(\hR)-$measurable function
      that is also Lipschitz continuous in $(y,z)$.

  \ss   Given  $\mu \in \cU^t$,  one can deduce from      \eqref{2l_growth},
  the version of \eqref{f_linear_growth} and   \eqref{f_Lip} without $\nu-$controls, H\"older's inequality, \eqref{eq:d251a} and \eqref{eq:esti_X_11}  that
   \beas
   && \hspace{-1cm}   E_{P^{t,x,\mu}} \neg  \bigg[   \underset{s \in [t,T]}{\sup}    \big| \ul{\cL}^{t,x}_{  s}\big|^q
   \neg +   \neg  \underset{s \in [t,T]}{\sup} \big| \ol{\cL}^{t,x}_{  s}\big|^q
  \neg+  \neg    \Big(\int_t^T  \neg \big| \hat{f}^{t,x}_\t   (s, 0,0 )   \big|  ds \Big)^q\bigg]
   \neg  \le  \neg  c_0    \neg + \neg  c_0 \,  E_{P^{t,x,\mu}} \neg \bigg[   \underset{s \in [t,T]}{\sup}   \big| B^{t,x}_s \big|^2
   \neg + \neg \int_t^T \dneg   | \hat{q}^t_s |^2        ds \bigg]   \\
  && \le \neg  c_0    \neg + \neg  c_0 \,  E_t \neg \bigg[   \underset{s \in [t,T]}{\sup}   \big| B^{t,x}_s \big(\cX^{t,x,\mu}\big) \big|^2
   \neg + \neg \int_t^T \dneg   \big| \hat{q}^t_s \big(\cX^{t,x,\mu}\big)  \big|^2        ds \bigg] \le
  \neg  c_0    \neg + \neg  c_0 \,  E_t \neg \bigg[   \underset{s \in [t,T]}{\sup}   \big|   \wt{X}^{t,x,\mu}_s  \big|^2
   \neg + \neg \int_t^T \dneg    | \mu_s     |^2        ds \bigg]   \neg < \neg  \infty .
  \eeas
 As $\ul{\cL}^{t,x}_T    \le h ( B^{t,x}_T  ) \le \ol{\cL}^{t,x}_T $,    it follows
 that   $  E_{P^{t,x,\mu}} \big[ | h ( B^{t,x}_T  ) |^q  \big]     < \infty $, i.e.  $\xi \in \hL^q \big( \cF^t_\t , {P^{t,x,\mu}} \big)$.
 Then Proposition \ref{prop_filtration_coincide} and Theorem 4.1  of  \cite{EHW_2011} shows  that
    the following Doubly reflected BSDE 
 on the probability space $ \big(\O^t, \cF^{P^{t,x,\mu}}_T, {P^{t,x,\mu}} \big)$
 \beas   
    \left\{\ba{l}
 \dis   \cY_s= h ( B^{t,x}_T  ) \neg + \neg  \int_s^T  \neg    \hat{f}^{t,x}  \big(r,    \cY_r, \cZ_r  \big)  \, dr  \neg+ \neg  \ul{\cK}_{\,T} \neg-\neg \ul{\cK}_{\,s}
  \neg- \neg \big( \ol{\cK}_T \neg - \neg \ol{\cK}_s \big)
   \neg-\neg \int^{P^{t,x,\mu}}_{[s,T]} \neg \cZ_r d W^{P^{t,x,\mu}}_r  , \q    s \in [t,T] ,   \vspace{1mm} \q \\
   \dis \ul{\cL}^{t,x}_{ s}   \le \cY_s
   \le  \ol{\cL}^{t,x}_{ s}  , \q    s \in [t,T] ,
   \q \hb{and} \q  \int_t^T \neg  \big(  \cY_s   -  \ul{\cL}^{t,x}_{ s}  \big)     d \ul{\cK}_{\,s}
  = \int_t^T \neg  \big(  \ol{\cL}^{t,x}_{ s} - \cY_s  \big)     d \ol{\cK}_s   = 0
     \ea \right.
    \eeas
    admits a unique solution    $\Big(\cY^{t,x, {P^{t,x,\mu}}} \neg \big(T,  h ( B^{t,x}_T  ) \big)  ,
 \cZ^{t,x, {P^{t,x,\mu}}}  \neg \big(T,  h ( B^{t,x}_T  ) \big)  ,
 \ul{\cK}^{t,x, {P^{t,x,\mu}}}  \neg \big(T,  h ( B^{t,x}_T  ) \big)  ,
   \ol{\cK}^{t,x, {P^{t,x,\mu}}}  \neg \big(T,  h ( B^{t,x}_T  ) \big)  \Big) \\ \in \hG^q_{\bF^{ P^{t,x,\mu }}}  \big([t,T], {P^{t,x,\mu}} \big)$.

\begin{prop} \label{prop_DRBSDE_transform}
    For any $(t,x)  \in [0,T] \times \hR^d$ and $\mu \in \cU^t$
  \beas
  P^t_0 \Big(  \Big(  \cY^{t,x,P^{t,x,\mu}}_s \big(T, h ( B^{t,x}_T  )\big) \Big)   \big(\cX^{t,x,\mu}\big)
  = Y^{t,x,\mu}_s \big(T,  h ( \wt{X}^{t,x,\mu}_T )     \big), ~ \fa s \in [t,T]\Big) = 1 .
  \eeas

 \end{prop}

  Let $ \wt{\cY}^{t,x,P^{t,x,\mu}}_s \neg \big(T, h ( B^{t,x}_T  )\big) $ be the $\bF^t-$version
 of $\cY^{t,x,P^{t,x,\mu}}_s \neg \big(T, h ( B^{t,x}_T  )\big) $.  For the constant $y^{t,x,\mu} \dfnn \wt{Y}^{t,x,\mu}_t \Big( T, h(\wt{X}^{t,x,\mu}_T ) \Big)  \\  \in \cF^t_t$, one can deduce from \eqref{eq:d265} that
   \beas
  1 = P^t_0 \Big\{ \Big(\cY^{t,x,P^{t,x,\mu}}_t \neg \big(T, h ( B^{t,x}_T  )\big) \Big) (\cX^{t,x,\mu}) 
   = y^{t,x,\mu} \Big\}
 = P^{t,x,\mu}  \Big\{ \cY^{t,x,P^{t,x,\mu}}_t \neg \big(T, h ( B^{t,x}_T  ) \big)  = y^{t,x,\mu} \Big\}.
  \eeas
  It follows that $ y^{t,x,\mu} = \wt{\cY}^{t,x,P^{t,x,\mu}}_t \neg \big(T, h ( B^{t,x}_T  ) \big) $. Hence,
  \bea  \label{eq:xxb011}
 w (t,x) = \underset{\mu \in \cU^t}{\sup} \wt{Y}^{t,x,\mu}_t \Big( T, h(\wt{X}^{t,x,\mu}_T ) \Big)
  = \underset{P \in \cQ^{t \neg ,x}_S }{\sup} \wt{\cY}^{t,x,P }_t \Big(T, h\big(B^{t,x}_T \big)\Big)  ,
 \eea
 which extended the value function of \cite{STZ_2011c} \big(see (5.9) therein\big)
 to the case of doubly reflected BSDEs based on more general forward SDEs.
 Thus our value function $w$ is closely related to the second-order doubly reflected BSDEs.
 On the other hand, when the generator $f \equiv 0$,   the right-hand-side of \eqref{eq:xxb011}
  is a doubly reflected version of the value function considered  in \cite{Nutz_2011}.

 \section{Proofs}

  \label{sec:Proofs}

\subsection{Proofs of Section \ref{sec:intro}  \& \ref{sec:zs_drgame}}

 \label{subsection:Proofs_S1_S2}

   \ss \no {\bf Proof of Lemma \ref{lem_countable_generate1}:}
     For any $s \in [t,T]  $,  it is clear that
$   \si \big(\cC^{t,T}_s \big)  \subset  \si   \Big\{ \big(B^{t,T}_r\big)^{-1} (\cE ) : r \in [t,s], \cE \in \sB(\hR^d) \Big\}  = \cF^{t,T}_s $.
 To see the reverse, we fix $r \in [t,s]$.  For any $x \in \hQ^d$ and $ \l \in \hQ_+$, let $\{s_j\}_{j \in \hN} \subset \hQ_{r,s} $
  with $\lmtd{j \to \infty} s_j =r $.
 Since $\O^{t,T}$ is the set of $\hR^d-$valued continuous functions on $[t,T]$ starting from $0$, we can deduce that
 \beas   
  \big(B^{t,T}_r\big)^{-1} \big(O_\l(x)\big) = \underset{n = \lceil \frac{2}{ \l} \rceil   }{\overset{\infty}{\cup}}
 \underset{m \in \hN }{\cup}  \underset{j \ge m}{\cap} \Big( \big(B^{t,T}_{s_j}\big)^{-1} \big(O_{\l-\frac{1}{n}}(x)\big) \Big)
 \in \si \big( \cC^{t,T}_s \big) .
 \eeas
  which implies that
  \beas
  \cO   \dfnn   \big\{ O_\l(x) :  x \in \hQ^d ,  \l \in \hQ_+  \big\} \subset \L_r \dfnn  \Big\{ \cE \subset \hR^d: \big(B^{t,T}_r\big)^{-1} \neg (\cE  ) \in \si \big( \cC^{t,T}_s \big) \Big\} .
  \eeas
Clearly,  $\cO$ generates $\sB(\hR^d)$ and $\L_r$ is a $\si -$field of $\hR^d$. Thus, one has $\sB(\hR^d)  \subset \L_r$.
Then  it follows that
 \beas
  \hspace{4.2cm}   \cF^{t,T}_s = \si \Big\{ \big(B^{t,T}_r\big)^{-1} (\cE ) : r \in [t,s], \cE \in \sB(\hR^d) \Big\}
   \subset \si \big( \cC^{t,T}_s \big) .
     \hspace{4.2cm} \hb{\qed}
 \eeas

    \ss \no {\bf Proof of Lemma \ref{lem_shift_inverse}:} For simplicity, let us  denote $\Pi^{T,S}_{t,s}$ by $\Pi$.
 We first show the continuity of $\Pi$.  Let $A$ be an open   subset  of $\O^{s,S}$. Given    $\o  \in \Pi^{-1}  (A)$,
   since $\Pi (\o) \in A$,  there exist a $\d > 0$ such that
  $  O_\d \big( \Pi (\o) \big) \dfnn \Big\{ \wt{\o} \in \O^{s,S} :  \underset{r \in [s,S]}{\sup}
   \big|\wt{\o}(r) - \big(\Pi (\o) \big)(r) \big| < \d \Big\} \subset A $.
  For any $\o' \in O_{  \d /2}( \o)   $,   one can deduce that
   \beas
      \underset{r \in [s,S]}{\sup}
   \big| \big(\Pi (\o')\big)(r) - \big(\Pi (\o) \big)(r) \big|
   \le \big|  \o' (s) -  \o (s) \big| + \underset{r \in [s,S]}{\sup}
   \big|  \o' (r) -  \o (r) \big|
  \le 2 \underset{r \in [t,T]}{\sup} |\o'(r) - \o (r)  | < \d ,
   \eeas
  which shows that  $ \Pi (\o') \in  O_\d \big( \Pi (\o) \big)  \subset A$
    or $\o' \in \Pi^{-1} (A)$.   Hence, $\Pi^{-1} (A)$   is  an open   subset   of $\O^{t,T}$.

  \ss   Now, let $r \in  [s,S]$.   For any $s' \in [s,r  ] $ and $ \cE \in \sB(\hR^d)$, one can deduce that
   \bea   \label{eq:g313}
       \Pi^{-1} \neg \Big( \big( B^{s,S}_{s'} \big)^{-1} \neg (\cE) \Big) \neg  =  \neg
    \Big\{\o \in \O^{t,T} \dneg : B^{s,S}_{s'} \big(\Pi (\o) \big)  \neg \in  \neg  \cE  \Big\}
    \neg  = \neg  \big\{\o \in \O^{t,T} \dneg : \o(s')  \neg -   \neg  \o(s) \neg  \in  \neg   \cE \big\}
      \neg = \neg   (B^{t,T}_{s'} \neg - \neg  B^{t,T}_s)^{-1} (\cE)  \neg \in \neg  \cF^{t,T}_r .
    \eea
    \if{0}
     Given $s' \in (s,r]$,
   since
   \bea  \label{eq:g301}
   \o \in \big( B^{t,T}_{s'}- B^{t,T}_s \big)^{-1} (\cE)  \hb{ is equivalent to }  \Pi (\o) \in \big( B^{s,S}_{s'}  \big)^{-1} (\cE),
   \eea
   we see that $ \Pi \Big( \big( B^{t,T}_{s'}- B^{t,T}_s \big)^{-1} (\cE) \Big)  \subset \big( B^{s,S}_{s'}  \big)^{-1} (\cE)$; On the
   other hand, for any $ \wt{\o} \in  \big( B^{s,S}_{s'}  \big)^{-1} (\cE)$, we set
   an $\o \in \O^{t,T}$ by $\o(r') \dfnn \b1_{\{r' \ge s\}}\wt{\o} (r') $,
   $\fa r' \in [t,T]$. Clearly, $  \Pi (\o) =\wt{\o} \in  \big( B^{s,S}_{s'}  \big)^{-1} (\cE) $,  the equivalence \eqref{eq:g301}
    again shows that $\o \in \big( B^{t,T}_{s'}- B^{t,T}_s \big)^{-1} (\cE) $. Then it follows that $ \wt{\o} =  \Pi (\o)
    \in \Pi \Big( \big( B^{t,T}_{s'}- B^{t,T}_s \big)^{-1} (\cE) \Big)  $,
     thus $ \Pi \Big( \big( B^{t,T}_{s'}- B^{t,T}_s \big)^{-1} (\cE) \Big) = \big( B^{s,S}_{s'}  \big)^{-1} (\cE) \in \cF^{s,S}_r$.
  \fi
  Thus all the generating sets of $\cF^{s,S}_r$ belong to $\L \dfnn \big\{ A \subset \O^{s,S}: \Pi^{-1}  (A) \in \cF^{t,T}_r \big\}$, which is clearly a $\si-$field of $\O^{s,S}$. It follows that $\cF^{s,S}_r \subset  \L$, i.e., $\Pi^{-1}  (A) \in \cF^{t,T}_r$ for
  any $A \in \cF^{s,S}_r$.

      Next, we show that the the induced probability $\wt{P} \dfnn P^{t,T}_0 \circ \Pi^{-1}$
    equals to $  P^{s,S}_0  $ on $\cF^{s,S}_S$:
   Since the Wiener measure on $\big( \O^{s,S}, \sB(\O^{s,S})\big)$ is unique (see e.g. Proposition I.3.3 of \cite{revuz_yor}),
    it suffices to show that  the canonical process $B^{s,S}$ 
 is a Brownian motion on $\O^{s,S}$ under $\wt{P}$:
  Let $ s \le r \le r ' \le S$. For any $\cE \in \sB(\hR^d)$,
 similar to \eqref{eq:g313}, one can deduce that
 \bea  \label{eq:g319}
   \Pi^{-1} \big(\big(B^{s,S}_{r'}-B^{s,S}_{r}\big)^{-1}(\cE)\big)
        =     (B^{t,T}_{r'} \neg - \neg  B^{t,T}_{r})^{-1} (\cE)  .
        \eea
         Thus, $  \wt{P}\Big(  \big(B^{s,S}_{r'} \neg - \neg B^{s,S}_{r}\big)^{-1}(\cE) \Big)
   =     P^{t,T}_0 \Big(   \Pi^{-1}
 \big(\big(B^{s,S}_{r'} \neg - \neg B^{s,S}_{r}\big)^{-1}(\cE)\big)  \Big)
 = P^{t,T}_0 \Big((B^{t,T}_{r'} \neg - \neg  B^{t,T}_{r})^{-1} (\cE) \Big)  $,
   which shows that   the distribution of $B^{s,S}_{r'} \neg - \neg B^{s,S}_{r}$ under $ \wt{P}$ is the same as that of
   $B^{t,T}_{r'} \neg - \neg  B^{t,T}_{r}$ under $P^{t,T}_0 $ (a $d-$dimensional normal distribution with mean $0$ and variance matrix
    $  (r'-r)I_{d \times d} $).

    \ss On the other hand,  for any $A \in  \cF^{s,S}_{r}$, since
      $ \Pi^{-1} (A) $ belongs to $  \cF^{t,T}_{r}$,   
         its independence from  $B^{t,T}_{r'} \neg - \neg  B^{t,T}_{r}$ under $P^{t,T}_0$  and \eqref{eq:g319} yield  that
         for any $ \cE \in \sB(\hR^d) $
           \beas
 && \hspace{-1cm} \wt{P}\Big(A \cap \big(B^{s,S}_{r'} \neg - \neg B^{s,S}_{r}\big)^{-1}(\cE) \Big)
   =     P^{t,T}_0 \Big(\Pi^{-1} (A) \cap  \Pi^{-1}
 \big(\big(B^{s,S}_{r'} \neg - \neg B^{s,S}_{r}\big)^{-1}(\cE)\big)  \Big) \\
 &  &     = P^{t,T}_0 \Big(\Pi^{-1} (A) \Big) \cd  P^{t,T}_0 \Big(   \Pi^{-1} \big(\big(B^{s,S}_{r'} \neg - \neg B^{s,S}_{r}\big)^{-1}(\cE)\big) \Big)
   = \wt{P} (A) \cd \wt{P}\Big( \big(B^{s,S}_{r'} \neg - \neg B^{s,S}_{r}\big)^{-1}(\cE) \Big) .
 \eeas
  Hence, $ B^{s,S}_{r'} \neg - \neg B^{s,S}_{r}$ is independent of $ \cF^{s,S}_{r} $ under $\wt{P}$.   \qed


  \ss \no {\bf Proof of Lemma \ref{lem_F_version}:}

  \ss \no (1)       First, let   $\xi \neg \in  \neg  \hL^1 \big( \cF^P_T, P\big)$ and $s  \neg \in \neg  [t,T]$.
 For any $A  \neg \in \neg   \cF^P_s$,
   there exists an $\wt{A}  \in \cF^t_s$   such that $ A \, \D \,  \wt{A}  \in \sN^{P}  $ (see e.g. Problem 2.7.3 of \cite{Kara_Shr_BMSC}).  
 Thus we have that
     $
       \int_A   \xi  d P   \neg  =  \neg    \int_{\wt{A}}  \xi  d P
 \neg =  \neg  \int_{\wt{A}}  E_P \big[ \xi \big|  \cF^t_s  \big]  d P
     \neg    =  \neg  \int_A   E_P \big[ \xi \big|  \cF^t_s  \big]   d P  
     $,
     which implies that
  \bea  \label{eq:xc053}
  E_{P} \big[ \xi \big| \cF^P_s  \big] =  E_P \big[ \xi \big|  \cF^t_s  \big]  ,  \q  \pas
  \eea
 Then it easily follows that any martingale $X$ with respect to $(\bF^t, P)$ is also a martingale with respect to $(\bF^P, P)$.

   \ss Next,   let $X \neg = \neg \{X_s\}_{s \in [t,T]}$  be a local martingale with respect to $(\bF^t, P)$.
  There exists    an increasing  sequence  $\{\t_n\}_{n \in \hN}$  of $\bF^t-$stopping times
 with  $P\Big( \lmtu{n \to \infty} \t_n = T \Big)=1$
         such that $  B^t_{\t_n \land\, \cd}  $\,, $n \in \hN$ are all martingales with respect to $(\bF^t, P)$.
For any $m \in \hN$,
         $   
         \nu_m \dfnn \inf\{s \in [t,T] : |X_s| > m \} \land T         $
             defines an $\bF^t-$stopping time.     In light of the Optional Sampling theorem,  $  X_{\t_n \land    \nu_m   \land \cd}  $
   is  a martingale  with respect to $(\bF^t, P)$. Thus, for any $  t \neg \le \neg  s  \neg < \neg  r  \neg \le \neg  T $,
    one has
    \bea \label{eq:a029}
       E_{P}  \Big[  X_{\t_n \land   \nu_m \land r}    \Big| \cF^P_s \Big]
  = E_{P}  \Big[  X_{\t_n \land   \nu_m \land r}    \Big| \cF^t_s \Big]
  = X_{\t_n \land   \nu_m \land s}    ,     \q     P-\hb{a.s.}
    \eea
              Since $P\Big( \lmtu{n \to \infty} \t_n = T \Big) =1$, when $n \to \infty$ in \eqref{eq:a029},
       the bounded convergence theorem implies that
     $  E_{P}  \Big[     X_{      \nu_m  \land  r }    \Big| \cF^P_s  \Big] \\ =  X_{     \nu_m \land s}  $,   \pas   ~
           Namely,   $X_{\nu_m  \land \cd } $ is a bounded  martingale with respect to $\big(  \bF^P, P  \big) $.
       Clearly, $\{ \nu_m \}_{k \in \hN}$ are $\bF^P-$stopping times with $  \lmtu{ k \to \infty}  \nu_m  = T$. Hence,
         $X$ is a local martingale with respect to $\big(  \bF^P, P  \big) $,  More general,
        any semi-martingale  with respect to $(\bF^t, P)$ is also a semi-martingale with respect to $\big(\bF^P, P \big)$.

  \ss \no (2) The uniqueness is obvious and it suffices to show  the existence for case $\hE = \hR $:
  Let  $\{X_s\}_{s \in [t,T]}$   be    a real-valued, $\bF^P-$adapted  continuous process.
   For each $s \in \hQ_{t,T}   $, we see from \eqref{eq:xc053} that
   \beas
 \wt{X}_s \dfnn E_P \big[ X_s \big|  \cF^t_s  \big] = E_{P} \big[ X_s \big| \cF^P_s  \big] =  X_s  , \q   \pas
   \eeas

   \ss   Set  $\cN \dfnn \big\{\o \in \O^t: $ the path $s \to X_s (\o)$ is not continuous \big\}$\, \cup \,
  \Big( \underset{s \in \hQ_{t,T} }{\cup } \{ X_s \ne \wt{X}_s \} \Big) \in \sN^P $.  Since
  $      X^n_s \dfnn  \sum^{1  +  \lfloor n (T  -  t) \rfloor }_{i=1}
   \wt{X}_{t+\frac{i-1}{n}} \b1_{\{s \in [t+ \frac{i-1}{n}, t+ \frac{i}{n} )\}} $, $   s \in [t,T]$
        is a real-valued, $\bF^t-$progressively measurable process for  any $n \in \hN$,  We see that
  $  \wt{X}_s \dfnn \Big( \lsup{n \to \infty} X^n_s \Big) \b1_{\{ \lsup{n \to \infty} X^n_s < \infty\}}  $
  also  defines a real-valued, $\bF^t-$progressively measurable process.

  Let $ \o \in \cN^c$ and $s \in [t,T]$. For any $n \in \hN$,  since
   $s \in [s_n, s_n+\frac{1}{n}) $ with $  \dis  s_n \dfnn t + \frac{\lfloor n(s-t) \rfloor}{n}$, one has
 $  X^n_s (\o) = \wt{X}_{s_n}(\o) = X_{s_n}(\o)$.
  Clearly, $\lmtu{n \to \infty} s_n =s$,  As $ n \to \infty$,   the continuity of $X$ shows that
 $   \lmt{n \to \infty} X^n_s (\o) =  \lmt{n \to \infty} X_{s_n}(\o) = X_s(\o) $,
 which implies that  $\cN^c \subset \big\{\o \in \O^t:  X_s(\o) = \wt{X}_s(\o) $, $\fa s \in [t,T]  \big\} $.
  Therefore, $\wt{X}$ is $P-$indistinguishable  from $X$, and it follows that  $\wt{X}$ also has $P-$a.s. continuous paths.

  Next, let $\{X_s\}_{s \in [t,T]}$ be a real-valued, $\bF^P-$progressively measurable process that is bounded.
   Since $ \sX_s  \dfnn \int_t^s X_r dr $, $s \in [t,T]$ defines an real-valued, $\bF^P-$adapted continuous
   process, we know from  part (1)   that $\sX$ has a unique $\bF^t-$version $\wt{\sX}$.
   For any $n \in \hN$,
   $ X^n_s  \dfnn n \big(\wt{\sX}_s -\wt{\sX}_{(s-1/n) \vee t} \big) $
   is clearly an  real-valued,  $\bF^t-$adapted continuous  process and thus an  $\bF^t-$progressively measurable  process.
   It follows that $  \wt{X}_s \dfnn \Big( \lsup{n \to \infty} X^n_s \Big) \b1_{\{ \lsup{n \to \infty} X^n_s < \infty\}}  $
   again  defines a real-valued, $\bF^t-$progressively measurable process.

     Set $ \wh{\cN} \dfnn \big\{\o \in \O^t: \sX_s (\o) \ne \wt{\sX}_s (\o) \hb{ for some }s \in [t, T] \big\} \in \sN^P $.
     For any $\o \in \wh{\cN}^c$,   one can deduce that
     \beas
     \lmt{n \to \infty} X^n_s (\o) =  \lmt{n \to \infty} n \big( \sX_s - \sX_{(s-1/n) \vee t} \big)
       =\lmt{n \to \infty} n \int_{(s-1/n) \vee t}^s X_r dr = X_s, \q  \hb{for a.e. }s \in [t,T] ,
     \eeas
      which implies that $ \wt{X}_s (\o) = X_s (\o)$ for $ds \times dP-$a.s. $(s,\o) \in [t,T] \times \O^t$.

  \ss Moreover, for general  real-valued, $\bF^P-$progressively measurable process $\{X_s\}_{s \in [t,T]}$, let $\wt{X}^m$
   be the $\bF^t-$version of $ \big\{X^m_s = (-m ) \vee \big( X_s \land m \big) \big\}_{s \in [t,T]}$ for any  $m \in \hN$.
   Then $  \wt{X}_s \dfnn \Big( \lsup{m \to \infty} \wt{X}^m_s \Big) \b1_{\{ \lsup{m \to \infty} \wt{X}^m_s < \infty\}}  $
     defines a real-valued, $\bF^t-$progressively measurable process. Let $\cD
      \dfnn \underset{m \in \hN }{\cup } \big\{(s,\o) \in [t,T] \times \O^t : X^m_s (\o)  \ne \wt{X}^m_s (\o)  \big\} $.
   Clearly, $ds \times dP (\cD) = 0$ and it holds for any $ (s,\o) \in ([t,T] \times \O^t) \backslash \cD$ that $ \wt{X}_s (\o) = X_s (\o)  $. \qed


    \begin{lemm} \label{lem_comp_DRBSDE}
  Given $t \in [0,T]$ and two $(t,q)-$parameter sets  $ \big(\xi_1,\ff_1,\ul{L}^1, \ol{L}^1 \big), \big(\xi_2,\ff_2,\ul{L}^2, \ol{L}^2 \big)$ with $  P^t_0 \big(\ul{L}^1_s \le \ul{L}^2_s, \, \ol{L}^1_s \le \ol{L}^2_s, \fa s \in [t,T] \big) =1 $,
 let  $ \big(Y^i ,Z^i ,\ul{K}^i , \ol{K}^i  \big) \in  \hC^q_{\ol{\bF}^t}([t,T]) \times
 \hH^{2,q}_{\ol{\bF}^t}([t,T],\hR^d) \times \hK_{\ol{\bF}^t}([t,T]) \times \hK_{\ol{\bF}^t}([t,T]) $,  $i=1,2$ be a solution
   of DRBSDE$ \big( P^t_0,\xi_i,\ff_i,\ul{L}^i, \ol{L}^i \big)$.
    For either $i=1$ or $i=2$,      if \, $\ff_i$ satisfies \eqref{ff_Lip}, then for any $\varpi \neg \in \neg (1,q]$
\beas
     E_t \bigg[ \underset{s \in [t,T]}{\sup} \big| \big(  Y^1_s \neg - \neg  Y^2_s \big)^+ \big|^\varpi  \bigg]
       \neg \le   \neg  C(T,\varpi, \g )   \Bigg\{ E_t \big[  \big|  ( \xi_1   \neg -  \neg  \xi_2)^+  \big|^\varpi \big]
           \neg +  \neg  E_t \bigg[  \bigg( \neg \int_t^T   \neg  \big(\ff_1 (r, Y^{3-i}_r ,Z^{3-i}_r )
 \neg - \neg  \ff_2 (r, Y^{3-i}_r,Z^{3-i}_r) \big)^+ dr    \bigg)^{\neg \varpi}  \bigg] \Bigg\}   .
  \eeas

\end{lemm}

    \ss \no {\bf Proof:}
 Without loss of generality, let $\ff_1$ satisfy \eqref{ff_Lip}. Fix $\varpi \neg \in \neg (1,q]$. We assume that
 \bea  \label{eq:s531}
  E_t \bigg[ \Big( \neg \int_t^T  \neg \big(\ff_1 (r, Y^2_r ,Z^2_r ) \neg - \neg \ff_2 (r, Y^2_r,Z^2_r) \big)^+  dr \Big)^{\neg \varpi} \, \bigg] \neg < \neg \infty ,
 \eea
 otherwise, the result  holds automatically.

  \ss  For $\fX = \xi,  Y,Z $,    we set $ \D \fX \dfnn \fX^1 -\fX^2 $.
  Applying Tanaka's formula to process $(\D Y)^+$ yields that
 \beas
 (\D Y_s)^+  & \tneg = & \tneg   (\D \xi   )^+ +   \int_s^T \b1_{\{ \D Y_r > 0 \}}
 \big( \ff_1 (r, Y^1_r,Z^1_r) - \ff_2 (r, Y^2_r,Z^2_r) \big) dr  - \frac12 \int_s^T  d \fL_r   \nonumber    \\
 & \tneg & \tneg  + \int_s^T \b1_{\{ \D Y_r > 0 \}}  \big(   d  \ul{K}^1_r- d  \ul{K}^2_r -  d \ol{K}^1_r + d   \ol{K}^2_r \big)
 - \int_s^T \b1_{\{ \D Y_r > 0 \}}  \D Z_r  d B^t_r
  ,   \q \fa t \le   s   \le T  ,
 \eeas
where $\fL$ is a real-valued, $\bF^t-$adapted, increasing and continuous process known as ``local time".
 Then we can deduce from   Lemma 2.1 of \cite{EHW_2011}   that
   \bea
  && \hspace{-1cm} \big| (\D Y_s)^+ \big|^\varpi  - \big| (\D Y_{s'})^+ \big|^\varpi
+ \frac{\varpi(\varpi-1)}{2} \int_s^{s'}  \b1_{\{ \D Y_r > 0 \}}\big| (\D Y_s)^+ \big|^{\varpi-2} |\D Z_r|^2 dr \nonumber \\
      && \hspace{-0.7cm} \le     \neg     \varpi \neg \int_s^{s'}  \dneg  \b1_{\{ \D Y_r \neg>  0 \}}
 \big| (\D Y_r)^+ \big|^{\varpi-1} \big( \ff_1 (r, Y^1_r,Z^1_r)  \neg - \neg  \ff_2 (r, Y^2_r,Z^2_r) \big) dr
    \neg+\neg \varpi  \neg  \int_s^{s'}  \dneg  \b1_{\{ \D Y_r > 0 \}}
 \big| (\D Y_r)^+ \big|^{\varpi-1}  \big(   d  \ul{K}^1_r \neg-\neg d  \ul{K}^2_r
 \neg-\neg  d \ol{K}^1_r \neg+\neg d   \ol{K}^2_r \big) \nonumber \\
 &&  \hspace{-0.4cm} -  \,   \varpi \int_s^{s'}  \neg  \b1_{\{ \D Y_r > 0 \}}
 \big| (\D Y_r)^+ \big|^{\varpi-1} \D Z_r d B^t_r
 \neg -  \neg  \frac{\varpi}{2}  \neg   \int_s^{s'}  \neg  \b1_{\{ \D Y_r > 0 \}} \big| (\D Y_r)^+ \big|^{\varpi-1} d \fL_r,   \q \fa   t \le s \le s' \le T   . \label{eq:n101a}
   \eea

  By the lower flat-off condition of DRBSDE$ \big( P^t_0,\xi_1,\ff_1,\ul{L}^1, \ol{L}^1 \big)$, it holds $P^t_0-$a.s. that
  $$
  0 \neg \le \neg \int_t^T  \neg \b1_{\{ \D Y_r > 0 \}} \big| (\D Y_r)^+ \big|^{\varpi-1}    d \ul{K}^1_r
    \neg  =   \neg   \int_t^T  \neg \b1_{\{ \ul{L}^1_r = Y^1_r > Y^2_r   \}}
     \big| (\ul{L}^1_r -  Y^2_r)^+ \big|^{\varpi-1}     d \ul{K}^1_r
   \neg \le \neg  \int_t^T  \neg   \b1_{\{ \ul{L}^1_r  > \ul{L}^2_r  \}}
    \big| (\ul{L}^1_r -  Y^2_r)^+ \big|^{\varpi-1}  d \ul{K}^1_r = 0.
  $$
 Similarly, the upper flat-off condition of DRBSDE$ \big( P^t_0,\xi_2,\ff_2,\ul{L}^2, \ol{L}^2 \big)$
implies that  $P^t_0-$a.s.
       $$
  0 \neg \le \neg  \int_t^T  \neg \b1_{\{ \D Y_r > 0 \}} \big| (\D Y_r)^+ \big|^{\varpi-1}    d \ol{K}^2_r
   \neg  = \neg    \int_t^T  \neg \b1_{\{  Y^1_r > Y^2_r = \ol{L}^2_r   \}}
    \big| (Y^1_r -\ol{L}^2_r )^+ \big|^{\varpi-1}   d \ol{K}^2_r
   \neg \le  \neg  \int_t^T  \neg  \b1_{\{ \ol{L}^1_r > \ol{L}^2_r \}}
    \big| (Y^1_r -\ol{L}^2_r )^+ \big|^{\varpi-1}  d \ol{K}^2_r = 0.
  $$
 Putting these two inequality back into \eqref{eq:n101a} and using Lipschitz continuity of $\ff_1$ in $(y, z)$,
 we obtain
 \bea
  && \hspace{-0.8cm} \big| (\D Y_s)^+ \big|^\varpi
+ \frac{\varpi(\varpi-1)}{2} \int_s^{s'}  \b1_{\{ \D Y_r > 0 \}}\big| (\D Y_s)^+ \big|^{\varpi-2} |\D Z_r|^2 dr \nonumber \\
      &&    \le   \big| (\D Y_{s'})^+ \big|^\varpi \neg + \neg     \varpi \neg \int_s^{s'}  \dneg  \b1_{\{ \D Y_r \neg>   0 \}}
 \big| (\D Y_r)^+ \big|^{\varpi-1} \Big[ \big( \g |\D Y_r| \neg + \neg  \g |\D Z_r|    \neg + \neg  \D \ff_+ \big) dr
          \neg  -  \neg  \D Z_r d B^t_r \Big] ,   ~ \;    \fa   t \neg \le  \neg  s \neg  \le  \neg s'  \neg \le \neg  T   ,
   \qq  \label{eq:n101}
   \eea
 where $\D \ff_+ \dfnn  \big( \ff_1 (r, Y^2_r,Z^2_r) - \ff_2 (r, Y^2_r,Z^2_r) \big)^+ $.

 \ss  Since  $ \hC^q_{\ol{\bF}^t }( [t,T] ) \subset \hC^\varpi_{\ol{\bF}^t }( [t,T] )$
  and $ \hH^{2,q}_{\ol{\bF}^t }([t,T], \hR^d ) \subset \hH^{2,\varpi}_{\ol{\bF}^t }([t,T], \hR^d ) $ by Jensen's inequality,  the Burkholder-Davis-Gundy inequality and H\"older's inequality imply that for some $c>0$
  \beas
 && \hspace{-1.5cm} E_t \Bigg[ \underset{ s \in [t,T]}{\sup} \bigg| \int_t^s  \neg \b1_{\{ \D Y_r \neg>   0 \}}
 \big| (\D Y_r)^+ \big|^{\varpi-1} \D Z_r d B^t_r \bigg| \, \Bigg]
  \le c E_t \Bigg[     \bigg( \int_t^T  \neg \big| (\D Y_r)^+ \big|^{2 \varpi-2}  \big| \D Z_r \big|^2 d r \bigg)^{1/2}  \Bigg]\\
 &&  \le c E_t \Bigg[   \underset{r \in [t,T]}{\sup} \big| \D Y_r \big|^{\varpi-1} \bigg( \int_t^T  \neg    \big| \D Z_r \big|^2 d r \bigg)^{1/2}  \Bigg]
  \le  c  \big\| \D Y \big\|^{\varpi-1}_{\hC^\varpi_{\ol{\bF}^t }( [t,T] )  }
  \big\|\D Z \big\|_{ \hH^{2,\varpi}_{\ol{\bF}^t }([t,T], \hR^d )} < \infty ,
  \eeas
  which shows that
  \bea   \label{eq:n121}
  \hb{ $\big\{ \int_t^s  \neg  \b1_{\{ \D Y_r \neg>   0 \}}
  \big| (\D Y_r)^+ \big|^{\varpi-1} \D Z_r d B^t_r \big\}_{s \in [t,T]}$
  is a uniformly integrable martingale with respect to $ \big( \ol{\bF}^t,P^t_0 \big) $. }
  \eea
   Then, letting $s=t$, $s' \neg= \neg T$ and taking the expectation $E_t$\ in \eqref{eq:n101},
   we can deduce from   H\"older's inequality, Young's inequality and \eqref{eq:s531}    that
    \beas
   && \hspace{-0.7cm}
  \frac{\varpi(\varpi \neg - \neg 1)}{2} E_t \int_t^T \b1_{\{ \D Y_r > 0 \}}\big| (\D Y_s)^+ \big|^{\varpi-2} |\D Z_r|^2 dr   \\
    && \hspace{-0.3cm} \le    E_t \big[ \, \big| \D \xi  \big|^\varpi \big]
   \neg + \neg \varpi    E_t    \bigg[ \bigg( \underset{r \in [t,T]}{\sup} \big| \D Y_r \big|^{\varpi-1} \bigg)
   \bigg(  \g (T \neg - \neg t) \underset{r \in [t,T]}{\sup}  | \D Y_r |
 + \g \sqrt{T \neg -\neg t} \Big(   \int_t^T  \dneg
     | \D Z_r |^2        dr  \Big)^{1/2}  +\int_t^T  \dneg      \D \ff_+      dr   \bigg) \bigg]   \\
 && \hspace{-0.3cm} \le  \varpi \big( 1 + \g (T \neg - \neg t) \big)
   \big\| \D Y \big\|^\varpi_{\hC^\varpi_{\ol{\bF}^t }( [t,T] )  }
   \neg +  \neg \g^\varpi (T \neg - \neg t)^{\varpi/2} \big\| \D Z \big\|^\varpi_{ \hH^{2,\varpi}_{\ol{\bF}^t }([t,T], \hR^d )}
  \neg  + \neg         E_t \neg \bigg[ \bigg( \neg  \int_t^T  \neg
    \D \ff_+  dr  \bigg)^{ \neg \varpi} \bigg]    \neg < \neg \infty .
   \eeas

   Hence, we can define an increasing sequence of $\bF^t-$stopping times
  \beas
  \t_n \dfnn \inf \Big\{ s \in [t,T]: \int_t^s  \neg  \b1_{\{ \D Y_r > 0\}}  \big| (\D Y_r)^+ \big|^{\varpi-2}  \big|\D Z_r\big|^2 dr > n \Big\} \land T, \q \fa n \in \hN
  \eeas
 such that $\lmtu{n \to \infty} \t_n =T$, $P^t_0-$a.s.
   Fix $n \in \hN$. Since
  $$
           \b1_{\{ \D Y_r > 0\}}    \big| (\D Y_r)^+ \big|^{\varpi-1}       |\D Z_r|
    \le      \frac{ \g }{\varpi \neg - \neg 1   }       \big| (\D Y_r)^+ \big|^\varpi
  \neg  + \neg  \frac{\varpi \neg - \neg 1}{4 \g}
  \b1_{\{ \D Y_r > 0\}}  \big| (\D Y_r)^+ \big|^{\varpi-2} \big| \D Z_r \big|^2 dr, ~  s \in [t,T] ,  
  $$
  letting       $s= \t_n \land s $ and $s' = \t_n $ in \eqref{eq:n101} yields that
    \bea
  && \hspace{-2cm} \big|  \big( \D Y_{\t_n \land s} \big)^+ \big|^\varpi +  \frac{\varpi(\varpi \neg - \neg 1)}{4}  \neg  \int_{\t_n \land s}^{\t_n}  \neg \big| (\D Y_r)^+ \big|^{\varpi-2} \b1_{\{ \D Y_r > 0\}} \big| \D Z_r \big|^2 dr
    \le    \big| ( \D Y_{\t_n} )^+  \big|^\varpi  \neg + \neg \varpi \int_{\t_n \land s}^{\t_n}  \neg  \big| (\D Y_r)^+ \big|^{\varpi-1}   \D \ff_+      dr    \nonumber \\
  &   &      +  \Big(\varpi \g  + \frac{\varpi \g^2}{\varpi \neg - \neg 1  } \Big)  \neg  \int_{\t_n \land s}^{\t_n}  \neg  \big| (\D Y_r)^+ \big|^\varpi dr  -   \varpi \int_{\t_n \land s}^{\t_n}  \neg
   \b1_{\{ \D Y_r > 0\}}    \big| (\D Y_r)^+ \big|^{\varpi-1}  \D Z_r d B^t_r ,
     \q    s \in [t, T]  .   \label{eq:n117}
   \eea
  Taking the expectation $E_t$,    we can deduce from Fubini's Theorem, \eqref{eq:n121} and Optional Sampling Theorem that
  \bea
  &     & \hspace{-2cm}    E_t \Big[ \big| \big( \D Y_{\t_n \land s} \big)^+ \big|^\varpi \Big] +
       \frac{\varpi(\varpi \neg - \neg 1)}{4}  E_t \neg  \int_{\t_n \land s}^{\t_n}  \neg
  \b1_{\{ \D Y_r > 0\}} \big| (\D Y_r)^+ \big|^{\varpi-2}      \big| \D Z_r \big|^2 dr    \nonumber   \\
    &&     \le   E_t \big[ \eta_n \big]     +  \Big(\varpi \g  + \frac{\varpi \g^2}{\varpi \neg - \neg 1  } \Big)  \neg
        \int_s^T \neg E_t   \Big[     \big| \big( \D Y_{\t_n \land r} \big)^+ \big|^\varpi \Big] dr  ,
     \q    s \in [t, T]  ,  \label{eq:n115}
 \eea
where $\eta_n \dfnn \big| (\D Y_{\t_n})^+   \big|^\varpi    +     \varpi \int_t^{\t_n}  \neg  \big| (\D Y_r)^+ \big|^{\varpi-1}
    \D \ff_+     dr$.

 \ss Let $C(T,\varpi,\g) $  denote a generic constant, depending on $ T,\varpi,\g$, whose form may vary from line to line.
 An application of   Gronwall's inequality to \eqref{eq:n115} yields that
   \bea  \label{eq:n119}
    E_t \Big[ \big| \big(\D Y_{\t_n \land s} \big)^+ \big|^\varpi \Big]  +
       \frac{\varpi(\varpi \neg - \neg 1)}{4}  E_t \neg  \int_{\t_n \land s}^{\t_n}  \neg
  \b1_{\{ \D Y_r > 0\}} \big| (\D Y_r)^+ \big|^{\varpi-2}  \big| \D Z_r \big|^2 dr
    \le C(T,\varpi,\g)  E_t \big[ \eta_n \big] ,    \q      s \in [t, T]  .
    \eea
   which together with  Fubini's Theorem  shows  that
    \beas
    E_t    \int_t^{\t_n} \neg | ( \D Y_s )^+ \big|^\varpi ds \le   E_t    \int_t^T \neg \big| \big( \D Y_{\t_n \land s} \big)^+ \big|^\varpi ds
       = \int_t^T \neg  E_t \Big[ \big| \big( \D Y_{\t_n \land s} \big)^+ \big|^\varpi \Big]   ds
        \le C(T,\varpi,\g)  E_t \big[ \eta_n \big]  .
    \eeas
    Then we can deduce from \eqref{eq:n117} that
           \beas
     E_t \bigg[ \underset{s \in [t,\t_n]}{\sup}  \big| ( \D Y_s )^+ \big|^\varpi  \bigg]
    \le  C(T,\varpi,\g)   E_t \big[ \eta_n \big]
          +   \varpi   E_t \bigg[ \underset{s \in [t,T]}{\sup} \bigg| \int_{\t_n \land s}^{\t_n}  \neg
           \b1_{\{ \D Y_r > 0\}} \big| (\D Y_r)^+ \big|^{\varpi-1} \D Z_r d B^t_r \bigg| \, \bigg]   .
   \eeas
       The Burkholder-Davis-Gundy inequality again implies that  for some $c>0$
  \beas
 && \hspace{-0.6 cm} 
    E_t \bigg[ \underset{s \in [t,T]}{\sup} \bigg| \int_s^T \b1_{\{r \le \t_n\}}
   \b1_{\{ \D Y_r > 0\}} \big| (\D Y_r)^+ \big|^{\varpi-1}   \D Z_r d B^t_r \bigg| \, \bigg]
   \le \neg c E_t \bigg[     \bigg( \int_t^T \neg \b1_{\{r \le \t_n\}}
  \b1_{\{ \D Y_r > 0\}} \big| (\D Y_r)^+ \big|^{2\varpi-2}  \big| \D Z_r \big|^2 d r \bigg)^{1/2}  \bigg]  \\
  &&  \le  \neg  c E_t \bigg[   \underset{r \in [t,\t_n]}{\sup} \big| (\D Y_r)^+ \big|^{\varpi/2}
 \bigg( \int_t^{\t_n}  \neg \b1_{\{ \D Y_r > 0\}} \big| (\D Y_r)^+ \big|^{\varpi-2}    \big| \D Z_r \big|^2 d r \bigg)^{ \neg  1/2}  \bigg]  \neg \le \neg  \frac{1}{2\varpi}    E_t \bigg[   \underset{r \in [t,\t_n]}{\sup} \big| (\D Y_r)^+ \big|^\varpi \bigg] \\
  &&  ~  + \frac{\varpi c^2}{2} E_t  \neg  \int_t^{\t_n}  \neg  \b1_{\{ \D Y_r > 0\}}
   \big| (\D Y_r)^+ \big|^{\varpi-2}   \big| \D Z_r \big|^2 d r
   \le  \neg  \frac{1}{2\varpi}    E_t   \bigg[   \underset{r \in [t,\t_n]}{\sup} \big| (\D Y_r)^+ \big|^\varpi \bigg]
   \neg + \neg  C(T,\varpi,\g)  E_t \big[ \eta_n \big] ,
  \eeas
  where   we  used   \eqref{eq:n119} with $s   =   t$ in the last inequality.
    Since $ E_t      \bigg[   \underset{r \in [t,\t_n]}{\sup} \big| (\D Y_r)^+ \big|^\varpi \bigg]
    \le  \big\| \D Y \big\|^\varpi_{\hC^\varpi_{\ol{\bF}^t }( [t,T] )  }   <  \infty$, it follows from   Young's inequality that
  \beas
   E_t \bigg[ \underset{s \in [t,\t_n]}{\sup}  \big| ( \D Y_s )^+ \big|^\varpi  \bigg]
  & \tneg \dneg \le & \tneg \dneg  C(T,\varpi,\g) \,  E_t \big[ \eta_n \big]
    \neg  \le  \neg   C(T,\varpi,\g) \Bigg\{ E_t \Big[ \big| \big(\D Y_{\t_n} \big)^+  \big|^\varpi  \Big]
 \neg  +  \neg  E_t \bigg[   \underset{r \in [t,\t_n]}{\sup} \big| (\D Y_r)^+ \big|^{\varpi-1}
  \dneg  \int_t^{\t_n}  \neg  \D \ff_+   dr \bigg] \Bigg\} \\
  & \tneg \dneg \le & \tneg \dneg   C(T,\varpi,\g) E_t \Big[ \big| \big(\D Y_{\t_n} \big)^+   \big|^\varpi  \Big]
  \neg   +  \neg  \frac{1}{2} E_t \bigg[   \underset{r \in [t,\t_n]}{\sup} \big| (\D Y_r)^+ \big|^\varpi \bigg]  \neg + \neg   C(T,\varpi,\g) E_t \bigg[    \bigg( \int_t^{\t_n}  \neg   \D \ff_+   dr \bigg)^{ \neg \varpi} \, \bigg].
  \eeas
   Hence, we have
  \bea  \label{eq:n125}
E_t \bigg[ \underset{s \in [t,\t_n]}{\sup}  \big| ( \D Y_s )^+ \big|^\varpi  \bigg]
    \le  C(T,\varpi,\g) \Bigg\{ E_t \Big[ \big| \big( \D Y_{\t_n} \big)^+  \big|^\varpi  \Big]
  +   E_t \bigg[    \bigg( \int_t^T   \neg    \D \ff_+   dr \bigg)^{\neg \varpi} \, \bigg] \Bigg\} .
  \eea
  Since $ \D Y \in   \hC^\varpi_{\ol{\bF}^t }( [t,T] )  $ and since $\lmtu{n \to \infty} \t_n =T$, $P^t_0-$a.s., the Dominated Convergence Theorem implies that
   $ \lmt{n \to \infty}  E_t \Big[ \big| \big( \D Y_{\t_n} \big)^+  \big|^\varpi  \Big] 
    = E_t \big[ \big| ( \D \xi )^+  \big|^\varpi  \big]$. Letting $ n \to \infty $ in \eqref{eq:n125} and applying the Monotone Convergence Theorem
    on its left-hand-side lead to that
    \bea   \label{eq:s541}
     E_t \bigg[ \underset{s \in [t,T]}{\sup} \big| \big(  Y^1_s \neg - \neg  Y^2_s \big)^+ \big|^\varpi  \bigg]
       \neg \le   \neg  C(T,\varpi, \g )   \Bigg\{ E_t \big[  \big|  ( \xi_1   \neg -  \neg  \xi_2)^+  \big|^\varpi \big]
           \neg +  \neg  E_t \bigg[  \bigg( \neg \int_t^T   \neg  \big(\ff_1 (r, Y^2_r ,Z^2_r )
 \neg - \neg  \ff_2 (r, Y^2_r,Z^2_r) \big)^+ dr    \bigg)^{\neg \varpi}  \bigg] \Bigg\}   .  \q
  \eea

   \ss \no {\bf Proof of Proposition \ref{prop_comp_DRBSDE}:}
  For either $i=1$ or $i=2$,      if \, $\ff_i$ satisfies \eqref{ff_Lip}, applying Lemma \ref{lem_comp_DRBSDE}
  with $\varpi = q $ yields that
  $  E_t \bigg[ \underset{s \in [t,T]}{\sup} \big| \big(  Y^1_s \neg - \neg  Y^2_s \big)^+ \big|^q  \bigg] = 0$.
   Hence,  it holds $P^t_0-$a.s. that   $  (\D Y_s)^+ = 0$, or $ Y^1_s \le Y^2_s$
  for any $s \in [t,T]$.    \qed

\ss \no {\bf Proof of Proposition \ref{prop_apriori_DRBSDE}:} For any $\varpi \in (1,q]$, it follows from Lemma \ref{lem_comp_DRBSDE} that
 \bea  \label{eq:s537}
     E_t \bigg[ \underset{s \in [t,T]}{\sup} \big| \big(  Y^1_s \neg - \neg  Y^2_s \big)^+ \big|^\varpi  \bigg]
       \neg \le   \neg  C(T,\varpi, \g )   \Bigg\{ E_t \big[  \big|  ( \xi_1   \neg -  \neg  \xi_2)^+  \big|^\varpi \big]
           \neg +  \neg  E_t \bigg[  \bigg( \neg \int_t^T   \neg  \big(\ff_1 (r, Y^2_r ,Z^2_r )
 \neg - \neg  \ff_2 (r, Y^2_r,Z^2_r) \big)^+ dr    \bigg)^{\neg \varpi}  \bigg] \Bigg\}   .
  \eea
 Exchanging the order of $Y^1$, $Y^2$ and applying Lemma \ref{lem_comp_DRBSDE} again give that
  \beas
     E_t \bigg[ \underset{s \in [t,T]}{\sup} \big| \big(  Y^2_s \neg - \neg  Y^1_s \big)^+ \big|^\varpi  \bigg]
       \neg \le   \neg  C(T,\varpi, \g )   \Bigg\{ E_t \big[  \big|  ( \xi_2   \neg -  \neg  \xi_1)^+  \big|^\varpi \big]
           \neg +  \neg  E_t \bigg[  \bigg( \neg \int_t^T   \neg  \big(\ff_2 (r, Y^2_r ,Z^2_r )
 \neg - \neg  \ff_1 (r, Y^2_r,Z^2_r) \big)^+ dr    \bigg)^{\neg \varpi}  \bigg] \Bigg\}   ,
  \eeas
  which together with \eqref{eq:s537} implies \eqref{eq:n211}.  \qed

 \ss \no {\bf Proof of Lemma \ref{lem_DX_varpi}:}  (1)  Set $\Th=(t,x,\mu,\nu)$ and fix $s \in [t,T]$. For any $s' \in [t,s]$, one can deduce from \eqref{FSDE}, \eqref{b_si_Lip} and \eqref{b_si_linear_growth} that
     \beas
   \q \underset{r \in [t,s']}{\sup}  \big| X^\Th_{r} - x  \big|   \le   \g  \neg   \int_t^{s'} \dneg
      \Big( 1 \neg + \neg  |x|  \neg + \neg  \underset{r \in [t, r]}{\sup} \big| X^\Th_{r} \neg - \neg x  \big|
   \neg + \neg  [\mu_r ]_{\overset{}{\hU}}   \neg + \neg   [\nu_r ]_{\overset{}{\hV}} \Big) dr
  + \neg \underset {r \in [t, s']}{\sup} \neg \left|\int_t^{r} \neg   \si \big( r, X^\Th_{r}, \mu_r, \nu_r \big)  dB^t_r \right| , \q P^t_0-a.s.
     \eeas
     Then H\"olders inequality,  Doob's martingale inequality, \eqref{b_si_Lip}, \eqref{b_si_linear_growth}
 and  Fubini's Theorem  imply  that
    \beas
         E_t \neg \left[  \underset{r \in [t,s']}{\sup}  \big| X^\Th_{r} - x  \big|^2 \right]
      & \dneg \dneg  \le & \dneg \dneg       c_0  E_t \neg \int_t^{s'}
   \neg \Big( 1 \neg + \neg |x|^2 +  \underset{r \in [t, r]}{\sup}   \big| X^\Th_{r} \neg - \neg x  \big|^2     \neg
    +   \neg   [\mu_r ]^2_{\overset{}{\hU}}    \neg + \neg   [\nu_r ]^2_{\overset{}{\hV}} \Big)  dr
  + c_0  E_t \neg       \int_t^{s'} \neg  \big| \si\big(r,    X^\Th_{r}, \mu_r, \nu_r \big) \big|^2  dr     \\
    & \dneg \dneg   \le & \dneg \dneg   c_0 (1 \neg + \neg |x|^2) (s-t)  +   c_0   \neg \int_t^{s'}
   \neg  E_t \neg \left[\underset{r \in [t, r]}{\sup}
    \big| X^\Th_{r} \neg - \neg x  \big|^2 \right] dr \neg
    + c_0  E_t \int_t^{s'} \neg \big( [\mu_r ]^2_{\overset{}{\hU}}
  \neg + \neg   [\nu_r ]^2_{\overset{}{\hV}} \big)  dr  , \q   s' \in [t,s].
         \eeas
 An application of  Gronwall's inequality yields that
      \beas
        E_t \left[  \underset{r \in [t,s']}{\sup}  \big| X^\Th_{r} - x  \big|^2 \right]
         \le   c_0 e^{ c_0 T   } (1 \neg + \neg |x|^2 ) (s-t)   + c_0 e^{ c_0 T   }
  E_t \int_t^{s'} \neg \big( [\mu_r ]^2_{\overset{}{\hU}}
  \neg + \neg   [\nu_r ]^2_{\overset{}{\hV}} \big)  dr    , \q s' \in [t,s].
      \eeas
      Taking $s'=s$ gives \eqref{eq:esti_X_2}.

 \ss \no (2)  Given $x' \in \hR^k$, we set $\D X_r \dfnn X^{t,x,\mu,\nu}_r \neg - \neg X^{t,x',\mu,\nu}_r$,
 $\fa r \in [t,T]$. By \eqref{b_si_Lip},
  \beas
  ~ \;    \underset{r \in [  t,s]}{\sup} \big| \D X_r \big|  \neg \le  \neg    |x  \neg - \neg  x'|
   \neg + \neg  \g  \neg \int_t^s  \neg  | \D X_r | dr
    \neg + \neg  \underset{r \in [t,s]}{\sup} \bigg|  \neg \int_t^r  \dneg   \Big( \si \big(r,X^{t,x,\mu,\nu}_r , \mu_r,\nu_r\big)
    \neg - \neg  \si \big(r,X^{t,x',\mu,\nu}_r , \mu_r,\nu_r\big) \Big) d B^t_r \bigg|, \q  \fa s  \neg \in \neg  [t,T] .
  \eeas
  Then  the Burkholder-Davis-Gundy inequality and \eqref{b_si_Lip}   imply that
    \beas
   E_t \bigg[  \underset{r \in [  t,s]}{\sup}  | \D X_r  |^\varpi \bigg] & \tneg \le&  \tneg  c_\varpi |x  \neg - \neg  x'|^\varpi
   + c_\varpi E_t \Bigg[\bigg( \neg \int_t^s | \D X_r | dr \bigg)^\varpi + \bigg(
   \int_t^s  \neg   | \D X_r  |^2  d  r \bigg)^{\neg \varpi/2}\Bigg]  \nonumber  \\
   & \tneg \le& \tneg  c_\varpi |x  \neg - \neg  x'|^\varpi
   + c_\varpi E_t \Bigg[\bigg( \neg \int_t^s | \D X_r | dr \bigg)^\varpi +
   \underset{r \in [  t,s]}{\sup}  | \D X_r  |^{\neg \varpi/2}  \bigg( \int_t^s  \neg   | \D X_r  |   d  r \bigg)^{\neg \varpi/2}\Bigg]  \nonumber \\
   & \tneg \le& \tneg  c_\varpi |x  \neg - \neg  x'|^\varpi
   + c_\varpi E_t \Bigg[\bigg( \neg \int_t^s | \D X_r | dr \bigg)^\varpi \Bigg] + \frac12 E_t \Bigg[
   \underset{r \in [  t,s]}{\sup}  | \D X_r  |^{\neg \varpi }   \Bigg]
    , \q  \fa s  \neg \in \neg  [t,T] .
  \eeas
    As $  E_t \bigg[  \underset{r \in [  t,T]}{\sup}  | \D X_r  |^\varpi \bigg]
     \le 1+ 2 E_t \bigg[  \underset{r \in [  t,T]}{\sup}  \big| X^{t,x,\mu,\nu}_r  \big|^2
     +  \underset{r \in [  t,T]}{\sup}  \big| X^{t,x',\mu,\nu}_r  \big|^2 \bigg] < \infty$ by \eqref{eq:esti_X_1},
    it follows from H\"older's inequality and Fubini's Theorem that
    \bea
      E_t \bigg[  \underset{r \in [  t,s]}{\sup}  | \D X_r  |^\varpi \bigg]   \neg  \le  \neg
      c_\varpi |x  \neg - \neg  x'|^\varpi
   \neg  + \neg  c_\varpi E_t    \neg \int_t^s \big| \D X_r \big|^\varpi dr
   \neg  \le  \neg     c_\varpi |x  \neg - \neg  x'|^\varpi
   \neg  +  \neg  c_\varpi   \neg \int_t^s \neg E_t \bigg[\underset{r' \in [t,r]}{\sup}\big| \D X_{r'} \big|^\varpi \bigg]  dr
    , ~  \fa s  \neg \in \neg  [t,T] .   \q   \label{eq:s018}
    \eea
   An application of    Gronwall's inequality yields \eqref{eq:esti_X_3}.

 \ss \no (3)  Next, Let us assume \eqref{b_si_Lip_u} for some $\l \in (0,1]$.  Given $    \mu'  \in    \cU^t$,
  we set $\D \cX_r \dfnn X^{t,x,\mu,\nu}_r \neg - \neg X^{t,x,\mu',\nu}_r$, $\fa r \in [t,T]$.   By  \eqref{b_si_Lip} and \eqref{b_si_Lip_u},
  \ \beas
       \underset{r \in [  t,s]}{\sup} \big| \D \cX_r \big|
  \neg \le  \neg  \g  \neg \int_t^s  \neg \Big( | \D \cX_r |
   \neg + \neg  \rho^\l_{\overset{}{\hU}} \big( \mu_r, \mu'_r \big) \Big) dr
    \neg + \dneg  \underset{r \in [t,s]}{\sup} \bigg|  \neg \int_t^r  \dneg   \Big( \si \big(r,X^{t,x,\mu,\nu}_r , \mu_r,\nu_r\big)
    \neg - \neg  \si \big(r,X^{t,x,\mu',\nu}_r , \mu'_r,\nu_r\big) \Big) d B^t_r \bigg|, ~  \fa s  \neg \in \neg  [t,T] .
  \eeas
  Then  one can deduce from
    the Burkholder-Davis-Gundy inequality, \eqref{b_si_Lip}, \eqref{b_si_Lip_u} and H\"older's inequality   that
   \beas
~ && \hspace{-1cm} E_t \bigg[  \underset{r \in [  t,s]}{\sup}  | \D \cX_r  |^\varpi \bigg]
  \le    c_\varpi E_t \Bigg[\bigg( \neg \int_t^s | \D \cX_r | dr \bigg)^\varpi  \neg + \neg  \bigg( \neg \int_t^s \rho^\l_{\overset{}{\hU}} \big( \mu_r, \mu'_r \big) dr \bigg)^\varpi \neg + \neg \bigg(
   \int_t^s  \neg   \Big( | \D \cX_r | \neg + \neg  \rho^\l_{\overset{}{\hU}} \big( \mu_r, \mu'_r \big) \Big)^2  d  r \bigg)^{\neg \varpi/2}\Bigg]     \\
  & &  \le    c_\varpi E_t \Bigg[\bigg( \neg \int_t^s | \D \cX_r | dr \bigg)^\varpi \neg + \neg \bigg( \neg \int_t^s \rho^{2\l}_{\overset{}{\hU}} \big( \mu_r, \mu'_r \big) dr \bigg)^{\neg \varpi/2} \neg + \neg
   \underset{r \in [  t,s]}{\sup}  | \D \cX_r  |^{\neg \varpi/2}  \bigg( \int_t^s  \neg   | \D \cX_r  |   d  r \bigg)^{\neg \varpi/2}\Bigg]     \\
   &&    \le     c_\varpi E_t \Bigg[\bigg( \neg \int_t^s | \D \cX_r | dr \bigg)^\varpi \neg + \neg \bigg( \neg \int_t^s \rho^{2 \l}_{\overset{}{\hU}} \big( \mu_r, \mu'_r \big) dr \bigg)^{\neg \varpi/2 } \Bigg]  \neg +  \neg  \frac12 E_t \Bigg[
   \underset{r \in [  t,s]}{\sup}  | \D \cX_r  |^{\neg \varpi }   \Bigg]
    , \q  \fa s  \neg \in \neg  [t,T] .
  \eeas
  Similar to \eqref{eq:s018}, it follows from  H\"older's inequality and  Fubini's Theorem that
        \beas
   ~ \;   E_t \bigg[  \underset{r \in [  t,s]}{\sup}  | \D \cX_r  |^\varpi \bigg]   \neg  \le  \neg
    c_\varpi     \neg \int_t^s E_t \bigg[\underset{r' \in [t,r]}{\sup}\big| \D \cX_{r'} \big|^\varpi \bigg]  dr   \neg +  \neg
     c_\varpi E_t \Bigg[\bigg( \neg \int_t^s \rho^{2 \l}_{\overset{}{\hU}} \big( \mu_r, \mu'_r \big) dr \bigg)^{\neg   \varpi/2 } \Bigg]
    , ~\;  \fa s  \neg \in \neg  [t,T] .
    \eeas
   Then    an application of    Gronwall's inequality yields \eqref{eq:s015}.
    Similarly, with \eqref{b_si_Lip_v} we can deduce \eqref{eq:s015b} for each  $    \nu \, '  \in    \cV^t$. \qed

  \begin{lemm}  \label{lem_dist_meas}
  Let $ \hM $ be a  separable metric space with metric $\rho_{\overset{}{\hM}}$. For any two $\hM-$valued, $\bF^t-$adapted
  \(resp. $\bF^t-$progressively measurable\) processes $ Y$, $Z$, the nonnegative-valued process $\rho_{\overset{}{\hM}}\big(Y,Z\big)$ is  also $\bF^t-$adapted \(resp.\; $\bF^t-$progressively measurable\).
  \end{lemm}

  \ss \no {\bf Proof:} Let $\{x_n\}_{n \in \hN}$ be the countable dense subset of $\hM$ and denote by $\sB(\hM)$
  the Borel$-\si-$field of $\hM$.  We first claim that for any $y,z \in \hM$ and $\l >0$,
   \bea \label{eq:s211}
   \rho_{\overset{}{\hM}}(y,z)  \neg < \neg  \l \hb{ if and only if there exist $n \neg  \in  \neg  \hN$ and $r  \neg \in  \neg  \hQ \cap (0,\l)$ such that $\rho_{\overset{}{\hM}}(y,x_n) \neg < \neg r$ and $\rho_{\overset{}{\hM}}(x_n,z)  \neg < \neg \l - r$. }
   \eea
 ``$\Longleftarrow$": This direction is obvious due to the triangle inequality. ``$\Longrightarrow$":
      If $\rho_{\overset{}{\hM}}(y,z)  \neg < \neg  \l $, we let $r$ be a positive rational number that is less than
       $ \frac12 \big(\l -\rho_{\overset{}{\hM}}(y,z)  \big)$. There exists an $n \in \hN$, such that $\rho_{\overset{}{\hM}}(y,x_n) < r$. By  the triangle inequality,
        \beas
         \rho_{\overset{}{\hM}}(x_n,z) \le \rho_{\overset{}{\hM}}(x_n,y) + \rho_{\overset{}{\hM}}(y,z)
         < r  +  \rho_{\overset{}{\hM}}(y,z)  < \l - r  .
        \eeas
        So we proved the claim \eqref{eq:s211}.

   Now, given  two $\hM-$valued, $\bF^t-$adapted processes $ Y$ and $Z$, for any $s \neg \in \neg [t,T]$
    and $\l \neg > \neg 0$,       \eqref{eq:s211} implies that
    \beas
    \big\{\o  \neg \in  \neg  \O^t  \neg :  \rho_{\overset{}{\hM}}\big(Y_s(\o),Z_s (\o)\big)  \neg   <  \neg  \l \big\}
   =   \underset{n \in \hN}{\cup}  \underset{r      \in       \hQ \cap (0,\l)}{\cup} \Big( \big\{\o \in \O^t  \neg  :
   Y_s(\o)  \neg \in  \neg  O_r(x_n)    \big\} \neg \cap  \neg
   \big\{\o  \neg \in  \neg  \O^t  \neg  :   Z_s(\o)  \neg \in \neg  O_{\l-r}(x_n)   \big\} \Big) \in \cF^t_s ,
    \eeas
   which shows $\rho_{\overset{}{\hM}}\big(Y,Z\big)$ is also $\bF^t-$adapted.

  If $ Y$ and $Z$ are further $\bF^t-$progressively measurable, then for any $s \neg \in \neg [t,T]$
    and $\l \neg > \neg 0$,   we see from   \eqref{eq:s211}  that
 \beas
   && \hspace{-0.8cm} \big\{(r,\o)   \neg \in  \neg  [t,s] \times \O^t  \neg :  \rho_{\overset{}{\hM}}\big(Y_r (\o),Z_r (\o)\big)  \neg   <  \neg  \l \big\}
  \\
  &&  =   \underset{n \in \hN}{\cup}  \underset{r      \in       \hQ \cap (0,\l)}{\cup}  \neg
   \Big( \big\{(r,\o)   \neg \in  \neg  [t,s]  \neg \times \neg  \O^t  \neg  :
   Y_r(\o)  \neg \in  \neg  O_r(x_n)    \big\} \neg \cap  \neg
   \big\{(r,\o)   \neg \in  \neg  [t,s]  \neg \times \neg  \O^t  \neg  :   Z_r(\o)  \neg \in \neg  O_{\l-r}(x_n)   \big\} \Big)
   \neg \in \neg  \sB([t,s])  \neg \otimes  \neg  \cF^t_s  .
    \eeas
    Namely, $\rho_{\overset{}{\hM}}\big(Y,Z\big)$ is $\bF^t-$progressively measurable as well. \qed

 \ss \no {\bf Proof of Lemma \ref{lem_estimate_Y}:}
   We set $\Th \dfnn (t,x, \mu, \nu)$.

        \ss  \no (1)  For any $x' \neg \in \neg  \hR^k $, let $\Th'  \neg \dfnn \neg  (t,x',\mu,\nu)$
         and $\D X  \neg \dfnn \neg   \wt{X}^{\Th'}  \neg - \neg  \wt{X}^\Th $.
   The measurability of $\big( f^\Th_T,  \D X   \big)$   and \eqref{f_Lip} show that
  \beas
  \q  f_\pm(s, \o,y,z)     \dfnn   \neg   f^\Th_T  (s,  \o, y,z    )
     \pm \g \big| \D X_s(\o)   \big|^{  2 / q} ,
     \q    \fa (s,\o,y,z  ) \in   [t,T] \times \O^t \times \hR \times \hR^d
    \eeas
    define two $\sP \big(  \bF^t  \big)      \otimes \sB(\hR)     \otimes \sB(\hR^d)/\sB(\hR)-$measurable functions
     that are Lipschitz continuous in $(y,z)$ with coefficient $\g$.
   We see from  H\"older inequality, \eqref{eq:s031} and \eqref{eq:esti_X_1}  that
\bea
      E_t    \bigg[            \Big(\int_t^T  \neg \big| f_\pm   (s, 0,0 )   \big|  ds \Big)^{\neg q} \bigg]
  \le   c_0 E_t      \bigg[    \Big(\int_t^T  \neg \big| f^\Th_T   (s, 0,0 )   \big|  ds \Big)^{\neg q}
         +      \underset{s \in [t,T]}{\sup}   \big| \D X_s \big|^2    \bigg]    <       \infty .   \label{eq:s041}
  \eea

    Fix $\e > 0$.   The  function
  $  \phi(  \fx) \dfnn \left( |\fx|^2+ \e   \right)^{1/q  }$,  $  \fx  \in   \hR^k$
   has the following derivatives: for any $i, j \in \{1, \cds \neg , k \}$
   \beas
    \pa_{i}  \, \phi(  \fx) \neg  = \neg  \frac{2}{q} \, \phi^{1-q} (  \fx) \, \fx_i  \q
\hb{and} \q  \pa^2_{i j}  \, \phi(  \fx)  \neg = \neg  \frac{2}{q} \, \phi^{1-q} (  \fx)
  \, \d_{i j} \neg + \neg \frac{4}{q^2} (1 \neg - \neg q)\phi^{1-2q} (  \fx) \, \fx_i \fx_j        .
   \eeas
  It is easy to estimate that
  \bea
 && |\fx|^{2/q}     \le \phi(  \fx) \le  |\fx|^{2/q} + \e^{1/q}, \q |D \phi(  \fx) |
  = \frac{2}{q} \, \phi^{1-q} (  \fx) \, |\fx|
  \le  \frac{2}{q} \, \big| \fx \big|^{\frac{2}{q}-1},  \q \fa \fx  \in \hR^k. \label{eq:s035a}
  \eea
  For any   $\fz \in \hR^{k \times d}$, since $trace \big(  D^2 \phi(  \fx) \fz \fz^T \big)
    = \frac{2}{q} \, \phi^{1-q} (  \fx) \,|\fz|^2
    \neg + \neg \frac{4}{q^2} (1 \neg - \neg q)\phi^{1-2q} (  \fx)  \,  \big| \fz^T \fx  \big|^2$, we also have
  \bea
         - \frac{2}{q} \,  |\fx|^{\frac{2}{q}-2}  \,|\fz|^2
    \le \frac{4}{q^2} (1 \neg - \neg q)|\fx|^{\frac{2}{q}-2}  \,|\fz|^2
    \le   trace \big(  D^2 \phi(  \fx) \fz \fz^T \big)
    \le \frac{2}{q} \,  |\fx|^{\frac{2}{q}-2}  \,|\fz|^2  ,  \q \fa \fx  \in \hR^k  .   \label{eq:s035b}
  \eea

 Let us define $\cF^t_T-$measurable random variables $\xi_\pm \neg \dfnn \neg h \big( \wt{X}^\Th_T \big) \pm \g \phi \big( \D X_T   \big)      $ as well as real$-$valued, $ \bF^t-$adapted  continuous processes
  \beas
  \ul{L}^\pm_s  \dfnn \ul{L}^\Th_s   \pm \g \phi \big(\D X_s  \big)
  \q  \hb{and}  \q
     \ol{L}^\pm_s  \dfnn \ol{L}^{\, \raisebox{-0.5ex}{\scriptsize $\Th$}}_s \pm \g \phi
      \big( \D X_s  \big),      \q   \fa    s \in [t,T] .
      \eeas
     Clearly,        $ \ul{L}^\pm_s   < \ol{L}^{\, \raisebox{-0.5ex}{\scriptsize $\pm$}}_s  $ for any $ s \in [t,T]$
     and $ \ul{L}^\pm_T  \neg \le \neg  \xi_\pm  \neg \le \neg  \ol{L}^\pm_T$, $P^t_0-$a.s.
      Since
      \bea   \label{eq:s037}
        E_t   \bigg[   \underset{s \in [t,T]}{\sup}    \big| \phi \big(\D X_s  \big) \big|^q       \bigg]
 \le c_0  E_t   \bigg[      \underset{s \in [t,T]}{\sup}   \big| \D X_s \big|^2  \bigg]  \neg + \neg  c_0 \e < \infty
 \eea
   by \eqref{eq:s035a} and    \eqref{eq:esti_X_1}, we see from
   \eqref{eq:s031} that
    $     \ul{L}^\pm ,\ol{L}^\pm \in \hC^q_{\bF^t} \big( [t, T] \big) $
  and   it follows    that
  $ \xi_\pm  \neg \in \neg  \hL^q(\cF^t_T)$. Then Theorem 4.1 of \cite{EHW_2011} shows that
   the DRBSDE$\big(P^t_0, \xi_\pm, f_\pm , \ul{L}^\pm ,    \ol{L}^{\, \raisebox{-0.5ex}{\scriptsize $\pm$}} \big)$
   admits a unique solution $\Big(Y^\pm,Z^\pm, \ul{K}^\pm, \ol{K}^\pm  \Big) \in \hG^q_{\ol{\bF}^t}([t,T])$.
 For any $(s,\o,y,z  ) \in   [t,T] \times \O^t \times \hR \times \hR^d $,   the $ 2/q-$H\"older continuity of $h$,
   \eqref{f_Lip} and \eqref{2l_growth}  imply that
   \if{0}
 \beas
    \big| f^{\Th'}_T (s,\o,y,z) - f^\Th_T (s,\o,y,z) \big|
    &=& \big| f  \big(s,\wt{X}^{\Th'}_s(\o),y,z, \mu_s(\o), \nu_s(\o)\big) - f  \big(s,\wt{X}^\Th_s(\o),y,z, \mu_s(\o), \nu_s(\o)\big)   \big| \\
    &\le& \g \big| \D X_s(\o)   \big|^{2/q} + \g  \big( \rho_{\overset{}{\hU}} \big( \mu_s (\o), \mu_s (\o) \big) \big)^{  2 / q}  \\
    \big| \ul{L}^{\Th'}_s(\o) -  \ul{L}^\Th_s(\o) \big| + \big|\ol{L}^{\, \raisebox{-0.5ex}{\scriptsize $\Th'$}}_s(\o) -  \ol{L}^{\, \raisebox{-0.5ex}{\scriptsize $\Th$}}_s(\o) \big| & \le & \g \big| \D X_s (\o) \big|^{2/q} \le
     \g \phi \big( \D X_s (\o) \big) \\
     \big| h \big( \wt{X}^{\Th'}_T (\o) \big) - h \big( \wt{X}^\Th_T (\o) \big) \big|  & \le & \g \big| \D X_T (\o) \big|^{2/q} \le
     \g \phi \big( \D X_T (\o) \big) .
    \eeas
    \fi
     \beas
     &&   \xi_- (\o) \le h \big( \wt{X}^{\Th'}_T (\o) \big) \le \xi_+ (\o) , \q  f_-  (s,\o,y,z) \le   f^{\Th'}_T (s,\o,y,z) \le f_+  (s,\o,y,z)  ,  \\
 && \qq \ul{L}^-_s(\o) \le   \ul{L}^{\Th'}_s(\o) \le \ul{L}^+_s(\o), \q \hb{and} \q \ol{L}^{\, \raisebox{-0.5ex}{\scriptsize $-$}}_s(\o)
   \le \ol{L}^{\, \raisebox{-0.5ex}{\scriptsize $\Th'$}}_s(\o) \le \ol{L}^{\, \raisebox{-0.5ex}{\scriptsize $+$}}_s(\o) .
     \eeas
     Clearly, there inequalities also hold if $ \Th' $ is replaced by $\Th$.
        Proposition \ref{prop_comp_DRBSDE} then yields that $P^t_0-$a.s.
      \bea \label{eq:s057}
       Y^-_s \le  Y^{\Th'}_s  \big(T, h\big(\wt{X}^{\Th'}_T \big) \big)  \le Y^+_s ,
       \q \hb{and} \q
       Y^-_s \le  Y^{\Th}_s  \big(T, h\big(\wt{X}^{\Th}_T \big) \big)  \le Y^+_s  , \q \fa s \in [t,T] .
      \eea

  \ss
  By \eqref{eq:s037}, the processes $ \wh{Y}^\pm_s \dfnn  Y^\pm_{s } \mp \g \phi \big(\D X_{s }  \big) $, $s \in [t,T]$ are of $  \hC^q_{\bF^t} \big( [t, T] \big)$.   Applying It\^o's formula      yields that
    \bea
     \wh{Y}^\pm_s
     &  \tneg  = &  \tneg  \xi_\pm \neg + \neg  \int_s^T  \neg  \Big[  f_\pm  (r,   Y^\pm_r, Z^\pm_r)  \pm \g \, \big(( D \phi )  \big( \D X_r \big) \big)^T       \D b_r  \pm \frac 12 \g \, trace \Big( D^2 \phi \big( \D X_r  \big)    \D \si_r  ( \D \si_r  )^T \Big) \Big]  \, dr   \nonumber \\
   &  \tneg &   \tneg       \neg+   \ul{K}^\pm_T  \neg-\neg \ul{K}^\pm_{ s  }
  \neg- \neg \big( \ol{K}^\pm_T \neg - \neg \ol{K}^\pm_{ s  } \big)
   \neg-\neg \int_s^T   \neg \wh{Z}^\pm_r d B^t_r  , \q    s \in [t,T] ,   \label{eq:s063}
    \eea
     where $\D b_r \dfnn b \big(r, \wt{X}^{\Th'}_r, \mu_r,  \nu_r \big)
     - b \big(r, \wt{X}^{ \Th }_r,  \mu_r,  \nu_r \big) $,
      $\D \si_r \dfnn \si \big(r, \wt{X}^{\Th'}_r, \mu_r,  \nu_r \big)
       - \si \big(r, \wt{X}^{ \Th }_r,  \mu_r,  \nu_r \big) $   
       and  $ \wh{Z}^\pm_r \dfnn   \Big( Z^\pm_r \mp \g \, \big(( D \phi )  \big( \D X_r \big) \big)^T
       \D \si_r \Big) $.
       To wit,
       $ \Big( \wh{Y}^\pm  , \wh{Z}^\pm  , \ul{K}^\pm, \ol{K}^\pm \Big)$
          solves the DRBSDE$\Big(P^t_0, \xi_\pm  , \hat{f}^\pm , \ul{L}^\Th ,
         \ol{L}^{\, \raisebox{-0.5ex}{\scriptsize $\Th$}} \Big)$
       with
     \beas
     \hat{f}^\pm (s, \o,y,z) & \dneg \dneg \dfnn & \dneg \dneg     f_\pm \Big(s, \o,y \pm \g \phi \big( \D X_s (\o)  \big),
      z \pm \g \, \big(( D \phi )   \big( \D X_s (\o) \big) \big)^T  \D \si_s (\o) \Big)
       \pm \g \, \big(( D \phi )  \big( \D X_r (\o) \big) \big)^T       \D b_r (\o) \\
  & \dneg \dneg & \dneg \dneg \qq \; \; \q \pm \frac 12 \g \, trace \Big( D^2 \phi \big( \D X_r (\o) \big)  \D \si_r (\o) ( \D \si_r (\o) )^T \Big)  ,     \q   \fa (s,\o,y,z  ) \in   [t,T] \times \O^t \times \hR \times \hR^d   .
     \eeas

 The measurability of
    $\big(  f_\pm,  b, \si ,  \wt{X}^{\Th'}, \wt{X}^\Th  ,  \mu, \nu  \big)$
   and the Lipschitz continuity of $f_\pm$ in $(y,z)$ imply that $ \hat{f}^\pm $ are also
     $\sP \big(  \bF^t  \big)   \neg   \otimes     \neg    \sB(\hR)    \neg   \otimes   \neg    \sB(\hR^d)/\sB(\hR)-$measurable functions
     that are Lipschitz continuous in $(y,z)$ with coefficient $\g$.
     Then we can deduce from \eqref{eq:s035b}, \eqref{eq:s035a}, \eqref{b_si_Lip}      that
           \bea
    \big| \hat{f}^\pm   (s, 0,0 )   \big|
   & \dneg \tneg \le & \tneg \dneg  \big| f_\pm   (s, 0,0 )   \big|  \neg  +  \neg  \g^2     \phi \big( \D X_s    \big)
     \neg + \neg   \big| \big(( D \phi )   \big( \D X_s   \big) \big) \big| \big( \g^2 | \D \si_s |
      \neg + \neg  \g | \D b_s | \big)
     \neg + \neg      \frac{ \g }{q} \,     \big|  \D X_r \big|^{\frac{2}{q}-2}  \,| \D \si_r |^2 \nonumber \\
    & \dneg \tneg  \le & \tneg \dneg \big| f_\pm   (s, 0,0 )   \big|     \neg + \neg  c_0 \big| \D X_s\big|^{2/q}
    \neg + \neg  c_0 \e^{1/q}  , \q  s \in [t,T] .  \label{eq:s043}
   \eea
   \if{0}
   It follows from \eqref{eq:s041} and \eqref{eq:esti_X_1}  that
\beas
    E_t   \bigg[         \Big(\int_t^T  \neg \big| \hat{f}^\pm   (s, 0,0 )   \big| \, ds \Big)^q\bigg]
      & \dneg  \dneg \le & \dneg  \dneg
        c_0  E_t   \bigg[         \Big(\int_t^T  \neg \big| f_\pm   (s, 0,0 )   \big| \, ds \Big)^q\bigg]
       \neg  + \neg  c_0   E_t   \bigg[      \underset{s \in [t,T]}{\sup}   \big| \D X_s \big|^2  \bigg]
        \neg  + \neg  c_0 (1  \neg  + \neg \e)  <  \infty .    
  \eeas
    Hence, the DRBSDE$\Big(P^t_0, \xi_\pm  , \hat{f}^\pm , \ul{L}^\Th ,
         \ol{L}^{\, \raisebox{-0.5ex}{\scriptsize $\Th$}} \Big)$ admits a unique solution in
          $ \hG^q_{\ol{\bF}^t}([t,T])$ by Theorem 4.1 of \cite{EHW_2011}.
   \fi
    Similarly, 
    one can deduce from \eqref{eq:s035a} \eqref{b_si_Lip},
      and \eqref{eq:esti_X_1}  that
  \beas
  && \hspace{-0.8 cm}   E_t   \bigg[         \Big(\int_t^T  \neg \big| \big(( D \phi )  \big( \D X_s \big) \big)^T
       \D \si_s   \big|^2 \, ds \Big)^{\neg q/2}\bigg]    \le
      c_0  E_t   \bigg[         \Big(\int_t^T  \neg \big|  \D X_s \big|^{4/q  }   \, ds \Big)^{\neg q/2}\bigg]
        \neg   \le   \neg    c_0   E_t   \bigg[      \underset{s \in [t,T]}{\sup}   \big| \D X_s \big|^2  \bigg]
           \neg < \neg  \infty ,
  \eeas
  which shows that $ \wh{Z}^\pm  \in \hH^{2,q}_\fF([t,T], \hR^d )$.
 Then  Proposition \ref{prop_apriori_DRBSDE} implies that
   \bea  \label{eq:s061}
     E_t \bigg[ \underset{s \in [t,T]}{\sup} \big|   \wh{Y}^+_s - \wh{Y}^-_s \big|^\varpi  \bigg]
        \le c_\varpi     \Bigg\{ E_t \big[  |   \xi_+ - \xi_-  |^\varpi \big]
         + E_t \bigg[    \bigg( \int_t^T   \neg   \big|   \hat{f}^+   (r, \wh{Y}^-_r ,  \wh{Z}^-_r  )
         - \hat{f}^-   (r, \wh{Y}^-_r ,  \wh{Z}^-_r  )  \big| dr \bigg)^{\neg \varpi} \, \bigg] \Bigg\}   .
  \eea
 Since
  \beas
  &&  \hspace{-0.8cm}  \hat{f}^+   (s, \wh{Y}^-_s ,  \wh{Z}^-_s  )
         - \hat{f}^-   (s, \wh{Y}^-_s ,  \wh{Z}^-_s  )   =   f^\Th_T \Big(s, \o,\wh{Y}^-_s  + \g \phi \big( \D X_s   \big),
      \wh{Z}^-_s + \g \, \big(( D \phi )   \big( \D X_s  \big) \big)^T  \D \si_s  \Big) \\
     &&  -      f^\Th_T \Big(s, \o,\wh{Y}^-_s  - \g \phi \big( \D X_s   \big),
  \wh{Z}^-_s - \g \, \big(( D \phi )   \big( \D X_s  \big) \big)^T  \D \si_s  \Big) + 2 \g \big| \D X_s \big|^{  2 / q} \\
  &  &     +  2 \g \, \big(( D \phi )  \big( \D X_s  \big) \big)^T       \D b_s
  +  \g \, trace \Big( D^2 \phi \big( \D X_s  \big)  \D \si_s  ( \D \si_s  )^T \Big)
   ,\q s \in [t,T]  ,
     \eeas
  similar to \eqref{eq:s043},    the Lipschitz continuity of $f^\Th_T$,
  \eqref{eq:s035b}, \eqref{eq:s035a} and \eqref{b_si_Lip} imply that
    \beas  
   \big| \hat{f}^+   (s, \wh{Y}^-_s ,  \wh{Z}^-_s  )
         - \hat{f}^-   (s, \wh{Y}^-_s ,  \wh{Z}^-_s  ) \big|    \le       c_0 \big| \D X_s\big|^{2/q}
    \neg + \neg  c_0 \e^{1/q}
   ,\q s \in [t,T]  .
     \eeas
 Putting this back into \eqref{eq:s061}, we can deduce from \eqref{eq:s057}, \eqref{eq:s035a} and \eqref{eq:esti_X_3} that
 \beas
  \qq  && \hspace{-2.5cm} E_t \bigg[ \underset{s \in [t,T]}{\sup} \big|   Y^{\Th'}_s  \big(T, h\big(\wt{X}^{\Th'}_T \big) \big)
    - Y^{\Th}_s  \big(T, h\big(\wt{X}^{\Th}_T \big) \big) \big|^\varpi  \bigg]
   \le  E_t \bigg[ \underset{s \in [t,T]}{\sup} \big|   Y^+_s - Y^-_s \big|^\varpi  \bigg]
   \le  c_\varpi E_t \bigg[ \underset{s \in [t,T]}{\sup} \big|   \wh{Y}^+_s - \wh{Y}^-_s \big|^\varpi  \bigg]  \\
    &&    + c_\varpi E_t \bigg[ \underset{s \in [t,T]}{\sup} \big|   \phi (\D X_s ) \big|^\varpi  \bigg]
     \le c_\varpi     \bigg\{
           E_t \bigg[    \underset{s \in [t,T]}{\sup} \big|   \D X_s \big|^{\neg \frac{2\varpi}{q}}
             \bigg] + \e^{ \frac{ \neg \varpi}{q} } \bigg\}
         \le c_\varpi \big( |x'-x|^{\frac{2\varpi}{q}} + \e^{ \frac{ \neg \varpi}{q} } \big) .
 \eeas
 Then letting $\e \to 0$ yields \eqref{eq:s025}.

 \ss \no (2)   Next, we assume that $\ul{l}, \ol{l}$ and $h$   satisfy \eqref{2l_growth2}, that
 $b, \si $ are $ \l - $H\"older continuous in $u$, and that  $f$ is $ 2 \l - $H\"older continuous in $u$
  for some $\l \in (0,1/q]$.
 For any $\mu' \in   \cU^t  $, let $\Th^* \dfnn (t,x, \mu', \nu)$ and   $ \D \cX \dfnn  \wt{X}^{\Th^*} - \wt{X}^\Th.$
  The measurability of $\big( f^\Th_T,  \D \cX ,  \mu', \mu  \big)$
      together with  Lemma \ref{lem_dist_meas} and \eqref{f_Lip} shows that
  \beas
  \q  \ff_\pm(s, \o,y,z)     \dfnn   \neg   f^\Th_T  (s,  \o, y,z    )
     \pm \g \big| \D \cX_s(\o)   \big|^{  2 / q}
    \neg   \pm \g   \,   \rho^{  2 / q}_{\overset{}{\hU}} \big(  \mu'_s (\o), \mu_s (\o) \big)     ,
    ~ \;\;   \fa (s,\o,y,z  ) \in   [t,T] \times \O^t \times \hR \times \hR^d
    \eeas
    define two $\sP \big(  \bF^t  \big)      \otimes \sB(\hR)     \otimes \sB(\hR^d)/\sB(\hR)-$measurable functions
     that are Lipschitz continuous in $(y,z)$ with coefficient $\g$.
   Similar to \eqref{eq:s041},  H\"older inequality, \eqref{eq:s031} and \eqref{eq:esti_X_1} imply that
      \beas
      E_t    \bigg[            \Big(\int_t^T  \neg \big| \ff_\pm   (s, 0,0 )   \big|  ds \Big)^{\neg q} \bigg]
  \le   c_0 E_t      \bigg[    \Big(\int_t^T  \neg \big| f^\Th_T   (s, 0,0 )   \big|  ds \Big)^{\neg q}
         +      \underset{s \in [t,T]}{\sup}   \big| \D \cX_s \big|^2   +    \int_t^T   \dneg  \Big( [ \mu'_s ]^2_{\overset{}{\hU}}    +    \big[\mu_s\big]^2_{\overset{}{\hU}}   \Big) ds   \bigg]       <       \infty .  
  \eeas

      The function $\ol{\psi}(\fx)= \psi  (|x| )$, $\fx  \in \hR^k$ has the following  derivatives: for any $i, j \in \{1, \cds \neg , k \}$
   \beas
    \pa_{i}  \, \ol{\psi}(  \fx) & \tneg  = & \tneg  \b1_{\{|\fx| < R_1 \}} \fx_i +
    \b1_{\{|\fx| \in [ R_1, R_2] \}} \psi '  (|\fx| ) |\fx|^{-1} \fx_i + \b1_{\{|\fx| > R_2 \}} \frac{2}{q}  |\fx|^{\frac{2}{q}-2}  \fx_i \\
\hb{and} \q  \pa^2_{i j}  \, \ol{\psi}(  \fx)  & \tneg  = & \tneg \b1_{\{|\fx| < R_1 \}} \d_{ij}   +
    \b1_{\{|\fx| \in [ R_1, R_2] \}} \big(\psi '  (|\fx| ) |\fx|^{-1} \d_{ij} -
    |\fx|^{-3}\fx_i\fx_j+ \psi ''  (|\fx| ) | \fx |^{-2} \fx_i \fx_j \big) \\
     & \tneg & + \b1_{\{|\fx| > R_2 \}} \big( \frac{2}{q}  |\fx|^{\frac{2}{q}-2}  \d_{ij} + \frac{4}{q^2} (1 \neg - \neg q)  |\fx|^{\frac{2}{q}-4} \, \fx_i \fx_j \big)    .
  \eeas
   We can estimate   that
  \bea
       \big|D \ol{\psi} (  \fx) \big|  & \dneg  \dneg =& \dneg  \dneg  \b1_{\{|\fx| < R_1 \}} |\fx| +
    \b1_{\{|\fx| \in [ R_1, R_2] \}} \psi '  (|\fx| )   + \b1_{\{|\fx| > R_2 \}} \frac{2}{q}  |\fx|^{\frac{2}{q}-1}
       \le \k_\psi + \frac{2}{q}  |\fx|^{\frac{2}{q}-1} ,  \q \fa  \fx \in \hR^k .  \label{eq:s135a}
       \eea
       For any  $\fz \in \hR^{k \times d}$, since
       \beas
      trace \big(  D^2 \ol{\psi} (  \fx) \fz \fz^T \big)
  & \dneg  \dneg =& \dneg  \dneg  \b1_{\{|\fx| < R_1 \}} |\fz|^2   +
    \b1_{\{|\fx| \in [ R_1, R_2] \}} \big(\psi '  (|\fx| ) |\fx|^{-1} |\fz|^2 -
    |\fx|^{-3} \big| \fz^T \fx  \big|^2 + \psi ''  (|\fx| ) | \fx |^{-2} \big| \fz^T \fx  \big|^2 \big) \qq \qq   \nonumber \\
     &  \dneg  \dneg  & + \b1_{\{|\fx| > R_2 \}} \big( \frac{2}{q}  |\fx|^{\frac{2}{q}-2}  |\fz|^2 + \frac{4}{q^2} (1 \neg - \neg q)  |\fx|^{\frac{2}{q}-4} \, \big| \fz^T \fx  \big|^2 \big)  ,
   \eeas
   similar to \eqref{eq:s035b}, one can show that
   \bea
  \big|trace \big(  D^2 \ol{\psi} (  \fx) \fz \fz^T \big)\big| & \dneg  \tneg  \le    &  \tneg  \dneg
  \b1_{\{|\fx| =0\}} |\fz|^2 \neg +  \neg \b1_{\{0 < |\fx| \le  R_2 \}}  \big( 1
  \neg  +  \neg  R^{-1}_1 \neg \underset{\l \in [R_1,R_2]}{\sup } 1 \vee \psi'(\l)   \neg + \neg \underset{\l \in [R_1,R_2]}{\sup } |\psi''(\l) | \big) |\fz|^2
     \neg + \neg  \b1_{\{|\fx| > R_2 \}}   \frac{2}{q}  |\fx|^{\frac{2}{q}-2}  |\fz|^2  \nonumber \\
     & \tneg  \dneg  \le  & \tneg  \dneg \b1_{\{|\fx| =0\}} |\fz|^2 + \b1_{\{|\fx| > 0\}} \k_\psi \big( 1 \land |\fx|^{\frac{2}{q}-2} \big) |\fz|^2  ,  \q \fa  \fx \in \hR^k   .   \label{eq:s135b}
  \eea

 Let us define $\cF^t_T-$measurable random variables $\eta_\pm \neg \dfnn \neg h \big( \wt{X}^\Th_T \big) \pm \g \ol{\psi} \big( \D \cX_T   \big)      $ as well as real$-$valued, $ \bF^t-$adapted  continuous processes
  \beas
  \ul{\sL}^\pm_s  \dfnn \ul{ L}^\Th_s   \pm \g \ol{\psi} \big(\D \cX_s  \big)
  \q  \hb{and}  \q
     \ol{\sL}^\pm_s  \dfnn \ol{ L}^{\, \raisebox{-0.5ex}{\scriptsize $\Th$}}_s \pm \g \ol{\psi}
      \big( \D \cX_s  \big),      \q   \fa    s \in [t,T] .
      \eeas
     Clearly,        $ \ul{\sL}^\pm_s   < \ol{\sL}^{\, \raisebox{-0.5ex}{\scriptsize $\pm$}}_s  $ for any $ s \in [t,T]$
     and $ \ul{\sL}^\pm_T  \neg \le \neg  \eta_\pm  \neg \le \neg  \ol{\sL}^\pm_T$, $P^t_0-$a.s.
      Since $\psi(\l) \le \l^{2/q}$ for any $\l \ge 0$, we see from      \eqref{eq:esti_X_1} that
      \bea   \label{eq:s037b}
        E_t   \bigg[   \underset{s \in [t,T]}{\sup}    \big| \ol{\psi} \big(\D \cX_s  \big) \big|^q       \bigg]
    \neg \le \neg  E_t   \bigg[      \underset{s \in [t,T]}{\sup}    |   \D \cX_s  |^2     \bigg]  < \infty .
 \eea
   Thus,    $     \ul{\sL}^\pm ,\ol{\sL}^\pm \in \hC^q_{\bF^t} \big( [t, T] \big) $ by \eqref{eq:s031},
  and   it follows    that
  $ \eta_\pm  \neg \in \neg  \hL^q(\cF^t_T)$. Then Theorem 4.1 of \cite{EHW_2011} shows that
   the DRBSDE$\big(P^t_0, \eta_\pm, \ff_\pm , \ul{\sL}^\pm ,    \ol{\sL}^{\, \raisebox{-0.5ex}{\scriptsize $\pm$}} \big)$
   admits a unique solution $\Big(\sY^\pm,\sZ^\pm, \ul{\sK}^\pm, \ol{\sK}^\pm  \Big) \in \hG^q_{\ol{\bF}^t}([t,T])$.
 For any $(s,\o,y,z  ) \in   [t,T] \times \O^t \times \hR \times \hR^d $,
 \eqref{2l_growth2}, \eqref{f_Lip} and \eqref{f_Lip_u}    imply that
        \beas
     &&   \eta_- (\o) \le h \big( \wt{X}^{\Th^*}_T (\o) \big) \le \eta_+ (\o) , \q  \ff_-  (s,\o,y,z) \le   f^{\Th^*}_T (s,\o,y,z) \le \ff_+  (s,\o,y,z)  ,  \\
 && \qq \ul{\sL}^-_s(\o) \le   \ul{\sL}^{\Th^*}_s(\o) \le \ul{\sL}^+_s(\o), \q \hb{and} \q \ol{\sL}^{\, \raisebox{-0.5ex}{\scriptsize $-$}}_s(\o)
   \le \ol{\sL}^{\, \raisebox{-0.5ex}{\scriptsize $\Th^*$}}_s(\o) \le \ol{\sL}^{\, \raisebox{-0.5ex}{\scriptsize $+$}}_s(\o) .
     \eeas
     Clearly, there inequalities also hold if $ \Th^* $ is replaced by $\Th$.
        Proposition \ref{prop_comp_DRBSDE} then yields that $P^t_0-$a.s.
      \bea \label{eq:s057b}
       \sY^-_s \le  \sY^{\Th^*}_s  \big(T, h\big(\wt{X}^{\Th^*}_T \big) \big)  \le \sY^+_s ,
       \q \hb{and} \q
       \sY^-_s \le  \sY^{\Th}_s  \big(T, h\big(\wt{X}^{\Th}_T \big) \big)  \le \sY^+_s  , \q \fa s \in [t,T] .
      \eea

  \ss
  By \eqref{eq:s037b}, the processes $ \wh{\sY}^\pm_s \neg \dfnn  \neg  \sY^\pm_{s }    \mp    \g \ol{\psi} \big(\D \cX_{s }  \big) $, $s  \neg \in \neg  [t,T]$ are of $  \hC^q_{\bF^t}  \neg  \big( [t, T] \big)$. Let
  \beas
 && \D \wt{b}_s \dfnn b \big(r, \wt{X}^{\Th^*}_s, \mu'_s,  \nu_s \big)
     - b \big(r, \wt{X}^{ \Th }_s,  \mu_s,  \nu_s \big) , \q
      \D \wt{\si}_s \dfnn \si \big(r, \wt{X}^{\Th^*}_s, \mu'_s,  \nu_s \big)
       - \si \big(r, \wt{X}^{ \Th }_s,  \mu_s,  \nu_s \big) \\
    &&  \qq  \hb{and }  \q   \wh{\sZ}^\pm_s \dfnn   \Big( \sZ^\pm_s \mp \g \, \big(( D \ol{\psi} )  \big( \D \cX_s \big) \big)^T      \D \wt{\si}_s \Big) , \q \fa s \in [t,T].
       \eeas
    Similar to \eqref{eq:s063},  It\^o's formula implies that  $ \Big( \wh{\sY}^\pm  , \wh{\sZ}^\pm  , \ul{\sK}^\pm, \ol{\sK}^\pm \Big)$
          solves the DRBSDE$\Big(P^t_0, \eta_\pm  , \hat{\ff}^\pm , \ul{\sL}^\Th ,
         \ol{\sL}^{\, \raisebox{-0.5ex}{\scriptsize $\Th$}} \Big)$       with
     \beas
     \hat{\ff}^\pm (s, \o,y,z) & \dneg \dneg \dfnn & \dneg \dneg     \ff_\pm \Big(s, \o,y \pm \g \ol{\psi} \big( \D \cX_s (\o)  \big),
      z \pm \g \, \big(( D \ol{\psi} )   \big( \D \cX_s (\o) \big) \big)^T  \D \wt{\si}_s (\o) \Big)
       \pm \g \, \big(( D \ol{\psi} )  \big( \D \cX_r (\o) \big) \big)^T       \D \wt{b}_r (\o) \\
  & \dneg \dneg & \dneg \dneg \qq \; \; \q \pm \frac 12 \g \, trace \Big( D^2 \ol{\psi} \big( \D \cX_r (\o) \big)  \D \wt{\si}_r (\o) ( \D \wt{\si}_r (\o) )^T \Big)   ,
     \q   \fa (s,\o,y,z  ) \in   [t,T] \times \O^t \times \hR \times \hR^d   .
     \eeas

    The measurability of
    $\big(  \ff_\pm,  b, \si ,  \wt{X}^{\Th^*}, \wt{X}^\Th  , \mu', \mu, \nu  \big)$
   and the Lipschitz continuity of $\ff_\pm$ in $(y,z)$ imply that $ \hat{\ff}^\pm $ are also
     $\sP \big(  \bF^t  \big)   \neg   \otimes     \neg    \sB(\hR)    \neg   \otimes   \neg    \sB(\hR^d)/\sB(\hR)-$measurable functions
     that are Lipschitz continuous in $(y,z)$ with coefficient $\g$.
     Then we can deduce from  \eqref{eq:s135a}, \eqref{b_si_Lip}, \eqref{b_si_Lip_u} and \eqref{eq:s135b}    that
           \bea
      \big| \hat{\ff}^\pm   (s, 0,0 )   \big|
    & \dneg \dneg  \le    & \dneg \dneg     \big| \ff_\pm   (s, 0,0 )   \big|   + \g^2     \ol{\psi} \big( \D \cX_s    \big)   +  \big| \big(( D \ol{\psi} )   \big( \D \cX_s   \big) \big) \big| \big( \g^2 | \D \wt{\si}_s | + \g | \D \wt{b}_s | \big)
  \neg + \neg \frac 12 \g \, trace \Big( D^2 \ol{\psi} \big( \D \cX_r   \big)  \D \wt{\si}_r ( \D \wt{\si}_r )^T \Big) \nonumber \\
          & \dneg \dneg   \le     & \dneg \dneg    \big| \ff_\pm   (s, 0,0 )   \big|      \neg + \neg  \g^2 \big| \D \cX_s\big|^{2/q}
     \neg + \neg  c_0 \Big( \k_\psi + \big| \D \cX_s\big|^{\frac{2}{q} - 1} \Big)  \Big(  \big|  \D \cX_s\big|
    \neg  +  \neg            \rho^{\l}_{\overset{}{\hU}}(\mu'_s , \mu_s)     \Big) \nonumber \\
    & \dneg \dneg   &    + \neg  \frac12 \g     \, \Big( \b1_{\{|\D \cX_s| =0\}}   + \b1_{\{|\D \cX_s| > 0\}} \k_\psi     \big( 1 \land  |\D \cX_s|^{\frac{2}{q}-2} \big) \Big) \Big(  \big|  \D \cX_s\big|^2
    \neg  +  \neg            \rho^{2 \l}_{\overset{}{\hU}}(\mu'_s , \mu_s)     \Big)  \nonumber  \\
         & \dneg \dneg   \le    & \dneg \dneg       \big| \ff_\pm   (s, 0,0 )   \big|
          \neg   +  \neg   c_0   \k_\psi    \Big(  \big|  \D \cX_s\big|
    \neg  +  \neg            \rho^{\l}_{\overset{}{\hU}}(\mu'_s , \mu_s)  + \big| \D \cX_s\big|^{2/q}
     \neg  + \neg   \rho^{2 \l}_{\overset{}{\hU}}(\mu'_s , \mu_s)  \Big)  ,  \q    s \neg \in \neg [t,T] , \q \label{eq:s043b}
   \eea
   where we used the Young's inequality in the last step:
   \beas
    \big| \D \cX_s\big|^{\frac{2}{q} - 1} \rho^{\l}_{\overset{}{\hU}}(\mu'_s , \mu_s) \le
   c_0 \Big(  \big| \D \cX_s\big|^{2/q} + \rho^{\frac{2\l}{q}}_{\overset{}{\hU}}(\mu'_s , \mu_s) \Big)
   \le  c_0 \Big(  \big| \D \cX_s\big|^{2/q} + \rho^{\l}_{\overset{}{\hU}}(\mu'_s , \mu_s)
   + \rho^{2 \l}_{\overset{}{\hU}}(\mu'_s , \mu_s)  \Big)   .
   \eeas
    Similarly, 
     one can deduce from \eqref{eq:s135a} \eqref{b_si_Lip}, \eqref{b_si_Lip_u},
    H\"older's inequality and \eqref{eq:esti_X_1}  that
  \beas
  && \hspace{-0.6 cm}   E_t   \bigg[         \Big(\int_t^T  \neg \big| \big(( D \ol{\psi} )  \big( \D \cX_s \big) \big)^T
       \D \wt{\si}_s   \big|^2 \, ds \Big)^{\neg q/2}\bigg]    \le
       c_0   E_t   \bigg[         \Big(\int_t^T  \neg  \big( \k^2_\psi +  \big| \D \cX_s\big|^{\frac{4}{q} - 2} \big)
          \big(  \big|  \D \cX_s\big|^2
    \neg  +  \neg            \rho^{2\l}_{\overset{}{\hU}}(\mu'_s , \mu_s)     \big) \, ds \Big)^{\neg q/2}\bigg] \\
    &&  \le  c_0   E_t   \bigg[         \Big(  \k^q_\psi + \underset{s \in [t,T]}{\sup}  \big| \D \cX_s\big|^{ 2 - q } \Big)
          \bigg(  \underset{s \in [t,T]}{\sup}  \big| \D \cX_s\big|^q
    \neg  +  \neg       \Big(  \int_t^T \neg  \rho^2_{\overset{}{\hU}}(\mu'_s , \mu_s)     \big) \, ds \Big)^{\frac{\l q}{2}} \bigg) \bigg] \\
   & & \le \neg   c_0 \k^q_\psi \Bigg( \bigg\{ E_t \bigg[   \underset{s \in [t,T]}{\sup}   \big| \D \cX_s \big|^2
    \bigg] \bigg\}^{q/2  }  \neg + \neg  \bigg\{ E_t \neg \int_t^T \neg \rho^2_{\overset{}{\hU}}(\mu'_s , \mu_s)  \big)  ds    \bigg\}^{\frac{\l q}{2}} \Bigg)
    \neg  + \neg  c_0 E_t \bigg[  1 \neg + \neg   \underset{s \in [t,T]}{\sup}   \big| \D \cX_s \big|^2
   \neg + \neg \int_t^T  \neg  \rho^2_{\overset{}{\hU}}(\mu'_s , \mu_s)     \big)  ds  \bigg]  \neg < \neg  \infty   ,
  \eeas
  where we used the Young's inequality in the last step:
   \beas
   \underset{s \in [t,T]}{\sup}  \big| \D \cX_s\big|^{ 2 - q } \Big(  \int_t^T   \rho^2_{\overset{}{\hU}}(\mu'_s , \mu_s)     \big) \, ds \Big)^{\frac{\l q}{2}} \le \frac{(1 - \l)q}{2}+ \frac{2-q}{2} \underset{s \in [t,T]}{\sup}  \big| \D \cX_s\big| +
     \frac{\l q}{2}    \int_t^T   \rho^2_{\overset{}{\hU}}(\mu'_s , \mu_s)     \big) \, ds   .
   \eeas
 Thus, $ \wh{\sZ}^\pm  \in \hH^{2,q}_\fF([t,T], \hR^d )$.
 Then  Proposition \ref{prop_apriori_DRBSDE} implies that
   \bea \label{eq:s103}
     E_t \bigg[ \underset{s \in [t,T]}{\sup} \big|   \wh{\sY}^+_s - \wh{\sY}^-_s \big|^\varpi  \bigg]
        \le c_\varpi     \Bigg\{ E_t \big[  |   \eta_+ - \eta_-  |^\varpi \big]
         + E_t \bigg[    \bigg( \int_t^T   \neg   \big|   \hat{\ff}^+   (r, \wh{\sY}^-_r ,  \wh{\sZ}^-_r  )
         - \hat{\ff}^-   (r, \wh{\sY}^-_r ,  \wh{\sZ}^-_r  )  \big| dr \bigg)^{\neg \varpi} \, \bigg] \Bigg\}   .
  \eea
 Similar to \eqref{eq:s043b},    the Lipschitz continuity of $f^\Th_T$,
  \eqref{eq:s135a}, \eqref{b_si_Lip}, \eqref{b_si_Lip_u} and \eqref{eq:s135b} imply that
       \beas
      \big| \hat{f}^+   (s, \wh{\sY}^-_s ,  \wh{\sZ}^-_s  )
         - \hat{f}^-   (s, \wh{\sY}^-_s ,  \wh{\sZ}^-_s  ) \big|  &\dneg \tneg\le & \dneg \tneg
          2 \g^2  \ol{\psi} \big( \D \cX_s   \big)  + 2 \g^2
      \Big| \big(( D \ol{\psi} )   \big( \D \cX_s  \big) \big)^T  \D \wt{\si}_s    \Big| + 2 \g \big| \D \cX_s   \big|^{  2 / q} \\
      & \dneg \tneg&   \dneg \tneg               + 2 \g    \rho^{  2 \l }_{\overset{}{\hU}} \big( \mu'_s , \mu_s  \big)
       +  2 \g \, \Big| \big(( D \ol{\psi} )  \big( \D \cX_r  \big) \big)^T       \D \wt{b}_r \Big|
   +  \g \, trace \Big( D^2 \ol{\psi} \big( \D \cX_r  \big)  \D \wt{\si}_r  ( \D \wt{\si}_r  )^T \Big)   \\
 & \dneg \tneg   \le & \dneg \tneg      c_0   \k_\psi    \Big(  \big|  \D \cX_s\big|
    \neg  +  \neg            \rho^{ \l}_{\overset{}{\hU}}(\mu'_s , \mu_s)  + \big| \D \cX_s\big|^{2 / q}
     \neg  + \neg   \rho^{2 \l}_{\overset{}{\hU}}(\mu'_s , \mu_s)  \Big)
         ,\q s \in [t,T]  .
     \eeas
      Putting this back into \eqref{eq:s103}, we can deduce from \eqref{eq:s057b},
      \eqref{eq:s015}  and H\"older's inequality that
 \beas
    && \hspace{-0.7cm} E_t \bigg[ \underset{s \in [t,T]}{\sup} \big|   Y^{\Th^*}_s  \big(T, h\big(\wt{X}^{\Th^*}_T \big) \big)
    - Y^{\Th}_s  \big(T, h\big(\wt{X}^{\Th}_T \big) \big) \big|^\varpi  \bigg]
   \le  E_t \bigg[ \underset{s \in [t,T]}{\sup} \big|   \sY^+_s - \sY^-_s \big|^\varpi  \bigg]     \\
    && \le    c_\varpi E_t \bigg[ \underset{s \in [t,T]}{\sup} \big|   \wh{\sY}^+_s - \wh{\sY}^-_s \big|^\varpi  \bigg]
   + c_\varpi   E_t \bigg[ \underset{s \in [t,T]}{\sup} \big|   \ol{\psi} (\D \cX_s ) \big|^\varpi  \bigg]   \\
  &&  \le   c_\varpi  \k^\varpi_\psi   \Bigg\{   E_t \bigg[ \underset{s \in [t,T]}{\sup} \big| \D \cX_s \big|^{  \varpi }\bigg]
   \neg + \neg  E_t \bigg[    \underset{s \in [t,T]}{\sup} \big|   \D \cX_s \big|^{\neg \frac{2\varpi}{q}} \, \bigg]
   \neg +  \neg   E_t \bigg[ \Big( \int_t^T \rho^{2 \l}_{\overset{}{\hU}}(\mu'_s , \mu_s)    ds  \Big)^{\frac{\varpi}{2}}\bigg]
   \neg +  \neg  E_t \bigg[ \Big( \int_t^T \rho^{2 \l}_{\overset{}{\hU}}(\mu'_s , \mu_s)  ds  \Big)^{\varpi}\bigg] \Bigg\} \\
     &&  \le    c_\varpi  \k^\varpi_\psi   \Bigg\{   E_t \bigg[ \Big( \int_t^T \rho^{2 \l}_{\overset{}{\hU}}(\mu'_s , \mu_s)    ds  \Big)^{\frac{\varpi}{2}}\bigg] + E_t \bigg[ \Big( \int_t^T \rho^{2 \l}_{\overset{}{\hU}}(\mu'_s , \mu_s)    ds  \Big)^{\frac{\varpi}{q}}\bigg] + E_t \bigg[ \Big( \int_t^T \rho^{2 \l}_{\overset{}{\hU}}(\mu'_s , \mu_s)    ds  \Big)^{ \varpi } \bigg] \Bigg\} .
 \eeas
   Hence, \eqref{eq_estimate_Y_u} holds. The proof of \eqref{eq_estimate_Y_v} is similar. \qed

\subsection{Proof of Dynamic Programming Principle}

\label{subsection:DPP}

     \begin{lemm} \label{lem_ocset_Omega}
  Let $0 \le t \le s \le T$.
 For any $\o \in \O^t$ and $\d > 0$,
   $ O^s_\d (\o) \dfnn \Big\{\o' \in \O^t: \underset{r \in [t,s]}{\sup}
   \big| \o' (r) -   \o (r) \big| < \d \Big\}$ is a $\cF^t_s-$measurable  open  subset  of $\O^t$.

  \end{lemm}

  \ss \no {\bf Proof:}
 \if{0}
 (1)    Applying Lemma \ref{lem_truncation_measurable} with $r=S=s$
   and using \eqref{eq:xxc023}  yield that
  $\Pi^{T,s}_{t,t} (A) \in \cF^{t,s}_s =\sB(\O^{t,s})$. Then in light of  Proposition 15.11 of \cite{Royden_real},
   we can find   an open subset $\wt{O}$ of $\O^{t,s}$
  and a  closed subset $\wt{F}$ of $\O^{t,s}$ such that $ \wt{F} \subset  \Pi^{T,s}_{t,t} (A) \subset \wt{O} $,
    $ P^{t,s}_0 \big(\wt{O} \big\backslash \Pi^{T,s}_{t,t} (A) \big)  < \e $ and
  $ P^{t,s}_0 \big(   \Pi^{T,s}_{t,t} (A) \big\backslash \wt{F} \big) < \e$. Applying Lemma \ref{lem_shift_inverse}
  with $(s,S,r)=(t,s,s)$ shows that $ O \dfnn \big(\Pi^{T,s}_{t,t} \big)^{-1}(\wt{O})$
  is an $\cF^{t,T}_s-$measurable,  open subset of $\O^{t,T}$ and  $ F \dfnn \big(\Pi^{T,s}_{t,t} \big)^{-1}(\wt{F})$
  is an $\cF^{t,T}_s-$measurable,  closed subset of $\O^{t,T}$.
  Lemma \ref{lem_shift_inverse} also shows that
   \beas
 &&   \e >  P^{t,s}_0 \big(\wt{O} \big\backslash \Pi^{T,s}_{t,t} (A) \big)
    = P^{t,T}_0 \Big( \big(\Pi^{T,s}_{t,t}\big)^{-1} \big(   \wt{O} \big\backslash \Pi^{T,s}_{t,t} (A) \big)\Big)
    = P^{t,T}_0 \Big(  O  \Big\backslash  \big(\Pi^{T,s}_{t,t}\big)^{-1} \big( \Pi^{T,s}_{t,t} (A) \big)\Big)   \qq   \\
         \hb{and similarly}   &&    \e >
    P^{t,s}_0 \big(   \Pi^{T,s}_{t,t} (A) \big\backslash \wt{F} \big)
    =  P^{t,T}_0 \Big( \big(\Pi^{T,s}_{t,t}\big)^{-1} \big(   \Pi^{T,s}_{t,t} (A)  \big) \Big\backslash F \Big)  .
   \eeas
  As $F \subset \big(\Pi^{T,s}_{t,t}\big)^{-1} \big(   \Pi^{T,s}_{t,t} (A)  \big) \subset O $, it remains that show that
  $ \big(\Pi^{T,s}_{t,t}\big)^{-1} \big(   \Pi^{T,s}_{t,t} (A)  \big) = A$:
   Clearly, $A \subset \big(\Pi^{T,s}_{t,t}\big)^{-1} \big(   \Pi^{T,s}_{t,t} (A)  \big)$.
   For any $\o \in \big(\Pi^{T,s}_{t,t}\big)^{-1} \big(   \Pi^{T,s}_{t,t} (A)  \big)$,
   as $\Pi^{T,s}_{t,t} (\o) \in \Pi^{T,s}_{t,t} (A)$, there exists a $\o' \in A$ such that $\Pi^{T,s}_{t,t} (\o) =
   \Pi^{T,s}_{t,t} (\o')$. i.e. $\o(r) = \o' (r)$, $\fa r \in [t,s]$. As $A  \in \cF^{t,T}_s$,
   Lemma \ref{lem_element} implies that $\o \in  \o' \otimes_s \O^{s,T} \in A $. Thus
   $ \big(\Pi^{T,s}_{t,t}\big)^{-1} \big(   \Pi^{T,s}_{t,t} (A)  \big) = A$.
  \fi
 Let $\o \in \O^t$ and $\d > 0$. Given $\o' \in O^s_{\d} (\o)$,
   for any $\o'' \in O_{\d'}(\o')    $ with $\d' \dfnn \d - \underset{r \in [t,s]}{\sup}
   \big| \o' (r) -   \o (r) \big| $, one has
    \beas
    \underset{r \in [t,s]}{\sup}  \big| \o'' (r) -   \o (r) \big|
    \le  \| \o''   -   \o'  \|_t + \underset{r \in [t,s]}{\sup}  \big| \o' (r) -   \o (r) \big|
    < \d'  + \underset{r \in [t,s]}{\sup}  \big| \o' (r) -   \o (r) \big| = \d .
    \eeas
  Thus $ O_{\d'}(\o') \subset O^s_{\d} (\o) $, which shows that $ O^s_{\d} (\o) $ is an open subset of $ \O^t$.
  Moreover, since $\O^t$ is the set of $\hR^d-$valued continuous functions on $[t,T]$ starting from $0$,
  we see that
 \beas
  \hspace{1.75cm}  O^s_{\d} (\o) =  \underset{r \in \hQ_{t,s}  }{\cap} \big\{ \o' \in \O^t :  | \o' (r) - \o (r) | < \d \big\}
   =  \underset{r \in \hQ_{t,s} }{\cap} \big\{ \o' \in \O^t : B^t_r(\o') \in O_\d  \big( \o(r) \big)  \big\} \in \cF^t_s .
  \hspace{1.75cm}    \hb{\qed}
 \eeas

\no  Given $ t \in [0,T]$ and $\beta \in \fB^t $, we define
 \beas
 I (t, x,\beta )  \dfnn     \underset{\mu \in \cU^t}{\sup} \; \wt{Y}^{t,x,\mu, \beta \lan \mu \ran }_t  \Big(T,h\big(\wt{X}^{t,x,\mu, \beta \lan \mu \ran}_T \big)  \Big), \q \fa x \in \hR^k .
 \eeas
  Taking supremum over $\mu \in \cU^t$ in \eqref{eq:j701} and \eqref{eq:k115} implies that
  \bea \label{eq:I_finite}
  \ul{l}(t,x) \le I (t, x, \beta ) \le \ol{l}(t,x) , \q   \fa  x \in   \hR^k ,
  \eea
  and that
 \bea  \label{eq:p571}
     \hb{   the function  $x \to I (t, x, \beta )   $ is continuous.}
 \eea

 Similarly, for any $\a \in \cA^t$,
  $   I (t, x,\a )  \dfnn     \underset{\nu \in \cV^t}{\sup} \; \wt{Y}^{t,x,  \a \lan \nu \ran, \nu }_t
   \Big(T,h\big(\wt{X}^{t,x,\a \lan \nu \ran, \nu}_T \big)  \Big) \in \big(\ul{l}(t,x) ,  \ol{l}(t,x) \big)  $
  is continuous in $ x \in \hR^k $.

\ss \no {\bf Proof of Theorem \ref{thm_DPP}:}
    Let $ \{t_n \}_{n \in \hN}$ denote  the countable set $\hQ_{t,T}  $ and
    let   $\{\t_{\mu,\beta} \neg : \mu \neg \in \neg \cU^{t}, \beta \neg \in  \neg \fB^{t} \}$  be a
    family of $\hQ_{t,T}  -$valued $\bF^{t}-$stopping times.

  \ss \no {\bf 1)} We fix    $\e >0$. For any  $n \in \hN$ and $x \in \hR^k$, since $ w_1(t_n,x) $
   is finite by \eqref{eq:n411},  there exists a  $ \beta^n_{x}  \in \fB^{t_n} $ 
     such that
   \bea  \label{eq:j231}
  w_1(t_n,x) \le  I(t_n,x,  \beta^n_{x})
   \le w_1(t_n,x)+ \frac13 \e.
   \eea
  The continuity of $  I(t_n,\cd,  \beta^n_{x}) $ by \eqref{eq:p571} and 
  Proposition  \ref{prop_w_conti}   assure that there exists  a $\d_n(x) > 0$ such that
   \bea   \label{eq:j233}
  \ba{c}
   \dis \big| I(t_n,x',  \beta^n_{x})  -  I(t_n,x,  \beta^n_{x}) \big| <     \frac13 \e ,
     \q \hb{and} \q    \big| w_1(t_n,x' ) - w_1(t_n,x)     \big| <  \frac13 \e  ,   \q \fa x' \in O_{\d_ n \neg (x)}(x).
     \ea
 \eea
 Given  $n \in \hN$,     Lindel\"of's covering theorem  (see e.g.  Theorem VIII.6.3 of  \cite{Dugundji_1966})
 shows that there exists    a sequence $\{x^n_i\}_{i \in \hN}$ of $\hR^k$ satisfying
   $ \dis \underset{i \in \hN}{\cup} \wt{O}^n_i =\hR^k$ with  $\wt{O}^n_i \dfnn O_{\d_n \neg (x^n_i) }(x^n_i)$.
  For any $i \in \hN$, let $\beta^n_i \dfnn  \beta^n_{x^n_i}$. We can deduce  from \eqref{eq:j231} and  \eqref{eq:j233} that
  for any $x' \in \wt{O}^n_i$
  \bea  \label{eq:j311}
  I(t_n,x',  \beta^n_i)  <  I(t_n,x^n_i,  \beta^n_i)  +  \frac13 \e  \le  w_1(t_n,x^n_i)+ \frac23 \e  < w_1(t_n,x' ) +\e .
  \eea

 \ss  Now,  fix $(\beta, \mu)  \in \fB^{t} \times \cU^{t}$. We  simply denote $\t_{\mu,\beta}$ by $\t$ and
 set $\Th = \big( t,x, \mu, \beta \lan \mu \ran  \big) $.
  For any   $ n, i \in \hN$,  define   
    \beas   
   A^n_i \dfnn \{\t=t_n\} \cap \big\{  \wt{X}^\Th_{t_n}  \in \wt{O}^n_i  \big\backslash  \underset{j < i }{\cup} \wt{O}^n_j \big\} \in  \cF^{t}_{t_n} \cap  \cF^{t}_\t       .
     \eeas
  Let  $m \in \hN$  and   $  A_m \dfnn \O^{t}   \Big\backslash   \underset{n,i =1 }{\overset{m}{\cup}} A^n_i \in \cF^t_\t$.
  Proposition   \ref{prop_paste_strategy} shows that
    \beas    
    \beta^m (r,\o,u )    \dfnn
        \begin{cases}
      \beta^n_i    \big( r, \Pi_{t,t_n}(\o),u \big)     , ~ & \hb{if $(r,\o ) \in
   \[\t,T\]_{A^n_i} = [t_n,T] \times  A^n_i$    for $ n,i \in \{1 \cds \neg , m\}$} ,   \\
     \beta (r,\o,u ),  & \hb{if $(r,\o )  \in \[t,\t\[ \,\cup \,\[\t, T\]_{A_m }  $   }
    \end{cases}        \q \fa u \in \hU
    \eeas
 defines a  $\fB^{t}-$strategy  such that it holds for $P^{t}_0-$a.s. $\o \in \O^{t}$ that
  for any $r \in \big[ \t(\o) , T \big] $
 \bea         \label{eq:j241}
    \big( \beta^m \lan \mu \ran   \big)^{\t,\,\o}_r
    = \begin{cases}
  \big( \beta^n_i \lan \mu^{t_n,\o}  \ran   \big)_r,\q & \hb{if $ \o   \in  A^n_i$  for $ n,i \in \{1 \cds \neg , m\}$} , \\
     \big( \beta \lan \mu \ran \big)^{\t,\,\o}_r ,  & \hb{if $ \o   \in A_m $.   }
    \end{cases}
      \eea

    Let $   \Th_m \neg \dfnn \neg   \big( t , x , \mu ,  \beta^m \lan \mu \ran \big)$.
   As    $ \beta^m  \lan \mu \ran   \neg = \neg   \beta \lan \mu \ran         $ on $\[t,\t\[$\,,
    we see from \eqref{eq:p611} that $  P^{t}_0-$a.s.
   \bea    \label{eq:q101}
   \wt{X}^{\Th_m}_s =X^{\Th_m  }_s = X^\Th_s  =  \wt{X}^\Th_s , \q \fa s \in [t,\t] .
   \eea
   Hence,  for any $\cF^{t}_\t -$measurable random variable $\xi$
   with $\ul{L}^{\Th_m}_\t =  \ul{L}^\Th_\t \le \xi \le \ol{L}^\Th_\t = \ol{L}^{\Th_m}_\t$, $  P^{t}_0-$a.s.,
      the DRBSDE$ \big(P^t, \xi , \\ f^{\Th_m}_\t,  \ul{L}^{\Th_m}_{\t \land \cd},\ol{L}^{\Th_m}_{\t \land \cd}   \big)$
      and the DRBSDE$ \big(P^t, \xi , f^\Th_\t,  \ul{L}^\Th_{\t \land \cd},\ol{L}^\Th_{\t \land \cd}   \big)$
      are essentially the same. To wit, we have
   \bea  \label{eq:q104}
   \big(Y^{\Th_m} (\t, \xi),Z^{\Th_m} (\t, \xi),\ul{K}^{\Th_m} (\t, \xi),
    \ol{K}^{\, \raisebox{-0.5ex}{\scriptsize $\Th_m$}} (\t, \xi) \big)
      =  \big(Y^\Th (\t, \xi),Z^\Th (\t, \xi),\ul{K}^{\Th}(\t, \xi),
       \ol{K}^{\, \raisebox{-0.5ex}{\scriptsize $\Th$}} (\t, \xi) \big)    .
    \eea

 Since  $    \wt{X}^{\Th_m}_\t 
   =    \wt{X}^\Th_\t $,   $  P^{t}_0-$a.s. by \eqref{eq:q101}, it holds for $ P^{t}_0-$a.s. $ \o \in \O^{t}$ that
   \beas     
   (\Th_m)^\o_\t  = \Big( \t(\o),\wt{X}^\Th_{\t(\o)}(\o), \mu^{\t,\o}, \big( \beta^m \lan \mu \ran   \big)^{\t,\,\o} \Big) .
   \eeas
   Then applying Proposition \ref{prop_DRBSDE_shift} and Proposition \ref{prop_FSDE_shift}
  with $\Th=\Th_m$  and $\xi = h \big(\wt{X}^{\Th_m}_T\big)$, we can deduce from \eqref{eq:j241},
        Proposition \ref{prop_dwarf_strategy} (1),    \eqref{eq:j701} and \eqref{eq:j311} that
         \bea
      && \hspace{-0.7 cm}   \Big( \wt{Y}^{\Th_m}_\t \big(T, h \big(\wt{X}^{\Th_m}_T\big)\big)\Big) (\o)
       =    \wt{Y}^{ (\Th_m)^\o_\t   }_{\t(\o)}
          \Big(T, h \Big( \big( \wt{X}^{\Th_m} \big)^{\t,\o}_T \Big) \Big)
           =   \wt{Y}^{ (\Th_m)^\o_\t   }_{\t(\o)}
          \Big(T, h \Big(     \wt{X}^{(\Th_m)^\o_\t}_T    \Big) \Big)  \nonumber \\
    &  &   =        \wt{Y}^{\t(\o),\wt{X}^\Th_{\t(\o)} (\o), \mu^{\t,\o}, (\beta \lan \mu \ran )^{\t,\o}}_{\t(\o)}
  \bigg(T, h\Big( \wt{X}^{\t(\o),\wt{X}^\Th_{\t(\o)} (\o), \mu^{\t,\o}, (\beta \lan \mu \ran )^{\t,\o}}_T \Big)\bigg)
  \nonumber \\
      &&   \le \b1_{\{\o \in A_m\}} \ol{l} \big( \t(\o),\wt{X}^\Th_{\t(\o)} (\o) \big)
     \neg + \neg \sum^m_{n,i =1 } \b1_{  \{\o \in A^n_i \}} \wt{Y}^{t_n,\wt{X}^\Th_{t_n}(\o),\mu^{t_n,\o}, \beta^n_i \lan \mu^{t_n,\o}  \ran }_{t_n}    \Big(T, h\Big(  \wt{X}^{t_n,\wt{X}^\Th_{t_n}(\o), \mu^{t_n,\o}, \beta^n_i \lan \mu^{t_n,\o}  \ran }_T \, \Big)\Big) \qq \q \label{eq:t023}  \\
      && \le \b1_{\{\o \in A_m\}} \ol{l} \big( \t(\o),\wt{X}^\Th_{\t(\o)} (\o) \big)
    + \sum^m_{n,i =1 } \b1_{\{\o \in A^n_i\}}   I \big(t_n, \wt{X}^\Th_{t_n}(\o), \beta^n_i \big)
      \le        \xi_m (\o)   ,  \q \hb{for $P^{t}_0-$a.s. $ \o \in \O^{t}$}  ,             \label{eq:h031}
      \eea
    where       $ \xi_m   \dfnn \b1_{ A_m } \ol{l} \big( \t ,\wt{X}^\Th_\t    \big)  +
     \b1_{   A^c_m  }  \Big( \big(  w_1\big(\t , \wt{X}^\Th_{\t }  \big) + \e \big)
     \, \land \, \ol{l}\big(\t , \wt{X}^\Th_{\t }  \big) \Big)  $.

 \ss   The continuity of $\ul{l}$ and \eqref{eq:n333} show that $ \ol{l} \big( \t ,\wt{X}^\Th_\t    \big) ,
   w_1\big(\t , \wt{X}^\Th_{\t }  \big)   \in  \cF^{t}_\t $,
    so is $ \xi_m     $.
   Also,  we see from     \eqref{eq:n411}    that
   \beas
   \q \;\;\;   \ul{L}^{\Th}_\t       =   \ul{l}\big(\t , \wt{X}^\Th_{\t }  \big)  \le  w_1\big(\t , \wt{X}^\Th_{\t }  \big)
     \le \ol{l}\big(\t , \wt{X}^\Th_{\t }  \big) = \ol{L}^{\Th}_\t   ~\; \hb{ and } ~\;
        \ul{L}^{\Th}_\t =   \ul{l}\big(\t , \wt{X}^\Th_{\t }  \big)
   \le  \xi_m  \le \ol{l}\big(\t , \wt{X}^\Th_{\t }  \big) = \ol{L}^{\Th}_\t
    ,   \q    P^{t}_0-  a.s.
   \eeas
       Thus both $ Y^{\Th}  (\t, \xi_m ) $  and $ Y^{\Th}  \big(\t, w_1\big(\t , \wt{X}^\Th_\t  \big) \big) $ is well-posed.
  Then  Proposition \ref{prop_apriori_DRBSDE} implies that
   \bea     \label{eq:j335}
   \big|  \wt{Y}^{\Th}_{t}  (\t,  \xi_m  )
   \neg  - \neg  \wt{Y}^\Th_{t}  \big(\t,    w_1\big(\t , \wt{X}^\Th_{\t }  \big)   \big)  \big|
     \neg \le \neg c_0  \big\|   \xi_m      \neg  -  \neg
        w_1\big(\t , \wt{X}^\Th_{\t }  \big)     \big\|_{ \hL^q  (\cF^{t}_\t  ) }
   \neg  \le  \neg   c_0 \big\|   \b1_{ A_m } \neg \big(
    \ol{l} \big( \t ,\wt{X}^\Th_\t    \big) \neg  -  \neg  \ul{l} \big(\t , \wt{X}^\Th_{\t }  \big)    \big)   \big\|_{ \hL^q  (\cF^{t}_\t  ) }  \neg + \neg  c_0  \,  \e  .    \q~\;\;
   \eea
     Applying  \eqref{eq:p677} with $(\z,\t,\xi) = \Big(\t, T,h \big( \wt{X}^{ \Th_m }_T \big) \Big)$,   applying Proposition \ref{prop_comp_DRBSDE} to \eqref{eq:h031}, and
   using  \eqref{eq:q104} for   $\xi =\xi_m$ yield that
    \bea
      \wt{Y}^{\Th_m}_{t} \Big(T,h \big( \wt{X}^{ \Th_m }_T \big) \Big)
 & =& \wt{Y}^{\Th_m}_{t} \Big(\t,
    \wt{Y}^{\Th_m}_\t \big(T, h \big( \wt{X}^{ \Th_m }_T \big) \big)   \Big)
   \le \wt{Y}^{\Th_m}_{t}  \big( \t,         \xi_m      \big)
 =  \wt{Y}^\Th_{t} \big(\t,         \xi_m     \big) \nonumber \\
 &\le& \wt{Y}^\Th_{t}  \big(\t,    w_1\big(\t , \wt{X}^\Th_{\t }  \big)   \big)
 +    c_0 \big\|   \b1_{ A_m } \neg \big(
    \ol{l} \big( \t ,\wt{X}^\Th_\t    \big) \neg  -  \neg  \ul{l} \big(\t , \wt{X}^\Th_{\t }  \big)    \big)   \big\|_{ \hL^q  (\cF^{t}_\t  ) } + c_0  \,  \e .  \label{eq:h041}
    \eea

Since $ \dis \underset{n,i \in \hN}{\cup} A^n_i   =\O^{t}$,
  one has $ \lmtd{m \to \infty}  A_m   =\es$. By \eqref{2l_growth} and H\"older's inequality,
 \bea
      \big\|       \ol{l} \big(\t , \wt{X}^\Th_{\t }  \big) \neg  -  \neg \ul{l} \big( \t ,\wt{X}^\Th_\t    \big)        \big\|_{ \hL^q  (\cF^{t}_\t  ) }
       \le     \underset{s \in [t,T]}{\sup}|\ol{l}(s,0)|
     +    \underset{s \in [t,T]}{\sup}|\ul{l}(s,0)|
     +     2 \g  \big\| \wt{X}^\Th   \big\|_{ \hC^2_{\bF^{t}  }  ([t,T],\hR^k  ) }
 < \infty .   \label{eq:t528}
 \eea
  Hence, the Dominated Convergence Theorem shows that  $ \lmtd{m \to \infty} \big\|   \b1_{ A_m } \neg \big(   \ol{l} \big(\t , \wt{X}^\Th_{\t }  \big) \neg  -  \neg \ul{l} \big( \t ,\wt{X}^\Th_\t    \big)    \big)   \big\|_{ \hL^q  (\cF^{t}_\t  ) } =0$.
   Letting  $m \to \infty$  in \eqref{eq:h041}  gives that
   \beas
  \wt{Y}^{t,x,\mu, \beta^m \lan \mu \ran}_{t} \Big(T,h\big(X^{ t,x,\mu, \beta^m \lan \mu \ran }_T\big)\Big) \le
  \wt{Y}^{t,x,\mu, \beta \lan \mu \ran}_{t} \Big( \t_{\mu,\beta}, \, w_1\big(\t_{\mu,\beta} , X^{t,x,\mu, \beta \lan \mu \ran }_{\t_{\mu,\beta}}  \big) \Big)+ c_0  \,  \e .
  \eeas
 Taking supremum over $\mu \in \cU^{t}$, we obtain
  \beas
  w_1(t,x) 
  \le  I \big( t,x , \beta^m \big) \le  \underset{\mu \in \cU^{t}}{\sup} \wt{Y}^{t,x,\mu, \beta \lan \mu \ran}_{t} \Big({\t_{\mu,\beta}}, \,  w_1\big({\t_{\mu,\beta}}, X^{t,x,\mu, \beta \lan \mu \ran }_{\t_{\mu,\beta}} \big) \Big)+ c_0  \,  \e .
  \eeas
 Then taking infimum over $\beta \in \fB^{t}$ and letting  $\e \to 0$ yield \eqref{eq:DPP0}. Similarly, one has \eqref{eq:DPP1}.

 \ms \no {\bf 2)} Next, assume ($\bV_\l$) for some $\l \in (0,1)$. We shall show the inverse of  \eqref{eq:DPP1}.

  \ss \no {\bf a)}  Fix  $(\beta,\mu) \in \wh{\fB}^{t} \times   \cU^{t}$.
  We simply denote $\t_{\mu,\beta}$ by $\t $.  For some $\k >0$ and some  non-negative  measurable process $\Psi$ on $(\O^t,\cF^t_T) $
   with $ C_\Psi \dfnn E_t \neg \int_t^T \neg   \Psi^2_r  \,  dr < \infty$, it holds
     $dr \times d P^t_0-$a.s. that
     \beas
       \big[ \beta(r,\o,u)  \big]_{\overset{}{\hV}} \le \Psi_r (\o) + \k [u]_{\overset{}{\hU}} ,
      \q  \fa u \in \hU  .
       \eeas
   Similar to   \eqref{eq:t315},   applying Proposition \ref{prop_rcpd_L1} with
   $ \xi =  \int_\t^T    \Psi^2_r  dr \in \hL^1 (\cF^t_T) $, we can deduce that
   for all $\o \in \O^t$ except on a $P^t_0-$null set $\cN_0$
     \bea     \label{eq:t325}
         E_{ \t(\o) } \bigg[ \int_{\t(\o)}^T \big( \Psi^{\t,\o}_r\big)^2   \, dr \bigg]
     =  E_t \bigg[   \int_\t^T \Psi^2_r \, dr \Big| \cF^t_\t  \bigg] (\o)
    \le E_t \bigg[   \int_t^T \Psi^2_r \, dr \Big| \cF^t_\t  \bigg] (\o)  <\infty .
     \eea
         For any $n   \in    \hN$, similar to \eqref{eq:t311} and the conclusion that follows,
        it holds for all $\o    \in    \O^t$ except on a   $P^t_0-$a.s. null set $\cN^1_n$ that
        for $dr \times dP^{t_n}_0-$a.s. $(r, \wt{\o}) \in [t_n,T] \times \O^{t_n}$
     \bea  \label{eq:t329}
     \big[ \beta^{t_n,\o} \big(r, \wt{\o},u\big) \big]_{\overset{}{\hV}}
      \le  \Psi^{t_n,\o}_r  (   \wt{\o}  ) + \k [u]_{\overset{}{\hU}} \, , \q \fa u \in \hU .
   \eea
   Also,     Proposition \ref{prop_dwarf_strategy} (2) shows that  there exists another $P^t_0-$null set $ \cN^2_n$ such that
 $\beta^{t_n,\o} \in \wh{\fB}^{t_n}$ for any $\o \in \big( \cN^2_n \big)^c$.

  Now, let us set $\varpi \dfnn \frac{\l+1}{2 \l}\land q >1$,  clearly, $\l \varpi \le  \frac{\l+1}{2 }  < 1$.
  Also, we  fix   $ \e \in \Big( 0,  \frac14 \land \big( 4  C_\Psi \big)^{  \frac{ 2 \l \varpi}{  \l \varpi - 1 } } \Big)$.
      For $ G^1_\e \dfnn \Big\{ \o \in \O^t: E_t \Big[   \int_t^T \Psi^2_r \, dr \big| \cF^t_\t  \Big] (\o) >
       \e^{  \frac{\l \varpi - 1 }{2 \l \varpi}}  \Big\}
         \in \cF^t_\t$,  one can   deduce that
       \bea   \label{eq:t321}
     P^t_0 ( G^1_\e ) \le \e^{  \frac{1-\l \varpi}{2 \l \varpi} }  \,  E_t \bigg[ \b1_{G^1_\e}  E_t \Big[ \int_t^T \Psi^2_r \, dr \Big| \cF^t_\t  \Big]
     \bigg]   \le \e^{  \frac{1-\l \varpi}{2 \l \varpi} } \,    E_t \int_t^T \Psi^2_r \, dr =    C_\Psi \e^{  \frac{1-\l \varpi}{2 \l \varpi} } <  \frac14 .
       \eea
   Let $  \fk(\e) \dfnn   \big\lceil 1 -\log_2 \e  \big\rceil   $. Given $\fk \in \hN$ with $ \fk \ge \fk(\e)$,
   There exist a $  \d (\fk) > 0$  and   a closed subset $F_\fk$ of $\O^t$ 
 with    $   P^t_0\big( F_\fk \big)  \neg > \neg  1 - 2^{-2\fk} $ such that for any $\o ,\o'    \in    F_\fk$
 with $\| \o-\o'\|_t < \d (\fk) $
  \bea  \label{eq:r011}
  \underset{r \in [t,T]}{\sup} \, \underset{u \in  \hU}{\sup} \; \rho_{\overset{}{\hV}} \big( \beta(r,\o,u) , \beta(r,\o',u) \big)
  < 2^{-2\fk}  .
   \eea
      As   $ G_\fk \dfnn \Big\{ \o \in \O^t: E_t \big[\b1_{F_\fk^c} \big|\cF^t_\t \big] (\o) >  2^{- \fk} \Big\}
 \in \cF^t_\t $, similar to \eqref{eq:t321}, one can deduce that
    $
     P^t_0 ( G_\fk )
      \le  2^\fk  P^t_0 \big( F_\fk^c \big) \le 2^{- \fk}$. Thus for $G^2_\e \dfnn
      \underset{\fk \ge  \fk(\e)  }{\cup} G_\fk  $, we have
      \beas
      P^t_0 \big( G^2_\e \big) \le 2^{1-\fk(\e)} \le \e .
      \eeas

   Using \eqref{lem_basic_complement}, \eqref{eq:n721} and applying Proposition \ref{prop_rcpd_L1}
   with $\xi \neg = \neg  \b1_{ F_\fk^c } $,
   we obtain that for all $\o  \neg \in \neg  \O^t$ except on a $P^t_0-$null set $\cN_\fk$
     \bea  \label{eq:t011b}
       P^{\t(\o)}_0 \big(     (F_\fk^{\t,\o})^c   \big) = P^{\t(\o)}_0 \big(     (F_\fk^c)^{\t,\o}   \big)
       = E_{\t(\o)}  \big[    \b1_{ (F_\fk^c)^{\t,\o}}   \big] =   E_{\t(\o)}  \Big[   \big( \b1_{ F_\fk^c } \big)^{\t,\o}  \Big]
        =    E_t \big[    \b1_{     F_\fk^c }  \big| \cF^t_\t    \big] (\o)  .
     \eea

 \ss \no {\bf b)}  As $  \cF^t_T  =  \sB(\O^t)   $ by \eqref{eq:xxc023},
    there exists an open set $\wt{O}_\e$ of $\O^t$   that  includes
 \beas
   \wt{G}_\e \dfnn G^1_\e  \cup  G^2_\e  \cup  \cN_0 \cup  \Big( \underset{n \in \hN}{\cup}  \cN^1_n  \Big) \cup  \Big( \underset{n \in \hN}{\cup} \cN^2_n  \Big)    \cup  \Big( \underset{\fk \ge  \fk(\e)  }{\cup} \cN_\fk \Big) ,
  \eeas
   and satisfies (see e.g. Proposition 15.11 of \cite{Royden_real})
    \bea  \label{eq:t341}
      P^t_0 \big( \wt{O}_\e \big) <  P^t_0 \big(G^1_\e \big) +  P^t_0 \big(  G^2_\e  \big)      + \e
         \le C_\Psi \e^{  \frac{1-\l \varpi}{2 \l \varpi} }   + 2 \, \e < \frac14 + 2 \, \e .
   \eea

      Fix $n \in \hN$ such that $\{\t =t_n \} \backslash \wt{O}_\e  \ne \es $, we let  $x \in \hR^k$
       and $\o \in  \{\t =t_n \} \backslash \wt{O}_\e $.  Since
      $I (t_n, x,\beta^{t_n,\o} )$ is finite by \eqref{eq:I_finite}, there exists   a  $\mu^n_{\o,x} \in \cU^{t_n} $ such that
    \beas   
    \wt{Y}^{t_n,x,\mu^n_{\o,x}, \beta^{t_n,\o}  \lan \mu^n_{\o,x} \ran }_{t_n}  \Big(T, h\Big(\wt{X}^{t_n,x,\mu^n_{\o,x}, \beta^{t_n,\o} \lan \mu^n_{\o,x} \ran}_T \Big)  \Big) \ge I (t_n, x,\beta^{t_n,\o} ) - \frac12 \e
     \ge  \wh{w}_1 (t_n, x  ) - \frac12 \e .
     \eeas
   Let $\wh{\fk} = \fk^n_{\o,x} \dfnn \fk(\e) \vee \Big\lceil \frac{\l \varpi}{1 - \l \varpi  }\log_2 \Big( \frac{8 \k^2}{\e}   E_{t_n} \neg   \int_{t_n}^T
       \big[  \big( \mu^n_{\o,x}  \big)_r   \big]^2_{\overset{}{\hU}}  dr \Big) \Big\rceil$.
 Given $x' \in O_\e (x)$ and $\o' \in \Big( \{\t =t_n \} \cap O^{\, \raisebox{0.5ex}{\scriptsize $t_n$} }_{ {\d(\wh{\fk}) }} (\o ) \Big) \Big\backslash \wt{O}_\e    $,
        applying \eqref{eq:s025} and \eqref{eq_estimate_Y_v} with $t =t_n$   and using
        H\"older's inequality, we obtain
      \bea
\qq && \hspace{-1.5cm}    \bigg| \wt{Y}^{t_n,x,\mu^n_{\o,x}, \beta^{t_n,\o}  \lan \mu^n_{\o,x} \ran}_{t_n}     \Big(T, h\Big(  \wt{X}^{t_n,x,\mu^n_{\o,x}, \beta^{t_n,\o}  \lan \mu^n_{\o,x} \ran}_T \Big)\Big)
    -   \wt{Y}^{t_n,x',\mu^n_{\o,x}, \beta^{t_n,\o'}  \lan \mu^n_{\o,x} \ran}_{t_n}
      \Big(T, h\Big(  \wt{X}^{t_n,x',\mu^n_{\o,x}, \beta^{t_n,\o'}  \lan \mu^n_{\o,x} \ran}_T \Big)\Big) \bigg|^{\varpi}
       \nonumber  \\
  &&   \le  c_0 \big|x-x'\big|^{\frac{2 \varpi}{q}}  +  c_\l (\k_\psi)^{\varpi}  E_{t_n} \Bigg[ \bigg( \int_{t_n}^T
     \rho^{2 }_{\overset{}{\hV}} \Big( \big( \beta^{t_n,\o}  \lan \mu^n_{\o,x} \ran \big)_r,
      \big( \beta^{t_n,\o'}  \lan \mu^n_{\o,x} \ran \big)_r \Big)   dr  \bigg)^{ \frac{\l \varpi}{2} } \Bigg] \nonumber \\
   &&  \q    +  c_\l (\k_\psi)^{\varpi}  E_{t_n} \Bigg[ \bigg( \int_{t_n}^T
     \rho^{2 }_{\overset{}{\hV}} \Big( \big( \beta^{t_n,\o}  \lan \mu^n_{\o,x} \ran \big)_r,
      \big( \beta^{t_n,\o'}  \lan \mu^n_{\o,x} \ran \big)_r \Big)   dr  \bigg)^{\l \varpi } \Bigg] .  \label{eq:t511}
      \eea
  For any $\wt{\o}  \neg \in  \neg F_{\wh{\fk}}^{t_n,\o ,\o' } \dfnn  F_{\wh{\fk}}^{t_n,\o  } \cap F_{\wh{\fk}}^{t_n,\o'  }  $, since $\o'  \neg \in  \neg   O^{\, \raisebox{0.5ex}{\scriptsize $t_n$} }_{ {\d(\wh{\fk} ) }} (\o )$,
   one has
   $   \underset{r \in [t,T]}{\sup}\big| ( \o  \otimes_{t_n} \wt{\o} ) (r) -
    ( \o ' \otimes_{t_n} \wt{\o} ) (r) \big|  \neg = \neg  \underset{r \in [t,t_n]}{\sup}\big| \o'   (r)  -   \o    (r) \big|< \d (\wh{\fk}) $.
        It then follows from \eqref{eq:r011} that for any $r \in [t_n,T]$
   \bea
    \rho_{\overset{}{\hV}} \Big( \big( \beta^{t_n,\o} \lan \mu^n_{\o,x} \ran \big)_r (\wt{\o}),
      \big( \beta^{t_n,\o'} \lan \mu^n_{\o,x} \ran \big)_r (\wt{\o}) \Big) & \tneg \dneg =& \tneg \dneg
    \rho_{\overset{}{\hV}} \Big(  \beta  \big( r ,     \o \otimes_{t_n} \wt{\o} , \big(\mu^n_{\o,x}\big)_r ( \o \otimes_{t_n} \wt{\o} )  \big) ,
  \beta \big( r, \o^n_i  \otimes_{t_n}  \wt{\o}, \big(\mu^n_{\o,x}\big)_r (\o \otimes_{t_n} \wt{\o})  \big)   \Big) \nonumber  \\
   & \tneg \dneg  <& \tneg \dneg    2^{-2 \wh{\fk}} < 2^{-\fk(\e)}   <   \e   .  \label{eq:t514}
    \eea
       Since both $\o$ and $\o^n_i$ belongs to $\{\t = t_n \} \cap G^c_\e \cap   G^{\raisebox{0.4ex}{\scriptsize $c$}}_{\wh{\fk}}  \cap   \cN^{\raisebox{0.4ex}{\scriptsize $c$}}_{\wh{\fk}}
       \cap \cN^{\raisebox{0.4ex}{\scriptsize $c$}}_0       \cap \big( \cN^1_n \big)^c \cap \big( \cN^2_n \big)^c   $,
   we see from   \eqref{eq:t325} and \eqref{eq:t329}  that
     \bea
  \q   && \hspace{-1.5 cm}   E_{t_n}     \int_{t_n}^T
     \rho^2_{\overset{}{\hV}} \Big( \big( \beta^{t_n,\o}  \lan \mu^n_{\o,x} \ran \big)_r,
      \big( \beta^{t_n,\o'}  \lan \mu^n_{\o,x} \ran \big)_r \Big)   dr     \nonumber     \\
  &&     \le 2   \int_{\O^{t_n}}  \neg   \int_{t_n}^T
     \bigg( \Big[ \beta^{t_n, \o }  \big( r, \wt{\o}, \big(\mu^n_{\o,x}\big)_r (\wt{\o})  \big)  \Big]^2_{\overset{}{\hV}}
     \neg  + \neg
     \Big[ \beta^{t_n, \o'}  \big( r, \wt{\o}, \big(\mu^n_{\o,x}\big)_r (\wt{\o})  \big)  \Big]^2_{\overset{}{\hV}} \bigg)  dr \, dP^{t_n}_0(\wt{\o}) \nonumber   \\
   &&  \le 4     E_{t_n} \neg  \int_{t_n}^T
     \big( \Psi^{t_n,\o}_r  \big)^2    dr
     + 4     E_{t_n} \neg   \int_{t_n}^T
     \big( \Psi^{t_n,\o'}_r  \big)^2    dr
     + 8 \k^2   E_{t_n} \neg   \int_{t_n}^T
       \big[  \big( \mu^n_{\o,x}  \big)_r   \big]^2_{\overset{}{\hU}} \,  dr  \nonumber  \\
  &&     \le     4   E_t  \bigg[   \int_t^T   \P^2_r \, dr \Big| \cF^t_\t\bigg] (\o)
  + 4 E_t  \bigg[   \int_t^T   \P^2_r \, dr \Big| \cF^t_\t\bigg] (\o')   + \e  \big( 2^{\wh{\fk}} \big)^{\frac{ 1-\l \varpi }{\l \varpi}}
  \le  8 \e^{  \frac{\l \varpi -1 }{2 \l \varpi}}  + \e  \big( 2^{\wh{\fk}} \big)^{\frac{ 1-\l \varpi }{\l \varpi}} ,  \label{eq:t517}
       \eea
    and \eqref{eq:t011b} shows that
     \beas   
   P^{t_n}_0 \Big( \big( F_{\wh{\fk}}^{t_n,\o ,\o' } \big)^c \Big)
   \neg  \le  \neg   P^{t_n}_0 \Big(  \big(F_{\wh{\fk}}^{t_n,\o }\big)^c
   \Big)  \neg + \neg  P^{t_n}_0 \Big(  \big(F_{\wh{\fk}}^{t_n,\o'}\big)^c   \Big)
      \neg = \neg  E_t \big[\b1_{F_{\wh{\fk}}^c} \big| \cF^t_\t \big] (\o )
       \neg + \neg  E_t \big[\b1_{F_{\wh{\fk}}^c} \big| \cF^t_\t \big] (\o')
        \neg \le   \neg   2^{1-\wh{\fk}}     .
     \eeas
     Putting this together with \eqref{eq:t514} and \eqref{eq:t517} back into \eqref{eq:t511}, one can deduce from H\"older's inequality   that
         \beas
 && \hspace{-0.7 cm}    \bigg| \wt{Y}^{t_n,x,\mu^n_{\o,x}, \beta^{t_n,\o}  \lan \mu^n_{\o,x} \ran}_{t_n}     \Big(T, h\Big(  \wt{X}^{t_n,x,\mu^n_{\o,x}, \beta^{t_n,\o}  \lan \mu^n_{\o,x} \ran}_T \Big)\Big)
    -   \wt{Y}^{t_n,x',\mu^n_{\o,x}, \beta^{t_n,\o'}  \lan \mu^n_{\o,x} \ran}_{t_n}
      \Big(T, h\Big(  \wt{X}^{t_n,x',\mu^n_{\o,x}, \beta^{t_n,\o'}  \lan \mu^n_{\o,x} \ran}_T \Big)\Big) \bigg|^{\varpi}
        \nonumber  \\
  &&  \hspace{-0.2 cm} \le  c_0 \e^{\frac{2 \varpi}{q}}  \neg + \neg  c_\l (\k_\psi)^{\varpi}   P^{t_n}_0 \Big(   F_{\wh{\fk}}^{t_n,\o ,\o' }  \Big)   \e^{\neg \l \varpi  }     \nonumber \\
 &&  +  c_\l (\k_\psi)^{\varpi} \big(E_{t_n}[1]\big)^{\frac12} \Big[ P^{t_n}_0 \Big( \big( F_{\wh{\fk}}^{t_n,\o ,\o' } \big)^c \Big) \Big]^{ \frac{ 1 -\l \varpi}{2} } \bigg(      E_{t_n}  \neg    \int_{t_n}^T
     \rho^2_{\overset{}{\hV}} \Big\{    \big( \beta^{t_n,\o}  \lan \mu^n_{\o,x} \ran \big)_r,
      \big( \beta^{t_n,\o'}  \lan \mu^n_{\o,x} \ran \big)_r \Big)   dr       \bigg)^{ \frac{\l \varpi}{2} }  \qq  ~  \\
  &&   +  c_\l (\k_\psi)^{\varpi} \bigg( P^{t_n}_0 \Big(   F_{\wh{\fk}}^{t_n,\o ,\o' }  \Big)     \e^{2 \l \varpi}
   \neg + \neg \Big[ P^{t_n}_0 \Big( \big( F_{\wh{\fk}}^{t_n,\o ,\o' } \big)^c \Big)
   \Big]^{ 1 - \l \varpi } \bigg(      E_{t_n}  \neg    \int_{t_n}^T
     \rho^2_{\overset{}{\hV}} \Big( \big( \beta^{t_n,\o}  \lan \mu^n_{\o,x} \ran \big)_r,
      \big( \beta^{t_n,\o'}  \lan \mu^n_{\o,x} \ran \big)_r \Big)   dr       \bigg)^{   \l \varpi }  \bigg\}    \\
      &&  \hspace{-0.2 cm} \le     c_\l (\k_\psi)^{\varpi} \e^{\l \varpi}        \neg
   +  \neg   c_\l  (\k_\psi)^{\varpi} \big(  \e^{  \frac{1-\l \varpi}{2 \l \varpi}} + \e  \big)^{ \frac{\l \varpi}{2} }
   +  \neg   c_\l  (\k_\psi)^{\varpi} \big(  \e^{  \frac{1-\l \varpi}{2 \l \varpi}} + \e  \big)^{ \l \varpi }
    \le c_\l (\k_\psi)^{\varpi}    \e^{\frac{1- \l \varpi}{4}}   ,
      \eeas
   where  We used $\l < 1 \le 2/q $  and  $ 2^{1-\wh{\fk}}    \le     2^{1 - \fk(\e) }    \le     \e$  in the second inequality.
   Moreover,   Proposition  \ref{prop_w_conti}     assures that there exists  a $\l_n  \in  ( 0, \e ) $ such that
   $\big|  \wh{w}_1(t_n, x  )  -  \wh{w}_1 (t_n,x' ) \big| <     \frac12 \e $  for any $x' \in O_{\l_n }(x)$.    
    Therefore, it holds for any $x' \in O_{\l_n }(x) $
    and $\o' \in \Big( \{\t =t_n \} \cap O^{\, \raisebox{0.5ex}{\scriptsize $t_n$} }_{ {\d(\wh{\fk}) }} (\o ) \Big) \Big\backslash \wt{O}_\e    $   that
     \bea   \label{eq:t525}
  \q  && \hspace{-1.2 cm} \wt{Y}^{t_n,x',\mu^n_{\o,x}, \beta^{t_n,\o'}  \lan \mu^n_{\o,x} \ran}_{t_n}
      \Big(T, h\Big(  \wt{X}^{t_n,x',\mu^n_{\o,x}, \beta^{t_n,\o'}  \lan \mu^n_{\o,x} \ran}_T \Big)\Big)
    \neg  \ge \neg \wt{Y}^{t_n,x,\mu^n_{\o,x}, \beta^{t_n,\o}  \lan \mu^n_{\o,x} \ran}_{t_n}     \Big(T, h\Big(  \wt{X}^{t_n,x,\mu^n_{\o,x}, \beta^{t_n,\o}  \lan \mu^n_{\o,x} \ran}_T \Big)\Big) \neg - \neg  c_0 \k_\psi \e^{\frac{1-\l \varpi}{4 \varpi}} \nonumber \\
   &&   \ge  \wh{w}_1 (t_n, x  ) - \frac12 \e  - c_0 \k_\psi \e^{\frac{1-\l \varpi}{4 \varpi}}
    >  \wh{w}_1 (t_n,x' )  - \e - c_0 \k_\psi \e^{\frac{1-\l \varpi}{4 \varpi}}
    >  \wh{w}_1 (t_n,x' ) - c_0 \k_\psi \e^{\frac{1-\l \varpi}{4 \varpi}} .
     \eea

   \ss \no {\bf c)}    For any $n \in \hN$, there exists a closed set $F^n_\e$ of $\O^t$   that is included in $\{\t = t_n \}$
   and satisfies $ P^t_0 (F^n_\e) \ge \big( P^t_0 \{\t = t_n \} - \frac{\e}{2^n} \big)^+ $.
   Set $\fN \dfnn \{n \in \hN: F^n_\e \backslash \wt{O}_\e  \ne \es\}$. We claim that $\fN \ne \es$: Assume not.
    Then  $ \underset{n \in \hN}{\cup} F^n_\e \subset \wt{O}_\e $   and it follows
     \beas
   \frac34 \ge \frac14 + 2 \, \e  > P^t_0 \big( \wt{O}_\e \big) \ge
        \sum_{n \in \hN}  P^t_0(F^n_\e  ) \ge  \sum_{n \in \hN} \Big( P^t_0 \{\t = t_n \}  \neg - \neg  \frac{\e}{2^{n   }} \Big)
        = 1 - \e \ge \frac34  .
     \eeas
   A contradiction appears.

     Let $n \in \fN$.       We see from    Lemma \ref{lem_ocset_Omega}  that
     $\Big\{ \fO^n_{\o,x}    \dfnn O^{\, t_n }_{\d(\fk^n_{\o,x} )}(\o) \times O_{\l_n}(x)  : \o  \neg  \in  \neg F^n_\e \backslash \wt{O}_\e , x \in \hR^k \Big\} $ together with
     $  \big( \O^t \big\backslash ( F^n_\e \backslash \wt{O}_\e ) \big) \times \hR^k $ form an open cover of $\O^t \times \hR^k$.
     Since the  canonical space $\O^t$ is separable,    the product space $ \O^t \times \hR^k  $ is still separable and thus Lindel\"of.     Then one can find
     a sequence $ \{(\o^n_i, x^n_i )\}_{i \in \hN} $ of $ \big( F^n_\e \backslash \wt{O}_\e \big) \times \hR^k $ such that
     $ \big( F^n_\e \backslash \wt{O}_\e \big) \times \hR^k \subset \underset{i \in \hN}{\cup}  \fO^n_i     $.
                with $\fO^n_i \dfnn \fO^n_{\o^n_i, x^n_i }$.
     We  assume without loss of generality that
       \beas
        \wt{\fO}^n_i \dfnn \fO^n_i
 \big\backslash \underset{j < i }{\cup}  \fO^n_j  \ne \es, \q \fa i \ge 2
    \eeas
    and set  $\wt{\fO}^n_1 \dfnn \fO^n_1$.
    Since   $\big\{O^{\, \raisebox{0.5ex}{\scriptsize $t_n$} }_{ {\d }} (\o):  \o \in \O^t$, $\d > 0$\big\} are all
       $\cF^t_{t_n}-$measurable  by Lemma \ref{lem_ocset_Omega}  again,    one can  inductively  show that each
    $\wt{\fO}^n_i$ is a disjoint union of finitely many measurable rectangular sets
    $\big\{ \fA^{n,i}_j \times \fE^{n,i}_j :  \,  \fA^{n,i}_j \in \cF^t_{t_n} , \fE^{n,i}_j \in \sB(\hR^k)\big\}^{J^n_i}_{j=1}$.
        For any $i \in \hN$, we set $ \big(\fk^n_i, \beta^n_i, \mu^n_i \big)
        = \Big( \fk^n_{\o^n_i, \, x^n_i},\beta^{t_n, \o^n_i}, \mu^n_{\o^n_i, \, x^n_i}  \Big)$.

 \ss \no {\bf d)} We let $\Th \dfnn \big( t,x, \mu, \beta \lan \mu \ran \big) $ and
     define disjoint sets:
     \beas   
   \q   A^{n,i}_j \dfnn \{\t=t_n\}   \cap    \fA^{n,i}_j
   \cap \big\{  \wt{X}^\Th_{t_n}  \in \fE^{n,i}_j \big\} \in  \cF^{t}_{t_n} \cap  \cF^{t}_\t ,
   \q \fa n \in \fN ,~  i  \in \hN \hb{ and }   j = 1 , \cds \neg , J^n_i  .
     \eeas
     Fix  $m \in \hN$ such that  $\fN_m \dfnn \fN \cap \{1, \cds \neg , m\} \ne \es$. Setting
     $\dis  A_m \dfnn  \O^{t} \Big\backslash \Big( \underset{n \in \fN_m  }{\cup} \underset{i  =1  }{\overset{m}{\cup}}
      \underset{ j =1  }{\overset{J^n_i}{\cup}}  A^{n,i}_j \Big) \in \cF^{t}_\t  $, we  see from
        Proposition   \ref{prop_paste_control}   that
      \bea  \label{eq:t031}
  \mu^m_r (\o)   \dfnn \begin{cases}
     \big(\mu^n_i\big)_r \big( \Pi_{t,t_n} (\o) \big) , & \hb{if $(r,\o ) \neg \in \neg
   \[\t,T\]_{A^{n,i}_j} \neg = \neg [t_n,T] \neg \times \neg   A^{n,i}_j$    for $ n \in \fN_m$,
    $i =1 \cds \neg, m $ and $j \neg= \neg 1 \cds \neg , J^n_i $} ,          \q        \\
     \mu_r(\o),  & \hb{if $(r,\o )   \neg \in \neg  \[t,\t\[ \,\cup \,\[\t, T\]_{A_m }      $   }
    \end{cases}
    \eea
    defines a  $\cU^{t}-$control such that for any $(r,\o) \in    \[ \t,T  \] $
    \bea \label{eq:j241b}
  \big(\mu^m \big)^{\t,\o}_r = \begin{cases}
    \big( \mu^n_i  \big)_r,\q & \hb{if $ \o   \in    A^{n,i}_j$    for  $ n \in \fN_m$,
     $i =1 \cds \neg, m $ and $j \neg= \neg 1 \cds \neg , J^n_i $} ,   \\
       \mu^{\t,\o}_r ,  & \hb{if $ \o   \in A_m $.   }
    \end{cases}
    \eea
   Let $\wh{\Th}_m \neg \dfnn \neg  \big(t,x,\mu^m, \beta \lan \mu^m \ran \big) $.  As  $\mu^m   \neg = \neg  \mu  $
    on\;$ \[t, \t \[ $\;and\;thus\;$   \beta \lan \mu^m \ran    \neg =  \neg   \beta \lan \mu \ran   $\;on\;$ \[t, \t \[ $,
      \eqref{eq:p611}\;shows\;that\;$  P^{t}_0-$a.s.
   \beas
   \wt{X}^{\wh{\Th}_m}_s =X^{\wh{\Th}_m  }_s = X^\Th_s  =  \wt{X}^\Th_s , \q \fa s \in [t,\t] .
   \eeas
   Thus similar to \eqref{eq:q104},   for any $\cF^{t}_\t -$measurable random variable $\xi$
   with $\ul{L}^{\wh{\Th}_m}_\t =  \ul{L}^\Th_\t \le \xi \le \ol{L}^\Th_\t = \ol{L}^{\wh{\Th}_m}_\t$, $  P^{t}_0-$a.s.,
        \bea
     \big(Y^{\wh{\Th}_m} (\t, \xi),Z^{\wh{\Th}_m} (\t, \xi),\ul{K}^{\wh{\Th}_m} (\t, \xi),
    \ol{K}^{\, \raisebox{-0.5ex}{\scriptsize $\wh{\Th}_m$}} (\t, \xi) \big)
      =  \big(Y^\Th (\t, \xi),Z^\Th (\t, \xi),\ul{K}^{\Th}(\t, \xi),
       \ol{K}^{\, \raisebox{-0.5ex}{\scriptsize $\Th$}} (\t, \xi) \big)  .        \label{eq:q104b}
       \eea

    Let $ n \in \fN_m$,  $i =1 \cds \neg, m $, $j \neg= \neg 1 \cds \neg , J^n_i $ and  $\o \in \wt{A}^{n,i}_j \dfnn A^{n,i}_j \backslash \big( \wt{O}_\e   \cup \wh{\cN}_m \big) $, where
       $\wh{\cN}_m$ is the $P^t_0-$null set such that $\big(\beta \lan \mu^m \ran\big)^{\t,\o} \in \cV^{\t(\o)}$
       for all $\o \in \wt{\cN}^c_m$   according to  Proposition \ref{prop_dwarf_strategy} (1).
       For  $(r,\wt{\o}) \in [t_n,T] \times \O^{t_n} $,  since $A^{n,i}_j \in \cF^t_{t_n}$, Lemma \ref{lem_element}
    shows that $ \o \otimes_{t_n} \wt{\o} \in A^{n,i}_j$. Then
   one can deduce from \eqref{eq:t031} that
     \bea
   \big( \beta \lan \mu^m \ran \big)^{t_n,\o}_r (\wt{\o}) & \tneg =&  \tneg
    \big( \beta \lan \mu^m \ran \big)_r \big( \o \otimes_{t_n} \wt{\o}  \big)
    = \beta \big( r , \o \otimes_{t_n} \wt{\o} , \mu^m_r ( \o \otimes_{t_n} \wt{\o} ) \big)  \nonumber   \\
     & \tneg =&  \tneg      \beta  \Big( r ,     \o \otimes_{t_n} \wt{\o} , \big(\mu^n_i\big)_r ( \wt{\o} )  \Big)
     = \beta^{t_n, \o }  \Big( r, \wt{\o}, \big(\mu^n_i\big)_r (\wt{\o})  \Big)
     = \big( \beta^{t_n, \o } \lan \mu^n_i \ran  \big) _r (\wt{\o}) .    \label{eq:t263a}
     \eea
    Clearly, $\o^n_i \in F^n_\e \backslash \wt{O}_\e \subset \{\t =t_n\} \backslash \wt{O}_\e$. As
    $\o \in \wt{A}^{n,i}_j $, we see that     $\o \in   \big( \{\t=t_n\}   \cap    \fA^{n,i}_j \big)  \backslash   \wt{O}_\e
     \subset  \Big( \{\t=t_n\}   \cap     O^{\, t_n }_{\d(\fk^n_i )}(\o^n_i) \Big)  \Big\backslash   \wt{O}_\e   $
     and that $\wt{X}^\Th_{t_n}(\o) \in \fE^{n,i}_j \subset O_{\l_n } (x^n_i)$. Applying
     \eqref{eq:t525} with $(\o,x,\o',x') = \big(\o^n_i,x^n_i, \o, \wt{X}^\Th_{t_n}(\o)\big)$,
     we see from \eqref{eq:j701} that
        \beas
    \wt{Y}^{t_n,\wt{X}^\Th_{t_n}(\o),\mu^n_i, \beta^{t_n,\o}  \lan \mu^n_i \ran}_{t_n}
      \Big(T, h\Big(  \wt{X}^{t_n,\wt{X}^\Th_{t_n}(\o),\mu^n_i, \beta^{t_n,\o}  \lan \mu^n_i \ran}_T \Big)\Big)
         \ge   \Big(     \wh{w}_1 \big(t_n,\wt{X}^\Th_{t_n}(\o) \big)  - c_0 \k_\psi \e^{\frac{1-\l \varpi}{4 \varpi}} \Big) \vee \ul{l} \big(
        t_n,\wt{X}^\Th_{t_n}(\o) \big) .
     \eeas

     \if{0}
   For any $\wt{\o} \in F^{\raisebox{0.5ex}{\scriptsize{$n,\o$}}}_{\fk^n_i} \dfnn  F_{\fk^n_i}^{\raisebox{0.5ex}{\scriptsize{$t_n,\o^n_i$}} } \cap F_{\fk^n_i}^{\raisebox{0.5ex}{\scriptsize{$t_n,\o$}}}  $,
    since $\o \in A^{n,i}_j \subset  O^{\, t_n }_{\d(\fk^n_i )}(\o^n_i)$,
   one has
   \beas
     \underset{r \in [t,T]}{\sup}\big| ( \o^n_i \otimes_{t_n} \wt{\o} ) (r)-
    ( \o \otimes_{t_n} \wt{\o} ) (r) \big| = \underset{r \in [t,t_n]}{\sup}\big| \o^n_i   (r) -  \o    (r) \big|
    < \d \big( \fk^n_i \big) .
    \eeas
      Then it follows from \eqref{eq:r011} that for any $r \in [t_n,T]$
   \beas
  ~  \rho_{\overset{}{\hV}} \Big( \big( \beta \lan \mu^m \ran \big)^{t_n,\o}_r (\wt{\o}),
      \big( \beta^n_i \lan \mu^n_i \ran \big)_r (\wt{\o}) \Big) \neg = \neg
    \rho_{\overset{}{\hV}} \Big(  \beta  \big( r ,     \o \otimes_{t_n} \wt{\o} , \big(\mu^n_i\big)_r ( \wt{\o} )  \big) ,
  \beta \big( r, \o^n_i  \otimes_{t_n}  \wt{\o}, \big(\mu^n_i\big)_r (\wt{\o})  \big)   \Big)  \neg  <  \neg  2^{-2 \fk^n_i}
   \neg < \neg  2^{-\fk(\e)}  \neg < \neg  \e  .
    \eeas
  Since both $\o$ and $\o^n_i$ belongs to $\{\t = t_n \} \cap G^c_\e \cap  \cN^c_0    \cap \big( \cN^1_n \big)^c
   \cap \big( \cN^2_n \big)^c \cap G^c_{\fk^n_i} \cap   \cN^c_{\fk^n_i}     $, similar to \eqref{eq:t517} and \eqref{eq:t519}, we can deduce from   \eqref{eq:t263a}, \eqref{eq:t263b} that  
     \beas
  \q   && \hspace{-1.5 cm}      E_{t_n}  \neg   \int_{t_n}^T
     \rho^2_{\overset{}{\hV}} \Big( \big( \beta \lan \mu^m \ran \big)^{t_n,\o}_r,
      \big( \beta^n_i \lan \mu^n_i \ran\big)_r \Big)   dr
              \le 2   \int_{\O^{t_n}}  \neg   \int_{t_n}^T
     \bigg( \Big[ \beta^{t_n, \o }  \big( r, \wt{\o}, \big(\mu^n_i\big)_r (\wt{\o})  \big)  \Big]^2_{\overset{}{\hV}}
     \neg  + \neg
     \Big[ \beta^{t_n, \o^n_i}  \big( r, \wt{\o}, \big(\mu^n_i\big)_r (\wt{\o})  \big)  \Big]^2_{\overset{}{\hV}} \bigg)  dr \, dP^{t_n}_0(\wt{\o}) \nonumber   \\
 &&  \le 4     E_{t_n} \neg  \int_{t_n}^T
     \big( \Psi^{t_n,\o}_r  \big)^2    dr
     + 4     E_{t_n} \neg   \int_{t_n}^T
     \big( \Psi^{t_n,\o^n_i}_r  \big)^2    dr
     + 8 \k^2   E_{t_n} \neg   \int_{t_n}^T
       \big[  \big( \mu^n_i  \big)_r   \big]^2_{\overset{}{\hU}}  dr  \nonumber  \\
  &&     \le     4   E_t  \bigg[   \int_t^T   \P^2_r \, dr \Big| \cF^t_\t\bigg] (\o)
        + 4 E_t  \bigg[   \int_t^T   \P^2_r \, dr \Big| \cF^t_\t\bigg] (\o^n_i)  + 8 \k^2 \fz
       \le 8 \fk     +       8 \k^2 \fz  ,
       \eeas
     shows that
     \beas
    P^{t_n}_0 \big( (F^{n,\o}_\fk  )^c \big) 
      \le P^{t_n}_0 \big(     (F_\fk^{t_n,\o^n_i})^c   \big) +  P^{t_n}_0 \big(     (F_\fk^{t_n,\o})^c   \big)
     = E_t \big[\b1_{F_\fk^c} \big| \cF^t_\t \big] (\o^n_i) + E_t \big[\b1_{F_\fk^c} \big| \cF^t_\t \big] (\o) \le 2 \sqrt{\e}  .
     \eeas
      Then applying \eqref{eq:s025} and  \eqref{eq_estimate_Y_v} with $t =t_n$ and $ \varpi  =    \frac{q+1}{2}  $     yields  that
      \beas
 && \hspace{-1.5cm}    \Bigg| \wt{Y}^{t_n,\wt{X}^\Th_{t_n}(\o),\mu^n_i, ({\beta \lan \mu^m \ran})^{t_n,\o}}_{t_n}     \bigg(T, h\Big(  \wt{X}^{t_n,\wt{X}^\Th_{t_n}(\o),\mu^n_i,({\beta \lan \mu^m \ran})^{t_n,\o}}_T \Big)\bigg)
    -   \wt{Y}^{t_n,\wt{X}^\Th_{t_n}(\o),\mu^n_i, \beta^n_i \lan \mu^n_i \ran}_{t_n}    \bigg(T, h\Big(  \wt{X}^{t_n,\wt{X}^\Th_{t_n}(\o),\mu^n_i,\beta^n_i \lan \mu^n_i \ran}_T \Big)\bigg) \Bigg|^{\frac{q+1}{2}} \\
  &&   \le    c_0  E_{t_n} \Bigg[ \bigg( \int_{t_n}^T
     \rho^2_{\overset{}{\hV}} \Big( \big( \beta \lan \mu^m \ran \big)^{t_n,\o}_r,
      \big( \beta^n_i \lan \mu^n_i \ran\big)_r \Big)   dr  \bigg)^{ \frac{q+1}{2q} } \Bigg]   \\
   &&   \le  c_0  \e^{\neg \frac{1+q}{q}}  P^{t_n}_0 \big(F^{n,\o}_\fk \big)
      +  c_0 \Big( P^{t_n}_0 \big( (F^{n,\o}_\fk  )^c \big) \Big)^{ \frac{ q-1}{2q} } \bigg(
      E_{t_n}  \neg   \int_{t_n}^T
     \rho^2_{\overset{}{\hV}} \Big( \big( \beta \lan \mu^m \ran \big)^{t_n,\o}_r,
      \big( \beta^n_i \lan \mu^n_i \ran\big)_r \Big)   dr
      \bigg)^{ \frac{q+1}{2q} }    \\
   &&   \le  c_0  \e^{\neg \frac{q+1}{q}}
      +  c_0 \e^{ \frac{ q-1}{4q} }\big(  \fk + \k^2 \fz  \big)^{ \frac{q+1}{2q} }     ,
     \eeas
     \fi
    Using similar arguments to those that lead to \eqref{eq:t023}, one can deduce from \eqref{eq:j241b}  and \eqref{eq:t263a}  that
               \bea
     &&  \hspace{-1cm}    \bigg( \wt{Y}^{\wh{\Th}_m}_\t \Big(T, h \Big(\wt{X}^{\wh{\Th}_m}_T\Big)\Big)\bigg) (\o)
        =        \wt{Y}^{\t(\o),\wt{X}^\Th_{\t(\o)} (\o),(\mu^m)^{\t,\o},({\beta \lan \mu^m \ran})^{\t,\o}}_{\t(\o)}
   \bigg(T, h\Big(  \wt{X}^{\t(\o),\wt{X}^\Th_{\t(\o)} (\o),(\mu^m)^{\t,\o}, ({\beta \lan \mu^m \ran})^{\t,\o}}_T \Big)\bigg)
    \nonumber  \\
   &&   \ge     \sum_{n \in \fN_m } \neg \sum^m_{i  =1 } \sum^{J^n_i}_{j=1} \b1_{\big\{\o \in \wt{A}^{n,i}_j \big\}} \wt{Y}^{t_n,\wt{X}^\Th_{t_n}(\o),\mu^n_i, \beta^{t_n, \o } \lan \mu^n_i \ran}_{t_n}    \bigg(T, h\Big(  \wt{X}^{t_n,\wt{X}^\Th_{t_n}(\o),\mu^n_i,\beta^{t_n, \o } \lan \mu^n_i \ran}_T \Big)\bigg)   \nonumber  \\
 &&  \q +
     \b1_{\{\o \in A_m \cup \wt{O}_\e   \}} \ul{l} \big( \t(\o),\wt{X}^\Th_{\t(\o)} (\o) \big)
    \ge   \wh{\xi}_m(\o) ,  \q \hb{for $P^{t}_0-$a.s. $ \o \in \O^{t}$,}   \label{eq:h031b}
  \eea
  where   $ \dis  \wh{\xi}_m \neg \dfnn    \neg
        \b1_{    A^c_m \cap \wt{O}^c_\e   }
   \Big(     \wh{w}_1 \big(\t,\wt{X}^\Th_\t(\o) \big)  - c_0 \k_\psi \e^{\frac{1-\l \varpi}{4 \varpi}} \Big)
      \vee    \ul{l} \big( \t ,\wt{X}^\Th_{\t }   \big)
     +   \b1_{    A_m \cup \wt{O}_\e  } \ul{l} \big( \t ,\wt{X}^\Th_{\t }   \big)      $.

\ss \no  {\bf e) }   Similar to \eqref{eq:j335}, we see from Proposition \ref{prop_apriori_DRBSDE} that
      \beas
   \big|  \wt{Y}^{\Th}_{t}  (\t,  \wh{\xi}_m  )
   \neg  - \neg  \wt{Y}^\Th_{t}  \big(\t,    \wh{w}_1\big(\t , \wt{X}^\Th_{\t }  \big)   \big)  \big|
     \neg \le \neg c_0  \big\|   \wh{\xi}_m      \neg  -  \neg
        \wh{w}_1\big(\t , \wt{X}^\Th_{\t }  \big)     \big\|_{ \hL^q  (\cF^{t}_\t  ) }
   \neg  \le  \neg   c_0 \big\|   \b1_{    A_m \cup \wt{O}_\e  } \neg \big(
    \ol{l} \big( \t ,\wt{X}^\Th_\t    \big) \neg  -  \neg  \ul{l} \big(\t , \wt{X}^\Th_{\t }  \big)    \big)   \big\|_{ \hL^q  (\cF^{t}_\t  ) }  \neg + \neg  c_0 \k_\psi \e^{\frac{1-\l \varpi}{4 \varpi}}  .    \q~\;\;
   \eeas
       Then applying  \eqref{eq:p677} with $(\z,\t,\xi) = \Big(\t, T, h \big( \wt{X}^{ \wh{\Th}_m }_T \big) \Big)$,
        applying  Proposition \ref{prop_comp_DRBSDE} to \eqref{eq:h031b}, and
   using  \eqref{eq:q104b} for   $\xi =\xi_m$ yield that
    \bea
    I(t,x, \beta) &\ge&  \wt{Y}^{\wh{\Th}_m}_{t} \Big(T,h \big( \wt{X}^{ \wh{\Th}_m }_T \big) \Big)
  =     \wt{Y}^{\wh{\Th}_m}_{t} \Big(\t,
    \wt{Y}^{\wh{\Th}_m}_\t \big(T, h \big( \wt{X}^{ \wh{\Th}_m }_T \big) \big)   \Big)
   \ge \wt{Y}^{\wh{\Th}_m}_{t}  \big( \t,         \wh{\xi}_m      \big)
 =  \wt{Y}^\Th_{t} \big(\t,         \wh{\xi}_m     \big) \nonumber \\
 &\ge& \wt{Y}^\Th_{t}  \big(\t,    \wh{w}_1\big(\t , \wt{X}^\Th_{\t }  \big)   \big)
 -    c_0 \big\|   \b1_{ A_m \cup \wt{O}_\e } \neg \big(
    \ol{l} \big( \t ,\wt{X}^\Th_\t    \big) \neg  -  \neg  \ul{l} \big(\t , \wt{X}^\Th_{\t }  \big)    \big)   \big\|_{ \hL^q  (\cF^{t}_\t  ) } -  c_0 \k_\psi \e^{\frac{1-\l \varpi}{4 \varpi}} .  \label{eq:h041b}
    \eea
   Let $A_\e \dfnn \lmtd{m \to \infty}   A_m = \underset{m \in \hN}{\cap} A_m $. As $m \to \infty$   in \eqref{eq:h041b},   
   the Dominated Convergence Theorem and \eqref{eq:t528}   show that
    \bea   \label{eq:t533}
      I(t,x, \beta)  \ge  \wt{Y}^\Th_{t}  \big(\t,    \wh{w}_1\big(\t , \wt{X}^\Th_{\t }  \big)   \big)
 -    c_0 \big\|   \b1_{ A_\e \cup \wt{O}_\e } \neg \big(
    \ol{l} \big( \t ,\wt{X}^\Th_\t    \big) \neg  -  \neg  \ul{l} \big(\t , \wt{X}^\Th_{\t }  \big)    \big)   \big\|_{ \hL^q  (\cF^{t}_\t  ) } -  c_0 \k_\psi \e^{\frac{1-\l \varpi}{4 \varpi}} .
    \eea

   Given $n \in \fN$ and $ \o \in F^n_\e \backslash \wt{O}_\e \subset \{\t = t_n\}$,
 since $\big( \o, \wt{X}_{t_n} (\o) \big) \in \big( F^n_\e \backslash \wt{O}_\e \big) \times \hR^k $, there exists
 $i \in \hN$, such that $\big( \o, \wt{X}_{t_n} (\o) \big) $ is in $ \wt{\fO}^n_i $ and thus further  belongs to
  some $\fA^{n,i}_j \times \fE^{n,i}_j$, $j =1, \cds \neg , J^n_i$. To wit,   $ \o \in \{\t = t_n\}
 \cap \fA^{n,i}_j \cap \big\{ \wt{X}_{t_n} \in \fE^{n,i}_j \big\} = A^{n,i}_j $.  It follows that
    \beas
 \Big( \underset{n \in \hN}{\cup}  F^n_\e  \Big)  \backslash \wt{O}_\e
  =\underset{n \in \hN}{\cup} \Big( F^n_\e \backslash \wt{O}_\e \Big)
 = \underset{n \in \fN}{\cup} \Big( F^n_\e \backslash \wt{O}_\e \Big)
 \subset  \underset{n \in \fN}{\cup} \underset{i  \in \hN}{\cup} \underset{j = 1  }{\overset{J^n_i}{\cup}}  A^{n,i}_j
 = \underset{m \in \hN}{\cup}   \underset{n \in \fN_m  }{\cup} \underset{i = 1  }{\overset{m}{\cup}}
   \underset{j = 1  }{\overset{J^n_i}{\cup}}  A^{n,i}_j   ,   
        \eeas
   which together with   \eqref{eq:t341} implies  that
   \beas
 P^t_0  \big( A_\e \cup \wt{O}_\e \big) &\neg =  & \neg  P^t_0  \bigg( \Big( \underset{m \in \hN}{\cup} \underset{n \in \fN_m  }{\cup} \underset{i = 1  }{\overset{m}{\cup}}
   \underset{j = 1  }{\overset{J^n_i}{\cup}}  A^{n,i}_j \Big)^c \cup \wt{O}_\e \bigg)
   \neg   \le    \neg     P^t_0 \Big( \Big( \underset{n \in \hN}{\cup} F^n_\e \Big)^c \Big) + P^t_0 \big( \wt{O}_\e \big)  \\
&\neg = & \neg  1  - \sum_{n \in \hN} P^t_0(F^n_\e)     +        P^t_0 \big(\wt{O}_\e  \big)
        \le          C_\Psi \e^{  \frac{1-\l \varpi}{2 \l \varpi} }      +     3 \, \e .
    \eeas
   Thus, letting $\e \to 0$ in \eqref{eq:t533}, we obtain
      \beas
        I(t,x, \beta)  \ge  \wt{Y}^{t,x, \mu, \beta \lan \mu \ran}_{t}
         \big(\t,    \wh{w}_1\big(\t , \wt{X}^{t,x, \mu, \beta \lan \mu \ran}_{\t }  \big)   \big)   .
      \eeas
   Eventually,  taking supremum over $\mu \in \cU^{t}$ on the right-hand-side and then  taking infimum over $\beta \in \wh{\fB}^{t}$ on both sides yield that
    \beas    
    \wh{w}_1(t,x)  \ge \underset{\beta \in \wh{\fB}^{\, t}}{\inf} \, \underset{\mu \in \cU^{t}}{\sup} \wt{Y}^{t,x,\mu, \beta \lan \mu \ran}_{t}
   \Big({\t_{\mu,\beta}}, \,  \wh{w}_1\big({\t_{\mu,\beta}}, \wt{X}^{t,x,\mu, \beta \lan \mu \ran }_{\t_{\mu,\beta}} \big) \Big)  .
   \eeas

 \ss \no   {\bf 3)} For any $\big( \ol{t},\ol{x} , \ol{y},\ol{z}, \ol{u},\ol{v} \big) \in [0,T] \times \hR^k \times \hR \times \hR^d \times \hU \times \hV$,  we   define
 \beas
 \ul{\fl} \big(\ol{t},\ol{x}\big) \dfnn  - \ol{l} \big(\ol{t},\ol{x}\big)  , \q
 \ol{\fl} \big(\ol{t},\ol{x}\big) \dfnn - \ul{l} \big(\ol{t},\ol{x}\big) , \q
  \fh \big( \ol{x}\big)    \dfnn    - h \big(\ol{x}\big)
  \q \hb{and} \q \ff \big( \ol{t},\ol{x} , \ol{y},\ol{z}, \ol{u},\ol{v} \big)
 \dfnn  -f \big( \ol{t},\ol{x} , -\ol{y},-\ol{z}, \ol{u},\ol{v} \big)   .
  \eeas
  Given $ \mu \in \cU^t$ and $ \nu \in \cV^t$, we still let $\Th$ stand for  $(t,x,\mu,\nu)$  and set  $\ul{\sL}^\Th_s \dfnn \ul{\fl}\big(s, \wt{X}^\Th_s \big) $
   and $\ol{\sL}^\Th_s \dfnn \ol{\fl}\big(s, \wt{X}^\Th_s \big) $, $s \in [t,T]$.
  For any $\bF^t-$stopping time $\t$ and any $\cF^t_\t-$measurable random variable $\xi$ with $ \ul{\sL}^\Th_\t \le  \xi \le \ol{\sL}^\Th_\t $, $P^t_0-$a.s.,   let
     $ \big(\sY^\Th (\t, \xi),\sZ^\Th (\t, \xi),\ul{\sK}^{\Th}(\t, \xi), \ol{\sK}^{\, \raisebox{-0.5ex}{\scriptsize $\Th$}} (\t, \xi) \big) $ denote the unique solution of the DRBSDE$\big(P^t_0,\xi,  \ff^\Th_\t,    \ul{\sL}^\Th_{\t \land \cd},\ol{\sL}^{\,\raisebox{-0.5ex}{\scriptsize $\Th$}}_{\t \land \cd}\big)$ in $ \hG^q_{\ol{\bF}^t } \big([t,T]\big) $, where
    \beas
     \ff^\Th_\t   (s,\o,y,z ) \dfnn \b1_{\{s < \t(\o)\}}  \ff  \Big(s,  \wt{X}^\Th_s(\o) , \, y, z, \mu_s (\o),\nu_s (\o) \Big) ,
    \q   \fa (s,\o,y,z  ) \in   [t,T] \times \O^t \times \hR \times \hR^d .
    \eeas
 Since $\ul{\sL}^\Th \neg =\neg - \ol{L}^{\,\raisebox{-0.5ex}{\scriptsize $\Th$}} $ and $\ol{\sL}^{\,\raisebox{-0.5ex}{\scriptsize $\Th$}} \neg =\neg  - \ul{L}^\Th $,    multiplying $-1$ in the DRBSDE$\big(P^t_0,\xi,  \ff^\Th_\t,    \ul{\sL}^\Th_{\t \land \cd},\ol{\sL}^{\,\raisebox{-0.5ex}{\scriptsize $\Th$}}_{\t \land \cd}\big)$ shows that $\big(\neg - \neg  \sY^\Th (\t, \xi),$ \\
  $ -\sZ^\Th (\t, \xi),-\ol{\sK}^{\, \raisebox{-0.5ex}{\scriptsize $\Th$}} (\t, \xi), -\ul{\sK}^{\Th}(\t, \xi)  \big) \neg \in \neg \hG^q_{\ol{\bF}^t } \big([t,T]\big) $ solves the DRBSDE$\big(P^t_0,-\xi,  f^\Th_\t, \ul{L}^\Th_{\t \land \cd},\ol{L}^{\,\raisebox{-0.5ex}{\scriptsize $\Th$}}_{\t \land \cd}\big)$. To wit
   \bea  \label{eq:s304}
    \big(-\sY^\Th (\t, \xi),-\sZ^\Th (\t, \xi),-\ol{\sK}^{\, \raisebox{-0.5ex}{\scriptsize $\Th$}} (\t, \xi), -\ul{\sK}^{\Th}(\t, \xi)  \big) = \big(Y^\Th (\t, -\xi), Z^\Th (\t, -\xi),\ol{K}^{\, \raisebox{-0.5ex}{\scriptsize $\Th$}} (\t, -\xi), \ul{K}^{\Th}(\t, -\xi)  \big) .
   \eea

Now let us consider the situation where    player II acts first by choosing
 a $\cV^t-$control to maximize
  $\wt{\sY}^{t,x,  \a \lan \nu \ran, \nu }_t   \Big(T, \\ \fh \Big(\wt{X}^{t,x,\a \lan \nu \ran, \nu}_T \Big)  \Big)$,\;where
  $\a \neg \in \neg \cA^t$\;is\;player\;I's\;response. So\;the\;priority\;value\;and\;intrinsic\;priority\;value\;of\;player\;II are
        \beas
\q  \fw_2 \big(t,x  \big) \dfnn \underset{\a \in \cA^t }{\inf} \; \underset{\nu \in \cV^t}{\sup} \; \wt{\sY}^{t,x,  \a \lan \nu \ran, \nu }_t   \Big(T,\fh \Big(\wt{X}^{t,x,\a \lan \nu \ran, \nu}_T \Big)  \Big) ~\; \hb{and } ~\;   \ol{\fw}_2 \big(t,x  \big) \dfnn \underset{\a \in \wh{\cA}^t }{\inf} \; \underset{\nu \in \cV^t}{\sup} \; \wt{\sY}^{t,x,\a \lan \nu \ran, \nu }_t   \Big(T,\fh \Big(\wt{X}^{t,x,\a \lan \nu \ran, \nu}_T \Big)  \Big) .
  \eeas
   For any   family $\{\t_{\nu,\a} \neg: \nu \neg\in \neg \cV^t, \a \neg\in \neg\cA^t \}$ of $\hQ_{t,T} -$valued, $\bF^t-$stopping times,   applying \eqref{eq:DPP0} yields that
   \bea     \label{eq:s307}
  \fw_2 (t,x) \le
    \underset{\a \in \cA^t}{\inf} \; \underset{\nu \in \cV^t}{\sup}    \;
    \wt{\sY}^{t,x,\a \lan \nu \ran, \nu }_t \Big(\t_{\nu,\a},
    \fw_2\big(\t_{\nu,\a},  \wt{X}^{t,x,\a \lan \nu \ran, \nu}_{\t_{\nu,\a}} \big) \Big)   .
   \eea
 For any $ \big(t,x  \big) \in [0,T] \times \hR^k$,  we see from   \eqref{eq:s304}  that
   \beas
  ~ \; \;  - \fw_2 \big(t,x\big) \neg = \neg  \underset{\a \in \cA^t }{\sup} \; \underset{\nu \in \cV^t}{\inf}  - \wt{\sY}^{t,x,  \a \lan \nu \ran, \nu }_t   \Big(T,\fh \Big(\wt{X}^{t,x,\a \lan \nu \ran, \nu}_T \Big)  \Big)
     \neg =  \neg  \underset{\a \in \cA^t }{\sup} \; \underset{\nu \in \cV^t}{\inf} \;   \wt{Y}^{t,x,  \a \lan \nu \ran, \nu }_t   \Big(T, h \Big(\wt{X}^{t,x,\a \lan \nu \ran, \nu}_T \Big)  \Big)  \neg = \neg  w_2 \big(t,x\big) .
    \eeas
    Putting it back into \eqref{eq:s307} and using \eqref{eq:s304}, we obtain \eqref{eq:DPP2}. Similarly, we
    have \eqref{eq:DPP3}  and its inverse holds under ($\bU_\l$) for some $\l \in (0,1)$.        \qed

 \subsection{Proofs of Section \ref{sec:PDE}}

 \label{subsection:Proof_S3}

  The proof of Theorem \ref{thm_viscosity} relies on the following comparison theorem
for generalized reflected BSDEs.

   \begin{prop}   \label{prop_comparison_RBSDE}

 Given $t \in [0,T]$  and   $i=1,2$, let  $\ff_i: [t,T] \times \O^t \times \hR \times \hR^d \to \hR $ be
  a    $\sP \big(  \ol{\bF}^t  \big)   \otimes \sB(\hR)  \otimes \sB(\hR^d)/\sB(\hR)-$measurable function.
 For some  $\xi_i \in \hL^q \big(\ol{\cF}^{\, t}_{\neg T} \big)$ and
   $L^i  \neg   \in  \neg  \hC^{+,q}_{\ol{\bF}^t} \big([t,T]\big)  $ \big(resp.\;$\hC^{-,q}_{\ol{\bF}^t} \big([t,T]\big)$\big)
   with $ \xi_i  \neg \ge $ \(resp.\;$\le$\) $ \neg  L^i_T   $, $P^t_0-$a.s.,
    let $ \big(Y^i ,Z^i ,V^i,K^i      \big)  \neg \in  \neg  \hC^q_{\ol{\bF}^t}( [t,T])
 \times  \hH^{2,q}_{\ol{\bF}^t}([t,T], \hR^d)
      \times \sV_{\ol{\bF}^t}( [t,T])   \times \hK_{\ol{\bF}^t}( [t,T]) $ be a solution of   the following
generalized reflected   backward stochastic differential equation with lower \(resp.\;upper\) obstacle  on the probability space $(\O^t,  \ol{\cF}^t_{\neg T}, P^t_0)$   \big($\ul{\hb{R}}$BSDE$ \big( P^t_0,\xi_i,  \ff_i,L^i  \big)$,
 resp.\;$\ol{\hb{R}}$BSDE$ \big( P^t_0,\xi_i,\ff_i,L^i  \big)$, for short\big):
  \bea   \label{RLUBSDE}
     \left\{\ba{l}
 \dis    L^i_s \le \(\hb{resp.}\, \ge \neg \) \, Y^i_s= \xi_i \neg + \neg  \int_s^T  \neg  \ff_i  (r,   Y^i_r, Z^i_r)  \, dr
  \neg+ \neg  V^i_{\,T} \neg-\neg V^i_{\,s} \neg+ \neg \(\hb{resp.}\, - \neg \) \, ( K^i_{\,T} \neg-\neg K^i_{\,s})
 \neg-\neg \int_s^T \neg Z^i_r d B^t_r  , \q    s \in [t,T] ,   \vspace{1mm} \q \\
   \dis     \int_t^T \neg  \big(  Y^i_s   -  L^i_s  \big)     d K^i_{\,s}   = 0 .
     \ea \right.
    \eea

If       $P^t_0(\xi_1   \le    \xi_2)    =    P^t_0 \big(L^1_s    \le    L^2_s ,
 \fa s    \in    [t,T] \big)    =    1 $, if
$V^1     -     V^2$ is a decreasing process, and if either of the following two conditions holds:

 \ss \no (\,i) $\ff_1$ satisfies \eqref{ff_Lip} and $  \ff_1 (s,  Y^2_s ,Z^2_s ) \le  \ff_2 (s, Y^2_s ,Z^2_s )  $,
     $ds \times dP^t_0-$a.s.,

 \ss \no (ii) $\ff_2$ satisfies   \eqref{ff_Lip} and   $  \ff_1 (s,  Y^1_s ,Z^1_s ) \le  \ff_2 (s, Y^1_s ,Z^1_s )  $,
     $ds \times dP^t_0-$a.s.;

 \ss \no  then
   $P^t_0 \big(  Y^1_s \le Y^2_s, \fa s \in [t,T] \big) =1 $.

   \end{prop}

 \ss  \no {\bf Proof:} We first show the comparison for
  $\ul{\hb{R}}$BSDE$ \big( P^t_0,\xi_i,\ff_i,L^i  \big)$, $i \neg = \neg 1,2$:
  Let     $ \D \fX \neg \dfnn  \neg \fX^1  \neg - \neg \fX^2 $ for $\fX   \neg =  \neg     Y,Z,V $.
  Similar to \eqref{eq:n101a}, applying Tanaka's formula to process $(\D Y)^+$
 \if{0}
 \beas
 (\D Y_s)^+ - (\D Y_{s'} )^+ & \tneg = & \tneg     \int_s^{s'} \b1_{\{ \D Y_r > 0 \}}
 \big( \ff_1 (r, Y^1_r,Z^1_r) - \ff_2 (r, Y^2_r,Z^2_r) \big) dr  - \frac12 \int_s^{s'} d \fL_r       \\
 & \tneg & \tneg  + \int_s^{s'} \b1_{\{ \D Y_r > 0 \}}  ( d \D V_r + d  K^1_r- d  K^2_r  )
 - \int_s^{s'} \b1_{\{ \D Y_r > 0 \}}  \D Z_r  d B^t_r
  ,   \q \fa t \le   s \le s' \le T  ,
 \eeas
where $\fL$ is a real-valued, $\bF^t-$adapted, increasing and continuous process known as ``local time".
\fi
 and using  Corollary 1 of  \cite{Hamadene_Popier_2008}, we obtain
     \bea
  && \hspace{-1cm} \big| (\D Y_s)^+ \big|^q  - \big| (\D Y_{s'})^+ \big|^q
+ \frac{q(q-1)}{2} \int_s^{s'}  \b1_{\{ \D Y_r > 0 \}}\big| (\D Y_s)^+ \big|^{q-2} |\D Z_r|^2 dr  \nonumber \\
 && \hspace{-0.6cm} \le       q \int_s^{s'}  \neg  \b1_{\{ \D Y_r > 0 \}}
 \big| (\D Y_r)^+ \big|^{q-1} \big( \ff_1 (r, Y^1_r,Z^1_r)  \neg - \neg  \ff_2 (r, Y^2_r,Z^2_r) \big) dr
  \neg +  \neg   q  \neg  \int_s^{s'}  \neg  \b1_{\{ \D Y_r > 0 \}} \big| (\D Y_r)^+ \big|^{q-1}
 ( d \D V_r  \neg + \neg  d  K^1_r  \neg - \neg  d  K^2_r )   \nonumber  \\
 && \hspace{-0.6cm} \q
     -    q     \int_s^{s'} \neg \b1_{\{ \D Y_r > 0 \}} \big| (\D Y_r)^+ \big|^{q-1}   \D Z_r    d B^t_r
       -     \frac{q}{2}  \neg   \int_s^{s'}  \neg  \b1_{\{ \D Y_r > 0 \}} \big| (\D Y_r)^+ \big|^{q-1} d \fL_r,
      \q  \fa t  \neg \le  \neg   s  \neg \le \neg  s'  \neg \le \neg  T    ,    \qq   \;\;     \label{eq:s411}
\eea
 where $\fL$ is a real-valued, $\bF^t-$adapted, increasing and continuous process known as ``local time".
   The flat-off condition of $(Y^1, Z^1, V^1, K^1)$ implies that $P^t_0-$a.s.
 \beas
   \q    0 \neg \le \neg \int_t^T  \neg  \b1_{\{ \D Y_r > 0 \}} \big| (\D Y_r)^+ \big|^{q-1}  d  K^1_r
  \neg =  \neg  \int_t^T  \neg  \b1_{\{ L^1_r = Y^1_r > Y^2_r    \}}
 \big| (L^1_r  \neg - \neg  Y^2_r)^+ \big|^{q-1}     d  K^1_r
  \neg  \le   \neg  \int_t^T  \neg  \b1_{\{ L^1_r > L^2_r    \}} \big| (L^1_r  \neg - \neg  L^2_r)^+  \big|^{q-1}  d  K^1_r
 \neg  = \neg  0 .
 \eeas
 Putting this back into \eqref{eq:s411} and using Lipschitz continuity of $\ff_1$ in $(y,z)$  yield that
 \beas
 & &  \hspace{-1cm} \big| (\D Y_s)^+ \big|^q  - \big| (\D Y_{s'})^+ \big|^q
  + \frac{q(q-1)}{2} \int_s^{s'} \b1_{\{ \D Y_r > 0 \}}\big| (\D Y_s)^+ \big|^{q-2} |\D Z_r|^2 dr     \\
 & & \hspace{-0.6cm} \le \neg      q  \neg \int_s^{s'}  \neg  \b1_{\{ \D Y_r > 0 \}}
 \big| (\D Y_r)^+ \big|^{q-1} \Big[ \Big(\g | \D Y_r|  \neg + \neg  \g |\D Z_r|  \neg + \neg \big( \ff_1 (r, Y^2_r,Z^2_r)
  \neg - \neg  \ff_2 (r, Y^2_r,Z^2_r) \big)^+ \Big) dr \neg -  \neg  \D Z_r d B^t_r  \Big]    ,
  ~   \fa t \neg  \le  \neg   s  \neg \le \neg  s'  \neg \le \neg  T   ,
  \eeas
   which is similar to \eqref{eq:n101} except that   $\varpi $ is specified   by $ q $.
   Then using  similar arguments to those
   that lead to \eqref{eq:s541}, we can deduce that
   \beas
  0  \le  E_t \bigg[ \underset{s \in [t,T]}{\sup} \big| (\D Y_s)^+ \big|^q \bigg]
 \le c_0  E_t \bigg[   \big| ( \xi_1 - \xi_2  )^+ \big|^q \bigg]
 + c_0 E_t \bigg[ \bigg( \int_t^T \big( \ff_1 (r, Y^1_r,Z^1_r)
 \neg - \neg \ff_2 (r, Y^2_r,Z^2_r) \big)^+ dr \bigg)^q \bigg] = 0 .
   \eeas
  Therefore, it holds $P^t_0-$a.s. that   $  (\D Y_s)^+ = 0$, or $ Y^1_s \le Y^2_s$
  for any $s \in [t,T]$.

 \ss  Next, we consider the case of reflected BSDEs with upper obstacles: For either $i=1$ or $i=2$, as
 $\big(Y^i ,Z^i ,V^i,K^i  \big)$ solves $\ol{\hb{R}}$BSDE$ \big( P^t_0,\xi_i,\ff_i,L^i  \big)$,
  the quadruplet $\big(\wh{Y}^i ,\wh{Z}^i ,\wh{V}^i,\wh{K}^i      \big)
 \dfnn \big(-Y^{3-i} ,-Z^{3-i} ,-V^{3-i},-K^{3-i}\big)$ solves the $\ul{\hb{R}}$BSDE $
 \big( P^t_0, \wh{\xi}_i,  \wh{\ff}_i,\wh{L}^i  \big)$  with $\wh{\xi}_i     \dfnn   - \xi_{3-i} $,  $   \wh{L}^i   \dfnn   - L^{3-i}_s  $ and the generator
  \beas
           \wh{\ff}_i (s,\o,y,z)    \dfnn    - \ff_{3-i} (s,\o,-y,-z),  \q
  \fa (s,\o,y,z)   \in     [t,T]    \times    \O^t    \times    \hR    \times    \hR^d   .
  \eeas

 It holds $P^t_0-$a.s. that   $ \wh{\xi}_1 - \wh{\xi}_2 = - \xi_2 + \xi_1 \le 0 $  and
 $ \wh{L}^1_s - \wh{L}^2_s = - L^2_s + L^1_s \le 0 $ for any $ s \in [t,T]$. Also, the process
$ \wh{V}^1 - \wh{V}^2 = - V^2 + V^1 $
 is  decreasing. For either $i=1$ or $i=2$, if   $\ff_i$ satisfies \eqref{ff_Lip} and $  \ff_1 (s,  Y^{3-i}_s ,Z^{3-i}_s )
 \le  \ff_2 (s, Y^{3-i}_s ,Z^{3-i}_s )  $, $ds \times dP^t_0-$a.s., then
   $\wh{\ff}_{3-i}$ satisfies \eqref{ff_Lip} and
 $  \wh{\ff}_1 (s, \wh{Y}^i_s, \wh{Z}^i_s) - \wh{\ff}_2 (s, \wh{Y}^i_s, \wh{Z}^i_s)  =
 -  \ff_2 (s, Y^{3-i}_s ,Z^{3-i}_s ) + \ff_1 (s,  Y^{3-i}_s ,Z^{3-i}_s ) \le 0 $, $ds \times dP^t_0-$a.s. Hence, all
 conditions to compare $\ul{\hb{R}}$BSDE$ \big( P^t_0, \wh{\xi}_1, \wh{\ff}_1,\wh{L}^1  \big)$
 with $\ul{\hb{R}}$BSDE$ \big( P^t_0, \wh{\xi}_2, \wh{\ff}_2,\wh{L}^2  \big)$  are satisfied. Then
 we can conclude that $P^t_0-$a.s., $\wh{Y}^1_s - \wh{Y}^2_s  = -  Y^2_s +  Y^1_s $ for any $s \in [t,T]$.   \qed

  \ss \no {\bf Proof of Theorem \ref{thm_viscosity}:}

 \ss \no {\bf  1)}
     We first show that
       $\ol{w}_1 $   is a viscosity subsolution    of   \eqref{eq:PDE} with   Hamiltonian $\ol{H}_1$
       when  $\hU_0 = \underset{i \in \hN}{\cup} F_i$ for closed subsets $\{F_i \}_{i \in \hN}$ of  $  \hU$.
 \if{0}
Clearly, $\cU^T = \hU_0$ and $\fB^T =  \wh{\fB}^T  $  collects all functions $\beta: \hU_0 \to \hV_0$ with  linear growth.
It follows that
  \beas
 \ol{w}_1(T,x) =    \wh{w}_1(T,x)=  \underset{\beta \in \wh{\fB}^T}{\inf} \, \underset{u \in \hU_0}{\sup}  \;
  \wt{Y}^{T,x,u, \beta ( u )}_T \Big(T,h\big( \wt{X}^{T,x,u, \beta ( u )}_T \big)  \Big)
  =  \underset{\beta \in \wh{\fB}^T}{\inf} \, \underset{u \in \hU_0}{\sup}  \;  \wt{Y}^{T,x,u, \beta ( u ) }_T \big(T,h(x)  \big)  = h(x), \q \fa x \in \hR^k
  \eeas
 \fi
  Let $(t_0,x_0, \vf) \in (0,T) \times \hR^k \times \hC^{1,2}\big([0,T] \times \hR^k\big)$
  be  such that $ \ol{w}_1  (t_0,x_0) = \vf (t_0,x_0)  $  and that  $\ol{w}_1-\vf$ attains a strict local   maximum  at $(t_0,x_0)$, i.e., for some   $\d_0 \in  \big(0,   t_0 \land ( T-t_0 ) \big) $
   \bea   \label{eq:a071}
     (\ol{w}_1 - \vf) (t,x) <  (\ol{w}_1 - \vf) (t_0,x_0) = 0,
   \q \fa (t,x) \in O_{\d_0}  (t_0,x_0) \big\backslash \big\{ (t_0,x_0) \big\} .
   \eea
   Let us simply  denote $\big(\vf(t_0   ,x_0),
     D_x \vf    (t_0  ,x_0),  D^2_x \vf    (t_0  ,x_0)\big) $ by $(y_0,z_0,\G_0 )$. Since $    \ul{l}(t_0,x_0) \le \vf (t_0,x_0) = \ol{w}_1 (t_0,x_0) \le \ol{l}(t_0,x_0) $  by \eqref{eq:n411}, it is clear that
 \beas
    \min   \bigg\{ \neg (\vf  \neg - \neg \ul{l}) (t_0   ,x_0),
 \max   \Big\{  \neg  -  \neg  \frac{\pa \vf }{\pa t}(t_0   ,x_0)
      \neg -  \neg     \ol{H}_{\neg 1} \neg \big(t_0   ,x_0, y_0,z_0,\G_0 \big)  ,
  (\vf  \neg - \neg \ol{l}) (t_0 ,x_0) \Big\}  \bigg\} \le 0
 \eeas
 if     $\ol{H}_1 \big(t_0,x_0, y_0,z_0,\G_0 \big) =  \infty$.

 \ss  To draw a contradiction,
 we assume that when $\ol{H}_1 \big(t_0,x_0, y_0,z_0,\G_0 \big) < \infty $,
 \bea   \label{eq:s637}
  \q  \varrho \dfnn \min   \bigg\{ \neg (\vf  \neg - \neg \ul{l}) (t_0   ,x_0),
 \max   \Big\{  \neg  -  \neg  \frac{\pa \vf }{\pa t}(t_0   ,x_0)
      \neg -  \neg     \ol{H}_{\neg 1} \neg \big(t_0   ,x_0, y_0,z_0,\G_0 \big)  ,
  (\vf  \neg - \neg \ol{l}) (t_0 ,x_0) \Big\}  \bigg\} >   0   .
 \eea
   Then   the continuity of $\vf$ and $\ol{l}$ implies that for some   $ \d_1 \in  (0, \d_0    )$
 \bea  \label{eq:s423}
    (\vf   \neg -  \neg \ul{l} )(t,x)  \neg \ge \neg  \frac34 \varrho , ~\;\; \fa (t,x)  \neg  \in  \neg  \ol{O}_{\d_1}  (t_0,x_0) . \q
 \eea
   As $  \vf (t_0,x_0)   \le \ol{l}(t_0,x_0)   $, we also see from \eqref{eq:s637}   that
    $  \neg  -  \neg  \frac{\pa \vf }{\pa t}(t_0,x_0) \neg -  \neg     \ol{H}_1 \neg \left(t_0,x_0,
 y_0,z_0,\G_0 \right) \ge   \varrho  $.
   Thus, one can find an  $m   \in \hN$    such that
  \bea \label{eq:s629}
  \neg  -  \neg  \frac{\pa \vf }{\pa t}(t_0,x_0) - \frac78 \, \varrho    \ge
             \underset{u \in \hU_0}{\sup} \;  \underset{v \in \sO^m_u}{\inf}  ~\, \lsup{\hU_0 \ni u' \to u}\; \underset{(t,x,y,z,\G)
   \in O_{1/m} (t_0,x_0,y_0,z_0,\G_0)}{\sup} \;  H (t ,x ,y ,z ,\G ,u',v) .
 \eea
As $\vf \in  \hC^{1,2}\big([0,T] \times \hR^k\big)$, there exists a $\d < \frac{1}{2m} \land \d_1$ such that for any $ (t,x) \in \ol{O}_\d (t_0,x_0)$
   \bea
 &&   \Big| \frac{\pa \vf }{\pa t}(t,x)   -    \frac{\pa \vf }{\pa t}(t_0,x_0) \Big|    \le    \frac18 \,\varrho
             ~\;\;   \label{eq:s631} \\
  \hb{and} &&
 \big| \vf  (t,x)     -    \vf (t_0,x_0) \big|    \vee    \Big| D_x \vf    (t  ,x)
   -    D_x \vf    (t_0  ,x_0)\Big|
    \vee    \Big| D^2_x \vf    (t  ,x)    -    D^2_x \vf    (t_0  ,x_0)\Big|    \le    \frac{1}{2m} \,, \nonumber
      \eea
 the latter of which together with \eqref{eq:s629} implies that
  \beas
  \neg  -  \neg  \frac{\pa \vf }{\pa t}(t_0,x_0) - \frac78 \, \varrho    \ge
             \underset{u \in \hU_0}{\sup} \;  \underset{v \in \sO^m_u}{\inf}  ~\, \lsup{\hU_0 \ni u' \to u}\;
 \underset{(t,x )     \in \ol{O}_\d  (t_0,x_0)}{\sup} \;  H (t ,x , \vf(t   ,x),
     D_x \vf    (t  ,x),  D^2_x \vf    (t  ,x)  ,u',v)   .
 \eeas

  Then  for any $ u \in    \hU_0 $,
 there exists a $ \fP (u) \in \sO^m_u$ such that
 \beas
   \neg  -  \neg  \frac{\pa \vf }{\pa t}(t_0,x_0) - \frac34 \, \varrho \ge
  \lsup{\hU_0 \ni u' \to u}\; \underset{(t,x )     \in \ol{O}_\d  (t_0,x_0)}{\sup} \;  H \big(t ,x , \vf(t   ,x),
     D_x \vf    (t  ,x),  D^2_x \vf    (t  ,x)  ,u', \fP(u) \big)
 \eeas
 and we can find   a $ \l(u) \in \big(0 , 1     \big)  $ such that
 for any $ u'   \in    \hU_0  \cap  O_{\l(u)}(u)  $
\bea \label{eq:s613}
   \neg  -  \neg  \frac{\pa \vf }{\pa t}(t_0,x_0) - \frac58 \, \varrho \ge
    \underset{(t,x )     \in \ol{O}_\d  (t_0,x_0)}{\sup} \;  H \big(t ,x , \vf(t   ,x),
     D_x \vf    (t  ,x),  D^2_x \vf    (t  ,x)  ,u', \fP(u)\big)    .
 \eea

 Let $\fO (u) \dfnn  O_{\l(u)}(u)$, $ u \in    \hU_0 $.  For any $i \in \hN$,  $\{ \fO (u) \}_{u \in F_i}$ together with $\hU \, \backslash \, F_i$
  form an open cover of $\hU$.      Since the separable metric space  $ \hU$ is   Lindel\"of,
 one can find a  sequence $\big\{ u^i_j \big\}_{j \in \hN}$ of $ F_i $  such that $F_i \subset \underset{j \in \hN}{\cup} \, \fO(u^i_j )$.     Let $\{ u_\ell \}_{\ell \in \hN} $
  represent the countable set  $ \big\{  u^i_j  \big\}_{i,j \in \hN} \subset \hU_0 $
  and let $\hat{v}$ be an arbitrary element of $\hV_0$.  It is  clear that
 \beas
  \fP(u) \dfnn \sum_{ \ell \in \hN} \b1_{\big\{ u \in \fO (u_\ell) \backslash \,
  \underset{\ell'<\ell}{\cup} \fO (u_{\ell'}) \big\}} \fP (u_\ell)
 + \b1_{\big\{ u \in \hU  \backslash \, \underset{\ell \in \hN}{\cup} \fO_\ell \big\}} \hat{v}  \in \hV_0  ,  \q  \fa     u  \in   \hU
 \eeas
   defines a $   \sB(\hU) /\sB(\hV)-$measurable function.

  \ss  For any $  u  \in    \hU_0$,   since   $ \hU_0 = \underset{i \in \hN}{\cup} F_i \subset
   \underset{i,j \in \hN}{\cup} \fO(u^i_j ) =  \underset{\ell \in \hN}{\cup} \fO(u_\ell ) $,
    there exists a $\ell \in \hN$ such that $u \in \fO (u_\ell) \backslash \,
  \underset{\ell'<\ell}{\cup} \fO (u_{\ell'})$.   We see from \eqref{eq:s613} that
 \beas
  \neg  -  \neg  \frac{\pa \vf }{\pa t}(t_0,x_0) - \frac58 \, \varrho &\ge&
    \underset{(t,x )     \in \ol{O}_\d  (t_0,x_0)}{\sup} \;  H \big(t ,x , \vf(t   ,x),
     D_x \vf    (t  ,x),  D^2_x \vf    (t  ,x)  ,u, \fP(u_\ell) \big)  \\
  &=&  \underset{(t,x )     \in \ol{O}_\d  (t_0,x_0)}{\sup} \;  H (t ,x , \vf(t   ,x),
     D_x \vf    (t  ,x),  D^2_x \vf    (t  ,x)  ,u, \fP(u) )  ,
 \eeas
which together with \eqref{eq:s631} implies that
 \bea  \label{eq:s617}
  \neg  -  \neg  \frac{\pa \vf }{\pa t}(t ,x ) - \frac12 \, \varrho \ge
          H (t ,x , \vf(t   ,x),
     D_x \vf    (t  ,x),  D^2_x \vf    (t  ,x)  ,u, \fP(u) )  , \q \fa (t,x) \in \ol{O}_\d(t_0,x_0) , ~ \fa u \in   \hU_0.
 \eea

 Let  $     \wp    \dfnn  \inf \big\{ (     \vf \neg - \neg \ol{w}_1) (t,x)  \neg : (t,x )  \neg \in \neg
     \ol{O}_\d  (t_0,x_0) \big\backslash O_{\frac{\d}{2}} (t_0,x_0) \big\}       $.
 Since the set $ \ol{O}_\d  (t_0,x_0) \big\backslash O_{\frac{\d}{2}} (t_0,x_0)$ is  compact, there exists
 a sequence $\{(t_n,x_n)\}_{n \in \hN} $ of $  \ol{O}_\d  (t_0,x_0) \big\backslash O_{\frac{\d}{2}} (t_0,x_0)$
 that  converges to some $(t_*,x_*) \in \ol{O}_\d  (t_0,x_0) \big\backslash O_{\frac{\d}{2}} (t_0,x_0) $
 and satisfies   $ \wp \neg = \neg  \lmtd{n \to \infty} ( \vf \neg - \neg \ol{w}_1 ) (t_n,x_n) $.
  The continuity of $\vf$ and the upper semicontinuity of $\ol{w}_1$ and imply that
 $\vf \neg - \neg \ol{w}_1$ is also lower semicontinuous.
 Thus, it follows that $ \wp \le (\vf \neg - \neg \ol{w}_1)(t_*,x_*) \le \lmtd { n \to \infty } (\vf \neg - \neg \ol{w}_1)(t_n,x_n ) = \wp  $,
 which together with \eqref{eq:a071} shows that
  \bea   \label{eq:s427}
  \wp  = \min \big\{ (\vf \neg - \neg \ol{w}_1) (t,x): (t,x ) \in
     \ol{O}_\d  (t_0,x_0) \big\backslash O_{\frac{\d}{2}} (t_0,x_0) \big\} = (\vf \neg - \neg \ol{w}_1)(t_*,x_*)  >  0 .
  \eea

  Then  we set     $\dis   \wt{\wp} \neg \dfnn  \neg   \frac{    \wp \land \varrho}{ 2 (1 \vee \g) T }  \neg   > \neg  0$
  and let $ \big\{(t_j, x_j)\big\}_{j \in \hN}$ be a sequence of $O_{\frac{\d}{4}} (t_0,x_0) $ such   that
     \beas   
 \lmt{j \to \infty} (t_j,x_j) = (t_0,x_0)   \q \hb{and} \q  \lmt{j \to \infty} \wh{w}_1 (t_j,x_j) = \ol{w}_1 (t_0,x_0) = \vf (t_0,x_0) .
     \eeas
    As $\lmt{j \to \infty} \big( \wh{w}_1 (t_j,x_j)- \vf (t_j,x_j) \big) = 0$, it holds for  some $j \in \hN$ that
    \bea  \label{eq:s437}
  \big| \wh{w}_1 (t_j,x_j)- \vf (t_j,x_j) \big| < \frac12 \wt{\wp}   t_0 .
   \eea

  In particular,   $ \fP (t,\o,u) \dfnn  \fP (u) $, $(t,\o,u) \in [t_j,T] \times \O^{t_j} \times \hU $
  is a $\hV_0-$valued, $ \sP(\bF^{t_j}) \otimes \sB(\hU) /\sB(\hV)-$measurable function.
   For any $(t,\o, u)  \in [t_j,T] \times \O^{t_j} \times \hU$, if
         $ u \in \fO_i \backslash \underset{j<i}{\cup} \fO_j $ for some $i \in \hN$, then
 \beas 
  \q  \big[ \fP (t, \o, u) \big]_{\overset{}{\hV}} = \big[\fP( u_i)\big]_{\overset{}{\hV}} \le   m  \neg  + \neg    m [u_i]_{\overset{}{\hU}}
  \neg  \le  \neg   m  \neg  + \neg   m [u]_{\overset{}{\hU}}
  \neg  + \neg   m \rho_{\overset{}{\hU}} (u, u_i)
 \neg  < \neg   m  \neg  + \neg   m [u]_{\overset{}{\hU}}  \neg  + \neg   m \l ( u_i)
  \neg  <  \neg      2m  \neg  + \neg   m [u]_{\overset{}{\hU}} \, ;
 \eeas
 otherwise, $\big[ \fP (t, \o, u) \big]_{\overset{}{\hV}} = \big[\hat{v}\big]_{\overset{}{\hV}}$\,.
 This shows that $ \fP$ satisfies \eqref{eq:r503} with $\Psi = 2 m + \big[\hat{v}\big]_{\overset{}{\hV}}$ and
 $\k =m$.
   Clearly,    \eqref{eq:r742} automatically holds for $\fP $. Hence, $\fP \in \wh{\fB}^{t_j}$.

  For any   $\mu \in   \cU^{t_j} $,
     we set $\Th_\mu \dfnn \big( t_j, x_j, \mu ,    \fP \lan  \mu \ran \big)$
  and  define  two  $\bF^{t_j}-$stopping times:
 \beas
  \q  \t_\mu   \dfnn  \inf \Big\{s  \neg \in \neg  ( t_j, T  ] \neg  : \big( s,\wt{X}^{\Th_\mu}_s \big)
  \neg  \notin   \neg   \ol{O}_{\frac34 \d} (t_0, x_0) \Big\}
 ~    \hb{and} ~   \z_\mu  \dfnn  \inf \Big\{s  \neg \in \neg  (  \t_\mu, T] \neg : \big( s, \wt{X}^{\Th_\mu}_s \big)
  \neg   \notin   \neg   \ol{O}_\d  (t_0,x_0) \big\backslash O_{\frac{\d}{2}} (t_0,x_0) \Big\} \land T .
 \eeas
    Since 
 $ \big| \big(T, \wt{X}^{\Th_\mu}_T\big) -(t_0,x_0) \big| \ge T -t_0 > \d_0  > \d_1 > \frac34 \d $,
one can deduce from the continuity of $ \wt{X}^{\Th_\mu} $ that
  \bea     \label{eq:a041}
   \t_\mu  < T  \q \hb{and} \q  \big( \t_\mu, \wt{X}^{\Th_\mu}_{\t_\mu} \big)  \in  \pa O_{\frac34 \d}(t_0, x_0) ,
 \q  P^{t_j}_0-a.s.
  \eea

 \ss     Given $n  \neg \in  \neg  \hN$,  we define
     $ q^n  \neg  (s)  \neg \dfnn \neg  \frac{\lceil 2^n s \rceil}{2^n}  \land  T$, $ s  \neg \in \neg  [0,T] $.
    Then  $ \t^n_\mu   \dfnn     q^n(\t_\mu)  \neg \land  \neg  \z_\mu$ is an  $\bF^{t_j}-$stopping time.
  Applying  \eqref{eq:p677} with $(\z,\t,\xi) \neg = \neg \Big(\t^n_\mu ,\, q^n(\t_\mu) ,  \,  \wh{w}_1   \big( q^n(\t_\mu),
  \wt{X}^{\Th_\mu}_{q^n(\t_\mu)}\big) \Big)$, we can deduce from Proposition \ref{prop_apriori_DRBSDE}
  and H\"older's inequality that
      \bea
   && \hspace{- 1.5 cm} \left|     \wt{Y}^{\Th_\mu}_{t_j}
   \Big( \t^n_\mu,   \wh{w}_1   \big( \t^n_\mu, \wt{X}^{\Th_\mu}_{\t^n_\mu} \big)     \Big)
   -\wt{Y}^{\Th_\mu}_{t_j}
  \Big(q^n(\t_\mu),   \wh{w}_1   \big( q^n(\t_\mu), \wt{X}^{\Th_\mu}_{q^n(\t_\mu)}\big)\Big)\right|^{\frac{q+1}{2}}
   \nonumber  \\
& = & \left|     \wt{Y}^{\Th_\mu}_{t_j}
   \Big( \t^n_\mu,   \wh{w}_1   \big( \t^n_\mu, \wt{X}^{\Th_\mu}_{\t^n_\mu} \big)     \Big)
   -\wt{Y}^{\Th_\mu}_{t_j}
  \bigg(\t^n_\mu,  \wt{Y}^{\Th_\mu}_{\t^n_\mu}
  \Big(q^n(\t_\mu), \wh{w}_1   \big( q^n(\t_\mu), \wt{X}^{\Th_\mu}_{q^n(\t_\mu)}\big)\Big)\bigg)\right|^{\frac{q+1}{2}}
  \nonumber  \\
  & \le &  c_0   E_{t_j} \bigg[ \Big| \wh{w}_1   \big( \t^n_\mu, \wt{X}^{\Th_\mu}_{\t^n_\mu} \big)
- \wt{Y}^{\Th_\mu}_{\t^n_\mu}
  \Big(q^n(\t_\mu), \wh{w}_1 \big( q^n(\t_\mu), \wt{X}^{\Th_\mu}_{q^n(\t_\mu)}\big)\Big) \Big|^{\frac{q+1}{2}} \bigg]
  \nonumber   \\
& = &  c_0   E_{t_j} \bigg[ \b1_{\{  q^n(\t_\mu) > \z_\mu \}}  \Big| \wh{w}_1   \big(   \z_\mu,
 \wt{X}^{\Th_\mu}_{  \z_\mu} \big)  -  \wt{Y}^{\Th_\mu}_{  \t^n_\mu }
  \Big(q^n(\t_\mu),   \wh{w}_1   \big( q^n(\t_\mu), \wt{X}^{\Th_\mu}_{q^n(\t_\mu)}\big)\Big)
\Big|^{\frac{q+1}{2}} \bigg]  \nonumber  \\
 &\le & c_0  \Big( P^t_0 \big( q^n(\t_\mu) > \z_\mu \big) \Big)^{\frac{q - 1}{2q}}
 \Bigg\{ E_{t_j} \bigg[    \Big| \wh{w}_1   \big(   \z_\mu,
 \wt{X}^{\Th_\mu}_{  \z_\mu} \big)  -  \wt{Y}^{\Th_\mu}_{  \t^n_\mu }
  \Big(q^n(\t_\mu),   \wh{w}_1   \big( q^n(\t_\mu), \wt{X}^{\Th_\mu}_{q^n(\t_\mu)}\big)\Big)
\Big|^q \bigg] \Bigg\}^{\frac{q+1}{2q}}   . \label{eq:s416}
    \eea
 By   \eqref{2l_growth},
 \bea \label{eq:s418}
  \Big| \wh{w}_1   \big(   \z_\mu, \wt{X}^{\Th_\mu}_{  \z_\mu} \big) \Big|
\le  \big( |\ul{l} | \vee |\ol{l} | \big)  \Big(\z_\mu, \wt{X}^{\Th_\mu}_{  \z_\mu}   \Big)
 \le  \ul{l}_*  + \ol{l}_* +  \g \underset{s \in [t_j,T]}{\sup} \big| \wt{X}^{\Th_\mu}_s \big|^{2/q} ,
 \eea
  where $ \ul{l}_* \dfnn \underset{s \in [t,T]}{\sup}   \big| \ul{l} (s, 0 ) \big|$
 and $ \ol{l}_* \dfnn \underset{s \in [t,T]}{\sup}   \big| \ol{l} (s, 0 ) \big|$.
 Similarly,
  \beas
  \q  \bigg| \wt{Y}^{\Th_\mu}_{  \t^n_\mu  }
  \Big(q^n(\t_\mu),   \wh{w}_1   \big( q^n(\t_\mu), \wt{X}^{\Th_\mu}_{q^n(\t_\mu)}\big)\Big) \bigg|
    \le  \Big( \big|\ul{L}^{\Th_\mu}_{\t^n_\mu} \big| \vee \big|\ol{L} ^{\Th_\mu}_{\t^n_\mu} \big| \Big)
    \neg  =  \neg  \big(  |\ul{l}  | \vee  |\ol{l}  | \big) \Big(\t^n_\mu, \wt{X}^{\Th_\mu}_{ \t^n_\mu }   \Big)
       \le     \ul{l}_*  \neg  +  \neg \ol{l}_*  \neg +  \neg  \g \underset{s \in [t_j,T]}{\sup} \big| \wt{X}^{\Th_\mu}_s \big|^{2/q} .
\eeas
 Putting it and \eqref{eq:s418} back into \eqref{eq:s416} gives that
  \bea
  && \hspace{-2.5cm} \left|     \wt{Y}^{\Th_\mu}_{t_j}
   \Big( \t^n_\mu,   \wh{w}_1   \big( \t^n_\mu, \wt{X}^{\Th_\mu}_{\t^n_\mu} \big)     \Big)
   -\wt{Y}^{\Th_\mu}_{t_j}
  \Big(q^n(\t_\mu),   \wh{w}_1   \big( q^n(\t_\mu), \wt{X}^{\Th_\mu}_{q^n(\t_\mu)}\big)\Big)\right|^{\frac{q+1}{2}}
  \nonumber \\
  &&   \le   c_0  \Big( P^t_0 \big(q^n( \t_\mu ) > \z_\mu \big) \Big)^{\frac{q - 1}{2q}}
  \Bigg\{ 1+ E_{t_j} \bigg[  \underset{s \in [t_j,T]}{\sup} \big| \wt{X}^{\Th_\mu}_s \big|^2  \bigg] \Bigg\}^{\frac{q+1}{2q}} .
  \label{eq:s420}
     \eea
 Since $ \t_\mu < \z_\mu $, $P^t_0-$a.s. by \eqref{eq:a041} and since $ \lmtd{n \to \infty} q^n(\t_\mu) = \t_\mu $,
   we see from \eqref{eq:esti_X_1} that the right-hand-side of \eqref{eq:s420} converges to $0$ as $n \to \infty$.
 Hence,   for some $ n_{\mu} \in \hN$
   \bea  \label{eq:s431}
  \left| \wt{Y}^{\Th_\mu}_{t_j}
   \Big(  \wh{\t}_\mu ,   \wh{w}_1   \big( \wh{\t}_\mu, \wt{X}^{\Th_\mu}_{\wh{\t}_\mu} \big)     \Big)
    -  \wt{Y}^{\Th_\mu}_{t_j}   \Big( q^{n_\mu}(\t_\mu),
 \wh{w}_1 \big( q^{n_\mu}(\t_\mu), \wt{X}^{\Th_\mu}_{q^{n_\mu}(\t_\mu)} \big) \Big) \right|
   <  \frac14 \wt{\wp} \,  t_0   ,
   \eea
where $ \wh{\t}_\mu \dfnn  \t^{n_\mu}_\mu = q^{n_\mu}(\t_\mu) \land \z_\mu $.

     As    $\wh{\t}_\mu$ is an $\bF^{t_j}-$stopping time,
      the continuity of $\vf$ and $\wt{X}^{\Th_\mu} $ show that
  $\sY^\mu_s \dfnn \vf \big( \wh{\t}_\mu \land s, \wt{X}^{\Th_\mu}_{\wh{\t}_\mu \land s} \big)
    - \wt{\wp} ( \wh{\t}_\mu \land s) $, $ s \in [t_j, T] $ defines a real-valued, $\bF^{t_j}-$adapted continuous process.
 Applying It\^o's formula to $ \sY^\mu $ yields that
    \bea
         \sY^\mu_s   & \tneg = &  \tneg   \sY^\mu_T     +   \int_s^T  \ff^\mu_r   dr
    -  \int_s^T  \sZ^\mu_r   d B^{t_j}_r ,   \q    s \in [t_j,T]  ,   \label{eq:s403}
    \eea
where $\sZ^\mu_r  \dfnn        \b1_{\{r < \wh{\t}_\mu\}}    D_x \vf    \big( r, \wt{X}^{\Th_\mu}_r \big)
  \cd   \si \big( r, \wt{X}^{\Th_\mu}_r, \mu_r,   \big(\fP \lan \mu \ran \big)_r \big)$ and
 $$
 \hspace{-3mm}   \ff^\mu_r   \neg   \dfnn   \neg       \b1_{\{ r \, < \wh{\t}_\mu  \}}
  \bigg\{ \wt{\wp}      -       \frac{\pa \vf }{\pa t} \big( r, \wt{X}^{\Th_\mu}_r \big)
      -       D_x \vf    \big( r, \wt{X}^{\Th_\mu}_r\big)   \cd    b \big( r, \wt{X}^{\Th_\mu}_r \neg , \mu_r,
  (\fP \lan \mu \ran )_r \big)
         -         \frac12  trace\Big( \si \si^T \big( r, \wt{X}^{\Th_\mu}_r \neg , \mu_r,
  (\fP \lan \mu \ran)_r \big) \cd    D^2_x \vf \big( r, \wt{X}^{\Th_\mu}_r\big) \Big) \neg \bigg\} .
     $$
 As $\vf  \neg \in  \neg   \hC^{1,2}\big([t,T] \neg  \times \neg  \hR^k\big)$,
  the measurability of $b$, $\si$ and $\fP$ show   that  both $\sZ^\mu$ and $\ff^\mu$ are
     $\bF^{t_j}-$progressively    measurable.

  \ss   Since it holds $P^t_0-$a.s. that
   \bea   \label{eq:a239}
  \big(\wh{\t}_\mu \land s , \wt{X}^{\Th_\mu}_{\wh{\t}_\mu \land s}   \big) \in   \ol{O}_\d  (t_0, x_0) ,
  \q \fa    s \in [t_j, T]  ,
   \eea
 we see from the continuity  of $\vf$ that $\sY^\mu$ is a bounded process, and
we can deduce from \eqref{b_si_linear_growth}, \eqref{b_si_Lip} as well as H\"older's inequality  that
  \bea
 && \hspace{-0.7cm}  E_{t_j} \left[  \bigg( \int_{t_j}^T \neg |\sZ^\mu_s|^2   \,  ds  \bigg)^{ q /2 } \right]
    =    (\g C_\vf)^q \,  E_{t_j} \left[  \bigg( \int_{t_j}^{\wh{\t}_\mu}
\Big( 1+ \big| \wt{X}^{\Th_\mu}_s \big|  + [\mu_s]_{\overset{}{\hU}}
  +   \big[ (\fP \lan \hat{\mu} \ran)_s \big]_{\overset{}{\hV}}   \Big)^2   ds  \bigg)^{ q /2 } \right]  \label{eq:s641}   \\
 &&    \le \neg   c_0  C^q_\vf  \bigg\{   \big(1  \neg + \neg  |x_0 |  \neg + \neg  \d    \big)^q
  \neg + \neg   \bigg( E_{t_j}   \neg   \int_{t_j}^T \dneg [\mu_s]^2_{\overset{}{\hU}} \,     ds    \bigg)^{ q /2 }
  \neg + \neg   \bigg( E_{t_j}    \neg  \int_{t_j}^T \dneg \big[ (\fP \lan \hat{\mu} \ran)_s \big]_{\overset{}{\hV}}^2 \,   ds
      \bigg)^{ q /2 } \bigg\}  \neg < \neg  \infty , ~ \hb{i.e. }
 \sZ^\mu \in \hH^{2,q}_{\bF^{t_j}  } \big( [t,T], \hR^d \big)   ,     \nonumber
  \eea
 where  $ C_\vf \dfnn \underset{(t,x) \in \ol{O}_\d (t_0,x_0) }{\sup} \big|D_x \vf (t,x)\big|$.  Moreover,
    \eqref{eq:a239} and \eqref{eq:s423} imply that $P^t_0-$a.s.
  \beas           
   \sY^\mu_s   \ge  \ul{l} ( \wh{\t}_\mu \land s, \wt{X}^{\Th_\mu}_{\wh{\t}_\mu \land s}) + \frac 34 \varrho - \wt{\wp} T
 > \ul{l} ( \wh{\t}_\mu \land s, \wt{X}^{\Th_\mu}_{\wh{\t}_\mu \land s}) = \ul{L}^{\Th_\mu}_{\wh{\t}_\mu \land s} ,
    \q  \fa    s \in [t_j, T]   ,
  \eeas
    which together with   \eqref{eq:s403} shows  that
 $\big\{\big( \sY^\mu_s, \sZ^\mu_s, 0, 0 \big) \big\}_{s \, \in  [t_j,T]}$
  solves the  $\ul{\hb{R}}$BSDE$\left( P^{t_j}_0 , \sY^\mu_T , \ff^\mu, \ul{L}^{\Th_\mu}_{\wh{\t}_\mu \land \cd} \right)$
 (see \eqref{RLUBSDE}).

 Since $\t_\mu \le \wh{\t}_\mu \le \z_\mu $ and since $ \ol{w}_1(t,x) \ge \wh{w}_1 (t,x)$
 for any $(t,x) \in [0,T] \times \hR^k$, we can deduce from  \eqref{eq:s427} that
      \beas
   \sY^\mu_T  \ge    \vf \big( \wh{\t}_\mu, \wt{X}^{\Th_\mu}_{\wh{\t}_\mu} \big)    -   \wt{\wp} T
 >  \vf \big( \wh{\t}_\mu, \wt{X}^{\Th_\mu}_{\wh{\t}_\mu} \big) - \wp \ge   \ol{w}_1  \big( \wh{\t}_\mu, \wt{X}^{\Th_\mu}_{\wh{\t}_\mu} \big)
 \ge \wh{w}_1 \big( \wh{\t}_\mu, \wt{X}^{\Th_\mu}_{\wh{\t}_\mu} \big) .
  \eeas
  Also,      \eqref{eq:a239}, \eqref{eq:s617}  and \eqref{f_Lip}  show   that for $ds \times d P^{t_j}_0-$a.s.
  $ (s,\o)  \in [t_j,T] \times \O^{t_j}$
     \beas
     \ff^\mu_s (\o)  & \tneg \ge &  \tneg  
 \b1_{\{ s \, < \wh{\t}_\mu (\o) \}} \left\{  \wt{\wp} + \frac12 \varrho   +    f \neg \left( s,\o, \wt{X}^{\Th_\mu}_s (\o),
       \sY^\mu_s (\o) -\wt{\wp}s,    \sZ^\mu_s (\o), \mu_s (\o),  ( \fP \lan \mu \ran )_s (\o) \right)  \right\}  \\
         &  \tneg \ge &  \tneg  \b1_{\{ s \, < \wh{\t}_\mu (\o) \}} \left\{  \wt{\wp} + \frac12 \varrho  - \g \wt{\wp} T
          +    f \neg \left(s, \o, \wt{X}^{\Th_\mu}_s (\o),
      \sY^\mu_s (\o),       \sZ^\mu_s  (\o), \mu_s (\o),  ( \fP \lan \mu \ran )_s  (\o) \right)  \right\}   \\
   &  \tneg  \ge  &  \tneg   f^{\Th_\mu }_{\wh{\t}_\mu} \big(s ,\o, \sY^\mu_s  (\o),  \sZ^\mu_s (\o) \big) .
       \eeas
   Clearly,
 $
  \left( \wt{Y}^{\Th_\mu} \Big( \wh{\t}_\mu,   \wh{w}_1  \big( \wh{\t}_\mu, \wt{X}^{\Th_\mu}_{\wh{\t}_\mu} \big)    \Big),
  \wt{Z}^{\Th_\mu} \Big( \wh{\t}_\mu,   \wh{w}_1  \big( \wh{\t}_\mu, \wt{X}^{\Th_\mu}_{\wh{\t}_\mu} \big)
   \Big),-\wt{\ol{K}}^{\,\raisebox{-0.5ex}{\scriptsize $\Th_\mu$}}
  \Big( \wh{\t}_\mu,   \wh{w}_1  \big( \wh{\t}_\mu, \wt{X}^{\Th_\mu}_{\wh{\t}_\mu} \big),
  \wt{\ul{K}}^{\Th_\mu} \Big( \wh{\t}_\mu,
  \wh{w}_1  \big( \wh{\t}_\mu, \wt{X}^{\Th_\mu}_{\wh{\t}_\mu} \big)    \Big)    \Big) \right)
 $
 solves $\ul{\hb{R}}$BSDE$\left( P^{t_j}_0 , \wh{w}_1  \big( \wh{\t}_\mu, \wt{X}^{\Th_\mu}_{\wh{\t}_\mu} \big) ,
 f^{\Th_\mu}_{\wh{\t}_\mu} , \ul{L}^{\Th_\mu}_{\wh{\t}_\mu \land \cd} \right)$.
 As   $f^{\Th_\mu}_{\wh{\t}_\mu}$ is  Lipschitz continuous in $(y, z)$, we know from Proposition
 \ref{prop_comparison_RBSDE} that $P^{t_j}_0-$a.s.
   \beas
   \sY^\mu_s   \ge \wt{Y}^{\Th_\mu}_s \left( \wh{\t}_\mu,
   \wh{w}_1  \big( \wh{\t}_\mu, \wt{X}^{\Th_\mu}_{\wh{\t}_\mu} \big)    \right)  , \q  \fa  s \in [t_j,T] .
   \eeas
  Letting $s=t_j$ and using \eqref{eq:s431}, we obtain
   \beas
       \vf   (t_j   , x_j  )  -   \frac34    \wt{\wp}   t_0
   &>&   \vf   (t_j   , x_j  )   -   \wt{\wp}   t_j   =   \sY^\mu_{t_j}
     \ge     \wt{Y}^{t_j, x_j, \mu ,    \fP \lan  \mu \ran}_{t_j} \Big( \wh{\t}_\mu,
   \wh{w}_1  \Big( \wh{\t}_\mu, \wt{X}^{t_j, x_j, \mu ,    \fP \lan  \mu \ran}_{\wh{\t}_\mu} \Big)    \Big)
 \nonumber  \\
      & > &       \wt{Y}^{t_j, x_j, \mu ,    \fP \lan  \mu \ran}_{t_j}   \Big( q^{n_\mu}(\t_\mu),
 \wh{w}_1 \Big( q^{n_\mu}(\t_\mu), \wt{X}^{t_j, x_j, \mu ,    \fP \lan  \mu \ran}_{q^{n_\mu}(\t_\mu)} \Big) \Big)
 - \frac14  \wt{\wp}   t_0    ,
      \eeas
 where we used the fact that $t_j \neg > \neg  t_0 -\frac14 \d  \neg > \neg  t_0 - \frac14 \d_0  \neg > \neg  \frac34 t_0 $.
 Taking supremum over $\mu \in \cU^{t_j} $ gives that
 \bea   \label{eq:s621}
   \vf   (t_j   , x_j  )  -   \frac34    \wt{\wp}   t_0 \ge  \underset{\mu \in \cU^{t_j}}{\sup}
  \wt{Y}^{t_j, x_j, \mu ,    \fP \lan  \mu \ran}_{t_j}   \Big( q^{n_\mu}(\t_\mu),
  \wh{w}_1 \Big( q^{n_\mu}(\t_\mu), \wt{X}^{t_j, x_j, \mu ,    \fP \lan  \mu \ran}_{q^{n_\mu}(\t_\mu)} \Big) \Big)
  - \frac14  \wt{\wp}   t_0 .
 \eea

 \ss Let $\{\t_{\mu,\beta} \neg : \mu \neg \in \neg \cU^{t_j}, \beta \neg \in  \neg \wh{\fB}^{t_j}\}$ be an arbitrary family
 of $\hQ_{t_j,T}  -$valued,
 $\bF^{t_j}-$stopping times such that $\t_{ \mu,\fP} = q^{n_\mu}(\t_\mu)   $ for any $\mu \in \cU^{t_j}$.
 Then \eqref{eq:s621}, \eqref{eq:s437}  and Theorem \ref{thm_DPP} imply that
  \beas
   \vf   (t_j   , x_j  )  -   \frac34    \wt{\wp}   t_0 &\ge&  \underset{\mu \in \cU^{t_j}}{\sup}
  \wt{Y}^{t_j, x_j, \mu ,    \fP \lan  \mu \ran}_{t_j}   \Big( q^{n_\mu}(\t_\mu),
  \wh{w}_1 \Big( q^{n_\mu}(\t_\mu), \wt{X}^{t_j, x_j, \mu ,    \fP \lan  \mu \ran}_{q^{n_\mu}(\t_\mu)} \Big) \Big)
  - \frac14  \wt{\wp}   t_0  \\
 &\ge &  \underset{\beta \in \wh{\fB}^{t_j}}{\inf} \underset{\mu \in \cU^{t_j}}{\sup}
  \wt{Y}^{t_j, x_j, \mu ,    \fP \lan  \mu \ran}_{t_j}   \Big( q^{n_\mu}(\t_\mu),
  \wh{w}_1 \Big( q^{n_\mu}(\t_\mu), \wt{X}^{t_j, x_j, \mu ,    \fP \lan  \mu \ran}_{q^{n_\mu}(\t_\mu)} \Big) \Big)
  - \frac14  \wt{\wp}   t_0 \\
& \ge &  \wh{w}_1 (t_j,x_j)  \neg - \neg  \frac14  \wt{\wp}   t_0
  \neg > \neg  \vf (t_j,x_j)  \neg - \neg  \frac34 \wt{\wp}   t_0 \,.
 \eeas
   A contradiction appears.
 Therefore,     $\ol{w}_1$  is a viscosity subsolution of   \eqref{eq:PDE} with Hamiltonian $\ol{H}_1$.

  Similarly, one can show that   $w^*_1$  is also a viscosity subsolution of   \eqref{eq:PDE} with Hamiltonian  $\ol{H}_1$
  when  $\hU_0 = \underset{i \in \hN}{\cup} F_i$ for closed subsets $\{F_i \}_{i \in \hN}$ of  $  \hU$.
  Using the transformation similar to that in the part (3) of proof of Theorem \ref{thm_DPP},
 one can show that if $ \hV_0 $ is a countable union of closed subsets of  $ \hV $, then
   $ \ul{w}_2 $ and $w^*_2$   are two   viscosity  supersolutions  of   \eqref{eq:PDE} with   Hamiltonian  $ \ul{H}_2 $.

\ms \no {\bf  2)} Next, by assuming   \($\bV_\l$\) for some $\l \in (0,1)$.
   we shall  show that $\ul{w}_1$  is a viscosity supersolution of   \eqref{eq:PDE} with Hamiltonian
 $\ul{H}_1$.
   Let $(t_0,x_0, \vf) \in (0,T) \times \hR^k \times \hC^{1,2}\big([0,T] \times \hR^k\big)$ be
   such that   $\ul{w}_1(t_0,x_0) = \vf (t_0,x_0)$
  and that  $\ul{w}_1-\vf$ attains a strict local   minimum  at $(t_0,x_0)$, i.e., for some
     $\d_0 \in  \big(0,    t_0 \land ( T-t_0 ) \big) $
   \beas   
     (\ul{w}_1 - \vf) (t,x) >  (\ul{w}_1 - \vf) (t_0,x_0) = 0,
   \q \fa (t,x) \in O_{\d_0} (t_0,x_0) \big\backslash \big\{ (t_0,x_0) \big\} .
   \eeas
 The continuity of $\ul{l}$, $\ol{l}$ and  \eqref{eq:n411} imply that
 $    \ul{l}(t,x) \neg \le  \neg  \ul{w}_1 (t,x)  \neg \le \neg  \ol{l}(t,x)  $, $ \fa (t,x)  \neg \in \neg  [0,T]
  \neg \times \neg  \hR^k    $.    So it suffices to show that
  \bea \label{eq:s611}
     &&    \max  \neg  \Big\{  \neg  - \dneg \frac{\pa }{\pa t} \vf (t_0,x_0)
      \neg - \neg  \ul{H}_1 \big(t_0,x_0, \vf (t_0,x_0), D_x \vf (t_0,x_0), D^2_x \vf (t_0,x_0)\big),
  (\vf  \neg -  \neg  \ol{l})(t_0,x_0)  \Big\}  \ge 0  ,  
     \eea
 which clearly holds if $\ul{H}_1 \big(t_0,x_0, \vf (t_0,x_0), D_x \vf (t_0,x_0), D^2_x \vf (t_0,x_0)\big) = - \infty$.

 \ss To make a contradiction,  we   assume    that when $\ul{H}_1 \big(t_0,x_0, \vf (t_0,x_0), D_x \vf (t_0,x_0), D^2_x \vf (t_0,x_0)\big) >  - \infty$, \eqref{eq:s611} does not hold, i.e.
  \beas
 \varrho \dfnn \min  \neg  \Big\{    \frac{\pa }{\pa t} \vf (t_0,x_0)
      \neg + \neg  \ul{H}_1 \big(t_0,x_0, \vf (t_0,x_0), D_x \vf (t_0,x_0), D^2_x \vf (t_0,x_0)\big),
  (\ol{l}   \neg -  \neg \vf )(t_0,x_0)  \Big\}  > 0 .
 \eeas
  As $\vf \in \hC^{1,2}\big([0,T] \times \hR^k\big) $, there exists a $\wh{u} \in \hU_0$ such that
 \beas
       \linf{(t,x)  \to  (t_0,x_0) }  \; \underset{v \in \hV_0}{\inf}
 \;    H (t,x, \vf (t,x), D_x \vf (t,x),   D^2_x \vf (t,x),  \wh{u},v)
 \ge \frac34 \varrho -  \frac{\pa }{\pa t} \vf (t_0,x_0)    .
 \eeas
 Then by the continuity of $\vf$ and $\ol{l}$, one can find  a $ \d \in  (0, \d_0    )$ such that
 \bea  \label{eq:s423b}
    \underset{v \in \hV_0}{\inf}
 \,  H (t,x, \vf (t,x), D_x \vf (t,x),   D^2_x \vf (t,x),  \wh{u},v)
 \neg \ge \neg  \frac12 \varrho  \neg - \neg   \frac{\pa }{\pa t} \vf (t,x)
 ~\; \hb{and} ~\; (\ol{l}  \neg -  \neg \vf)(t,x)  \neg \ge \neg  \frac34 \varrho , ~\; \fa (t,x)  \neg  \in  \neg  \ol{O}_{\d}  (t_0,x_0) . \q
 \eea

 Similar to \eqref{eq:s427} we set
 $        \wp    \dfnn   \min   \big\{ (\ul{w}_1 \neg -    \vf) (t,x)  \neg : (t,x )  \neg \in \neg
     \ol{O}_\d  (t_0,x_0) \big\backslash O_{\frac{\d}{2}} (t_0,x_0) \big\}       $
   \if{0}
 Since the set $ \ol{O}_\d  (t_0,x_0) \big\backslash O_{\frac{\d}{2}} (t_0,x_0)$ is  compact, there exists
 a sequence $\{(t_n,x_n)\}_{n \in \hN} $ of $  \ol{O}_\d  (t_0,x_0) \big\backslash O_{\frac{\d}{2}} (t_0,x_0)$
 that  converges to some $(t_*,x_*) \in \ol{O}_\d  (t_0,x_0) \big\backslash O_{\frac{\d}{2}} (t_0,x_0) $
 and satisfies   $ \wp \neg = \neg  \lmtd{n \to \infty} (\ul{w}_1  \neg - \neg  \vf) (t_n,x_n) $.
 The lower semicontinuity of $\ul{w}_1$ and the continuity of $\vf$ imply that  $\ul{w}_1  \neg - \neg  \vf$ is also lower semicontinuous.
 Thus, it follows that $ \wp \le (\ul{w}_1-\vf)(t_*,x_*) \le \lmtd { n \to \infty } (\ul{w}_1-\vf)(t_n,x_n ) = \wp  $,
 which together with \eqref{eq:a071b} shows that
  \bea   \label{eq:s427b}
  \wp  = \min \big\{ (\ul{w}_1 - \vf) (t,x): (t,x ) \in
     \ol{O}_\d  (t_0,x_0) \big\backslash O_{\frac{\d}{2}} (t_0,x_0) \big\} = (\ul{w}_1-\vf)(t_*,x_*)  >  0 .
  \eea
  \fi
  and     $\dis   \wt{\wp} \neg \dfnn  \neg   \frac{    \wp \land \varrho}{ 2 (1 \vee \g) T }  \neg   > \neg  0$.
    Let $ \big\{(t_j, x_j)\big\}_{j \in \hN}$ be a sequence of $O_{\frac{\d}{4}} (t_0,x_0) $ such   that
     \beas   
 \lmt{j \to \infty} (t_j,x_j) = (t_0,x_0)   \q \hb{and} \q  \lmt{j \to \infty} \wh{w}_1 (t_j,x_j) = \ul{w}_1 (t_0,x_0) = \vf (t_0,x_0) .
     \eeas
    As $\lmt{j \to \infty} \big( \wh{w}_1 (t_j,x_j)- \vf (t_j,x_j) \big) = 0$, it holds for  some $j \in \hN$ that
    \bea  \label{eq:s437b}
  \big| \wh{w}_1 (t_j,x_j)- \vf (t_j,x_j) \big| < \frac12 \wt{\wp}   t_0 .
   \eea

 Clearly,  $ \wh{\mu}_s \dfnn \wh{u}$, $s \in [t_j, T]$  is  a constant   $  \cU^{t_j}-$control.
  Fix $\beta \in \wh{\fB}^{t_j}$.  We set $\Th_\beta \dfnn \big( t_j, x_j, \wh{\mu} ,    \beta \lan  \wh{\mu} \ran \big)$
  and  define  two  $\bF^{t_j}-$stopping times:
 \beas
  \q  \t_\beta   \dfnn  \inf \Big\{s  \neg \in \neg  ( t_j, T  ] \neg  : \big( s,\wt{X}^{\Th_\beta}_s \big)
  \neg  \notin   \neg   \ol{O}_{\frac34 \d} (t_0, x_0) \Big\}
 ~    \hb{and} ~   \z_\beta  \dfnn  \inf \Big\{s  \neg \in \neg  (  \t_\beta, T] \neg : \big( s, \wt{X}^{\Th_\beta}_s \big)
  \neg   \notin   \neg   \ol{O}_\d  (t_0,x_0) \big\backslash O_{\frac{\d}{2}} (t_0,x_0) \Big\} \land T .
 \eeas
    Since 
 $ \big| \big(T, \wt{X}^{\Th_\beta}_T\big) -(t_0,x_0) \big| \ge T -t_0 > \d_0 > \frac34 \d $,
one can deduce from the continuity of $ \wt{X}^{\Th_\beta} $ that
  \bea     \label{eq:a041b}
   \t_\beta  < T  \q \hb{and} \q  \big( \t_\beta, \wt{X}^{\Th_\beta}_{\t_\beta} \big)  \in  \pa O_{\frac34 \d}(t_0, x_0) ,
 \q  P^{t_j}_0-a.s.
  \eea

 \ss     Given $n  \neg \in  \neg  \hN$,  we define
     $ q^n  \neg  (s)  \neg \dfnn \neg  \frac{\lceil 2^n s \rceil}{2^n} \neg \land \neg T$, $ s  \neg \in \neg  [0,T] $.
    Then  $ \t^n_\beta  \neg \dfnn  \neg   q^n(\t_\beta)  \neg \land  \neg  \z_\beta$ is an  $\bF^{t_j}-$stopping time.
  Similar to \eqref{eq:s420}, applying  \eqref{eq:p677} with $(\z,\t,\xi) \neg = \neg \Big(\t^n_\beta ,\, q^n(\t_\beta) ,  \,  \wh{w}_1   \big( q^n(\t_\beta),
  \wt{X}^{\Th_\beta}_{q^n(\t_\beta)}\big) \Big)$, we can deduce from Proposition \ref{prop_apriori_DRBSDE},
    H\"older's inequality and \eqref{2l_growth} that
  \if{0}
      \bea
   && \hspace{- 1.5 cm} \left|     \wt{Y}^{\Th_\beta}_{t_j}
   \Big( \t^n_\beta,   \wh{w}_1   \big( \t^n_\beta, \wt{X}^{\Th_\beta}_{\t^n_\beta} \big)     \Big)
   -\wt{Y}^{\Th_\beta}_{t_j}
  \Big(q^n(\t_\beta),   \wh{w}_1   \big( q^n(\t_\beta), \wt{X}^{\Th_\beta}_{q^n(\t_\beta)}\big)\Big)\right|^{\frac{q+1}{2}}
   \nonumber  \\
& = & \left|     \wt{Y}^{\Th_\beta}_{t_j}
   \Big( \t^n_\beta,   \wh{w}_1   \big( \t^n_\beta, \wt{X}^{\Th_\beta}_{\t^n_\beta} \big)     \Big)
   -\wt{Y}^{\Th_\beta}_{t_j}
  \bigg(\t^n_\beta,  \wt{Y}^{\Th_\beta}_{\t^n_\beta}
  \Big(q^n(\t_\beta), \wh{w}_1   \big( q^n(\t_\beta), \wt{X}^{\Th_\beta}_{q^n(\t_\beta)}\big)\Big)\bigg)\right|^{\frac{q+1}{2}}
  \nonumber  \\
  & \le &  c_0   E_{t_j} \bigg[ \Big| \wh{w}_1   \big( \t^n_\beta, \wt{X}^{\Th_\beta}_{\t^n_\beta} \big)
- \wt{Y}^{\Th_\beta}_{\t^n_\beta}
  \Big(q^n(\t_\beta), \wh{w}_1 \big( q^n(\t_\beta), \wt{X}^{\Th_\beta}_{q^n(\t_\beta)}\big)\Big) \Big|^{\frac{q+1}{2}} \bigg]
  \nonumber   \\
& = &  c_0   E_{t_j} \bigg[ \b1_{\{  q^n(\t_\beta) > \z_\beta \}}  \Big| \wh{w}_1   \big(   \z_\beta,
 \wt{X}^{\Th_\beta}_{  \z_\beta} \big)  -  \wt{Y}^{\Th_\beta}_{  \t^n_\beta }
  \Big(q^n(\t_\beta),   \wh{w}_1   \big( q^n(\t_\beta), \wt{X}^{\Th_\beta}_{q^n(\t_\beta)}\big)\Big)
\Big|^{\frac{q+1}{2}} \bigg]  \nonumber  \\
 &\le & c_0  \Big( P^t_0 \big( q^n(\t_\beta) > \z_\beta \big) \Big)^{\frac{q - 1}{2q}}
 \Bigg\{ E_{t_j} \bigg[    \Big| \wh{w}_1   \big(   \z_\beta,
 \wt{X}^{\Th_\beta}_{  \z_\beta} \big)  -  \wt{Y}^{\Th_\beta}_{  \t^n_\beta }
  \Big(q^n(\t_\beta),   \wh{w}_1   \big( q^n(\t_\beta), \wt{X}^{\Th_\beta}_{q^n(\t_\beta)}\big)\Big)
\Big|^q \bigg] \Bigg\}^{\frac{q+1}{2q}}   . \label{eq:s416b}
    \eea
 By   \eqref{2l_growth},
 \bea \label{eq:s418b}
  \Big| \wh{w}_1   \big(   \z_\beta, \wt{X}^{\Th_\beta}_{  \z_\beta} \big) \Big|
\le  \big( |\ul{l} | \vee |\ol{l} | \big)  \Big(\z_\beta, \wt{X}^{\Th_\beta}_{  \z_\beta}   \Big)
 \le  \ul{l}_*  + \ol{l}_* +  \g \underset{s \in [t_j,T]}{\sup} \big| \wt{X}^{\Th_\beta}_s \big|^{2/q} ,
 \eea
  where $ \ul{l}_* \dfnn \underset{s \in [t,T]}{\sup}   \big| \ul{l} (s, 0 ) \big|$
 and $ \ol{l}_* \dfnn \underset{s \in [t,T]}{\sup}   \big| \ol{l} (s, 0 ) \big|$.
 Similarly,
  \beas
  \q  \bigg| \wt{Y}^{\Th_\beta}_{  \t^n_\beta  }
  \Big(q^n(\t_\beta),   \wh{w}_1   \big( q^n(\t_\beta), \wt{X}^{\Th_\beta}_{q^n(\t_\beta)}\big)\Big) \bigg|
    \le  \Big( \big|\ul{L}^{\Th_\beta}_{\t^n_\beta} \big| \vee \big|\ol{L} ^{\Th_\beta}_{\t^n_\beta} \big| \Big)
    \neg  =  \neg  \big(  |\ul{l}  | \vee  |\ol{l}  | \big) \Big(\t^n_\beta, \wt{X}^{\Th_\beta}_{ \t^n_\beta }   \Big)
       \le     \ul{l}_*  \neg  +  \neg \ol{l}_*  \neg +  \neg  \g \underset{s \in [t_j,T]}{\sup} \big| \wt{X}^{\Th_\beta}_s \big|^{2/q} .
\eeas
 Putting it and \eqref{eq:s418b} back into \eqref{eq:s416b} gives that
 \fi
  \bea
  && \hspace{-2.5cm} \left|     \wt{Y}^{\Th_\beta}_{t_j}
   \Big( \t^n_\beta,   \wh{w}_1   \big( \t^n_\beta, \wt{X}^{\Th_\beta}_{\t^n_\beta} \big)     \Big)
   -\wt{Y}^{\Th_\beta}_{t_j}
  \Big(q^n(\t_\beta),   \wh{w}_1   \big( q^n(\t_\beta), \wt{X}^{\Th_\beta}_{q^n(\t_\beta)}\big)\Big)\right|^{\frac{q+1}{2}}
  \nonumber \\
  &&   \le   c_0  \Big( P^t_0 \big(q^n( \t_\beta ) > \z_\beta \big) \Big)^{\frac{q - 1}{2q}}
  \Bigg\{ 1+ E_{t_j} \bigg[  \underset{s \in [t_j,T]}{\sup} \big| \wt{X}^{\Th_\beta}_s \big|^2  \bigg] \Bigg\}^{\frac{q+1}{2q}} .
  \label{eq:s420b}
     \eea
 Since $ \t_\beta < \z_\beta $, $P^t_0-$a.s. by \eqref{eq:a041b} and since $ \lmtd{n \to \infty} q^n(\t_\beta) = \t_\beta $,
   we see from \eqref{eq:esti_X_1} that the right-hand-side of \eqref{eq:s420b} converges to $0$ as $n \to \infty$.
 Hence,   for some $ n_{\overset{}{\beta}} \in \hN$
   \bea  \label{eq:s431b}
  \left| \wt{Y}^{\Th_\beta}_{t_j}
   \Big(  \wh{\t}_\beta ,   \wh{w}_1   \big( \wh{\t}_\beta, \wt{X}^{\Th_\beta}_{\wh{\t}_\beta} \big)     \Big)
    -  \wt{Y}^{\Th_\beta}_{t_j}   \Big( q^{n_\beta}(\t_\beta),
 \wh{w}_1 \big( q^{n_\beta}(\t_\beta), \wt{X}^{\Th_\beta}_{q^{n_\beta}(\t_\beta)} \big) \Big) \right|
   <  \frac14 \wt{\wp} \,  t_0   ,
   \eea
where $ \wh{\t}_\beta \dfnn  \t^{n_\beta}_\beta = q^{n_\beta}(\t_\beta) \land \z_\beta $.

     As    $\wh{\t}_\beta$ is an $\bF^{t_j}-$stopping time,
      the continuity of $\vf$ and $\wt{X}^{\Th_\beta} $ show that
  $\sY^\beta_s \dfnn \vf \big( \wh{\t}_\beta \land s, \wt{X}^{\Th_\beta}_{\wh{\t}_\beta \land s} \big)
    + \wt{\wp} ( \wh{\t}_\beta \land s) $, $ s \in [t_j, T] $ defines a real-valued, $\bF^{t_j}-$adapted continuous process.
 Applying It\^o's formula to $ \sY^\beta $ yields that
    \bea
         \sY^\beta_s   & \tneg = &  \tneg   \sY^\beta_T     +   \int_s^T  \ff^\beta_r   dr
    -  \int_s^T  \sZ^\beta_r   d B^{t_j}_r ,   \q    s \in [t_j,T]  ,   \label{eq:s403b}
    \eea
where $\sZ^\beta_r  \dfnn        \b1_{\{r < \wh{\t}_\beta\}}    D_x \vf    \big( r, \wt{X}^{\Th_\beta}_r \big)
  \cd   \si \big( r, \wt{X}^{\Th_\beta}_r, \wh{u},   \big(\beta \lan \wh{\mu} \ran \big)_r \big)$ and
 $$
 \hspace{-3mm}   \ff^\beta_r   \neg   \dfnn   \neg    -  \b1_{\{ r \, < \wh{\t}_\beta  \}}
  \bigg\{ \wt{\wp}      +       \frac{\pa \vf }{\pa t} \big( r, \wt{X}^{\Th_\beta}_r \big)
      +       D_x \vf    \big( r, \wt{X}^{\Th_\beta}_r\big)   \cd    b \big( r, \wt{X}^{\Th_\beta}_r \neg , \wh{u},
  (\beta \lan \wh{\mu} \ran )_r \big)
         +         \frac12  trace\Big( \si \si^T \big( r, \wt{X}^{\Th_\beta}_r \neg , \wh{u},
  (\beta \lan \wh{\mu} \ran)_r \big) \cd    D^2_x \vf \big( r, \wt{X}^{\Th_\beta}_r\big) \Big) \neg \bigg\} .
     $$
 As $\vf  \neg \in  \neg   \hC^{1,2}\big([t,T]  \neg \times \neg  \hR^k\big)$,
  the measurability of $b$, $\si$ and $\beta$ show   that  both $\sZ^\beta$ and $\ff^\beta$ are
     $\bF^{t_j}-$progressively    measurable.

  \ss   Since it holds $P^t_0-$a.s. that
   \bea   \label{eq:a239b}
  \big(\wh{\t}_\beta \land s , \wt{X}^{\Th_\beta}_{\wh{\t}_\beta \land s}   \big) \in   \ol{O}_\d  (t_0, x_0) ,
  \q \fa    s \in [t_j, T]  ,
   \eea
  we see from the continuity  of $\vf$ that $\sY^\beta$ is a bounded process. And similar to \eqref{eq:s641},
  we can deduce from \eqref{b_si_linear_growth}, \eqref{b_si_Lip} as well as H\"older's inequality  that
  $\sZ^\beta \in \hH^{2,q}_{\bF^{t_j}  } \big( [t,T], \hR^d \big) $.
  \if{0}
  \beas
 && \hspace{-1.5cm}  E_{t_j} \left[  \bigg( \int_{t_j}^T \neg |\sZ^\beta_s|^2   \,  ds  \bigg)^{ q /2 } \right]
    =    (\g C_\vf)^q \,  E_{t_j} \left[  \bigg( \int_{t_j}^{\wh{\t}_\beta}
\Big( 1+ \big| \wt{X}^{\Th_\beta}_s \big|  + [\wh{u}]_{\overset{}{\hU}}
  +   \big[ (\beta \lan \hat{\mu} \ran)_s \big]_{\overset{}{\hV}}   \Big)^2   ds  \bigg)^{ q /2 } \right] \\
 &&    \le    c_0  C^q_\vf  \bigg\{   \big(1+ |x_0 | + \d + [\wh{u}]_{\overset{}{\hU}} \big)^q
+  \bigg( E_{t_j}    \int_{t_j}^T \neg \big[ (\beta \lan \hat{\mu} \ran)_s \big]_{\overset{}{\hV}}^2
  \,  ds    \bigg)^{ q /2 } \bigg\} < \infty , \q \hb{i.e. }
 \sZ^\beta \in \hH^{2,q}_{\bF^{t_j}  } \big( [t,T], \hR^d \big)   ,
  \eeas
 where  $ C_\vf \dfnn \underset{(t,x) \in \ol{O}_\d (t_0,x_0) }{\sup} \big|D_x \vf (t,x)\big|$.
 \fi
 Moreover,
    \eqref{eq:a239b} and \eqref{eq:s423b} imply that $P^t_0-$a.s.
  \beas           
   \sY^\beta_s   \le  \ol{l} ( \wh{\t}_\beta \land s, \wt{X}^{\Th_\beta}_{\wh{\t}_\beta \land s}) -\frac 34 \varrho + \wt{\wp} T
 < \ol{l} ( \wh{\t}_\beta \land s, \wt{X}^{\Th_\beta}_{\wh{\t}_\beta \land s}) = \ol{L}^{\Th_\beta}_{\wh{\t}_\beta \land s} ,
    \q  \fa    s \in [t_j, T]   ,
  \eeas
    which together with   \eqref{eq:s403b} shows  that
 $\big\{\big( \sY^\beta_s, \sZ^\beta_s, 0, 0 \big) \big\}_{s \, \in  [t_j,T]}$
  solves the  $\ol{\hb{R}}$BSDE$\left( P^{t_j}_0 , \sY^\beta_T , \ff^\beta, \ol{L}^{\Th_\beta}_{\wh{\t}_\beta \land \cd} \right)$.

 Since $\t_\beta \le \wh{\t}_\beta \le \z_\beta $ and since $ \ul{w}_1(t,x) \le \wh{w}_1 (t,x)$
 for any $(t,x) \in [0,T] \times \hR^k$, we can deduce that 
      \beas
   \sY^\beta_T  \le    \vf \big( \wh{\t}_\beta, \wt{X}^{\Th_\beta}_{\wh{\t}_\beta} \big)    + \wt{\wp} T
 <  \vf \big( \wh{\t}_\beta, \wt{X}^{\Th_\beta}_{\wh{\t}_\beta} \big) +\wp \le   \ul{w}_1  \big( \wh{\t}_\beta, \wt{X}^{\Th_\beta}_{\wh{\t}_\beta} \big)
 \le \wh{w}_1 \big( \wh{\t}_\beta, \wt{X}^{\Th_\beta}_{\wh{\t}_\beta} \big) .
  \eeas
  Also,      \eqref{eq:a239b}, \eqref{eq:s423b}  and \eqref{f_Lip}  show   that for $ds \times d P^{t_j}_0-$a.s.
  $ (s,\o)  \in [t_j,T] \times \O^{t_j}$
     \beas
     \ff^\beta_s (\o) & \tneg \le &  \tneg  
 \b1_{\{ s \, < \wh{\t}_\beta (\o) \}} \left\{ - \wt{\wp} - \frac12 \varrho   +    f \neg \left( s,\o, \wt{X}^{\Th_\beta}_s(\o),
       \sY^\beta_s (\o)-\wt{\wp}s,    \sZ^\beta_s (\o), \wh{u},  ( \beta \lan \wh{\mu} \ran )_s (\o) \right)  \right\}  \\
         &  \tneg \le &  \tneg  \b1_{\{ s \, < \wh{\t}_\beta (\o) \}} \left\{ - \wt{\wp} - \frac12 \varrho  + \g \wt{\wp} T
          +    f \neg \left(s, \o, \wt{X}^{\Th_\beta}_s (\o),
      \sY^\beta_s (\o),       \sZ^\beta_s (\o),\wh{u},  ( \beta \lan \wh{\mu} \ran )_s (\o) \right)  \right\}   \\
 &  \tneg \le &  \tneg   f^{\Th_\beta }_{\wh{\t}_\beta} \big(s , \o, \sY^\beta_s (\o),  \sZ^\beta_s (\o) \big)    .
       \eeas
   Clearly,
 $
 \left( \wt{Y}^{\Th_\beta} \Big( \wh{\t}_\beta,   \wh{w}_1  \big( \wh{\t}_\beta, \wt{X}^{\Th_\beta}_{\wh{\t}_\beta} \big)    \Big), \wt{Z}^{\Th_\beta} \Big( \wh{\t}_\beta,   \wh{w}_1  \big( \wh{\t}_\beta, \wt{X}^{\Th_\beta}_{\wh{\t}_\beta} \big)    \Big),
\wt{\ul{K}}^{\Th_\beta} \Big( \wh{\t}_\beta,   \wh{w}_1  \big( \wh{\t}_\beta, \wt{X}^{\Th_\beta}_{\wh{\t}_\beta} \big)    \Big),\wt{\ol{K}}^{\,\raisebox{-0.5ex}{\scriptsize $\Th_\beta$}}  \Big( \wh{\t}_\beta,   \wh{w}_1  \big( \wh{\t}_\beta, \wt{X}^{\Th_\beta}_{\wh{\t}_\beta} \big)    \Big) \right)
 $
 solves $\ol{\hb{R}}$BSDE$\left( P^{t_j}_0 , \wh{w}_1  \big( \wh{\t}_\beta, \wt{X}^{\Th_\beta}_{\wh{\t}_\beta} \big) ,
 f^{\Th_\beta}_{\wh{\t}_\beta} , \ol{L}^{\Th_\beta}_{\wh{\t}_\beta \land \cd} \right)$.
 As   $f^{\Th_\beta}_{\wh{\t}_\beta}$ is  Lipschitz continuous in $(y, z)$, we know from Proposition
 \ref{prop_comparison_RBSDE} that $P^{t_j}_0-$a.s.
   \beas
   \sY^\beta_s   \le \wt{Y}^{\Th_\beta}_s \left( \wh{\t}_\beta,
   \wh{w}_1  \big( \wh{\t}_\beta, \wt{X}^{\Th_\beta}_{\wh{\t}_\beta} \big)    \right)  , \q  \fa  s \in [t_j,T] .
   \eeas
  Letting $s=t_j$ and using \eqref{eq:s431b}, we obtain
\bea
       \vf   (t_j   , x_j  )  +   \frac34    \wt{\wp}   t_0
   &<&   \vf   (t_j   , x_j  )   +   \wt{\wp}   t_j   =   \sY^\beta_{t_j}
     \le     \wt{Y}^{t_j, x_j, \wh{\mu} ,    \beta \lan  \wh{\mu} \ran}_{t_j} \Big( \wh{\t}_\beta,
   \wh{w}_1  \Big( \wh{\t}_\beta, \wt{X}^{t_j, x_j, \wh{\mu} ,    \beta \lan  \wh{\mu} \ran}_{\wh{\t}_\beta} \Big)    \Big)
 \nonumber  \\
      & <&       \wt{Y}^{t_j, x_j, \wh{\mu} ,    \beta \lan  \wh{\mu} \ran}_{t_j}   \Big( q^{n_\beta}(\t_\beta),
 \wh{w}_1 \Big( q^{n_\beta}(\t_\beta), \wt{X}^{t_j, x_j, \wh{\mu} ,    \beta \lan  \wh{\mu} \ran}_{q^{n_\beta}(\t_\beta)} \Big) \Big)
 + \frac14  \wt{\wp}   t_0    , \label{eq:s433}
      \eea
 where we used the fact that $t_j \neg > \neg  t_0 -\frac14 \d  \neg > \neg  t_0 - \frac14 \d_0  \neg > \neg  \frac34 t_0 $.

  \ss Let $\{\t_{\mu,\beta} \neg : \mu \neg \in \neg \cU^{t_j}, \beta \neg \in  \neg \wh{\fB}^{t_j}\}$ be an arbitrary family
 of $\hQ_{t_j,T}  -$valued,
 $\bF^{t_j}-$stopping times such that $\t_{ \wh{\mu},\beta} = q^{n_\beta}(\t_\beta)   $ for any $\beta \in \wh{\fB}^{t_j}$.
 By   \eqref{eq:s433},   it holds for any $\beta \in \wh{\fB}^{t_j}$ that
    \beas
 \vf   (t_j   , x_j  )  +   \frac34    \wt{\wp}   t_0
  & <&       \wt{Y}^{t_j, x_j, \wh{\mu} ,    \beta \lan  \wh{\mu} \ran}_{t_j}   \Big( \t_{ \wh{\mu},\beta},
 \wh{w}_1 \Big( \t_{ \wh{\mu},\beta}, \wt{X}^{t_j, x_j, \wh{\mu} ,
  \beta \lan  \wh{\mu} \ran}_{\t_{ \wh{\mu},\beta}} \Big) \Big)     + \frac14  \wt{\wp}   t_0   \\
  &\le &   \underset{\mu \in \cU^{t_j}}{\sup}
  \wt{Y}^{t_j, x_j,  \mu, \beta \lan \mu \ran}_{t_j} \left( \t_{\mu,\beta},   \wh{w}_1  \big( \t_{\mu,\beta},
  \wt{X}^{t_j, x_j, \mu, \beta \lan \mu \ran}_{\t_{\mu,\beta}} \big)    \right) + \frac14  \wt{\wp}   t_0     .
  \eeas
 Then  taking infimum over $\beta \in \wh{\fB}^{t_j} $, we see from \eqref{eq:s437b}  and Theorem \ref{thm_DPP} that
   \beas
   \vf  \big(t_j   , x_j     \big)  \neg + \neg  \frac34 \wt{\wp}   t_0
  \neg \le \neg   \underset{\beta \in \wh{\fB}^{t_j}}{\inf} \underset{\mu \in \cU^{t_j}}{\sup}
  \wt{Y}^{t_j, x_j,  \mu, \beta \lan \mu \ran}_{t_j} \Big( \t_{\mu,\beta},   \wh{w}_1  \Big( \t_{\mu,\beta},
  \wt{X}^{t_j, x_j, \mu, \beta \lan \mu \ran}_{\t_{\mu,\beta}} \Big)    \Big)  \neg + \neg  \frac14  \wt{\wp}   t_0
   \neg   \le   \neg   \wh{w}_1 (t_j,x_j)  \neg + \neg  \frac14  \wt{\wp}   t_0
  \neg < \neg  \vf (t_j,x_j)  \neg + \neg  \frac34 \wt{\wp}   t_0 \,.
 \eeas
   A contradiction appears.
 Therefore, $\ul{w}_1$  is a viscosity supersolution of   \eqref{eq:PDE} with Hamiltonian  $\ul{H}_1$ under \($\bV_\l$\).

  Using the transformation similar to that in the part (3) of proof of Theorem \ref{thm_DPP},
 one can show that
   $ \ol{w}_2 $   is a  viscosity  subsolution  of   \eqref{eq:PDE} with   Hamiltonian  $ \ol{H}_2 $ given  \($\bU_\l$\)
 for some $\l \in (0,1)$.  \qed

\subsection{Proofs of Section \ref{sec:shift_prob}}

 \label{subsection:Proofs_S4}

   \ss \no {\bf Proof of Lemma \ref{lem_element}:}
     Set  $   \L    \dfnn \Big\{ A \subset \O^t:
 A =  \underset{\o \in A}{\cup} \big(\o \otimes_s \O^s \big)  \Big\} $. For any     $A \in \L$,    we claim that
  \bea   \label{eq:f221}
    \o \otimes_s \O^s   \subset A^c   \hb{ for any }  \o \in A^c .
 \eea
  Assume not, there is an $\o \in A^c$ and an $ \wt{\o} \in \O^s$ such that $\o \otimes_s  \wt{\o}  \in  A $,
 thus $\big(\o \otimes_s  \wt{\o}\big) \otimes_s  \O^s \subset  A$. Then
   $\o \in \o \otimes_s  \O^s =\big(\o \otimes_s  \wt{\o}\big) \otimes_s  \O^s \subset  A$. A contradiction appear.

   \ss For any $r \in [t,s]$ and  $ \cE \in \sB(\hR^d)$,   if $ \o \in   \big( B^t_r \big)^{-1}   \big( \cE \big)$, then for any $\wt{\o} \in \O^s$,  $ \big(\o \otimes_s \wt{\o}\big)(r) = \o(r) \in \cE$, i.e., $ \o \otimes_s \wt{\o} \in \big( B^t_r \big)^{-1}   \big( \cE \big)$. Thus $  \o \otimes_s \O^s \subset \big( B^t_r \big)^{-1}   \big( \cE \big)$, which implies that $\big( B^t_r \big)^{-1}   \big( \cE \big) \in \L$.
  In particular, $\es \in \L$ and $ \O^t \in \L$.  For any  $A \in \L$, \eqref{eq:f221} implies that $A^c \in \L $.
    For any $\{A_n\}_{n \in \hN} \subset \L$,
  $  \underset{n \in \hN}{\cup} A_n =  \underset{n \in \hN}{\cup} \Big( \underset{\o \in A_n}{\cup} \big(\o \otimes_s \O^s \big) \Big)  =  \underset{\o \in \underset{n \in \hN}{\cup} A_n}{\cup} \big(\o \otimes_s \O^s \big)   $,
  namely, $ \underset{n \in \hN}{\cup} A_n \in \L$. Thus,   $\L$ is   a $\si-$field of $\O^t$
  containing all   generating sets of $\cF^t_s$.
   It then follows that $ \cF^t_s   \subset \L $, proving the lemma.  \qed

     \ss \no {\bf Proof of Lemma \ref{lem_concatenation}:}
  Let $A  $ be an open subset of $\O^t$. Given $ \wt{\o} \in A^{s,\o}$, there exists a $\d$ such that
  $ O_\d \big(\o \otimes_s \wt{\o}\big) \subset A$. For any $\wt{\o}' \in O_\d(\wt{\o})$, one can deduce that
   $   \underset{r \in[t,T]}{\sup} \big| (\o \otimes_s \wt{\o}')(r) - (\o \otimes_s \wt{\o})(r) \big|
   =\underset{r \in[s,T]}{\sup} | \wt{\o}' (r) - \wt{\o} (r) | < \d $,
       which shows that $\o \otimes_s \wt{\o}' \in O_\d ( \o \otimes_s \wt{\o} ) \subset A $, i.e. $\wt{\o}' \in A^{s,\o} $.
     Hence,  $ A^{s,\o}$ is an open subset of $\O^s$.  If $A $ is a closed subset of $\O^t$, then $(A^c)^{s,\o}$ is an open subset
     of $\O^s$ and it    follows from
     \eqref{lem_basic_complement} that $A^{s,\o} = \big( (A^c)^{s,\o}\big)^c$ is a closed subset of $\O^s$.

  \ss  Next, let $ r \in [s, T]$.  For any $t' \in [t,r] $ and $ \cE \in \sB(\hR^d)$, we can deduce that
       \beas
    \Big( \big( B^t_{t'} \big)^{-1} (\cE) \Big)^{s,\o}   = \begin{cases}
    \O^s, & \hb{if  $t' \in [t,s)$ and $\o(t') \in \cE$}; \\
    \es, & \hb{if  $t' \in [t,s)$ and $\o(t') \notin \cE$}; \\
    \big\{ \wt{\o} \in \O^s  \neg :
     \o(s)+ \wt{\o} (t_i) \in \cE     \big\}
       =           \big( B^s_{t_i} \big)^{-1}    ( \cE'    ) \in \cF^s_r   , \q  & \hb{if  $t' \in [s,r]$} ,
   \end{cases}
   \eeas
    where $\cE' = \{ x - \o(s) : x \in \cE \} \in \sB(\hR^d) $.
      Thus all the generating sets of $\cF^t_r$ belong to   $ \L   \dfnn \Big\{ A \subset \O^t:  A^{s,\o}    \in  \cF^s_r \Big\} $.
      In particular,    $\es, \O^t \in \L$.    For any $A \in \L$ and $\{A_n\}_{n \in \hN} \subset \L$,  we see from  \eqref{lem_basic_complement} and \eqref{lem_basic_union} that
  $      \big(A^c\big)^{s,\o}   =   \big( A^{s,\o} \big)^c     \in \cF^s_r $  and
        $  \Big( \underset{n \in \hN}{\cup} A_n \Big)^{s,\o}
        =  \underset{n \in \hN}{\cup}   A^{s,\o}_n   \in \cF^s_r   $, i.e. $A^c, \underset{n \in \hN}{\cup} A_n \in \L$.
       So $\L$ is   a $\si-$field of $\O^t$, it follows that  $ \cF^t_r   \subset \L $,
        i.e.,     $A^{s,\o} \in \cF^s_r$ for any  $ A \in \cF^t_r$.

  On the other hand, since the continuity of paths in $\O^t$  shows that
   \bea   \label{eq:t131}
    \o \otimes_s  \O^s
        =     \Big\{ \o' \in \O^t  \neg : \o'(t')
     \neg  = \neg  \o(t'),  \fa t'\in \hQ \cap  [t,s)    \Big\}
         =      \underset{t' \in \hQ \cap [t,s) }{\cap} \big( B^t_{t'} \big)^{-1} \neg \big( \o(t')  \big)
          \in \cF^t_s  .
   \eea
  For any  $ \wt{A} \in \cF^s_r$, applying Lemma \ref{lem_shift_inverse} with $S=T$ gives that
 $\Pi^{-1}_{t,s} \big(\wt{A}\big) \in \cF^t_r $, which together with \eqref{eq:t131} shows that
   $\o \otimes_s \wt{A} = \Pi^{-1}_{t,s} \big(\wt{A}\big) \cap \big( \o \otimes_s \O^s \big) \in \cF^t_r $. \qed

 \if{0}

   \ss  On the other hand, for any $s' \in [s,r] $ and $ \cE \in \sB(\hR^d)$,  the continuity of paths in $\O^t$ shows that
  \beas
    \o \otimes_s   \big( B^s_{s'} \big)^{-1}  ( \cE )
    & \tneg   =&  \tneg    \Big\{ \o' \in \O^t  \neg : \o'(t')
     \neg  = \neg  \o(t'),  \fa t'\in \hQ \cap  [t,s) \hb{ and }  \o'(s') \in \wt{\cE}      \Big\}\\
     & \tneg   =&  \tneg   \bigg( \underset{t' \in \hQ \cap [t,s) }{\cap} \big( B^t_{t'} \big)^{-1} \neg \big( \o(t')  \big) \bigg)     \bigcap \Big(  \big( B^t_{s'} \big)^{-1} \neg \big( \wt{\cE} \big) \Big) \in \cF^t_r ,
   \eeas
  where $\wt{\cE} = \{ x + \o(s) : x \in \cE \} \in \sB(\hR^d) $. Thus all the generating sets of $\cF^s_r$ belong to $  \wt{\L}   \dfnn \Big\{ \wt{A} \subset \O^s:  \o \otimes_s \wt{A}   \in \cF^t_r \Big\} $. In particular, $\es, \O^t \in \wt{\L}$.  For any $\wt{A} \in \wt{\L}$ and $\{\wt{A}_n\}_{n \in \hN} \subset \wt{\L}$, we can deduce that
   $\o \otimes_s \wt{A}^c=  \big( \o \otimes_s \O^s \big) \, \big\backslash \big( \o \otimes_s \wt{A} \big) \in \cF^t_r$  and $ \o \otimes_s \Big( \underset{n \in \hN}{\cup} \wt{A}_n \Big)=
       \underset{n \in \hN}{\cup}\big( \o \otimes_s  \wt{A}_n \big) \in \cF^t_r $.
               So $\wt{\L}$ is   a $\si-$field of $\O^s$, it follows that $ \cF^s_r   \subset \wt{\L} $,
                i.e.,     $\o \otimes_s \wt{A} \in \cF^t_r$ for any  $ \wt{A} \in \cF^s_r$.
      \qed

 \fi

   \ss \no {\bf Proof of Proposition \ref{prop_shift1}:} Let $\xi $ be a  $\cF^t_r -$measurable random variable for some $r \in [s,T]$. For any $ \cM \in \sB(\hM)$, since $\xi^{-1}(\cM) \in \cF^t_r$, Lemma \ref{lem_concatenation} shows that
   \bea   \label{eq:f211}
    \big(\xi^{s,\o}\big)^{-1}(\cM)
    = \big\{\wt{\o} \in \O^s:    \xi ( \o \otimes_s  \wt{\o} ) \in \cM \big\}
    = \big\{\wt{\o} \in \O^s:      \o \otimes_s  \wt{\o}   \in \xi^{-1}(\cM) \big\}
    =\big( \xi^{-1}(\cM) \big)^{s,\o} \in \cF^s_r .
    \eea
 Thus $\xi^{s,\o}$ is $\cF^s_r -$measurable. Next,  consider a  $\hM-$valued,  $\bF^t-$adapted   process $  \{X_r \}_{r \in [t, T]}$. For any $r \in [s,T]$
  and $ \cM \in \sB(\hM)$, similar to \eqref{eq:f211}, one can deduce that
  $   \big(X^{s,\o}_r\big)^{-1}(\cM)    =\big( X^{-1}_r(\cM) \big)^{s,\o} \in \cF^s_r $, which shows that
   $ \big\{X^{s,\o}_r \big\}_{  r \in [s,T]}$ is $\bF^s-$adapted.  \qed

   \ss \no {\bf Proof of Proposition \ref{prop_shift2}:}  
     For  any $\cE \in \sB \big([t,T_0]\big)$ and $A \in \cF^t_{T_0}$, Lemma \ref{lem_concatenation} shows that
    \beas 
    \big(\cE \times A\big)^{s,\o} = \big\{\big(r, \wt{\o}\big) \in [s, T_0] \times \O^s: \big(r, \o \otimes_s \wt{\o}\big) \in \cE \times A \big\}
    =  \big( \cE \cap [s, T_0] \big) \times A^{s,\o} \in \sB\big([s,T_0]\big) \otimes \cF^s_{T_0}.
    \eeas
    Hence, the rectangular measurable set $ \cE \times A \in \L_{T_0} \dfnn \big\{ \cD \subset [t,T_0] \times \O^t : \cD^{s,\o} \in \sB\big([s,T_0]\big) \otimes \cF^s_{T_0} \big\}$. In particular, $\es \times \es \in \L_{T_0}$ and $[t,T_0] \times \O^t \in \L_{T_0}  $.
              For any $\cD \in \L_{T_0}$ and $\{\cD_n\}_{n \in \hN} \subset \L_{T_0}$, similar to \eqref{lem_basic2}, we can deduce that
      $ \big( ( [t,T_0] \times \O^t) \backslash \cD \big)^{s,\o}
    =     \big(\cD^{s,\o}\big)^c     \in \sB\big([s,T_0]\big) \otimes \cF^s_{T_0} $\,,
    and that $      \Big( \underset{n \in \hN}{\cup} \cD_n \Big)^{s,\o}
 =  \underset{n \in \hN}{\cup}   \cD_n^{s,\o}   \in \sB\big([s,T_0]\big) \otimes \cF^s_{T_0}  $.
     Thus    $\L_{T_0}$ is   a $\si-$field of $[t,T_0] \times \O^t$. As the product $\si-$field $ \sB([t,T_0]) \otimes \cF^t_{T_0}$ is generated by all rectangular measurable sets $ \big\{\cE \times A:    \cE \in \sB \big([t,T_0]\big), \,    A \in \cF^t_{T_0} \big\}$,
 we can deduce that $ \sB([t,T_0]) \otimes \cF^t_{T_0} \subset \L_{T_0}$, i.e., $ \cD^{s,\o} \in \sB\big([s,T_0]\big) \otimes \cF^s_{T_0}$ for any $ \cD \in \sB([t,T_0]) \otimes \cF^t_{T_0}$.

  \ss Let   $  \{X_r \}_{r \in [t, T]}$  be an $\hM-$valued,  measurable process
   on $ \big(\O^t,\cF^t_T\big)$. For any $ \cM \in \sB(\hM)$, since $X^{-1}(\cM) \in \sB\big([t,T]\big) \otimes \cF^t_T$,
   applying the above result with $T_0=T$ yields that
     \beas
    \big(X^{s,\o}\big)^{-1}(\cM) &\tneg \neg =& \tneg  \neg \{ (r,\wt{\o}) \in [s,T] \times \O^{s}  \neg  :X^{s,\o}_r (\wt{\o})  \neg \in \neg  \cM  \}
    = \{ (r,\wt{\o}) \in [s,T] \times \O^{s}  \neg : X_r (\o \otimes_s \wt{\o}) \neg  \in  \neg \cM  \} \nonumber \\
    &\tneg \neg =&\tneg  \neg \big\{ (r,\wt{\o}) \in [s,T] \times \O^{s}  \neg  : (r, \o \otimes_s \wt{\o}) \neg  \in
     \neg  X^{-1}(\cM)  \big\}
      =  \big( X^{-1}(\cM) \big)^{s,\o}        \in     \sB\big([s,T]\big) \otimes \cF^s_T,  
    \eeas
    thus $ \big\{X^{s,\o}_r \big\}_{  r \in [s,T]}$ is  a   measurable process on $ \big(\O^s,\cF^s_T\big)$.

   \ss Next, we consider    an $\hM-$valued, $\bF^t-$progressively measurable process
 $  \{X_r \}_{r \in [t, T]}$. For any  $T_0 \in [s,T] $ and    $\wt{\cM} \in \sB(\hM) $, the $\bF^t-$progressive  measurability of $X$
   assures that $\wt{\cD} \dfnn \big\{ (r, \o' ) \in [t,T_0] \times \O^t: \, X_r ( \o'   ) \in \wt{\cM} \big\}
    \in \sB\big([t,T_0]\big) \otimes \cF^t_{T_0}  $, thus $\wt{\cD}^{s,\o} \in \sB\big([s,T_0]\big) \otimes \cF^s_{T_0}$.
     It follows that
         \beas
      && \hspace{-2cm}   \big\{ (r, \wt{\o}) \in [s,T_0] \times \O^s: \, X^{s,\o}_{r}(\wt{\o}) \in \wt{\cM} \big\}
       =   \big\{ (r, \wt{\o}) \in [s,T_0] \times \O^s: \, X_{r}\big( \o \otimes_s \wt{\o}\big) \in \wt{\cM} \big\} \nonumber \\
     &&  =   \Big\{ (r, \wt{\o}) \in [s,T_0] \times \O^s: \,\big(r, \o \otimes_s \wt{\o}\big)
      \in  \wt{\cD} \, \Big\} = \wt{\cD}^{s,\o} \in \sB\big([s,T_0]\big) \otimes \cF^s_{T_0} , 
             \eeas
     which shows that   $ \big\{X^{s,\o}_r \big\}_{  r \in [s,T]}$ is $\bF^s-$progressively measurable.   


     \ss Moreover, for any $ \cD \in \sP\big(\bF^t\big)$, since
          $\b1_{\cD}=\big\{ \b1_{\cD}(r, \o'):   r \in [t,T],\, \o' \in \O^t \big\}$
    is an $\bF^t-$progressively measurable process, $ \big(\b1_{\cD}\big)^{s,\o}= \b1_{\cD^{s,\o}} $
    is  an $\bF^s-$progressively measurable process, where we used the fact that
      \beas
       \big(\b1_{\cD}\big)^{s,\o}\big(r, \wt{\o}\big)= \b1_{\cD} \big(r, \o \otimes_s \wt{\o}\big)
       = \b1_{\cD^{s,\o}} \big(r,   \wt{\o}\big) , \q \fa r \in [s,T],~ \wt{\o} \in \O^s .
      \eeas
      Thus,  $ \cD^{s,\o} \in   \sP\big(\bF^s\big)$.      \qed


  \ss \no {\bf Proof of Corollary \ref{cor_shift2}:}  
   Similar to the proof of Proposition \ref{prop_shift2},  one can show that
      for  any $\cD \in \sP\big(\bF^t\big)$ and $\cM \in \sB(\hM)$,
      $ 
          \big(\cD \times \cM\big)^{s,\o}
     = \cD^{s,\o} \times \cM   \in \sP\big(\bF^s\big) \otimes \sB(\hM)
     $,
     and  that $\L  \dfnn \big\{ \cJ \subset [t,T] \times \O^t \times \hM: \cJ^{s,\o} \in \sP\big(\bF^s\big) \otimes \sB(\hM) \big\}$ forms a $\si-$field of $[t,T] \times \O^t \times \hM$. Thus it follows that
     $ \sP\big(\bF^t\big) \otimes \sB(\hM) \subset \L $, i.e., $ \cJ^{s,\o}
   \in \sP\big(\bF^s\big) \otimes \sB(\hM)$ for any $  \cJ    \in    \sP\big(\bF^t\big) \otimes \sB(\hM)$.

    \ss Next,     let $f :  [t,T] \times \O^t \times \hM \rightarrow  \wt{\hM} $
   be a    $\sP\big(\bF^t\big) \otimes \sB(\hM)/\sB\big(\hM\big)-$measurable function.
    For any       $ \cE  \in \sB\big( \wt{\hM} \big) $, the    measurability of $f$
   assures that $ \wt{\cJ}  \dfnn \big\{ (r, \o',x ) \in [t,T ] \times \O^t \times \hM: \, f  (r, \o',x   ) \in  \cE  \big\}
    \in \sP\big(\bF^t\big) \otimes \sB(\hM)  $. Thus, $\wt{\cJ}^{s,\o}
    \neg =\neg  \big\{ \big(r, \wt{\o},x\big) \neg \in \neg  [s, T] \times \O^s \times \hM \neg  : f^{s,\o}\big(r, \wt{\o},x\big) \neg \in \neg  \cE \big\}   \neg \in \neg \sP\big(\bF^s\big) \neg \otimes \neg \sB(\hM)$,
    which gives the measurability of $f^{s,\o}$.   \qed


\ss \no {\bf Proof of Lemma \ref{lem_bundle}:}
Set $      A \dfnn  \{\o' \in \O^t \neg : \t(\o') = \t(\o)\}     $.
  For any  $c \in \hR^d$,     $ \xi^{-1}(c)=\{\o' \in \O^t \neg : \xi (\o')   =c \}  \in  \cF^t_\t  $,
    thus $\xi^{-1}(c) \cap  A   \in \cF^t_{\t(\o)}   $.
 Since $ \o \in  A     = \underset{c \in \hR^d}{\cup} \big( \xi^{-1}(c) \cap  A    \big) $,
 we can find some     $ \wt{c} \in \hR^d$ such that $ \o \in \xi^{-1}(\wt{c}) \cap  A   \in \cF^t_{\t(\o)}    $.
   Then lemma \ref{lem_element} implies that $\o \otimes_\t \O^{\t(\o)} \subset \xi^{-1}(\wt{c}) \cap  A $.
    It follows that
  $   \xi^{\t,\o}(\wt{\o})= \xi \big(\o \otimes_\t \wt{\o} \big)= \xi (\o ) =\wt{c}  $ for any  $   \wt{\o} \in \O^{\t(\o)} $.
 It is clear that  $\t \in \cF^t_\t$, thus $ \t\big(\o \otimes_\t \wt{\o} \big)= \t ^{\t,\o} (\wt{\o}) =\t(\o)$
 for any $\wt{\o} \in \O^{\t(\o)}$. \qed

  \if{0}

   \ss \no {\bf Proof of Corollary \ref{cor_concatenation_inverse}:} For any     $  P^s_0 - $null set $A$, there exists
    an $\ol{A} \in   \cF^s_T$ with $   P^s_0   (\ol{A}) =0$ that contains $A$.
      Lemma \ref{lem_concatenation_inverse} and Proposition \ref{prop_concatenation_inverse} show that
    $ \Pi^{-1}_{t,s} (\ol{A}) \in   \cF^t_T   $
    with $  P^t_0   \big( \Pi^{-1}_{t,s} (\ol{A}) \big) =0$. As $ \Pi^{-1}_{t,s} (A) \subset
     \Pi^{-1}_{t,s} (\ol{A}) $, we see that $  \Pi^{-1}_{t,s} (A) $ is a $ P^t_0 - $null set.

     On the other hand, for any
             $dr \times dP^s_0 - $null set $\cD$,  we can find
    some $\ol{\cD} \in \sB\big([s,T]  \big) \otimes \cF^s_T$ with $ \big(dr \times dP^s_0\big) (\ol{\cD}) =0$ that includes
    $\cD$.  Applying Lemma \ref{lem_concatenation_inverse} and Proposition \ref{prop_concatenation_inverse} again yields  that
    $ \wh{\Pi}^{-1}_{t,s} (\ol{\cD}) \in \sB\big([s,T]  \big) \otimes \cF^t_T \subset \sB\big([t,T]  \big) \otimes \cF^t_T $
    with $ \big(dr \times dP^t_0\big) \big(\wh{\Pi}^{-1}_{t,s} (\ol{\cD}) \big) =0$. Since $\wh{\Pi}^{-1}_{t,s} (\cD) \subset
    \wh{\Pi}^{-1}_{t,s} (\ol{\cD}) $, it follows    that $ \wh{\Pi}^{-1}_{t,s} (\cD) $ is a $dr \times dP^t_0 - $null set.
            \qed

   \fi

  \ss \no {\bf Proof of Proposition \ref{prop_rcpd_L1}:} Fix $s \in [t,T]$. Let us  first show that
     \bea   \label{eq:n715}
    E_s \big[ \xi^{s,\o} \big]= E_t \big[\xi\big| \cF^t_s \big](\o) \in  \hR , \q \hb{for } P^t_0-a.s. ~   \o \in \O^t .
    \eea
     In virtue of Theorem 1.3.4 and (1.3.15) of  \cite{Stroock_Varadhan},
   $P^t_0$ has a regular conditional probability distribution  with respect to    $\cF^t_s$, i.e.
      a family $\{P^\o_s \}_{\o \in \O^t } \subset \cP^t$ satisfying:

\ss \no   (\,\,i)  For any $A \in \cF^t_T$, the mapping $\o \to P^\o_s (A)$ is $\cF^t_s-$measurable;

\vspace{-7mm}
  \bea    \label{rcpd_1}
 \hb{(\,ii)    For any } \xi \in \hL^1 \big( \cF^t_T \big),~  E_{P^\o_s}  [\xi]   =   E_t \big[   \xi  \big| \cF^t_s \big] (\o) ,          ~   \hb{for  } P^t_0-a.s.~ \o \in \O^t ;\hspace{6 cm}
  \eea

      \vspace{-9mm}
   \bea    \label{rcpd_2}
    \hb{(iii)   For any }  \o \in \O^t  ,  ~    P^\o_s  \big( \o \otimes_s \O^s \big) =1 . \hspace{10.8cm}
      \eea
    Given $\o \in \O^t$, since $\o \otimes_s \wt{A} \in \cF^t_T$ for any $\wt{A} \in \cF^s_T$
 by Lemma \ref{lem_concatenation}, one  can deduce from    \eqref{rcpd_2}  that
  \beas   
  P^{s,\o}\big(\wt{A}\big) \dfnn  P^\o_s \big( \o \otimes_{s } \wt{A} \big), \q  \fa  \wt{A} \in \cF^{s}_T
  \eeas
  defines   a probability measure  on $\big(\O^{s},\cF^{s}_T \big)$.

   \ss   For any  $  \wt{A} \in \cF^s_T$,   \eqref{rcpd_2} and   \eqref{rcpd_1} implies  that  for $P^t_0-$a.s. $ \o \in \O^t$
    \beas
  P^{s,\o} \big(\wt{A}\big) =  P^\o_s \big(\o \otimes_s \wt{A}\big) =  P^\o_s \big( ( \o \otimes_s \O^s ) \cap \Pi^{-1}_{t,s} (\wt{A}) \big)
    = P^\o_s  \big(   \Pi^{-1}_{t,s} (\wt{A})  \big)
    =   E_t \Big[ \b1_{\Pi^{-1}_{t,s} (\wt{A})} \big|\cF^t_s \Big] (\o) .
    \eeas
    It is easy to see that $ \Pi^{-1}_{t,s}(\cF^s_T) = \si \big( B^t_r -B^t_s; r \in [s,T] \big) $. Thus
    $\Pi^{-1}_{t,s} (\wt{A})$  is independent of $\cF^t_s$. Applying
      Lemma \ref{lem_shift_inverse}  with $S=T$  yield that for $P^t_0-$a.s. $ \o \in \O^t$
    \beas
   P^{s,\o} \big(\wt{A}\big)  =  E_t \Big[ \b1_{\Pi^{-1}_{t,s} (\wt{A})} \big|\cF^t_s\Big] (\o)  =
  E_t \Big[\b1_{\Pi^{-1}_{t,s} (\wt{A})} \Big] = P^t_0 \big( \Pi^{-1}_{t,s} (\wt{A})\big)  =  P^s_0(\wt{A})   .
 \eeas
   Since   $\cC^s_T$ is a countable set by Lemma \ref{lem_countable_generate1}, we can find
  a $P^t_0-$null set $ \cN   $ such that for any $\o \in \cN^c$,
 $    P^{s,\o} \big(\wt{A}\big)  = P^s_0(\wt{A})  $ holds  for each $\wt{A} \in \cC^s_T  \cup \{\O^s\}$.
  To wit,
  $  \cC^s_T \cup \{\O^s\}  \subset  \L \dfnn \Big\{ \wt{A} \in  \cF^s_T  :  P^{s,\o} \big(\wt{A}\big)  = P^s_0(\wt{A})   \hb{ for any } \o \in  \cN^c  \Big\} $.
      It is easy to see that  $\L$  is   a  Dynkin system. As $\cC^s_T $ is closed under intersection,
   Lemma \ref{lem_countable_generate1} and Dynkin System Theorem    show  that
     $   \cF^s_T  = \si \big(     \cC^s_T      \big)  \subset \L  $. Namely,    it holds  for any $\o \in   \cN^c$ that
      \bea  \label{eq:n711}
      P^{s,\o} \big(\wt{A}\big)  = P^s_0(\wt{A}) , \q \fa \wt{A} \in  \cF^s_T    .
      \eea

  Let $A \in \cF^t_T$. For any $\o \in \O^t$, we have
 \bea   \label{eq:n721}
   (\b1_A)^{s,\o}(\wt{\o}) 
  = \b1_{\{ \o \otimes_s \wt{\o} \in A \}}
  =  \b1_{\{ \wt{\o} \in A^{s,\o}   \}} = \b1_{ A^{s,\o}} (\wt{\o}) , \q \fa          \wt{\o} \in \O^s   .
  \eea
 By   \eqref{rcpd_1}, there exists a $P^t_0-$null set $\cN(A)$ such that
  $ P^\o_s  (   A  )       =  E_t \big[ \b1_A |\cF^t_s\big] (\o)  \in \hR$.
  Given $\o \in \big( \cN^c \cup \cN(A) \big)^c$,
    since  $ A^{s,\o} \in \cF^s_T$ by  Lemma \ref{lem_concatenation},
     we see from \eqref{eq:n721}, \eqref{eq:n711} and \eqref{rcpd_2}   that
   \beas
 E_s \big[ (\b1_A)^{s,\o} \big] &=&  E_s \big[\b1_{A^{s,\o}} \big] =  P^s_0(A^{s,\o}) = P^{s,\o} (A^{s,\o})
  = P^\o_s \big(\o \otimes_s A^{s,\o}\big)   \\
    &=&  P^\o_s \big( ( \o \otimes_s \O^s ) \cap A \big)
    = P^\o_s  (   A  )       =  E_t \big[ \b1_A |\cF^t_s\big] (\o)  \in \hR  .
   \eeas
  Then it follows that   \eqref{eq:n715} holds for each  simple $\cF^t_T-$measurable random variable.

 \if{0}
 \ss \no (2) Let $\eta $ be a simple $\cF^P_T-$measurable random variable,
 i.e, $\eta= \underset{i=1}{\overset{N}{\sum}}\l_i \b1_{A_i}$ for $\l_i \in \hR$ and $A_i \in  \cF^P_T$.
 Since  $  E_{\ol{P}}\big[ \eta \big| \cF^t_\t \big]
  =   E_{\ol{P}}\bigg[\underset{i=1}{\overset{N}{\sum}}\l_i \b1_{A_i} \Big| \cF^t_\t \bigg]
  = \underset{i=1}{\overset{N}{\sum}} \l_i E_{\ol{P}}\big[ \b1_{A_i} \big| \cF^t_\t \big] $, \ppas,
   one can deduce from   part (1) that for \ppas ~$\o \in \O^t$,
    $\eta^{\t,\o} \in \hL^1 \big( \cF^{P^{\t(\o)}_0}_T, \ol{P^{\t(\o)}_0}\big) $ and
  \beas
 \qq E_{\ol{P^{\t(\o)}_0}}\big[ \eta^{\t,\o} \big]
  =  E_{\ol{P^{\t(\o)}_0}}\bigg[ \underset{i=1}{\overset{N}{\sum}}\l_i \big(\b1_{A_i}\big)^{\t,\o} \bigg]
  = \underset{i=1}{\overset{N}{\sum}}\l_i \, E_{\ol{P^{\t(\o)}_0}}\Big[  \big(\b1_{A_i}\big)^{\t,\o} \Big]
  = \underset{i=1}{\overset{N}{\sum}}\l_i \, E_{\ol{P}}\big[ \b1_{A_i} \big| \cF^t_\t\big](\o)
  = E_{\ol{P}}\big[ \eta \big| \cF^t_\t \big](\o) .
  \eeas
  \fi

    \ss   Now, for any $\xi \in   \hL^1 \big( \cF^t_T \big)$,
        $\xi^+$ can be approximated from below by a sequence of positive simple $\cF^t_T-$measurable random variables:
       $\xi_n = \underset{i=1}{\overset{n^2-1}{\sum}} \frac{i}{n} \b1_{A_n}$ with $A_n \dfnn \big\{\xi^+ \in \big[\frac{i}{n},\frac{i+1}{n}\big)\big\} \in \cF^t_T$. It is clear that
       $\lmtu{n \to \infty}   \xi^{s,\o}_n  = \big(\xi^+\big)^{s,\o}  $ for any $\o \in \O^t$. Then Monotone Convergence Theorem
   implies that  for $P^t_0-$a.s. $ \o \in \O^t$
    \bea   \label{eq:n733}
  E_s \Big[ \big(\xi^+\big)^{s,\o}  \Big] =  \lmtu{n \to \infty} E_s \big[\xi^{s,\o}_n   \big] =  \lmtu{n \to \infty} E_t \big[\xi_n \big| \cF^t_s \big] (\o)
   = E_t \big[\xi^+ \big| \cF^t_s \big] (\o)    < \infty .
    \eea
  Similarly, $E_s \Big[ \big(\xi^-\big)^{s,\o}  \Big]   = E_t \big[\xi^- \big| \cF^t_s \big] (\o)  < \infty $
   for $P^t_0-$a.s. $ \o \in \O^t$, which together with  \eqref{eq:n733} shows that \eqref{eq:n715} holds for both
      $|\xi|$ and $\xi$.
   Let $\t$ take values in $\{t_n\}_{n \in \hN} \subset [t,T]$, it then   follows that for $P^t_0-$a.s. $ \o \in \O^t$
        \beas
          E_{\t(\o)} \big[ \big|\xi^{\t,\o}\big| \big]
          & = &  E_{\t(\o)} \big[ |\xi|^{\t,\o} \big]
           =  \sum_{n \in \hN} \b1_{ \{ \t(\o)=t_n \} } E_{t_n} \big[ |\xi|^{t_n,\o} \big]
           =  \sum_{n \in \hN} \b1_{ \{ \t(\o)=t_n \} } E_t \big[ |\xi| |\cF^t_{t_n} \big] (\o)  \\
          & = &  \sum_{n \in \hN} \b1_{ \{ \t(\o)=t_n \} } E_t \big[ |\xi| |\cF^t_\t \big] (\o)
           = E_t \big[ |\xi| \big| \cF^t_\t\big](\o) \in  [0, \infty)
        \eeas
  Similarly, we see that \eqref{eq:f475} holds for $P^t_0-$a.s. $ \o \in \O^t$.   \qed

     \ss \no {\bf Proof of Corollary \ref{cor_rcpd11}:}
  Given a $P^t_0-$null set    $  \cN  $,   there exists an $A \in \cF^t_T$ with $P^t_0(A)=0$ such that $\cN \subset A$.
    For any  $\o \in \O^t$, applying   \eqref{eq:n721} with $s = \t(\o)$ gives that  $ (\b1_A)^{\t,\o} = \b1_{A^{\t,\o}}  $;
    Also,  \eqref{lem_basic_subset} and  Lemma \ref{lem_concatenation}  show that   $  \cN^{\t,\o}   \subset       A^{\t,\o} \in  \cF^{\t(\o)}_T$. Then  \eqref{eq:f475}
         imply  that for      $P^t_0-$a.s.   $\o  \in \O^t  $
     \beas
     P^{\t(\o)}_0   \big( A^{\t,\o} \big) = E_{\t(\o)} \big[\b1_{A^{\t,\o}}\big] =E_{\t(\o)} \big[ (\b1_A)^{\t,\o} \big]
         =        E_t \neg \big[ \b1_A \big| \cF^t_\t \big] (\o)       =      0 ,
        \eeas
    thus  $   \cN^{\t,\o} \in   \sN ^{P^{\t(\o)}_0} $.
   Next let $\xi_1$ and $ \xi_2$ be two real-valued random variables with  $\xi_1 \le \xi_2$,   $P^t_0-$a.s.
   Since  $   \cN  \dfnn \{ \o \in \O^t: \xi_1 (\o) > \xi_2 (\o) \} \in   \sN^{P^t_0} $,
             it holds  for  $P^t_0-$a.s. ~  $\o  \in \O^t  $ that
  \beas
     \hspace{0.9cm}
     0 = P^{\t(\o)}_0 \big( \cN^{\t,\o} \big) \neg = \neg P^{\t(\o)}_0 \big\{ \wt{\o} \neg \in \neg \O^{\t(\o)} \neg :
          \xi_1 \big(\o \neg  \otimes_\t \neg  \wt{\o} \big) > \xi_2 \big(\o \neg  \otimes_\t \neg  \wt{\o} \big)  \big\}
               \neg = \neg  P^{\t(\o)}_0  \big\{ \wt{\o} \neg \in \neg \O^{\t(\o)} \neg : \xi^{\t,\o}_1 (\wt{\o})
               > \xi^{\t,\o}_2 (\wt{\o})    \big\}  .  \hspace{0.9cm}  \hb{\qed}
  \eeas

\ss \no {\bf Proof of Proposition \ref{prop_rcpd3}:}
  For each  $\o  \in \O^t  $, applying    Proposition \ref{prop_shift2} with $s = \t(\o)$  and $T_0=T$  shows that
    $ \big\{X^{\t,\o}_r \big\}_{  r \in [\t(\o),T]}$  is a  measurable process
    on $ \big(\O^{\t(\o)},\cF^{\t(\o)}_T\big)$.
    Since $E_t \Big[\big( \neg \int_\t^T \neg \big|X_r  \big|^p  dr\big)^{ \wh{p}/p }\Big]    < \infty$,
       the integral    $   \int_t^T \b1_{\{r \ge \t(\o)\}} \big| X_r(\o) \big|^p  dr $ is well-defined for all $\o \in \O^t$
except on an $P^t_0-$null set $\cN  $.  Let  $  \xi \neg \dfnn \neg \b1_{ \cN^c} \neg  \int_\t^T  \neg | X_r |^p  dr  $,
 so $ \xi^{\wh{p}/p}  \neg \in  \neg  \hL^1\big(\cF^t_T  \big) $.
 Given $\o  \neg \in \neg  \O^t$,     \eqref{lem_basic_complement} and  Lemma \ref{lem_bundle} show that
    \bea  \label{eq:f479}
   \q         \xi^{\t,\o} (\wt{\o})
= \b1_{\big( \cN^c\big)^{\t,\o}} (\wt{\o})
  \cd \int_{\t(\o \otimes_\t \wt{\o})}^T \big| X_r \big(\o \otimes_\t \wt{\o}\big) \big|^p dr
= \b1_{\big( \cN^{\t,\o}\big)^c}  (\wt{\o}) \cd \int_{\t(\o)}^T \big| X^{\t,\o}_r (\wt{\o}) \big|^p dr , \q \fa \wt{\o} \in \O^{\t(\o)} .     \q
    \eea
  By Corollary \ref{cor_rcpd11}, it holds for $P^t_0-$a.s.~  $\o \in \O^t$ that   $\cN^{\t,\o} \in \sN^{P^{\t(\o)}_0}$.
  Then  \eqref{eq:f479} and \eqref{eq:f475} imply  that for $P^t_0-$a.s. $\o \in \O^t$
\beas
   \hspace{0.5cm}   E_{ \t(\o) }\bigg[ \Big( \int_{\t(\o)}^T \big| X^{\t,\o}_r \big|^p  dr \Big)^{\wh{p}/p} \bigg]
 = E_{ \t(\o) } \Big[   \big( \xi^{\wh{p}/p}\big)^{\t,\o} \Big]
=  E_t\Big[   \xi^{\wh{p}/p} \big| \cF^t_\t\Big](\o)
=  E_t\bigg[ \Big(  \int_\t^T  |  X_r |^p  dr \Big)^{\wh{p}/p} \Big| \cF^t_\t\bigg](\o) < \infty. \hspace{0.5cm} \hb{\qed}
\eeas

      \ss \no {\bf Proof of Corollary \ref{cor_rcpd3}:} Firstly, let  $\{X_r\}_{r \in [t,T]} \neg  \in  \neg  \hH^{p,\wh{p}}_{\bF^t}([t,T],\hE) $.
 For each $\o \neg  \in \neg  \O^t$, applying Proposition \ref{prop_shift2} with $s = \t (\o)$  shows that
  $\{X^{\t,\o}_r\}_{r \in [\t(\o),T]}$ is  $\bF^{\t(\o)}-$progressively measurable.
     Then    Proposition \ref{prop_rcpd3} implies that for   $P^t_0-$a.s.~$\o \in \O^t$,
 $\{X^{\t,\o}_r\}_{r \in [\t(\o),T]} \in \hH^{p,\wh{p}}_{\bF^{\t(\o)}}\big([\t(\o),T],\hE,P^{\t(\o)}_0\big) $.

 \ss  Next, let $\{X_r\}_{r \in [t,T]} \neg  \in \neg  \hC^p_{\bF^t}([t,T], \hE)$
 with continuous paths except on an $P^t_0-$null set $\cN  $.  Given $\o  \neg \in \neg  \O^t$, Proposition  \ref{prop_shift1} shows  that
   $\{X^{\t,\o}_r\}_{r \in [\t(\o),T]}$ is $\bF^{\t(\o)}-$adapted, and     the path $ [\t(\o),T] \ni  r \to X^{\t,\o}_r (\wt{\o}) =
  X_r \big( \o \otimes_\t \wt{\o}\big) $ is  continuous for any  $\wt{\o} \in \big( \cN^c \big)^{\t,\o}  $.
  By  \eqref{lem_basic_complement} and Corollary \ref{cor_rcpd11}, it holds for $P^t_0-$a.s. ~ $\o \in \O^t$ that
   $ P^{\t(\o)}_0 \big(  ( \cN^c  )^{\t,\o} \big) = P^{\t(\o)}_0 \big(  ( \cN^{\t,\o}  )^c \big) =1   $.
    Moreover, applying Proposition \ref{prop_rcpd_L1}
  with $\xi =\underset{r \in [t,T]}{\sup} |X_r|  \in \hL^p \big( \cF^t_T \big)$   yields that   for $P^t_0-$a.s. ~ $\o \in \O^t$,
  $     E_{ \t(\o) } \bigg[ \, \underset{r \in [\t(\o),T]}{\sup} |X^{\t,\o}_r|^p  \bigg]
   \le  E_{ \t(\o) } \big[   |\xi^{\t,\o}|^p  \big] <\infty$. \qed

\ss \no {\bf Proof of Proposition \ref{prop_rcpd5}:}
 We set $ \cD_r \dfnn \{ \o \in \O^t: (r,\o) \in \cD \}   $,     $\fa r \in [t,T]$.
   Fubini Theorem shows that
   \beas  
    0 \neg = \neg  (dr \times dP^t_0) \, (\cD \cap \[\t,T\])
  \neg   =\neg  \int_{ \O^t} \dneg \bigg( \int_{\t(\o)}^T \b1_{\cD_r} (\o) \, dr  \bigg) dP^t_0  (\o)
     \neg  = \neg  E_t\bigg[\int_\t^T    \b1_{\cD_r}dr \bigg].
   \eeas
 Thus $ \int_\t^T    \b1_{\cD_r}dr \in \hL^1(\cF^t_T)$ is equal to $0$, $P^t_0-$a.s., which together with
 \eqref{eq:f475} implies that
    \bea  \label{eq:f339}
     E_{\t(\o)} \bigg[ \Big( \int_{ \t }^T \b1_{\cD_r}dr \Big)^{\t,\o} \bigg]
   = E_t \bigg[ \int_\t^T    \b1_{\cD_r}dr \Big|\cF^t_\t \bigg] (\o)=0.
  \eea
 holds for any $\o  \in \O^t $ except on an $P^t_0-$null set $\cN  $. Given $\o \in \cN^c$,
 applying Proposition \ref{prop_shift2} with $s=\t(\o)$ and $  T_0 =T $ yields that
  $\cD^{\t,\o}  \in   \sB\big([\t(\o),T]\big) \otimes \cF^{\t(\o)}_T $.   Since
  \beas
    (\cD^{\t,\o})_r
    = \big\{\wt{\o} \in \O^{\t(\o)}: \big(r,  \wt{\o}\big) \in \cD^{\t,\o} \big\}
    = \big\{\wt{\o} \in \O^{\t(\o)}: \big(r,  \o \otimes_\t \wt{\o}\big) \in \cD   \big\}
    = \big\{\wt{\o} \in \O^{\t(\o)}: \o \otimes_\t \wt{\o} \in \cD_r\big\}
            \eeas
     for any $ r \in [\t(\o),T]$,  we can deduce from
       Fubini Theorem, Lemma \ref{lem_bundle}   and \eqref{eq:f339}  that
     \beas
   \hspace{0.6cm}   \big( dr \neg  \times \neg d P^{\t(\o)}_0 \big) \big( \cD^{\t,\o} \big)
                 & \tneg = & \tneg    \neg \int_{\O^{\t(\o)}}\neg \left(   \int_{ \t(\o)}^T \neg  \b1_{\cD^{\t,\o}_r} (\wt{\o}) dr \right) d P^{\t(\o)}_0 (\wt{\o})
         =    \neg \int_{\O^{\t(\o)}}\neg  \left(   \int_{ \t(\o \otimes_\t \wt{\o})}^T \neg  \b1_{\cD_r} \big(  \o \otimes_\t \wt{\o}\big) dr \right) d P^{\t(\o)}_0 (\wt{\o}) \nonumber \\
   &   \tneg       = & \tneg    \neg
    \int_{\O^{\t(\o)}} \neg \Big(   \int_{ \t  }^T \neg  \b1_{\cD_r}  dr \Big)^{\t,\o} \neg  (  \wt{\o} ) \, d P^{\t(\o)}_0 (\wt{\o}) =   E_{\t(\o)} \bigg[ \Big( \int_{ \t }^T \b1_{\cD_r}dr \Big)^{\t,\o} \bigg] = 0 .
      \hspace{2.6cm} \hb{\qed}  
      \eeas

      \ss \no {\bf Proof of Proposition \ref{prop_dwarf_strategy}:} (1)   Let $\mu \in \cU^t$. Given $\o \in \O^t $,
    we see from Proposition \ref{prop_shift2}   that   $\mu^{\t,\o} = \{\mu^{\t,\o}_r\}_{r \in [\t(\o),T]}$ is an
     $\bF^{\t(\o)}-$progressively measurable process, and    Lemma \ref{lem_bundle} shows that
 \bea \label{eq:t313}
  \bigg( \int_\t^T \big[ \mu_r\big]^2_{\overset{}{\hU}} \, dr \bigg)^{\t, \o} (\wt{\o})
  = \int_{\t(\o \otimes_\t \wt{\o})}^T \big[ \mu_r (\o \otimes_\t \wt{\o}) \big]^2_{\overset{}{\hU}} \, dr
  = \int_{\t(\o )}^T \big[ \mu^{\t,\o}_r (  \wt{\o}) \big]^2_{\overset{}{\hU}} \, dr , \q \fa \wt{\o} \in \O^{\t(\o)} .
 \eea
 As $\cD \dfnn \{(r,\o) \in [t,T] \times \O^t: \mu_r (\o) \in \hU \, \backslash \hU_0 \}$ has zero $dr \times dP^t_0-$measure,
 we see from Proposition  \ref{prop_rcpd5} that for $P^t_0-$a.s. $\o \in \O^t$,
  \beas
 0 = \big( dr \times d P^{\t(\o)}_0 \big) \big(\cD^{\t,\o}\big) =
 \big( dr \times d P^{\t(\o)}_0 \big) \big( \big\{ (r,\wt{\o}) \in [\t(\o),T] \times \O^{\t(\o)}: \mu^{\t,\o}_r ( \wt{\o})
 \in \hU \, \backslash \hU_0 \big\}  \big) .
  \eeas
 On the other hand,    applying Proposition \ref{prop_rcpd_L1}
  with $\xi = \int_\t^T \big[ \mu_r\big]^2_{\overset{}{\hU}} \, dr  \in \hL^1 \big( \cF^t_T \big)$
   and using \eqref{eq:t313} yield that   for $P^t_0-$a.s. ~ $\o \in \O^t$,
  \bea \label{eq:t315}
       E_{ \t(\o) } \bigg[ \int_{\t(\o)}^T \big[ \mu^{\t,\o}_r\big]^2_{\overset{}{\hU}} \, dr \bigg]
     =  E_t \bigg[   \int_\t^T \big[ \mu_r\big]^2_{\overset{}{\hU}} \, dr \Big| \cF^t_\t  \bigg] (\o)
    \le E_t \bigg[   \int_t^T \big[ \mu_r\big]^2_{\overset{}{\hU}} \, dr \Big| \cF^t_\t  \bigg] (\o)
    <\infty ,
  \eea
 thus  $  \mu^{\t,\o} \in \cU^{\t(\o)}$.
 Similarly, one can show that for any $\nu \in \cV^t$,   it holds for $P^t_0-$a.s. $\o \in \O^t$
  that  $\nu^{\t,\o} \in \cV^{\t(\o)}$.

  \ss \no (2)   Let  $\a \in  \cA^t$.
  There exists a  $dr \times d P^t_0-$a.s. null set $\cD \subset [t,T] \times \O^t$ such that
  $\a(r, \o, \hV_0) \neg \subset \neg  \hU_0  $ and $\eqref{eq:r503}$ holds
   for all $(r,\o) \in ( [s,T] \times \O^t ) \backslash \cD$.
     Given $\o \in \O^t $,    Corollary \ref{cor_shift2} shows  that   $\a^{s,\o}(r, \wt{\o},v) =
    \a (r, \o \otimes_s \, \wt{\o} , v) $, $\fa  (r, \wt{\o},v) \in [s,T] \times \O^s \times \hV$ is
    $ \sP\big(\bF^s\big) \otimes \sB(\hV)/\sB\big(\hU \big)-$measurable, and we can deduce that
     for any $(r,\wt{\o}) \in ([s,T] \times \O^s) \backslash \cD^{s,\o} $
     \big(or equivalent, $(r,\o \otimes_s \wt{\o}) \in ( [s,T] \times \O^t ) \backslash \cD$\big),
     $\a^{s,\o} \big(r, \wt{\o},\hV_0 \big)    =  \a  \big(r, \o \otimes_s \wt{\o}, \hV_0 \big) \subset \hU_0$  and
     \bea  \label{eq:t311}
  \big[ \a^{s,\o} \big(r, \wt{\o},v\big) \big]_{\overset{}{\hU}}
  =\big[ \a  \big(r, \o \otimes_s \wt{\o}, v\big) \big]_{\overset{}{\hU}}
    \le \Psi_r \big( \o \otimes_s \wt{\o} \big) + \k [v]_{\overset{}{\hV}}
    = \Psi^{s,\o}_r  (   \wt{\o}  ) + \k [v]_{\overset{}{\hV}} \, , \q \fa v \in \hV .
     \eea
      In light of Proposition \ref{prop_rcpd5} and Proposition \ref{prop_rcpd3},  it holds for $P^t_0-$a.s. ~ $\o \in \O^t$ that
      $\cD^{s,\o}$ is a $dr \times d P^s_0-$a.s. null set and that
          $ \Psi^{s,\o}  $  is a non-negative  measurable process on   $(\O^s,\cF^s_T) $ with
    $E_s \neg \int_s^T \neg  ( \Psi^{s,\o}_r )^2  \,  dr < \infty$.
     Hence,   $\a^{s,\o} \in \cA^s$ for $P^t_0-$a.s. ~ $\o \in \O^t$.

 Moreover,  assuming that $\a $ additionally satisfies \eqref{eq:r742}, i.e. $\a \in \wh{\cA}^t$, we shall
  show that $ \a^{s,\o} $ also satisfies \eqref{eq:r742} for $P^t_0-$a.s. ~ $\o \in \O^t$.
  Given    $n \in \hN$,  there exist a $\d_n$ and    a    closed subset $F_n     $ of $\O^t$
     with $  P^t_0  (  F_n  ) > 1- \frac{1}{n} $ such that
      for any   $\o ,\o' \in F_n$ with $ \|\o - \o' \|_t < \d_n$
 \bea  \label{eq:t457}
 \underset{r \in [t,T]}{\sup} \, \underset{v \in \hV}{\sup} \; \rho_{\overset{}{\hU}} \big( \a(r,\o,v) ,  \a(r,\o',v) \big) < \frac{1}{n}      .
 \eea
 For any $\o \in \O^t$,   Lemma \ref{lem_concatenation} shows that $   F^{s,\o}_n   $ is a
    closed subset of $\O^s$.
       Applying Proposition \ref{prop_rcpd_L1} with
   $ \xi =   \b1_{F_n}$ and using \eqref{eq:n721} show that  for all $\o \in \O^t$ except on a $P^t_0-$null set $\cN_n$
     \bea  \label{eq:t011}
       P^s_0 \big(     F^{s,\o}_n   \big)
    =       E_s \big[     \b1_{    F^{s,\o}_n }     \big]
 =       E_s \big[    \big( \b1_{     F_n }  \big)^{  s,\o}      \big]
  =    E_t \big[    \b1_{     F_n }  \big| \cF^t_s    \big] (\o) < \infty .
     \eea
    As  $ \lmt{n \to \infty} E_t  \Big[ 1-  E_t \big[   \b1_{  F_n }  \big|\cF^t_s    \big] \Big]
   = 1 - \lmt{n \to \infty}  P^t_0 \big( F_n \big) = 0 $, there exists a subsequence
   $\{n_i\}_{i \in \hN}$ of $\hN$ such that $\dis \lmt{i \to \infty}   E_t \big[   \b1_{  F_{n_{\neg i}} }  \big|\cF^t_s \big] \\
    = 1 $ holds except on a $P^t_0-$null set $\cN$.
   Thus, taking $n=n_i$ in \eqref{eq:t011} and letting $i \to \infty$ yield that
   $ \lmt{i \to \infty}    P^s_0 \big(  F^{ s,\o}_{n_{ i}} \big) = 1  $
   for all $ \o \in \O^t$ except on  the $P^t_0-$null set
    $\wt{\cN} \dfnn \big( \underset{i \in \hN}{\cup} \cN_{n_i} \big) \cup \cN  $.

    Now,   fix $\o \in \wt{\cN}^c$.
    For any $\e  > 0$, there exists an  $i \in \hN$ such that $n_i > \frac{1}{\e}$ and
    that $    P^s_0   \big(  F^{  s,\o}_{n_{  i}}  \big)  > 1-\e $.
    Let $\wt{\o}   \in F^{s,\o}_{ n_{  i}}$.  For any $  \wt{\o}' \in F^{s,\o}_{ n_{  i}}      $ with
    $ \big\| \wt{\o} - \wt{\o}' \big\|_s    \neg  <  \neg  \d_{n_i} $,
    since $  \underset{r \in [t,T]}{\sup} \big| ( \o  \neg \otimes_s     \wt{\o})(r)
       -     ( \o  \neg \otimes_s    \wt{\o}' ) (r) \big|
     \neg = \neg  \underset{r \in [s,T]}{\sup} \big|  \wt{\o} (r)    -    \wt{\o}' (r) \big|
      \neg < \neg   \d_{n_i} $,
     we see from \eqref{eq:t457} that
         \beas
     \underset{ r \in [t, T] }{\sup} \, \underset{v \in \hV}{\sup} \; \rho_{\overset{}{\hU}} \big(\a^{s,\o}(r, \wt{\o},v), \a^{s,\o}(r, \wt{\o}',v) \big)
   = \underset{ r \in [t, T] }{\sup} \, \underset{v \in \hV}{\sup} \; \rho_{\overset{}{\hU}} \big(\a (r, \o \otimes_s \wt{\o},v), \a (r, \o \otimes_s \wt{\o}',v) \big)      < \frac{1}{n_i} < \e       .
         \eeas
   Hence,   $\a^{s,\o}$ satisfies \eqref{eq:r742}.

 \ss   Similarly, one can show that for any $\beta \in  \fB^t$ \big(resp.\;$\beta \in \wh{\fB}^t$\big),   it holds for $P^t_0-$a.s. $\o \in \O^t$  that  $\beta^{s,\o}   \in  \fB^s$ \big(resp.\;$\beta^{s,\o}   \in \wh{\fB}^s$\big).      \qed

 \ss \no {\bf Proof of Lemma \ref{lem_shift_converge_proba}:} Let $ \{\xi_i\}_{i \in \hN} $ be a sequence of $ \hL^1(\cF^t_T)$
  that converges to 0 in probability $P^t_0$, i.e.
   \bea  \label{eq:p011}
   \lmtd{i \to \infty} E_t  \big[\b1_{ \{ |\xi_i| > 1/n    \}} \big]
     =  \lmtd{i \to \infty} P^t_0 \big( |\xi_i| > 1 / n   \big) = 0 , \q    \fa  n \in \hN .
    \eea
   In particular, $\lmtd{i \to \infty} E_t  \big[\b1_{ \{ |\xi_i| > 1  \}} \big]    =   0 $
   allows us to extract a subsequence $S_1 = \big\{ \xi^1_i \big\}_{i \in \hN}$
  from $ \{\xi_i\}_{i \in \hN}$ such that $ \lmt{i \to \infty} \b1_{\{|\xi^1_i| > 1\}}  = 0$, $P^t_0-$a.s.
  Clearly, $S_1$ also satisfies \eqref{eq:p011}. Then by $\lmtd{i \to \infty} E_t  \big[\b1_{ \{ |\xi^1_i| > 1/2  \}} \big]
     =   0$, we  can find    a subsequence $S_2 = \big\{ \xi^2_i \big\}_{i \in \hN}$ of $S_1$
    such that $ \lmt{i \to \infty} \b1_{\{|\xi^2_i| > 1/2 \}}  = 0$, $P^t_0-$a.s.  Inductively, for each $n \in \hN$ we
    can  select a subsequence $S_{n+1} = \{\xi^{n+1}_i\}_{i \in \hN}$ of $ S_n = \{\xi^n_i\}_{i \in \hN}$
     such that $ \lmt{i \to \infty} \b1_{\big\{ |\xi^{n+1}_i| > \frac{1}{n+1} \big\}} = 0$, $P^t_0-$a.s.

 For any $ i \in \hN $, we set $\wt{\xi}_i \dfnn \xi^i_i$, which belongs to $S_n$ for  $ n =1,\cds, i$. Given $n \in \hN$,
 since $\{\wt{\xi}_i\}^\infty_{i = n} \subset S_n$,  it holds $P^t_0-$a.s. that $ \lmt{i \to \infty} \b1_{\big\{|\wt{\xi}_i| > \frac{1}{n}\big\}} = 0$. Then   Bound Convergence Theorem and  Proposition \ref{prop_rcpd_L1} imply that
  \bea  \label{eq:p015}
   0= \lmt{i \to \infty} E_t \Big[ \b1_{ \{|\wt{\xi}_i| > 1/n  \}} \big|\cF^t_\t\Big](\o)=
    \lmt{i \to \infty} E_{\t(\o)} \Big[ \big( \b1_{ \{|\wt{\xi}_i| > 1/n  \}} \big)^{\t,\o} \Big]
  \eea
  holds for all $\o \in \O^t$ except on a $P^t_0-$null set $\cN_n$. Let $\o \in \Big(\underset{n \in \hN}{\cup} \cN_n\Big)^c$.
  For any $n \in \hN$, one can deduce that
  \beas
   \big( \b1_{ \{|\wt{\xi}_i| > 1/n  \}} \big)^{\t,\o} (\wt{\o})
    =   \b1_{ \big\{|\wt{\xi}_i(\o \otimes_\t \wt{\o}) | > 1/n  \big\} }
    =     \b1_{ \big\{ \big|\wt{\xi}^{\,\t,\o}_{\, i} (\wt{\o}) \big| > 1/n  \big\}}
    =   \Big( \b1_{ \big\{|\wt{\xi}^{\,\t,\o}_{\, i}| > 1/n  \big\}} \Big) (\wt{\o}) , \q \fa \wt{\o} \in \O^{\t(\o)} .
   \eeas
   which together with \eqref{eq:p015} leads to that
   $  \lmt{i \to \infty} P^{\t(\o)}_0 \Big(   |\wt{\xi}^{\,\t,\o}_{\,i}| > 1 / n    \Big)
   = \lmt{i \to \infty} E_{\t(\o)} \Big[ \big( \b1_{ \{|\wt{\xi}_i| > 1/n  \}} \big)^{\t,\o} \Big] = 0 $.   \qed

 \ss \no {\bf Proof of Proposition \ref{prop_FSDE_shift}:}
 Since $\wt{X}^\Th \neg \in \neg  \hC^2_{ \bF^t }  \neg  ([t,T], \hR^k)$, Corollary \ref{cor_rcpd3} shows that
  for $P^t_0-$a.s. $\o  \neg \in \neg  \O^t$,
 $ \big\{  (\wt{X}^\Th  )^{\t,\o}_r \big\}_{r \in [\t(\o), T]} \\ \in \hC^2_{ \bF^{\t(\o)} }\big([\t(\o),T], \hR^k,P^{\t(\o)}_0 \big)$.
 Thus, it suffices to show that for $P^t_0-$a.s. $\o \in \O^t$, $\big(\wt{X}^\Th \big)^{\t,\o}$ solves \eqref{eq:p201}.

 As the $\bF^t-$version of $X^\Th$, $\wt{X}^\Th$ also satisfies \eqref{FSDE}.
    Thus, one can deduce that   except on  a $P^t_0-$null set $\cN$
     \bea \label{eq:p205}
    \wt{X}^{\Th}_{\t \vee s}    \neg - \neg   \wt{X}^{\Th}_\t
     =  \neg  \int_{\t \vee s}^\t  \neg  b \big(r, \wt{X}^{\Th}_r,\mu_r,\nu_r \big)   dr
     \neg + \neg  \int_{\t \vee s}^\t  \neg  \si \big(r, \wt{X}^{\Th}_r,\mu_r,\nu_r \big)   dB^t_r
    =     \neg   \int_t^s  \neg  b^{\,\Th}_{\,r} dr
    \neg  +  \neg  \int_t^s  \neg  \si^\Th_r   dB^t_r ,  ~            s  \in [t , T]  ,
 \eea
    where $b^{\,\Th}_{\,r} \neg \dfnn \neg  \b1_{\{r > \t \}} b \big(r, \wt{X}^{\Th}_r,\mu_r,\nu_r \big) $
     and $\si^\Th_r  \neg \dfnn \neg  \b1_{\{r > \t \}} \si \big(r, \wt{X}^{\Th}_r,\mu_r,\nu_r \big)$.  In light of \eqref{b_si_Lip}, \eqref{eq:p311} and \eqref{eq:esti_X_1}, $M^\Th_s \dfnn \int_t^s   \si^\Th_r dB^t_r$, $s \in [t , T]$
 is a square-integrable martingale with respect to $ \big(  \ol{\bF}^t,P^t_0 \big) $.

  \ss  Since $ \wt{X}^{\Th}_\t     \in   \bF^t_\t  $,
    Lemma \ref{lem_bundle} shows that for any $\o \in \O^t$ and $\wt{\o} \in \O^{\t(\o)}$
 \bea      \label{eq:h111}
      \big(\wt{X}^{\Th}_\t  \big)  (\o \otimes_\t \wt{\o})    = \big(\wt{X}^{\Th}_\t  \big)(\o)
    = \wt{X}^{\Th}_{\t(\o)   } (\o)   \q \hb{and} \q   \t \big(\o \otimes_\t \wt{\o}\big) =\t(\o)   .
 \eea
 Fix $\o \in \O^t$. It easily follows that    for any   $   \wt{\o} \in \O^{\t(\o)}  $
 \bea
 \big( \wt{X}^{\Th}_{\t   \vee s} \big) (\o \otimes_\t \wt{\o})
  & \tneg  \dneg =& \tneg  \dneg  \wt{X}^{\Th}_{\t(\o \otimes_\t \wt{\o}) \vee s}   (\o \otimes_\t \wt{\o})
  \neg = \neg  \wt{X}^{\Th}_{\t(\o) \vee s}   (\o \otimes_\t \wt{\o})
    \neg  = \neg  \wt{X}^{\Th}_s   (\o \otimes_\t \wt{\o})
   \neg = \neg   \big( \wt{X}^{\Th} \big)^{\t,\o}_s (\wt{\o}),  ~\;    \fa s \in [\t(\o),T]  ,  \qq \q  \label{eq:p411x}  \\
        b^{\,\Th}_{\,r}       (\o \otimes_\t \wt{\o})
  & \tneg  \dneg =& \tneg  \dneg      \b1_{\{r >  \t ( \o \otimes_\t \wt{\o})  \}} b \big(r, \wt{X}^{\Th}_r( \o \otimes_\t \wt{\o}),\mu_r( \o \otimes_\t \wt{\o}),\nu_r( \o \otimes_\t \wt{\o}) \big) \nonumber  \\
    & \dneg  \dneg =& \dneg  \dneg
   \b1_{\{r >  \t ( \o  )  \}}  b \big(r, \big(\wt{X}^{\Th} \big)^{\t,\o}_r (  \wt{\o}), \mu^{\t,\o}_r( \wt{\o}),\nu^{\t,\o}_r(  \wt{\o}) \big) ,\q     \fa r \in [t,T]  , \label{eq:p411a}  \\
       \hb{and similarly}    \;\;\;   \q & \tneg  \dneg  & \tneg  \dneg \nonumber  \\
  \si^\Th_r     (\o \otimes_\t \wt{\o})   & \dneg  \dneg =& \dneg  \dneg
               \b1_{\{r >  \t ( \o  )  \}}   \si \big(r, \big(\wt{X}^{\Th}  \big)^{\t,\o}_r (  \wt{\o}),
   \mu^{\t,\o}_r(   \wt{\o}),\nu^{\t,\o}_r(   \wt{\o}) \big)  ,\q     \fa r \in [t,T]  .  \label{eq:p411b}
  \eea
 Then   for any  $\wt{\o}  \neg \in \neg  \big( \cN^c \big)^{\t,\o}$,
   applying  \eqref{eq:p205} to the path $\o  \neg \otimes_\t \neg  \wt{\o}$ over period $\big[\t(\o), T \big]$ and using
     \eqref{eq:h111}-\eqref{eq:p411a} yield
         \beas
   \big(\wt{X}^{\Th}\big)^{\t,\o}_s \neg  ( \wt{\o})
    & \tneg \dneg = &  \tneg \tneg
     \wt{X}^{\Th}_{\t(\o)   } (\o)
     +      \int_{\t(\o)}^s   b \big(r, \big(\wt{X}^{\Th} \big)^{\t,\o}_r (  \wt{\o}), \mu^{\t,\o}_r( \wt{\o}),\nu^{\t,\o}_r(  \wt{\o}) \big) dr    +  \big( M^\Th \big)^{\t,\o}_s (  \wt{\o}),
      \q   s  \in    [\t(\o) , T] .
 \eeas
   By  \eqref{lem_basic_complement} and Corollary \ref{cor_rcpd11}, it holds for $P^t_0-$a.s. ~ $\o \in \O^t$ that
  $ P^{\t(\o)}_0 \big(  ( \cN^c  )^{\t,\o} \big) = P^{\t(\o)}_0 \big(  ( \cN^{\t,\o}  )^c \big) =1   $. Hence, it remains to show
  that   for $P^t_0-$a.s.~$\o \in \O^t$, it holds   $ P^{\t(\o)}_0-$a.s. that
  \beas
     \big( M^\Th  \big)^{\t,\o}_s
         =       \int_{\t(\o)  }^s \si \big(r, \big(\wt{X}^{\Th}\big)^{\t,\o}_r ,\mu^{\t,\o}_r ,\nu^{\t,\o}_r  \big) \, dB^{\t(\o)}_r   , \q s \in [\t(\o),T] .
  \eeas

     Since $M^\Th $  is a square-integrable martingale with respect to $ \big(  \ol{\bF}^t,P^t_0 \big) $,
     we know that (see e.g. Problem 3.2.27 of    \cite{Kara_Shr_BMSC})
    there is a sequence of $ \hR^{k \times d}-$valued, $ \bF^t -$simple processes
   $  \Big\{\F^n_s = \sum^{\ell_n}_{i=1} \xi^n_i \, \b1_{ \big\{s \in (t^n_i, t^n_{i+1}] \big\} } ,
  \,   s \in [t,T] \Big\}_{n \in \hN}$ \big(where $t=t^n_1< \cds< t^n_{\ell_n+1}=T$
   and $\xi^n_i  \in   \cF^t_{  t^n_i}$  for $i=1, \cds, \ell_n$\big) such that
      \beas
     P^t_0 \neg  - \neg  \lmt{n \to \infty} \int_t^T     trace\Big\{ \big(  \F^n_r -  \si^\Th_r  \big)
     \big(    \F^n_r   -   \si^\Th_r  \big)^T  \Big\} \,   ds =0  \q
    \hb{and}   \q
       P^t_0  \neg - \neg   \lmt{n \to \infty} \, \underset{s \in [t,T]}{\sup} \left|  M^n_s  - M^\Th_s  \right| =0   .
    \eeas
   where $  M^n_s \dfnn \int_t^s \F^n_r dB^t_s = \sum^{\ell_n}_{i=1}  \xi^n_i  \big(  B^t_{s \land t^n_{i+1}} -  B^t_{s \land t^n_i} \big)  $. Then it directly follows that
         \beas
     P^t_0 \neg  - \neg  \lmt{n \to \infty} \int_\t^T     trace\Big\{ \big(  \F^n_r -  \si^\Th_r  \big)
     \big(    \F^n_r   -   \si^\Th_r  \big)^T  \Big\} \,   ds =0 ~
    \q  \hb{and}  \q
       P^t_0  \neg - \neg   \lmt{n \to \infty} \, \underset{s \in [\t,T]}{\sup} \left| M^n_s - M^\Th_s  \right| =0   .
    \eeas

    By Lemma \ref{lem_shift_converge_proba},  $    \{ \F^n \}_{n \in \hN}$ has a subsequence  $  \Big\{ \wh{\F}^n_s = \sum^{\wh{\ell}_n}_{i=1} \wh{\xi}^{\,n}_{\,i} \b1_{ \big\{s \in (\wh{t}^{\,n}_{\,i}, \wh{t}^{\,n}_{\, i+1}] \big\} } ,
  \,   s \in [t,T] \Big\}_{n \in \hN}$ such that   except on a $P^t_0-$null set  $\wh{\cN}$
   \bea
  0 &=& P^{\t(\o)}_0 \neg  - \neg  \lmt{n \to \infty} \bigg( \int_\t^T     trace\Big\{ \big(  \wh{\F}^n_r - \si^\Th_r  \big)
     \big(    \wh{\F}^n_r   -  \si^\Th_r \big)^T  \Big\} ds
      \bigg)^{\t,\o}       \nonumber   \\
       &=&  P^{\t(\o)}_0 \neg  - \neg  \lmt{n \to \infty}   \int_{\t(\o)}^T
        trace\Big\{ \Big( \big( \wh{\F}^n \big)^{\t,\o}_r - \big( \si^\Th \big)^{\t,\o}_r  \Big)
     \Big( \big( \wh{\F}^n \big)^{\t,\o}_r -  \big( \si^\Th \big)^{\t,\o}_r  \Big)^T  \Big\}  \,  ds   \label{eq:p415a} \\
    \hb{and}\q 0 &= &   P^{\t(\o)}_0  \neg - \neg   \lmt{n \to \infty} \bigg( \underset{s \in [\t , T]}{\sup}
        \big|  \wh{M}^n_s  - M^\Th_s  \big| \,  \bigg)^{\t,\o}
        =     P^{\t(\o)}_0  \neg - \neg   \lmt{n \to \infty}  \,  \underset{s \in [\t(\o),T]}{\sup}
        \Big| \big( \wh{M}^n  \big)^{\t, \o}_s  - \big( M^\Th \big)^{\t, \o}_s  \Big|    ,  \qq  \label{eq:p415b}
    \eea
    where    $ \wh{M}^n_s \dfnn \int_t^s \wh{\F}^n_r dB^t_s =   \sum^{\wh{\ell}_n}_{i=1}  \wh{\xi}^{\,n}_{\,i}  \Big(  B^t_{s \land \wh{t}^{\,n}_{\, i+1}} -  B^t_{s \land \wh{t}^{\,n}_{\,i}} \Big)  $ and we made similar deductions to \eqref{eq:p411a}.

 \ss Fix    $\o \in \wh{\cN}^c$. For any $n \in \hN$ and $i=1,\cds \neg, \wh{\ell}_n$, we set  $\wh{t}^{\,n}_{\,i} (\o) \dfnn \wh{t}^{\,n}_{\,i} \vee \t (\o)$ and Proposition \ref{prop_shift1} implies that
  $ \big( \wh{\xi}^{\,n}_{\,i} \big)^{\t,\o} \in \cF^{\t(\o)}_{\wh{t}^{\,n}_{\,i} (\o) }$  since $ \wh{\xi}^{\,n}_{\,i} \in \cF^t_{\wh{t}^{\,n}_{\,i}} \subset \cF^t_{\wh{t}^{\,n}_{\,i} (\o) } $.
   It holds for any $s \in [\t(\o),T]$ and $\wt{\o} \in \O^{\t(\o)}$
  \beas
   \big( \wh{\F}^n \big)^{\t,\o}_s (\wt{\o}) =   \wh{\F}^n_s (\o \otimes_\t \wt{\o})
   =  \sum^{\wh{\ell}_n}_{i=1} \wh{\xi}^{\,n}_{\,i} (\o \otimes_\t \wt{\o}) \, \b1_{ \big\{s \in (\wh{t}^{\,n}_{\,i}, \wh{t}^{\,n}_{\, i+1}] \big\} } =    \sum^{\wh{\ell}_n}_{i=1} \big( \wh{\xi}^{\,n}_{\,i} \big)^{\t,\o} (  \wt{\o}) \, \b1_{ \big\{s \in \big( \wh{t}^{\,n}_{\,i}(\o), \wh{t}^{\,n}_{\, i+1}(\o) \big] \big\}}   ,
   \eeas
   so $ \big( \wh{\F}^n \big)^{\t,\o}_s $    is an $ \hR^{k \times d}-$valued, $ \bF^{\t(\o)} -$simple process.
    Applying Proposition 3.2.26 of    \cite{Kara_Shr_BMSC}, we see from    \eqref{eq:p415a} that
    \bea   \label{eq:p417}
      0  =     P^{\t(\o)}_0  \neg - \neg   \lmt{n \to \infty}  \,  \underset{s \in [\t(\o),T]}{\sup}
        \Bigg|  \int_{\t(\o)}^s \big( \wh{\F}^n \big)^{\t,\o}_r dB^{\t(\o)}_r
         - \int_{\t(\o)  }^s \big( \si^\Th \big)^{\t,\o}_r dB^{\t(\o)}_r  \Bigg|    .
    \eea
 For any $\wt{\o} \in \O^{\t(\o)}$, one can deduce that
   \beas
   \big( \wh{M}^n  \big)^{\t, \o}_s (\wt{\o}) & \tneg =&  \tneg   
      \sum^{\wh{\ell}_n}_{i=1}   \wh{\xi}^{\,n}_{\,i} (\o  \neg  \otimes_\t  \neg  \wt{\o})
     \Big((\o  \neg  \otimes_\t  \neg  \wt{\o}) \big(s  \neg \land \neg  \wh{t}^{\,n}_{\, i+1} \big)
       \neg - \neg  (\o  \neg \otimes_\t \neg  \wt{\o}) \big(s  \neg \land \neg  \wh{t}^{\,n}_{\, i} \big) \Big)
      \neg  =  \neg  \sum^{\wh{\ell}_n}_{i=1}   \big( \wh{\xi}^{\,n}_{\,i} \big)^{\t,\o} (  \wt{\o})
     \Big(  \wt{\o} \big(s  \neg \land \neg  \wh{t}^{\,n}_{\, i+1}(\o) \big)
      \neg  -  \neg   \wt{\o}  \big(s  \neg \land \neg  \wh{t}^{\,n}_{\, i}(\o) \big) \Big)   \\
    & \tneg =& \tneg   \sum^{\wh{\ell}_n}_{i=1} \big( \wh{\xi}^{\,n}_{\,i} \big)^{\t,\o} (  \wt{\o}) \Big(B^{\t(\o)}_{s \land \wh{t}^{\,n}_{\, i+1}(\o)}   - B^{\t(\o)}_{s \land \wh{t}^{\,n}_{\, i}(\o)} \Big) (  \wt{\o})
   = \bigg( \int_{\t(\o)}^s \big( \wh{\F}^n \big)^{\t,\o}_r dB^{\t(\o)}_r  \bigg) (\wt{\o})   , \q s \in [\t(\o), T]  ,
   \eeas
   which together with \eqref{eq:p415b}, \eqref{eq:p417} and \eqref{eq:p411b} shows that $P^{\t(\o)}_0 -$a.s.
   \bea
  \hspace{1.1cm}  \big( M^\Th \big)^{\t, \o}_s
   = \int_{\t(\o)  }^s \big( \si^\Th \big)^{\t,\o}_r dB^{\t(\o)}_r
    =  \int_{\t(\o)  }^s \si \big(r, \big(\wt{X}^{\Th}\big)^{\t,\o}_r ,\mu^{\t,\o}_r ,\nu^{\t,\o}_r  \big) \, dB^{\t(\o)}_r  ,
   \q  s \in [\t(\o),T] .  \qq \hb{\qed}  \hspace{1  cm}          \label{eq:p515}
   \eea

 \ss \no {\bf Proof of Proposition \ref{prop_DRBSDE_shift}:}
  Since $ \Big(\wt{Y}^\Th  ,  \wt{Z}^\Th ,  \wt{\ul{K}}^{\,\raisebox{-0.5ex}{\scriptsize $\Th$}} ,
    \wt{\ol{K}}^{\, \raisebox{-0.8ex}{\scriptsize $\Th$}}   \Big) \dfnn \Big(\wt{Y}^\Th(T, \xi) ,  \wt{Z}^\Th(T, \xi),  \wt{\ul{K}}^{\,\raisebox{-0.5ex}{\scriptsize $\Th$}}(T, \xi),
    \wt{\ol{K}}^{\, \raisebox{-0.8ex}{\scriptsize $\Th$}} (T, \xi) \Big) \in \hG^q_{\bF^t} \big([t,T]\big)$,
  Corollary \ref{cor_rcpd3} shows that for $P^t_0-$a.s. $\o  \neg \in \neg  \O^t$, the shifted processes
   $ \Big\{  \Big(\big( \wt{Y}^\Th   \big)^{\t,\o}_r ,
  \big( Z^\Th   \big)^{\t,\o}_r , \big( \wt{\ul{K}}^{\,\raisebox{-0.5ex}{\scriptsize $\Th$}}   \big)^{\t,\o}_r ,
  \big( \wt{\ol{K}}^{\, \raisebox{-0.8ex}{\scriptsize $\Th$}}   \big)^{\t,\o}_r  \Big) \Big\}_{r \in [\t(\o),T]}  $

 \ss \no $ \in \hG^q_{\bF^{\t(\o)}} \big([\t(\o),T]\big)$.
   Thus, it suffices to show that for $P^t_0-$a.s. $\o \in \O^t$,
   \bea  \label{eq:p541}
   \hb{  $\Big(\big( \wt{Y}^\Th   \big)^{\t,\o}, \big( Z^\Th   \big)^{\t,\o}, \big( \ul{K}^\Th   \big)^{\t,\o},
  \Big( \wt{\ol{K}}^{\, \raisebox{-0.8ex}{\scriptsize $\Th$}}   \Big)^{\t,\o}  \Big)$ solves
   DRBSDE$\Big( P^{\t(\o)}_0, f^{\Th^\o_\t}_T, \ul{L}^{\Th^\o_\t} , \ol{L}^{\Th^\o_\t} \Big)$. }
  \eea

 As the $\bF^t-$version of $\Big(Y^\Th(T, \xi) ,  Z^\Th(T, \xi),  \ul{K}^{\Th}(T, \xi),
    \ol{K}^{\, \raisebox{-0.5ex}{\scriptsize $\Th$}} (T, \xi) \Big)$,
     $\Big(\wt{Y}^\Th   ,  \wt{Z}^\Th ,  \wt{\ul{K}}^{\,\raisebox{-0.5ex}{\scriptsize $\Th$}} ,
    \wt{\ol{K}}^{\, \raisebox{-0.8ex}{\scriptsize $\Th$}}   \Big)$ also satisfies
     DRBSDE$\Big( P^t_0,  \\  f^\Th_T, \ul{L}^\Th , \ol{L}^\Th \Big)$.
    Thus, it holds  except on  a $P^t_0-$null set $\cN$ that
       \bea   \label{eq:p511}
    \left\{\ba{l}
 \dis   \wt{Y}^{\Th}_{\t \vee s} - \xi + \wt{\ol{K}}^{\, \raisebox{-0.8ex}{\scriptsize $\Th$}}_T
     -   \wt{\ol{K}}^{\, \raisebox{-0.8ex}{\scriptsize $\Th$}}_{\,\t \vee s}-  \wt{\ul{K}}^{\,\raisebox{-0.5ex}{\scriptsize $\Th$}}_{\,T}
      +  \wt{\ul{K}}^{\,\raisebox{-0.5ex}{\scriptsize $\Th$}}_{\,\t \vee s}
 =    \int_{\t \vee s}^T    f^\Th_T  \big(r,   \wt{Y}^{\Th}_r, \wt{Z}^\Th_r \big)  \, dr
             \neg-\neg \int_{\t \vee s}^T \neg \wt{Z}^\Th_r d B^t_r \ss \\
  \dis  \hspace{2.7cm} =
   \int_s^T  \neg  \b1_{\{r > \t\}} f   \Big( r,\wt{X}^{\Th}_r, \wt{Y}^{\Th}_r, \wt{Z}^\Th_r , \mu_r , \nu_r  \Big)  \, dr
  \neg-\neg M^\Th_T + M^\Th_s  ,  \q     s \in [t,T] ,   \vspace{1mm} \q \\
   \dis \ul{l} \big(s, \wt{X}^\Th_s \big)   \le     \wt{Y}^{\Th}_s   \neg \le   \neg   \ol{l} \big(s, \wt{X}^\Th_s \big) ,
    ~  s \in [t,T] , ~ \hb{and} ~
     \int_\t^T \dneg  \Big(  \wt{Y}^{\Th}_s  \neg   -  \neg  \ul{l} \big(s, \wt{X}^\Th_s \big)  \Big) d \wt{\ul{K}}^{\,\raisebox{-0.5ex}{\scriptsize $\Th$}}_{\,s}
       = \int_\t^T \neg  \Big( \ol{l} \big(s, \wt{X}^\Th_s \big) \neg - \neg \wt{Y}^{\Th}_s  \Big) d \wt{\ol{K}}^{\, \raisebox{-0.8ex}{\scriptsize $\Th$}}_s   = 0 ,
            \ea \right.
    \eea
  where $M^\Th_s \dfnn \int_t^s  \neg \b1_{\{r > \t\}}   \wt{Z}^\Th_r d B^t_r $, $s \in [t,T]$.

     \if{0}
  Fix $\o \neg \in \neg  \O^t$. For any $   \wt{\o}  \neg \in \neg  \O^{\t(\o)}  $,
  since $\t \big(\o  \neg \otimes_\t \neg  \wt{\o}\big)  \neg = \neg \t(\o)$
  by Lemma \ref{lem_bundle}, similar to \eqref{eq:p411x} and \eqref{eq:p411a} one can deduce that
 \bea
 &&  \hspace{-0.8cm}\big( \wt{Y}^{\Th}_{\t   \vee s} \big) (\o \neg \otimes_\t \neg  \wt{\o})
   \neg   =   \neg     \big( \wt{Y}^{\Th} \big)^{\t,\o}_s (\wt{\o}),
      ~\big( \wt{\ul{K}}^{\,\raisebox{-0.5ex}{\scriptsize $\Th$}}_{\t   \vee s} \big) (\o  \neg \otimes_\t \neg  \wt{\o})
        \neg =  \neg    \big( \wt{\ul{K}}^{\,\raisebox{-0.5ex}{\scriptsize $\Th$}} \big)^{\t,\o}_s (\wt{\o}) ,
   ~\big( \wt{\ol{K}}^{\, \raisebox{-0.8ex}{\scriptsize $\Th$}}_{\t   \vee s} \big) (\o  \neg \otimes_\t \neg  \wt{\o})
    \neg  =   \neg  \big( \wt{\ol{K}}^{\, \raisebox{-0.8ex}{\scriptsize $\Th$}} \big)^{\t,\o}_s (\wt{\o}),~
      \fa s \in [\t(\o),T],  \qq  \label{eq:p537a} \\
  &&  \hspace{0.2cm}  \Big(\b1_{\{r > \t\}} f   \big( r,\wt{X}^{\Th}_r, \wt{Y}^{\Th}_r, \wt{Z}^\Th_r , \mu_r , \nu_r  \big) \Big)
       (\o \neg \otimes_\t  \neg  \wt{\o}) \nonumber \\
  &&  \hspace{1.5cm} \neg  =   \neg   \b1_{\{r >  \t ( \o  )  \}}  f   \Big( r,\big(\wt{X}^{\Th}\big)^{\t,\o}_r (\wt{\o}), \big(\wt{Y}^{\Th}\big)^{\t,\o}_r(\wt{\o}), \big(\wt{Z}^\Th\big)^{\t,\o}_r(\wt{\o}) , \mu^{\t,\o}_r(\wt{\o}) , \nu^{\t,\o}_r(\wt{\o})  \Big)  \, dr ,\q   \fa r \in [t,T] ,   \label{eq:p537b} \\
  &&  
   \hspace{0.2cm}   \Big(\b1_{\{r > \t\}}  \wt{Z}^{\Th}_r \Big) (\o \otimes_\t \wt{\o}) =
         \b1_{\{r > \t(\o)\}} \big( \wt{Z}^{\Th} \big)^{\t,\o}_r (\wt{\o}) ,\q   \fa r \in [t,T] .   \label{eq:p537c}
  \eea
  \fi

 \ss  Fix $\o \neg \in \neg  \O^t$.   For any  $\wt{\o}  \neg \in \neg  \big( \cN^c \big)^{\t,\o}$,
   applying  \eqref{eq:p511} to the path $\o  \neg \otimes_\t \neg  \wt{\o}$ over period $\big[\t(\o), T \big]$
    as well as using  deductions similar to \eqref{eq:p411x} and \eqref{eq:p411a}, we obtain that
  \beas
     \left\{\ba{l}
 \dis   \big( \wt{Y}^{\Th} \big)^{\t,\o}_s (\wt{\o}) \neg - \neg  \xi^{\t,\o} (\wt{\o})  \neg + \neg  \Big( \wt{\ol{K}}^{\, \raisebox{-0.8ex}{\scriptsize $\Th$}} \Big)^{\t,\o}_T (\wt{\o})
      \neg - \neg   \Big( \wt{\ol{K}}^{\, \raisebox{-0.8ex}{\scriptsize $\Th$}}\Big)^{\t,\o}_s (\wt{\o})
       \neg - \neg   \Big(\wt{\ul{K}}^{\,\raisebox{-0.5ex}{\scriptsize $\Th$}}\Big)^{\t,\o}_T (\wt{\o})
       \neg + \neg  \Big( \wt{\ul{K}}^{\,\raisebox{-0.5ex}{\scriptsize $\Th$}}\Big)^{\t,\o}_s (\wt{\o})
      \neg + \neg   \big( M^\Th\big)^{\t,\o}_T(\wt{\o}) \neg - \neg \big( M^\Th\big)^{\t,\o}_s (\wt{\o})  \ss \\
  \dis = \neg \int_{\t(\o)}^T   \neg    f   \Big( r,\big(\wt{X}^{\Th}\big)^{\t,\o}_r(\wt{\o}), \big(\wt{Y}^{\Th}\big)^{\t,\o}_r(\wt{\o}), \big(\wt{Z}^\Th\big)^{\t,\o}_r(\wt{\o}) , \mu^{\t,\o}_r(\wt{\o}) , \nu^{\t,\o}_r(\wt{\o})  \Big)   dr
   ,  ~  s \in [\t(\o),T] ,\ss \\
   \dis \ul{l} \big(s, \big( \wt{X}^\Th \big)^{\t,\o}_s (\wt{\o}) \big) \neg \le  \neg \big( \wt{Y}^{\Th}_s \big)^{\t,\o} (\wt{\o})  \le       \ol{l} \big(s, \big( \wt{X}^\Th  \big)^{\t,\o}_s (\wt{\o})  \big) ,
    ~  s \in \big[\t(\o),T\big] , \q \hb{and} ~  \ss  \\
   \dis  \int_{\t(\o)}^T \dneg  \Big( \big( \wt{Y}^{\Th }\big)^{\t,\o}_s (\wt{\o})
      \neg - \neg  \ul{l} \big(s, \big( \wt{X}^\Th \big)^{\t,\o}_s (\wt{\o}) \big)  \Big) d \Big(\wt{\ul{K}}^{\,\raisebox{-0.5ex}{\scriptsize $\Th$}} \Big)^{\t,\o}_s (\wt{\o})
      \neg = \neg \int_{\t(\o)}^T \dneg  \Big( \ol{l} \big(s, \big( \wt{X}^\Th  \big)^{\t,\o}_s (\wt{\o})  \big)
        \neg - \neg  \big( \wt{Y}^{\Th }\big)^{\t,\o}_s (\wt{\o})  \Big)
         d \Big( \wt{\ol{K}}^{\, \raisebox{-0.8ex}{\scriptsize $\Th$}} \Big)^{\t,\o}_s (\wt{\o})    = 0 .
            \ea \right.
    \eeas
  By  \eqref{lem_basic_complement} and Corollary \ref{cor_rcpd11}, it holds for $P^t_0-$a.s. ~ $\o \in \O^t$ that
  $ P^{\t(\o)}_0 \big(  ( \cN^c  )^{\t,\o} \big) \neg = \neg  P^{\t(\o)}_0 \big(  ( \cN^{\t,\o}  )^c \big)
   \neg = \neg 1   $.
  Hence, one can deduce from the above system of equations and Proposition \ref{prop_FSDE_shift}   that for $P^t_0-$a.s. $\o  \neg \in \neg  \O^t$,  it holds $P^{\t(\o)}_0-$a.s. that
           \bea   \label{eq:p527}
    \left\{\ba{l}
 \dis   \big( \wt{Y}^{\Th} \big)^{\t,\o}_s   - \xi^{\t,\o}   + \Big( \wt{\ol{K}}^{\, \raisebox{-0.8ex}{\scriptsize $\Th$}} \Big)^{\t,\o}_T     -  \Big( \wt{\ol{K}}^{\, \raisebox{-0.8ex}{\scriptsize $\Th$}}\Big)^{\t,\o}_s
      -  \Big(\wt{\ul{K}}^{\,\raisebox{-0.5ex}{\scriptsize $\Th$}}\Big)^{\t,\o}_T    + \Big( \wt{\ul{K}}^{\,\raisebox{-0.5ex}{\scriptsize $\Th$}}\Big)^{\t,\o}_s
      \neg + \neg   \big( M^\Th\big)^{\t,\o}_T  \neg - \neg \big( M^\Th\big)^{\t,\o}_s   \ss \\
     \dis = \neg  \int_{\t(\o)}^T  \neg    f    \Big( r, \big(\wt{X}^{\Th}\big)^{\t,\o}_r , \big(\wt{Y}^{\Th}\big)^{\t,\o}_r , \big(\wt{Z}^\Th\big)^{\t,\o}_r ,  \mu^{\t,\o}_r  , \nu^{\t,\o}_r   \Big)   dr
     \neg = \neg \int_{\t(\o)}^T  \neg  f^{\Th^\o_\t}_T   \Big( r,  \big(\wt{Y}^{\Th}\big)^{\t,\o}_r , \big(\wt{Z}^\Th\big)^{\t,\o}_r   \Big)   dr ,  ~  s \in [\t(\o),T] ,   \vspace{1mm}  \\
   \dis \ul{L}^{\Th^\o_\t}_s   = \ul{l} \Big(s,   \wt{X}^{\Th^\o_\t}_s    \Big) \neg \le  \neg \big( \wt{Y}^{\Th}_s \big)^{\t,\o} \neg \le   \neg   \ol{l} \Big(s, \wt{X}^{\Th^\o_\t}_s  \Big) = \ol{L}^{\Th^\o_\t}_s  ,
    ~  s \in \big[\t(\o),T\big] , \q \hb{and} ~  \ss  \\
    \dis  \int_{\t(\o)}^T \neg  \Big( \big( \wt{Y}^{\Th }\big)^{\t,\o}_s
      \neg - \neg   \ul{L}^{\Th^\o_\t}_s    \Big) d \Big(\wt{\ul{K}}^{\,\raisebox{-0.5ex}{\scriptsize $\Th$}} \Big)^{\t,\o}_s
       = \int_{\t(\o)}^T \neg  \Big( \ol{L}^{\Th^\o_\t}_s
        \neg - \neg  \big( \wt{Y}^{\Th }\big)^{\t,\o}_s    \Big)
         d \Big( \wt{\ol{K}}^{\, \raisebox{-0.8ex}{\scriptsize $\Th$}} \Big)^{\t,\o}_s      = 0 .
            \ea \right.
    \eea

   Similar to \eqref{eq:p515}, we can deduce from Proposition 3.2.26 and Problem 3.2.27 of    \cite{Kara_Shr_BMSC}
     (both work for continuous local martingales) that   for $P^t_0-$a.s.~$\o \in \O^t$, it holds   $ P^{\t(\o)}_0-$a.s. that
  \beas
   \q   \big( M^\Th  \big)^{\t,\o}_s
     \neg=\neg \int_{\t(\o)  }^s \neg \Big(\b1_{\{r > \t\}}  \wt{Z}^{\Th}  \Big)^{\t,\o}_r dB^{\t(\o)}_r
     \neg=\neg \int_{\t(\o)  }^s  \neg  \b1_{\{r > \t(\o)\}} \big( \wt{Z}^{\Th} \big)^{\t,\o}_r (\wt{\o}) dB^{\t(\o)}_r
         \neg=\neg       \int_{\t(\o)  }^s  \neg  \big(\wt{Z}^\Th\big)^{\t,\o}_r  \, dB^{\t(\o)}_r   , ~  s \in [\t(\o),T] ,
  \eeas
 which together with \eqref{eq:p527} gives \eqref{eq:p541}. Therefore, it holds  for $P^t_0-$a.s.~$\o \in \O^t$ that 
  \beas
  \wt{Y}^{\Th^\o_\t}_s \big(T, \xi^{\t,\o}\big) (\wt{\o}) = \big( \wt{Y}^\Th (T,\xi) \big)^{\t,\o}_s (\wt{\o})
  = \wt{Y}^\Th_s (T,\xi)   ( \o \otimes_\t \wt{\o}) , \q  \fa    (s,\wt{\o}) \in [\t(\o),T] \times \O^{\t(\o)} .
  \eeas
 Taking $(s,\wt{\o})=\big(\t(\o), \Pi_{t,\t(\o)}(\o)\big)$ gives that
 $ \wt{Y}^{\Th^\o_\t}_{\t(\o)} \big(T, \xi^{\t,\o}\big)   = \big( \wt{Y}^\Th_\t (T,\xi) \big) ( \o ) $
 for $P^t_0-$a.s.~$\o \in \O^t$.   \qed

   \ss \no {\bf Proof of Lemma \ref{lem_shift_inverse2}:}
   For any $\wt{\cE} \in \sB([s,r])$ and $A \in \cF^s_r   $, applying    Lemma \ref{lem_shift_inverse} with $S=T$ yields that
      \beas
       \wh{\Pi}^{\,-1}_{\,t,s} \big( \wt{\cE} \times A  \big) =
       \big\{ (r,\o) \in [s,T] \times \O^t: \big(r, \Pi_{t,s}(\o)\big) \in \wt{\cE} \times A \big\}
       \neg  =  \neg \wt{\cE} \times \Pi^{\,-1}_{\,t,s} (A)
      \in  \sB([s,r])  \otimes \cF^t_r  ,
      \eeas
       which shows that all rectangular measurable sets of  $\sB([s,r]) \otimes \cF^s_r$
       belongs to $\wt{\L} \dfnn \big\{ \cD   \subset [s,r] \times \O^s: \wh{\Pi}^{\,-1}_{\,t,s}  (\cD) \in \sB([s,r])  \otimes \cF^t_r \big\}$. Clearly, $\wt{\L}$ is a $\si-$field of $[s,r] \times \O^s$. Thus it follows that $\sB([s,r]) \otimes \cF^s_r \subset  \wt{\L}$, i.e.,
        \bea \label{eq:xxc041}
        \wh{\Pi}^{\,-1}_{\,t,s}  (\cD) \in  \sB([s,r]) \otimes \cF^t_r, \q \fa  \cD \in \sB([s,r]) \otimes \cF^s_r.
        \eea

   \ss   Next, we show that  $ \big( dr \times dP^t_0 \big)  \circ \wh{\Pi}^{-1}_{t,s}
    = \big( dr \times dP^s_0 \big)$ on $\sB([s,T])\otimes \cF^s_T$: For any $ \cE \in \sB \big([s,T]  \big) $ and $ A \in  \cF^s_T $, Lemma \ref{lem_shift_inverse} with $S=T$
   again implies that
   \beas
      \big( dr \times dP^t_0 \big) \big( \wh{\Pi}^{-1}_{t,s} ( \cE \times A )\big)
   =  \big( dr \times dP^t_0 \big) \big( \cE \times  \Pi^{-1}_{t,s} (  A )\big)
   = | \cE |  P^t_0 \big(  \Pi^{-1}_{t,s} (  A )\big)
   = | \cE |   P^s_0  (    A )
   = \big( dr \times dP^s_0 \big) ( \cE \times A )  ,
   \eeas
   where $|\cE|$ denotes the Lesbegue measure of $\cE$. Thus the
   collection $\cC_s$ of all   rectangular measurable sets of $\sB\big([s,T]  \big) \otimes \cF^s_T$
       is contained in $ \L  \dfnn \big\{ \cD   \subset [s,T] \times \O^s \neg : \, \big( dr \times dP^s_0 \big) (\cD) = \big( dr \times dP^t_0 \big) \big( \wh{\Pi}^{-1}_{t,s} (\cD)\big) \big\}$. In particular, $\es \times \es \in  \L$
       and $ [s,T] \times \O^s  \in  \L $. For any $\cD \in \L$, one can deduce that
       \beas
      \hspace{-3mm}      \big( dr \neg \times \neg  dP^s_0 \big) \big(( [s,T]  \neg \times \neg  \O^s) \backslash \cD \big)
        & \tneg  \tneg =& \tneg  \tneg  \big( dr  \neg \times  \neg  dP^s_0 \big) \big([s,T]  \neg \times  \neg  \O^s \big)
          \neg - \neg  \big( dr  \neg \times \neg  dP^s_0 \big) (\cD)
          \neg =\neg \big( dr  \neg \times \neg  dP^t_0 \big) \big( \wh{\Pi}^{-1}_{t,s} \big([s,T] \neg  \times \neg  \O^s \big) \big)
       \neg -\neg \big( dr \neg  \times \neg  dP^t_0 \big) \big( \wh{\Pi}^{-1}_{t,s} (\cD)\big)  \\
       &  \tneg  \tneg    =& \tneg  \tneg  \big( dr  \neg \times \neg  dP^t_0 \big) \big( \wh{\Pi}^{-1}_{t,s} \big([s,T]  \neg \times \neg  \O^s \big)
       - \wh{\Pi}^{-1}_{t,s} (\cD)\big) \neg = \neg  \big( dr  \neg \times  \neg dP^t_0 \big) \big( \wh{\Pi}^{-1}_{t,s}
        \big(( [s,T] \times \O^s ) \backslash \cD \big)\big) .
       \eeas
      On the other hand, for any pairwisely-disjoint sequence $\{\cD_n\}_{n \in \hN} $ of $\L$ (i.e.   $\cD_m \cap \cD_n = \es$ given $m \ne n$), it is clear that $\big\{ \wh{\Pi}^{-1}_{t,s} (\cD_n) \big\}_{n \in \hN} $ is also a pairwisely-disjoint
      sequence. It follows that
       \beas
         \big( dr \neg \times \neg  dP^s_0 \big) \Big( \underset{n \in \hN}{\cup}\cD_n\Big)
         &=& \sum_{n \in \hN}  \big( dr \neg \times \neg  dP^s_0 \big) \big( \cD_n\big)
         =  \sum_{n \in \hN}  \big( dr \neg \times \neg  dP^t_0 \big) \big( \wh{\Pi}^{-1}_{t,s} (\cD_n)\big)    \\
         &=&   \big( dr \neg \times \neg  dP^t_0 \big) \Big( \underset{n \in \hN}{\cup} \wh{\Pi}^{-1}_{t,s} (\cD_n)\Big)
         =  \big( dr \neg \times \neg  dP^t_0 \big) \Big(  \wh{\Pi}^{-1}_{t,s} \big(\underset{n \in \hN}{\cup}\cD_n \big)\Big) .
       \eeas
       Hence, $\L$ is a Dynkin system. Since $\cC_s$ is closed under intersection, the Dynkin System Theorem shows that
        $      \sB\big([s,T]  \big) \otimes \cF^s_T = \si( \cC_s ) \subset \L$, i.e. $ \big( dr \times dP^t_0 \big)  \circ \wh{\Pi}^{-1}_{t,s}  = \big( dr \times dP^s_0 \big)$ on $\sB([s,T])\otimes \cF^s_T$.

 \ss Finally, let us discuss the $ \sP(\bF^t)/\sP(\bF^s)-$measurability of  $ \wh{\Pi}_{t,s}$:
 Let $\cD \in \sP(\bF^s) $.
 For any $r \in [s,T]$, since $ \cD \cap \big([s,r] \times \O^s\big) \in \sB([s,r]) \otimes \cF^s_r $,
   \eqref{eq:xxc041} implies that
  \bea
 \big( \wh{\Pi}^{-1}_{t,s}  (\cD) \big) \cap \big([t,r] \times \O^t\big)
   &=& \big( \wh{\Pi}^{-1}_{t,s}  (\cD) \big) \cap \big([s,r] \times \O^t\big)
   =\wh{\Pi}^{-1}_{t,s}  (\cD) \cap \big( \wh{\Pi}^{-1}_{t,s}  ([s,r] \times \O^s ) \big) \nonumber \\
   &=& \wh{\Pi}^{-1}_{t,s} \Big( \cD \cap \big([s,r] \times \O^s\big) \Big) \in \sB([s,r]) \otimes \cF^t_r  \subset
   \sB([t,r]) \otimes \cF^t_r \,.  \q  \label{eq:r201}
  \eea
  On the other hand, it is clear that  $ \big( \wh{\Pi}^{-1}_{t,s}  (\cD) \big) \cap \big([t,r] \times \O^t\big) = \es $
for any $ r \in [t,s)$, which together with \eqref{eq:r201} implies that $ \wh{\Pi}^{-1}_{t,s} (\cD) \in \sP(\bF^t) $.
As $\wh{\Pi}^{-1}_{t,s} (\cD)  \subset [s,T] \times \O^t$, we see that  $ \wh{\Pi}^{-1}_{t,s} (\cD) \in \sP_s (\bF^t) $.
  \qed

 \ss \no {\bf Proof of Proposition \ref{prop_paste_control}:} (1) Let us first discuss  the $\bF^t-$progressive measurability of $\hU-$valued process $\wh{\mu}$.
  Given $s \in [t,T]$ and $U \in \sB \big(\hU  \big)  $, we have to show that
   $\big\{ (r,  \o)  \in   \neg  [t,s] \times \O^t  \neg  :  \wh{\mu}_r (\o) \in  \neg  U \big\} \in \sB\big([t,s]\big) \otimes \cF^t_s$\,:      The $\bF^t-$progressive measurability of process $\mu$ implies that for any $\cD \subset \sB\big([t,s]\big) \otimes \cF^t_s$
    \bea  \label{mu_progressive}
    \big\{ (r,  \o) \in \cD  : \,  \mu_r (\o) \in U \big\}
      =  \big\{ (r,  \o) \in [t,s] \times \O^t: \,  \mu_r (\o) \in U \big\} \cap \cD
       \in \sB\big([t,s]\big) \otimes \cF^t_s ,
    \eea
    which together with \eqref{paste_control} leads to that
    \beas
        \big\{ (r,  \o) \in [t,s] \times \{\t > s \}: \, \wh{\mu}_r (\o) \in U \big\}
       = \big\{ (r,  \o) \in [t,s] \times \{\t > s \}: \,  \mu_r (\o) \in U \big\}    \in \sB\big([t,s]\big) \otimes \cF^t_s .
       \eeas
      Thus we only need  to     show that  $ \big\{ (r,  \o) \in [t,s] \times \{\t =t_n  \}:
      \, \wh{\mu}_r (\o) \in U \big\} \in \sB\big([t,s]\big) \otimes \cF^t_s$ for each $ t_n \in [t,s]$:

   \ss \no  (\,i) For  $n > N$ with $ t_n \le s $, \eqref{paste_control} and    \eqref{mu_progressive} imply that
       \beas
          \big\{ (r,  \o) \in [t,s] \times \{\t =t_n \}: \, \wh{\mu}_r (\o) \in U \big\}
       = \big\{ (r,  \o) \in [t,s] \times \{\t = t_n \}: \,  \mu_r (\o) \in U \big\}
       \in \sB\big([t,s]\big) \otimes \cF^t_s .
       \eeas


 \ss \no  (ii) For  $n \le   N$  with $ t_n \le s $,  let $A^n_0 \dfnn \{\t=t_n\} \Big\backslash \Big(\overset{\ell_n}{\underset{i=1}{\cup}} A^n_i \Big) \in \cF^t_{t_n} $. 
    One can deduce from \eqref{paste_control} and    \eqref{mu_progressive} that
       \bea \label{eq:k351}
       \big\{ (r,  \o) \in [t,s] \times   A^n_0   : \, \wh{\mu}_r (\o) \in U \big\}
        =  \big\{ (r,  \o) \in [t,s] \times   A^n_0   : \,  \mu_r (\o) \in U \big\}
        \in \sB\big([t,s]\big) \otimes \cF^t_s .
     \eea
    For any $i =1, \cds \neg , \ell_n$, since $A^n_i $ is an $\cF^t_{t_n}-$measurable subset of $ \{\t=t_n\}$, we see from \eqref{paste_control} that
     \beas
  \q    \big\{ (r,  \o) \neg \in  \neg  [t,s]  \neg \times  \neg   A^n_i  \neg   :   \wh{\mu}_r (\o)  \neg \in \neg  U \big\}
  \neg  =  \neg  \big\{ (r,  \o)  \neg \in \neg  [t,t_n)  \neg \times  \neg  A^n_i  \neg  :   \mu_r (\o)  \neg \in \neg  U \big\}
       \neg  \cup  \neg  \big\{ (r,  \o)  \neg \in \neg  [t_n,s]  \neg \times \neg   A^n_i  \neg  :   \big(\mu^n_i\big)_r\big(\Pi_{t,t_n} (\o)\big)  \neg \in \neg  U \big\} .
     \eeas
 Clearly, $ \big\{ (r, \o) \neg \in \neg  [t,t_n)  \neg \times  \neg  A^n_i  \neg :  \mu_r (\o)  \neg \in \neg  U \big\}
  \neg  \in \neg  \sB\big([t,s]\big) \neg \otimes  \neg  \cF^t_s $ by \eqref{mu_progressive}.
      Since $  \cD^n_i  \neg \dfnn \neg  \big\{ (r, \wt{\o})  \neg \in \neg  [t_n,s]  \neg \times  \neg \O^{t_n} \neg :
  \big(\mu^n_i\big)_r(\wt{\o})  \neg \in \neg  U \big\}
  \neg \in \neg  \sB\big([t_n,s]\big)  \neg \otimes  \neg  \cF^{t_n}_s$
 by  the $\bF^{t_n}-$progressive measurability of process $\mu^n_i  $, one can deduce
 from Lemma \ref{lem_shift_inverse2} that
     \beas
   && \hspace{-2cm}  \big\{ (r,  \o) \in [t_n,s] \times  A^n_i  : \, \big(\mu^n_i\big)_r\big(\Pi_{t,t_n} (\o)\big) \in U \big\}
   =  \big\{ (r,  \o) \in [t_n,s] \times  A^n_i  : \,  \big(r, \Pi_{t,t_n} (\o)\big) \in \cD^n_i  \big\}    \\
   &&=  \big\{ (r,  \o) \in [t_n,T] \times \O^t    : \,  \wh{\Pi}_{t,t_n} (r,\o)  \in \cD^n_i  \big\}
   \cap \big( [t_n,s] \times  A^n_i \big)  \\
   &&=   \wh{\Pi}^{\,-1}_{t,t_n} \big( \cD^n_i \big)  \cap \big( [t_n,s] \times  A^n_i \big)   \in
      \sB\big([t_n,s]\big)  \neg \otimes  \neg  \cF^t_s \subset  \sB\big([t,s]\big)  \neg \otimes  \neg  \cF^t_s .
     \eeas
 It follows that $ \big\{ (r,  \o)   \in     [t,s]    \times      A^n_i  \neg   :   \wh{\mu}_r (\o)    \in    U \big\} \in  \sB\big([t,s]\big)    \otimes \cF^t_s$.
 Then taking union over $i =1,\cds, \ell_n$ and combining with \eqref{eq:k351} yield that
   $ \big\{ (r,  \o)   \in     [t,s]    \times   \{\t = t_n \}  :  \,  \wh{\mu}_r (\o)    \in    U \big\}
    \in  \sB\big([t,s]\big)    \otimes \cF^t_s$.

 \ss \no (2) For $n =1, \cds, N$ and $ i =1,\cds, \ell_n$, since
  \beas
  \q  \big\{ (r,  \o) \neg \in \neg  [t_n, T]  \neg \times \neg  A^n_i \neg : \,
   \wh{\mu}_r (\o)  \neg \in \neg  \hU \, \backslash \hU_0 \big\}
    =     \big( [t_n, T]  \neg \times \neg  A^n_i \big) \cap  \wh{\Pi}^{-1}_{t,t_n} \big( \big\{ (r,  \wt{\o})
      \neg \in \neg  [t_n,T]  \neg \times  \neg  \O^{t_n}  \neg  :
    \,  \big(\mu^n_i\big)_r (\wt{\o})  \neg \in \neg  \hU \, \backslash \hU_0 \big\} ,
  \eeas
    Lemma \ref{lem_shift_inverse2} again implies that
  \bea
    (dr \neg \times \neg  d P^t_0 ) \big( \big\{ (r,  \o)  \neg \in \neg  [t_n, T]  \neg \times \neg  A^n_i  \neg  : \,
     \wh{\mu}_r (\o)  \neg \in \neg  \hU \, \backslash \hU_0 \big\}\big)
       & \tneg  \dneg \le &  \tneg  \dneg  (dr  \neg \times \neg  dP^t_0) \Big( \wh{\Pi}^{-1}_{t,t_n} \big( \big\{ (r,  \wt{\o})
        \neg \in  \neg  [t_n,T]  \neg \times  \neg  \O^{t_n}  \neg  :
    \,  \big(\mu^n_i\big)_r (\wt{\o})  \neg \in \neg  \hU \, \backslash \hU_0 \big\} \big) \Big) \nonumber \\
    &  \tneg  \dneg   & \hspace{-5cm} =  (dr  \neg \times \neg  d P^{t_n}_0 ) \big( \big\{ (r,  \wt{\o})  \neg \in  \neg [t_n,T]
     \neg \times  \neg   \O^{t_n}  \neg  :
    \,  \big(\mu^n_i\big)_r (\wt{\o})  \neg \in \neg  \hU \, \backslash \hU_0 \big\} \big) =0 .  \label{eq:xb011}
  \eea
    Clearly,
   $    (dr \neg \times \neg  d P^t_0 ) \big( \big\{ (r,  \o)  \neg \in \neg  \[t,\t\[ \, \cup \, \[\t,T\]_{A_0}  \neg  : \,
     \wh{\mu}_r (\o)  \neg \in \neg  \hU \, \backslash \hU_0 \big\}\big)
    \le (dr \neg \times \neg  d P^t_0 ) \big( \big\{ (r,  \o)  \neg \in \neg  [t,T] \times  \O^t   \neg  : \,
      \mu_r (\o)  \neg \in \neg  \hU \, \backslash \hU_0 \big\}\big) =0 $, which together with
      \eqref{eq:xb011} shows that $ \wh{\mu}_r \in \hU_0 $, $dr \times   d P^t_0 -$a.s.

  \ss \no (3) Next, we show that $ E_t \int_t^T \big[ \wh{\mu}_r  \big]^2_{\overset{}{\hU}} \, dr   < \infty$:
  By \eqref{paste_control},
   \beas
  && \hspace{-1cm}   E_t \int_t^T \big[ \wh{\mu}_r  \big]^2_{\overset{}{\hU}} \, dr
   =  \int_{ \o \in \O^t}  \int_t^T \big[ \wh{\mu}_r (\o) \big]^2_{\overset{}{\hU}} \, dr \, d P^t_0 (\o) \\
  && =  \left( \int_{ \o \in \O^t}  \int_t^{\t(\o)}
  + \int_{ \o \in A_0}  \int_{\t(\o)}^T \right) \big[  \mu_r (\o) \big]^2_{\overset{}{\hU}} \, dr \, d P^t_0 (\o)
  + \sum^N_{n=1}\sum^{\ell_n}_{i=1} \int_{ \o \in A^n_i }  \int_{t_n}^T \Big[ \big(\mu^n_i\big)_r\big(\Pi_{t,t_n} (\o)\big) \Big]^2_{\overset{}{\hU}} \, dr \, d P^t_0 (\o)  \\
  && \le   \int_{ \o \in \O^t}  \int_t^T  \big[  \mu_r (\o) \big]^2_{\overset{}{\hU}} \, dr \, d P^t_0 (\o)
  + \sum^N_{n=1}\sum^{\ell_n}_{i=1} \int_{ \o \in \O^t  }  \int_{t_n}^T \Big[ \big(\mu^n_i\big)_r\big(\Pi_{t,t_n} (\o)\big) \Big]^2_{\overset{}{\hU}} \, dr \, d P^t_0 (\o) .
  \eeas
  For any $ n =1, \cds \neg , N$ and $i= 1, \cds \neg , \ell_n $, applying  Lemma \ref{lem_shift_inverse} with $(s,S)=(t_n,T)$ yields that
  \beas
 && \hspace{-2cm} \int_{ \o \in \O^t }  \int_{t_n}^T \Big[ \big(\mu^n_i\big)_r\big(\Pi_{t,t_n} (\o)\big) \Big]^2_{\overset{}{\hU}}
   \, dr \, d P^t_0 (\o)
    =     \int_{ \wt{\o} \in \O^{t_n}  }  \int_{t_n}^T \Big[ \big(\mu^n_i\big)_r ( \wt{\o}) \Big]^2_{\overset{}{\hU}} \, dr \, d P^t_0 \big( \Pi^{-1}_{t,t_n} ( \wt{\o}) \big)  \\
 && =    \int_{ \wt{\o} \in \O^{t_n}  }  \int_{t_n}^T \Big[ \big(\mu^n_i\big)_r ( \wt{\o}) \Big]^2_{\overset{}{\hU}} \, dr \, d P^{t_n}_0   ( \wt{\o})
 + E_{t_n} \int_{t_n}^T \Big[ \big(\mu^n_i\big)_r  \Big]^2_{\overset{}{\hU}} \, dr   < \infty .
   \eeas
  Thus it follows that
   \bea   \label{eq:r473}
  E_t \int_t^T \big[ \wh{\mu}_r  \big]^2_{\overset{}{\hU}} \, dr
    \le  E_t \int_t^T \big[  \mu_r  \big]^2_{\overset{}{\hU}} \, dr
  + \sum^N_{n=1}\sum^{\ell_n}_{i=1}  E_{t_n} \int_{t_n}^T \Big[ \big(\mu^n_i\big)_r \Big]^2_{\overset{}{\hU}} \, dr < \infty ,
  \eea
  which together with part (1) shows that $ \wh{\mu} \in \cU^t$.

 \ss \no (4) Let   $(r,\o) \in    \[ \t,T  \]_{A^n_i}   = [t_n,T] \times  A^n_i$
    for some $ n=1 \cds \neg , N$ and
    $ i=1,\cds \neg , \ell_n $. For any  $\wt{\o} \in \O^{t_n}$, since
     $   \o \otimes_{t_n} \wt{\o} \in A^n_i  $ by Lemma \ref{lem_element},
  it follows from    \eqref{paste_control}    that
        $  \wh{\mu}^{t_n,\o}_r (\wt{\o})     =    \wh{\mu}_r \big(\o \otimes_{t_n} \wt{\o} \big)
   = \big(\mu^n_i\big)_r \Big(\Pi_{t,t_n} \big(\o \otimes_{t_n} \wt{\o} \big) \Big)
   = \big(\mu^n_i\big)_r ( \wt{\o} ) $.

     On the other hand, we consider $(r,\o) \in    \[ \t,T  \]_{A_0}  $.
   For any  $\wt{\o} \in \O^{\t(\o)} $,    we claim that
  \bea  \label{eq:r461}
  \o \otimes_\t \wt{\o} \in A_0 .
  \eea
 Assume not, i.e. $\o \otimes_\t \wt{\o} \in A^n_i $ for some $ n=1 \cds \neg , N$ and
    $ i=1,\cds \neg , \ell_n $.   By    Lemma \ref{lem_bundle},
     $\t(\o) = \t \big(\o \otimes_\t \wt{\o} \big) = t_n $. So
     $\o \otimes_{t_n} \wt{\o} = \o \otimes_\t \wt{\o} \in A^n_i   $ and
     Lemma \ref{lem_element}   shows that  $ \o \neg \in \neg A^n_i $, a contradiction appears.
     Thus $ \o \otimes_\t \wt{\o} \in A_0 $.
     As $r \ge \t(\o)= \t \big(\o \otimes_\t \wt{\o} \big) $, we see that $\big(r,\o \otimes_\t \wt{\o} \big) \in    \[ \t,T  \]_{A_0}$. Then      \eqref{paste_control} yields   that
      $   \wh{\mu}^{\t,\o}_r (\wt{\o}) =   \wh{\mu}_r \big(\o \otimes_\t \wt{\o} \big)
       = \mu_r \big(\o \otimes_\t \wt{\o} \big) =\mu^{\t, \o}_r (\wt{\o}) $.  \qed

  \if{0}
     \begin{lemm} \label{lem_truncation_measurable}
Let $0 \le t \le S \le T $.  For any $r \in [t,S]$ and $A \in \cF^{t,T}_r$, we have $ \Pi^{T,S}_{t,t} (A) \in \cF^{t,S}_r$.

  \end{lemm}

  \ss \no {\bf Proof:} For\;simplicity, let\;us\;denote\;$\Pi^{T,S}_{t,t} $ by $\Pi$. Fix\;$r  \neg \in \neg  [t,S]$.
  For\;any\;$r'  \neg \in \neg  [t, r] $\;and\;$\cE  \neg \in \neg   \sB (\hR^d)$, we\;can\;deduce\;that
      \bea   \label{eq:t171}
   \Pi \Big( \big( B^{t,T}_{r'} \big)^{-1} (\cE) \Big) = \Pi \Big(\big\{ \o \in \O^{t,T}:  \o(r') \in \cE \big\}\Big)
   \subset  \big\{ \wt{\o} \in \O^{t,S}:  \wt{\o} (r') \in \cE \big\} = \big( B^{t,S}_{r'} \big)^{-1} (\cE) \in \cF^{t,S}_r .
   \eea
     For any $\wt{\o} \in \O^{t,S}$ with $ \wt{\o} (r') \in \cE$,
    $ \o(s) \dfnn  \wt{\o}(s \land S) $, $\fa s \in [t,T]$  defines a path  of $\O^{t,T} $ with $\o(r') = \wt{\o}(r') \in \cE $.
    As $\wt{\o} = \Pi(\o)$, we see that the  $\subset$ in \eqref{eq:t171} is actually an equality. Thus,
    all the generating sets of $\cF^{t,T}_r $ belong to
    $ \L_r \dfnn \big\{ A \in \cF^{t,T}_r: \Pi(A ) \in \cF^{t,S}_r \big\}$. In particular, $\es, \O^{t,T} \in \L_r$.

    Let $A \in \L_r$.
  we claim that $\Pi(A) \cap \Pi(A^c) = \es$. Assume not, there exist $\o_1 \in A$ and $\o_2 \in A^c$ such that
  $\Pi(\o_1)= \Pi(\o_2)$, i.e. $\o_1(\l) = \o_2(\l)$, $\fa \l \in [t,S]$. As $A  \in \cF^{t,T}_r$,
   Lemma \ref{lem_element} implies that $\o_2 \in  \o_1 \otimes_r \O^{r,T} \in A $, a contradiction appears. Thus,
   $\Pi(A) \cap \Pi(A^c) = \es $, which together with the fact that $ \Pi(A) \cup \Pi(A^c) 
   = \Pi(\O^{t,T} ) = \O^{t,S}  $ shows that $ \Pi(A^c) =  \big( \Pi(A) \big)^c \in  \cF^{t,S}_r$.
    As $ A \in  \cF^{t,T}_r$,  we  clearly have $A^c \in \cF^{t,T}_r$. Thus $ A^c \in \L_r$.

On the other hand, for $\{A_n\}_{n \in \hN} \subset \L_r$, one can easily deduce that
 $   \underset{n \in \hN}{\cup} A_n  \in \cF^{t,T}_r  $  and
 $  \Pi \big( \underset{n \in \hN}{\cup} A_n \big)   = \underset{n \in \hN}{\cup} \Pi  ( A_n  ) \in  \cF^{t,S}_r $.
 Namely, $ \underset{n \in \hN}{\cup} A_n \in \L_r $. Thus we see that $ \L_r$ is a $\si-$field of $\O^{t,T} $.
 It follows that $ \L_r  =  \cF^{t,T}_r $, i.e. $ \Pi(A) \in \cF^{t,S}_r $ for any $  A \in \cF^{t,T}_r $. \qed
   \fi

 \ss \no {\bf Proof of Proposition \ref{prop_paste_strategy}:}
  (1)  We first assume that $\a \in  \cA^t$ and $\{ \a^n_i \}^{\ell_n}_{i=1}   \subset  \cA^{ t_n} $
   for $n =1,\cds \neg ,  N$.
   Then \eqref{eq:r503} holds for  all $(r,  \o) \in [t,T] \times \O^t$ except on a  $dr  \neg \times \neg  d P^t_0-$null set $\cD $. And for $n \neg = \neg 1, \cds \neg , N$ and $i \neg = \neg 1,\cds \neg ,\ell_n$,
   $\a^n_i$ satisfy \eqref{eq:r503} for some $\k^n_i >0$ and some  non-negative $\bF^{t_n}-$measurable process $\Psi^{n,i}$  with $ E_{t_n} \neg \int_{t_n}^T \neg \big( \Psi^{n,i}_s \big)^2   ds < \infty$: i.e.,
   it holds all $(r, \wt{\o}) \in [t_n,T] \times \O^{t_n}$ except on a  $dr  \neg \times \neg  d P^{t_n}_0-$null set $\cD^n_i$ that
  \beas
 \big[ \a^n_i(r,\wt{\o},v) \big]_{\overset{}{\hU}} \le \Psi^{n,i}_r (\wt{\o}) + \k^n_i [v]_{\overset{}{\hV}} \,, \q \fa v \in \hV .
   \eeas

  \ss Fix  $n \neg = \neg 1, \cds \neg , N$ and $i \neg = \neg 1,\cds \neg ,\ell_n$.
      The $ \sP_{t_n}(\bF^t)/\sP(\bF^{t_n})-$measurability of
      the mapping $\wh{\Pi}_{t,t_n} \neg : [t_n,T] \times \O^t \to   [t_n,T] \times \O^{t_n} $
     by Lemma \ref{lem_shift_inverse2} shows that  the function
     $ \a^n_i \big(   \wh{\Pi}_{t,t_n}   \neg   (r,\o),v \big) =   \a^n_i \big( r, \Pi_{t,t_n} \neg (\o),v \big) $,
   $ \fa (r,\o, v) \in [t_n,T] \times \O^t \times \hV$
   is   $\sP_{t_n}(\bF^t)   \otimes \sB(\hV)    \big/    \sB(\hU) -$measurable, which together with
  the fact  $ \[\t,T\]_{A^n_i}  \neg = \neg [t_n,T]     \times     A^n_i \neg \in \neg \sP(\bF^t)$ implies
   that the function $\wh{\a}$ is $\sP \big(  \bF^t  \big)   \otimes \sB(\hV)    \big/    \sB(\hU) -$measurable.

     \if{0}
   Given $U \in \sB(\hU)$,
 \beas
  \{ (r,\o,v) \in [t,T] \times \O^t \times \hV: \wh{\a} (r,\o,v) \in U\}
 & = & \underset{n=1}{\overset{N}{\cup}}  \underset{i=1}{\overset{\ell_n}{\cup}}
 \big( [t_n,T]  \neg   \times  \neg   A^n_i  \neg \times \neg \hV  \big) \cap \{ (r,\o,v) \in [t_n,T] \times \O^t \times \hV:  \a^n_i (r,\o,v) \in U\} \\
 && \bigcup \big(\[t,\t\[ \, \cup \, \[\t,T\]_{ A_0 }\big) \cap \{ (r,\o,v) \in [t,T] \times \O^t \times \hV: \a (r,\o,v) \in U\}
 \eeas
 \fi

 \ss  Similar to \eqref{eq:xb011}, one can deduce from Lemma \ref{lem_shift_inverse2} that
 \beas
    (dr \neg \times \neg  d P^t_0 ) \big( \big\{ (r,  \o)  \neg \in \neg  [t_n, T]  \neg \times \neg  A^n_i  \dneg  :
     \wh{\a}(r, \o, \hV_0)
    \backslash \hU_0 \neg \ne \neg \es \big\}\big)
       \neg  \le \neg   (dr  \neg \times \neg  d P^{t_n}_0 ) \big( \big\{ (r,  \wt{\o})  \neg \in  \neg [t_n,T]
     \neg \times  \neg   \O^{t_n}  \dneg  :
     \a^n_i (r,  \wt{\o} , \hV_0)   \backslash \hU_0 \neg  \ne \neg \es \big\} \big) \neg = \neg 0 .
  \eeas
   As   $    (dr \neg \times \neg  d P^t_0 ) \big( \big\{ (r,  \o)  \neg \in \neg  \[t,\t\[ \, \cup \, \[\t,T\]_{A_0}  \neg  : \,
     \wh{\a}(r, \o, \hV_0)    \backslash \hU_0  \ne \es \big\}\big)
    \le (dr \neg \times \neg  d P^t_0 ) \big( \big\{ (r,  \o)  \neg \in \neg  [t,T] \times  \O^t   \neg  : \,
      \a(r, \o, \hV_0)   \backslash \hU_0  \ne \es \big\}\big) =0 $, we see that
      $ \wh{\a}(r,  \hV_0)    \subset \hU_0 $, $dr \times   d P^t_0 -$a.s.

        Since the mapping   $ \wh{\Pi}_{t,t_n}$ is also
  $\sB([t_n,T])\otimes \cF^t_T / \sB([t_n,T])\otimes \cF^{t_n}_T-$measurable by Lemma \ref{lem_shift_inverse2} again,
       the process
     $ \Psi^{n,i} \big( \wh{\Pi}_{t,t_n}   \neg   (r,\o) \big)  \neg = \neg  \Psi^{n,i} \big( r, \Pi_{t,t_n} \neg (\o)  \big) $,
   $ \fa (r,\o ) \neg \in \neg [t_n,T] \neg \times \neg \O^t  $      is   $\sB([t_n,T]) \neg \otimes  \neg  \cF^t_T    \big/    \sB(\hR) -$measurable,
  which together with the fact  $ \[\t,T\]_{A^n_i}  \neg = \neg [t_n,T]     \times     A^n_i \neg \in \neg \sB \big( [t,T] \big) \otimes \cF^t_T$ gives rise to a non-negative
   measurable process on $(\O^t,\cF^t_T) $:
           \beas
  \q \wh{\Psi}_r  ( \o )   \dfnn        \begin{cases}
      \Psi^{n,i}_r \big(  \Pi_{t,t_n} \neg (\o)  \big)  ,
       & \hb{if $(r,\o ) \neg \in \neg    \[\t,T\]_{A^n_i}  \neg = \neg
      [t_n,T]  \neg \times \neg   A^n_i$    for $ n  \neg = \neg  1 \cds \neg , N$ and
    $ i \neg = \neg 1,\cds \neg , \ell_n $} ,   \\
     \Psi_r ( \o ),  & \hb{if $(r,\o) \in \[t,\t\[ \, \cup \, \[\t,T\]_{ A_0 }$.     }
    \end{cases}
    \eeas
    Similar to \eqref{eq:r473}, one can show that
    $   E_t \int_t^T   \wh{\Psi}^2_r     \, dr
    \le  E_t \int_t^T    \Psi^2_r    \, dr
    + \sum^N_{n=1}\sum^{\ell_n}_{i=1}  E_{t_n} \int_{t_n}^T   \big( \Psi^{n,i}_r \big)^2   \, dr < \infty $.

   \ss  Let $\wh{\k}  \dfnn  \k  \vee \, \max\{ \k^n_i: n = 1, \cds, N \hb{ and } i=1,\cds, \ell_n \} > 0$.
      For any $(r, \o) \in ( \[t,\t\[ \, \cup \, \[\t,T\]_{ A_0 } ) \big\backslash \cD $, one has
 \beas
 \big[ \wh{\a} (r,\o,v)  \big]_{\overset{}{\hU}}
    =   \big[ \a(r,\o,v)  \big]_{\overset{}{\hU}} \le \Psi_r (\o) + \wh{\k} [v]_{\overset{}{\hV}}
    = \wh{\Psi}_r (\o) + \wh{\k} [v]_{\overset{}{\hV}} \, , \q \fa v \in \hV ;
 \eeas
  On the other hand,   if $(r, \o) \in [t_n,T]  \neg \times \neg   A^n_i  \big\backslash \wh{\Pi}^{-1}_{t,t_n} ( \cD^n_i ) $
    for some $n \neg = \neg 1, \cds \neg , N$ and $i \neg = \neg 1,\cds \neg ,\ell_n$, then $  \big(r, \Pi_{t,t_n} (\o) \big) = \wh{\Pi}_{t,t_n}(r,\o) \in ( [t_n,T]  \neg \times \neg \O^{t_n} )   \big\backslash \cD^n_i $ and it follows that
 \beas
  \big[ \wh{\a} (r,\o,v)  \big]_{\overset{}{\hU}}
    =   \big[ \a^n_i \big(   r, \Pi_{t,t_n} (\o),v \big) \big]_{\overset{}{\hU}}
     \le \Psi^{n,i}_r \big( \Pi_{t,t_n} (\o) \big) + \wh{\k} [v]_{\overset{}{\hV}}
    = \wh{\Psi}_r   (\o)   + \wh{\k} [v]_{\overset{}{\hV}} \, , \q \fa v \in \hV .
 \eeas
  Since $ \wt{\cD} \dfnn  \cD \cup \Big( \underset{n=1}{\overset{N}{\cup}}
  \underset{i=1}{\overset{\ell_n}{\cup}}   \,   \wh{\Pi}^{-1}_{t,t_n} ( \cD^n_i )  \Big) $  is
    a $dr    \times    d P^t_0-$null set by Lemma \ref{lem_shift_inverse2}, we see that
    $\wh{\a}$ satisfies \eqref{eq:r503} $dr \times d P^t_0-$a.s. Therefore, $\wh{\a}$ is an $\cA^t-$strategy.

 \ss \no (2)    Next, let us verify \eqref{eq:r401}:  Fix $\nu \in \cV^t$.  For $n \neg = \neg 1, \cds \neg , N$,
    Proposition \ref{prop_dwarf_strategy} (1) shows that $\nu^{\, t_n,\o} 
    \in \cV^{t_n}$ for all $\o \in \O^t$ except on a $P^t_0-$null set $\cN_n$. Let $\o \in \underset{n=1}{\overset{N}{\cap}} \cN^c_n$ and $r \in [\t(\o),T]$.
     If $\o \in A^n_i$ for some $n \neg = \neg 1, \cds \neg , N$ and $i \neg = \neg 1,\cds \neg ,\ell_n$, then    $ r \ge \t(\o)= t_n$. For any $\wt{\o} \in \O^{t_n}$, since
     $   \o \otimes_{t_n} \wt{\o} \in A^n_i  $ by Lemma \ref{lem_element},
  it follows from    \eqref{paste_strategy}    that
       \beas
   \big( \wh{\a} \lan \nu \ran \big)^{t_n,\o}_r (\wt{\o})
  =\big( \wh{\a} \lan \nu \ran \big)_r \big(\o \otimes_{t_n} \wt{\o} \big)
   = \wh{\a} \big( r, \o \otimes_{t_n} \wt{\o}, \nu_r(\o \otimes_{t_n} \wt{\o}) \big)
    =      \a^n_i \big(r, \wt{\o}, \nu^{t_n ,\o}_r (\wt{\o}) \big)
         =  \big( \a^n_i \lan \nu^{t_n ,\o} \ran \big)_r  ( \wt{\o}  ) .
   \eeas
  Otherwise, if $ \o  \in     A_0   $, we have seen from \eqref{eq:r461} that $ \o \otimes_\t \O^{\t(\o)} \subset A_0$.
   For any  $\wt{\o} \in \O^{\t(\o)} $,
     since $r \ge \t(\o)= \t \big(\o \otimes_\t \wt{\o} \big) $ by Lemma \ref{lem_bundle}, we see that $\big(r,\o \otimes_\t \wt{\o} \big) \in    \[ \t,T  \]_{A_0}$. Then      \eqref{eq:r461} leads to     that
      \beas
         \big( \wh{\a} \lan \nu \ran \big)^{\t,\o}_r (\wt{\o})
          =   \big( \wh{\a} \lan \nu \ran \big)_r \big(\o \otimes_\t \wt{\o} \big)
          =  \wh{\a} \big( r, \o \otimes_\t \wt{\o}, \nu_r(\o \otimes_\t \wt{\o}) \big)
          =   \a \big( r, \o \otimes_\t \wt{\o}, \nu_r(\o \otimes_\t \wt{\o}) \big)
          = \big(  \a \lan \nu \ran \big)^{\t,\o}_r (\wt{\o})  .
       \eeas

    \ss \no (3) Now,  let  $\a \in \wh{\cA}^t$ and $\{ \a^n_i \}^{\ell_n}_{i=1}  \subset \wh{\cA}^{\, t_n} $ for $n =1,\cds \neg ,N$.  We shall show that     $\wh{\a}$ satisfies \eqref{eq:r742}, and it thus belongs to $  \wh{\cA}^t $:
 Fix $\e > 0$.  There exist  $\d>0$ and a closed subset $F$ of $\O^t$ 
 with    $   P^t_0\big( F \big)  \neg > \neg  1-\frac{\e}{2}$ such that for any $\o ,\o'   \in    F$ with $ \|\o-\o' \|_t < \d $
\bea   \label{eq:r637}
\underset{r \in [t,T]}{\sup} \, \underset{v \in  \hV}{\sup} \; \rho_{\overset{}{\hU}} \big( \a(r,\o,v) ,  \a(r,\o',v) \big)
 < \frac{\e}{2}    .
 \eea
 As  $A_0 \in  \cF^t_T = \sB(\O^t) $ by \eqref{eq:xxc023},
 we can find   a closed subset $ F_0 $ of $\O^t$ that is included in $A_0$
 and satisfies    $\dis  P^t_0 \big( A_0   \backslash F_0 \big) < \frac{\e}{8}$
 (see e.g. Proposition 15.11 of \cite{Royden_real}).

 Let $ \ell_*  \dfnn \sum^N_{n=1}  \ell_n   $. Given  $n=1,\cds \neg , N$ and $i=1,\cds \neg, \ell_n$,
 there exist  $\d^n_i>0$ and a closed subset $\wt{F}^n_i$ of $\O^{t_n}$
  with    $   P^{t_n}_0\big( \wt{F}^n_i \big)  \neg > \neg  1- \frac{\e}{4 \ell_*}$ such that
  for any $\wt{\o}, \wt{\o}' \in  \wt{F}^n_i$
    with $\| \wt{\o} - \wt{\o}' \|_{t_n} < \d^n_i $ 
\bea   \label{eq:r637b}
\underset{r \in [t_n , T]}{\sup} \, \underset{v \in \hV}{\sup} \; \rho_{\overset{}{\hU}} \big( \a^n_i (r,\wt{\o},v) ,  \a^n_i (r,\wt{\o}',v) \big)  < \frac{\e}{4 \ell_*}  \,   .
 \eea
      Applying    Lemma \ref{lem_shift_inverse} with $(s,S)=(t_n,T)$ shows that $\Pi^{-1}_{t,t_n} \big( \wt{F}^n_i \big)$
 is a   closed subset of $\O^t$     and that
      \bea   \label{eq:r705}
   P^t_0\big(   \Pi^{-1}_{t,t_n}  ( \wt{F}^n_i  ) \big)
     =   P^{t_n}_0 \big(   \wt{F}^n_i \big) > 1- \frac{\e}{4 \ell_*}    .
 \eea
  Similar to $F_0$, one can find a  closed subset $F^n_i$ of $\O^t$   that is included in
 $  A^n_i$ and satisfies $\dis  P^t_0 \big( A^n_i \backslash F^n_i \big) < \frac{\e}{8\ell_*}$. Then
  \bea \label{eq:r723}
 P^t_0 ( F_0 ) +  \sum^N_{n = 1} \sum^{\ell_n}_{i=1} P^t_0 \big(      F^n_i    \big)
 =P^t_0 ( A_0 ) - P^t_0 ( A_0 \backslash F_0 ) +  \sum^N_{n = 1} \sum^{\ell_n}_{i=1}   P^t_0 \big(A^n_i \big) - \sum^N_{n = 1} \sum^{\ell_n}_{i=1} P^t_0 \big( A^n_i \backslash F^n_i \big)  > 1-\frac{\e}{4} .
  \eea

Since $\wh{F}_0 \dfnn F \cap F_0$ and  $\wh{F}^n_i \dfnn  F \cap F^n_i \cap \Pi^{-1}_{t,t_n}  \big( \wt{F}^n_i \big) $,
 $n=1,\cds \neg , N$ and $i=1,\cds \neg, \ell_n$ are disjoint closed subsets of $\O^t$,
 we let $ \ul{\d} >0$ stand for the minimal distance between any two of them.   Let
   \beas
   \wh{F} \dfnn   \wh{F}_0 \cup \big(\cup \{\wh{F}^n_i  : n=1,\cds \neg , N \hb{ and } i=1,\cds \neg, \ell_n \} \big)
  ~  \hb{ and } ~  \wh{\d} \dfnn \big(  \hb{$\frac12$} \ul{\d} \, \big) \land \d \land \hb{$\frac12$} \min\{\d^n_i, n=1,\cds \neg , N \hb{ and } i=1,\cds \neg, \ell_n \} .
   \eeas

   For any $A_1 , A_2 \in \cF^t_T$, one has
   \bea   \label{eq:r719}
  P^t_0( A_1 \cap A_2 ) =  P^t_0( A_1) +  P^t_0( A_2 ) - P^t_0( A_1 \cup A_2 ) \ge P^t_0( A_1) +  P^t_0( A_2 ) - 1  .
   \eea
    Taking $A_1 =  F  $ and $A_2 = \Big( \underset{n =1 }{\overset{N}{ \cup } }
  \underset{i=1}{\overset{\ell_n}{\cup}} F^n_i \Big)  \cup  F_0  $, we can deduce from \eqref{eq:r723}   that
 \bea   \label{eq:r727}
        P^t_0 ( F \cap F_0   )  +  \sum^N_{n = 1} \sum^{\ell_n}_{i=1} P^t_0 (F \cap  F^n_i   )
     > 1-\frac{3}{4} \e .
 \eea
 Also, for  $n=1,\cds \neg , N$ and $i  \neg = \neg  1, \cds \neg , \ell_n$, letting $A_1  \neg =  \neg F \cap  F^n_i  $
  and $A_2  \neg = \neg  \Pi^{-1}_{t,t_n} \big( \wt{F}^n_i  \big)$  in \eqref{eq:r719}, we see from  \eqref{eq:r705} that
 \beas
    P^t_0  \big( \wh{F}^n_i \big) =   P^t_0 \Big(  F \cap F^n_i \cap \Pi^{-1}_{t,t_n}  \big( \wt{F}^n_i \big)   \Big)
  >    P^t_0  \big( F \cap F^n_i \big) -   \frac{\e}{4 \ell_*} ,
 \eeas
 which together with \eqref{eq:r727} leads to that
   \beas
   P^t_0 \big(\wh{F}\big) =   P^t_0 \big(\wh{F}_0\big)
  + \sum^N_{n = 1} \sum^{\ell_n}_{i=1}  P^t_0  \big( \wh{F}^n_i \big)  > 1-  \e .
 \eeas

   Now let $\o,\o' \in \wh{F}$ with  $ \|\o - \o'\|_t   <   \wh{\d}       $.
  If $\o \in \wh{F}_0  $, so is $\o'$ since $ \wh{\d} \le \frac12 \ul{\d}$.
  As $\wh{F}_0 \subset F \cap A_0$ and $\wh{\d} \le \d$, it follows from \eqref{paste_strategy} and \eqref{eq:r637} that
 \beas
 \underset{r \in [t,T]}{\sup} \, \underset{v \in  \hV}{\sup} \; \rho_{\overset{}{\hU}} \big( \wh{\a}(r,\o,v) ,  \wh{\a}(r,\o',v) \big)
= \underset{r \in [t,T]}{\sup} \, \underset{v \in  \hV}{\sup} \; \rho_{\overset{}{\hU}} \big( \a(r,\o,v) ,  \a(r,\o',v) \big)
 < \frac{\e}{2}.
 \eeas
 On the other hand, if $\o \neg \in \neg \wh{F}^n_i  $ for some $n=1,\cds \neg , N$ and $i  \neg = \neg  1, \cds \neg , \ell_n$, so is
$\o'$.  Since $\wh{F}^n_i  \neg \subset \neg  F  \neg \cap \neg  A^n_i    \cap    \Pi^{-1}_{t,t_n}  \big( \wt{F}^n_i \big)$  and
 \beas
  \underset{r \in [t_n,T]}{\sup}  \big| \big(\Pi_{t,t_n}  \neg  (\o')\big)(r)
   \neg - \neg \big(\Pi_{t,t_n}  \neg  (\o)\big)(r)  \big|
  \neg \le  \neg  | \o'(t_n)  \neg - \neg \o(t_n) |
    + \neg  \underset{r \in [t_n,T]}{\sup} |\o'(r) \neg - \neg \o(r)  |
      \neg  \le  \neg  2 \underset{r \in [t,T]}{\sup} \, |\o'(r)-\o(r)  | \neg < \neg 2 \wh{\d}  \le  \neg  \d^n_i  ,
 \eeas
 we can deduce from  \eqref{paste_strategy} and \eqref{eq:r637b}  that
    \beas
   &&  \hspace{-2cm} \underset{r \in [t,T]}{\sup} \, \underset{v \in \hV}{\sup} \; \rho_{\overset{}{\hU}} \big( \wh{\a}(r,\o,v) ,    \wh{\a}(r,\o',v) \big)
 \le \underset{r \in [t, t_n)}{\sup} \, \underset{v \in \hV}{\sup} \; \rho_{\overset{}{\hU}} \big( \a(r,\o,v), \a(r,\o',v)\big) \\
 && +  \underset{r \in [t_n, T]}{\sup} \, \underset{v \in \hV}{\sup} \; \rho_{\overset{}{\hU}} \big(  \a^n_i \big(r,\Pi_{t,t_n}(\o),v \big) ,
     \a^n_i \big(r,\Pi_{t,t_n}(\o'),v \big) \big) < \frac{\e}{2} + \frac{\e}{4 \ell_*}  <  \e   .
    \eeas
  Therefore, $\wh{\a}   $ satisfies \eqref{eq:r742}, to wit, $\wh{\a} \in \wh{\cA}^t $.    \qed

     \subsection{Proofs of Section \ref{sec:one_control}}

  \label{subsection:Proofs_S5}

  \ss \no {\bf Proof of Lemma \ref{lem_hSd_property}:}
 As a subspace of $\hR^{d \times d}$, $\hS_d$ is also a normed vector space (We can regard the restriction to $\hS_d$ of
 the Euclidean norm $|\cd| $ on $\hR^{d \times d}$   as the Euclidean norm on $\hS_d$.)
 Since each $\hR^{d \times d} -$valued symmetric     matrix is uniquely determined by its lower (or upper) triangle, we see
 that
 \beas
 \phi(\G)\dfnn (g_{11},g_{21},g_{22},g_{31}, g_{32},g_{33}, \cds , g_{d1},\cds , g_{dd}) ,  \q  \fa \G = (g_{ij})^d_{i,j=1} \in \hS_d
 \eeas
 defines a bijection between $ \hS_d $ and $ \hR^{\frac{d( 1 \neg + \neg d  )}{2}}$.
 Clearly, $ |\phi(\G) | \le  |\G| \le \sqrt{2} |\phi(\G) | $, $\fa \G \in \hS_d$, thus
$\phi$ is a  homeomorphism. Then the separability of $\hR^{\frac{d( 1 \neg + \neg d  )}{2}}$ leads to
that of $\hS_d$.

 Moreover,  as   $det(\cd)$  is a continuous function on $\hR^{d \times d}$, its restriction on $\hS_d$
 is also continuous w.r.t. the relative Euclidean topology on $\hS_d$.
 \qed

   \ss \no {\bf Proof of Lemma \ref{lem_pqv}:}
  Let $i,j \neg \in \neg \{1,\cds \neg ,d\}$. For any $n \neg \in \neg \hN$,
  we set $\t^{n,i}_0   \neg =  \neg  t$ and recursively   define
   $\bF^t-$stopping times
 \beas
 \t^{n,i}_\ell \dfnn \inf \Big\{ s \in \big[\t^{n,i}_{\ell-1} ,T \big] :  \;
   \Big| \, B^i_s - B^i_{\t^{n,i}_{\ell-1}}   \Big| > 2^{-n} \Big\} \land T,
 \q  \fa  \,  \ell \in \hN .
 \eeas
   Clearly,  $ \dis \sT^{n,i,j}_s \dfnn \linf{m \to \infty} \sum^m_{\ell = 1  } B^{t,i}_{\t^{n,i}_{\ell-1} \land s }
  \Big( B^{t,j}_{\t^{n,i}_{\ell} \land s } -  B^{t,j}_{\t^{n,i}_{\ell-1} \land s } \Big)$, $s \in [t,T]$
  is an $ \hR \cup \{ - \infty\}-$valued, $\bF^t-$progressively measurable  process, so is
  $ \sT^{i,j}  \dfnn  \linf{n \to \infty} \sT^{n,i,j} $.

 For any   $P \in \cQ^t$, as Lemma \ref{lem_F_version} (1) shows that
 $B^t  $ is also a continuous semi-martingale with respect to $(\bF^P, P)$,
     we know from Theorem 2 of  \cite{Karandikar_1995} that
   $
  \lmt{n \to 0} \, \underset{s \in [t,T]}{\sup} \,
  \Big| \sT^{n,i,j}_s - \int^P_{[t,s]} B^{t,i}_r d B^{t,j}_r \Big| = 0$,   \pas ~
 Thus it holds      \pas ~  that
 \bea \label{eq:xa011}
 \sT^{i,j}_s   = \int^P_{[t,s]}  B^{t,i}_r d B^{t,j}_r , \q \fa   s \in [t,T]  .
 \eea
 This gives rise to a pathwise definition of
 the $(i,j)-$th cross variance of $B^t$ as well as its density:
  \beas
 \q \big\lan B^{t,i}, B^{t,j} \big\ran_s   \dfnn B^{t,i} B^{t,j} -  \sT^{i,j}_s  -  \sT^{j,i}_s   ~\;  \hb{and} ~\;
\hat{a}^{t,i,j}_s \dfnn    \lsup{m \to \infty}  \,  m  \Big(  \big\lan B^{t,i}, B^{t,j} \big\ran_s
 - \big\lan B^{t,i}, B^{t,j} \big\ran_{(s-1\neg /m)^{\neg {}^+} }  \Big)    , \q  s  \in [t,T] ,
  \eeas
  both of which  are   $\hR \cup \{\infty\}-$valued, $\bF^t-$progressively measurable  processes.

 \ss  For any $P \in \cQ^t$, we see from    \eqref{eq:xa011} that  \pas ~
   \beas
  \big\lan B^{t,i}, B^{t,j} \big\ran_s   =   B^{t,i} B^{t,j} -  \int^P_{[t,s]} B^{t,i}_r d B^{t,j}_r
  - \int^P_{[t,s]} B^{t,j}_r d B^{t,i}_r =  \big\lan B^{t,i}, B^{t,j} \big\ran^P_s  ,\q  \fa    s  \in [t,T]   .
  \eeas
  Then  \eqref{eq:a075} easily follows.   \qed

  \ss \no {\bf Proof of Lemma \ref{lem_hat_a_half}:}
 Let $t \in [0,T]$. We see from Lemma \ref{lem_pqv} that $|\hat{a}^t| $ is a $[0,\infty]-$valued,
 $\bF^t-$progressively measurable process. It follows that
 $ \dis \b1_{\{|\hat{a}^t| \in (0,\infty) \}} \frac{ 1 }{|\hat{a}^t|} $
 is a $[0,\infty)-$valued, $\bF^t-$progressively measurable process and  thus that
 $\fn^t \dis  \dfnn  \b1_{\{|\hat{a}^t| \in (0,\infty) \}} \frac{ \hat{a}^t }{|\hat{a}^t|} $ is  an $\hS_d-$valued, $\bF^t-$progressively measurable process.
 Since the determinant $det(\cd)$ is continuous on
$\hS_d$ by Lemma \ref{lem_hSd_property},
 $    \hat{\fn}^t \dfnn \b1_{\{det(\fn^t)>0\}} \fn^t
+  \b1_{\{ det(\fn^t) \le 0 \}} I_{d \times d}     $
 defines  an $\hS^{> 0}_d-$valued, $\bF^t-$progressively measurable process.

 For any $j \in \hN$, let $ \dis c_j \dfnn - \frac{1 \times 3 \times \cds \times (2j-3)  }{2^j \; j!}$,
  which is the $j-$th coefficient of the power series of $\sqrt{1-x}$, $x \in [-1,1]$.
 When a   $\G \in \hS^{> 0}_d $ has   $|\G| \le 1$, we know (see e.g. Theorem VI.9 of \cite{Reed_Simon_1972}) that
 $\vs \dfnn I_{d \times d} + \sum_{j \in \hN} c_j (I_{d \times d} - \G)^j$ is
 the unique   element in $ \hS^{> 0}_d $  such that $\vs^2 =\vs \cd \vs = \G $.
 Consequently,   $ q^t   \dfnn I_{d \times d} + \sum_{j \in \hN} c_j (I_{d \times d} - \hat{\fn}^t )^j $
 is the  unique $\hS^{> 0}_d-$valued, $\bF^t-$progressively measurable process such that
 \beas
 (q^t)^2  
 = \hat{\fn}^t = \b1_{\{|\hat{a}^t| \in (0,\infty) \}}
 \b1_{\{det(\hat{a}^t)>0\}} \frac{ \hat{a}^t }{|\hat{a}^t|}
 + \big( \b1_{\{|\hat{a}^t| \in (0,\infty) \}}
 \b1_{\{det(\hat{a}^t) \le 0\} } + \b1_{\{|\hat{a}^t| =0  \,  \hb{\small or} \,  \infty \}} \big) \, I_{d \times d}  .
 \eeas
 It follows that
 $  \hat{q}^t \dfnn  q^t  \big( \b1_{\{|\hat{a}^t| \in (0,\infty) \}}  \sqrt{|\hat{a}^t|}
 + \b1_{\{|\hat{a}^t| =0  \,  \hb{\small or} \,  \infty \}} \big)   $ is the
   unique $\hS^{> 0}_d-$valued, $\bF^t-$progressively measurable process satisfying
 \bea \label{eq:xd015}
 (\hat{q}^t)^2
  = \b1_{\{|\hat{a}^t| \in (0,\infty) \}} \b1_{\{det(\hat{a}^t)>0\}}   \hat{a}^t
 + \b1_{\{|\hat{a}^t| \in (0,\infty) \}} \b1_{\{det(\hat{a}^t) \le 0\} } \, |\hat{a}^t|  I_{d \times d}
 + \b1_{\{|\hat{a}^t| =0  \,  \hb{\small or} \,  \infty \}}   \, I_{d \times d}  .
  \eea
 Given $P \in \cQ^t_W$, since $|\G| \in (0,\infty)$ for each $ \G \in \hS^{>0}_d$, we can deduce from the second part of \eqref{eq:xd011} and \eqref{eq:xd015}  that \pas,
 $ (\hat{q}^t)^2_s = \hat{q}^t_s \cd \hat{q}^t_s = \hat{a}^t_s $ for a.e. $s \in [t,T]$.    \qed

    \ss \no {\bf Proof of Lemma \ref{lem_X_alpha}:} Let $s \in [t,T]$.
 For any   $r \in [t,s]$ and $\cE \in   \sB(\hR^d)$,  similar to \eqref{eq:a247}, one can deduce from  the $\bF^t-$adaptness of  $ \wt{X}^{t,x,\mu} $    that
        \beas
    \big(\cX^{t,x,\mu}\big)^{-1} \Big( \big(  B^t_r \big)^{-1}(\cE)\Big)
 \neg  =  \neg  \Big\{\o \in \O^t  \neg :   \cX^{t,x,\mu}_{r} (\o)      \in \cE \Big\}
    \neg  =  \neg  \left\{
      \ba{ll}
    \hb{} \tneg    \cN^{\,t,x}_\mu \neg \cup \neg  \Big( (\cN^{\,t,x}_\mu)^c  \neg \cap \neg  \big\{\o \in \O^t  \neg : \wt{X}^{t,x,\mu}_{r}  (\o)    \neg   \in  \neg  \cE \big\}  \Big)
    \neg  \in  \neg   \ol{\cF}^t_{\neg s} ,     & \hb{if }0 \in \cE, \ms \\
  \hb{} \tneg  (\cN^{\,t,x}_\mu)^c  \neg \cap \neg  \big\{\o \in \O^t \neg :
 \wt{X}^{t,x,\mu}_{r}  (\o)  \neg  \in  \neg  \cE \big\}
   \neg   \in   \neg   \ol{\cF}^t_{\neg s}    ,   & \hb{if }0 \notin \cE .
       \ea
       \right.
   \eeas
   where $\cE_x = \{x+x': x' \in \cE\} \in \sB(\hR^d)$.
    Thus $\big(  B^t_r \big)^{-1}(\cE)  \in \L^t_s \dfnn \Big\{A \subset \O^t: \big( \cX^{t,x,\mu}\big)^{-1}(A)
   \in  \ol{\cF}^t_{\neg s}    \Big\}$.
    Clearly, $\L^t_s$  is     a $\si-$field of $\O^t$.
    So it follows that
       \bea \label{eq:d031}
    \cF^t_s  =\si \left( \big(B^t_r \big)^{-1} (\cE) ; \,  r \in [t,s] , ~ \cE   \in   \sB(\hR^d)  \right)  \subset \L^t_s   \,.
       \eea

       For any $ \cN \in \sN^{P^{t,x,\mu}}$, it is contained in some $A \in \cF^t_T$ with
  $P^{t,x,\mu}(A)=0$.    By \eqref{eq:a151} and
       \eqref{eq:a155},  $\big(\cX^{t,x,\mu}\big)^{-1}( A )  \in \cF^t_T$ and        $        P^t_0 \Big( \big(\cX^{t,x,\mu}\big)^{-1}( A ) \Big) = P^{t,x,\mu} (A) =0 $.
                  Then, as a subset of         $   \big(\cX^{t,x,\mu}\big)^{-1}( A ) $,
             \bea   \label{eq:d261}
              \big(\cX^{t,x,\mu}\big)^{-1} \big(\cN\big)   \in  \sN^{P^t_0} \subset  \ol{\cF}^t_{\neg s} .
              \eea
           Thus,    $\sN^{P^{t,x,\mu}} \dneg  \subset  \neg  \L^t_s $, which together with  \eqref{eq:d031} yields
           $    \cF^{P^{t,x,\mu}}_s   \dneg = \neg   \si \big(  \cF^t_s     \cup     \sN^{P^{t,x,\mu}}  \big)
     \neg  \subset  \neg  \L^t_s  $, i.e.
           $\big(\cX^{t,x,\mu}\big)^{-1} \big(\cF^{P^{t,x,\mu}}_s   \big) \neg \subset \neg \ol{\cF}^t_{\neg s}  $.

    \ss    For any $ A \in \cF^{P^{t,x,\mu}}_T $, we know (see e.g. Proposition 11.4 of \cite{Royden_real})
  that $A  = \wt{A} \,\cup\,  \cN $ for some  $\wt{A} \in \cF^t_T$
    and $ \cN   \neg \in \neg  \sN^{P^{t,x,\mu}} $. 
 Since $(\cX^{t,x,\mu})^{-1}  \big(\wt{A}\, \big) \neg \in \neg  \cF^t_T $ by \eqref{eq:a151}
    and since $(\cX^{t,x,\mu})^{-1}  \big( \cN \big)  \neg \in \neg  \sN^{P^t_0}$ by \eqref{eq:d261}, one can deduce that
         \beas
   ~   P^t_0 \circ (\cX^{t,x,\mu})^{-1}  (A) = P^t_0 \Big( (\cX^{t,x,\mu})^{-1}  \big(\wt{A}\, \big) \cup (\cX^{t,x,\mu})^{-1}  \big(\cN \big) \Big)
            = P^t_0 \Big( (\cX^{t,x,\mu})^{-1}  \big(\wt{A}\, \big) \Big)
                =P^{t,x,\mu} \big(\wt{A}\,\big)    = P^{t,x,\mu}(A)   .  \neg  \qq  \hb{\qed}
     \eeas

   \begin{lemm}  \label{lem_Y_transform}
Let  $(t,x)  \in [0,T] \times \hR^d$ and $\mu \in \cU^t$.
  If $\cY$ is an $\hM-$valued, $\bF^{P^{t,x,\mu}}-$adapted process, then $\cY (\cX^{t,x,\mu}) $ is
  $\ol{\bF}^t  -$adapted. Moreover,  if $\cY  \in \hC^p_{\bF^{P^{t,x,\mu}}}
  \big([t,T], \hE, P^{t,x,\mu}  \big) $ \Big(resp.\;$\hK^p_{\bF^{P^{t,x,\mu}}}
  \big([t,T],   P^{t,x,\mu}  \big) $\Big) for some $p  \in [1,\infty)$, then  $\cY (\cX^{t,x,\mu})
  \in \hC^{p }_{\ol{\bF}^t} \big([t,T], \hE  \big)$ \Big(resp.\;$\hK^p_{\ol{\bF}^t}
  \big([t,T]   \big) $\Big).

\end{lemm}

  \ss \no {\bf Proof:}  The first conclusion directly follows from Lemma \ref{lem_X_alpha}.   If
  $\cY  \in \hC^p_{\bF^{P^{t,x,\mu}}}  \big([t,T], \hE, P^{t,x,\mu}  \big) $ \Big(resp.\;$\hK^p_{\bF^{P^{t,x,\mu}}}
  \big([t,T],  \\ P^{t,x,\mu}  \big) $\Big) for some $p  \in [1,\infty)$, let
  $A \dfnn \{\o \in \O^t:    \hb{the path $s  \neg \to \neg  \cY_s     (\o)$   is continuous (resp.\;increasing)}\}$.
     Then we see  from  \eqref{eq:d265}    that
   \bea   \label{eq:xd511}
    1 = P^{t,x,\mu} (A) = P^t_0  \circ  (\cX^{t,x,\mu})^{-1} ( A )  = P^t_0
   \big( \big\{ \o \in \O^t: \cX^{t,x,\mu} (\o) \in A  \big\} \big) .
   \eea
    Namely, $\cY \big( \cX^{t,x,\mu} \big)$ has $P^t_0-$a.s. continuous (resp.\;increasing) paths.
      Applying \eqref{eq:d265} again yields that
    \beas
      E_t \bigg[ \underset{s \in [t,T]}{\sup} \big| \cY_s \big( \cX^{t,x,\mu} \big)  \big|^p \bigg]
     = E_{P^{t,x,\mu}} \bigg[ \underset{s \in [t,T]}{\sup} \big| \cY_s   \big|^p \bigg] < \infty  ~\;
  \Big(\hb{resp.} ~ E_t \Big[   \big| \cY_T \big( \cX^{t,x,\mu} \big)  \big|^p \Big]
     = E_{P^{t,x,\mu}} \Big[  \big| \cY_T   \big|^p \Big] < \infty   \Big).
    \eeas
Thus, $\cY (\cX^{t,x,\mu})  \in \hC^{p }_{\ol{\bF}^t} \big([t,T], \hE  \big)$ \Big(resp.\;$\hK^p_{\ol{\bF}^t}
  \big([t,T]   \big) $\Big).  \qed

  \begin{lemm}  \label{lem_Z_transform}
  Let  $(t,x)  \in [0,T] \times \hR^d$ and $\mu \in \cU^t$.
  If $\cZ$ is an $\hM-$valued, $\bF^{P^{t,x,\mu}}-$progressively measurable process, then $\cZ (\cX^{t,x,\mu}) $ is
  $\ol{\bF}^t  -$progressively measurable. Consequently, for any $p \in [1,\infty)$ if $\cZ  \in \hH^{p, loc}_{\bF^{P^{t,x,\mu}}}
  \big([t,T], \hE, P^{t,x,\mu}  \big) $ \big(resp.\;$\hH^{p, \wh{p}}_{\bF^{P^{t,x,\mu}}}
   ([t,T], \hE, P^{t,x,\mu}   )$ for some $\wh{p} \in [1,\infty)$\big), then  $\cZ (\cX^{t,x,\mu})
  \in \hH^{p, loc}_{\ol{\bF}^t} \big([t,T], \hE  \big)$ \big(resp.\;$\hH^{p, \wh{p}}_{\ol{\bF}^t}  ([t,T], \hE   )$\big).

  \end{lemm}

  \ss \no {\bf Proof:}
  Let $\cZ$ be an $\hM-$valued, $\bF^{P^{t,x,\mu}}-$progressively measurable process.
  Given $s \in [t,T]$, we define $\Pi^{x,\mu}_{t,s} (r,\o) \dfnn \big(r, \cX^{t,x,\mu}(\o) \big)$,
  $\fa  (r,\o) \in [t,s] \times \O^t$. For any $\cE \in \sB([t,s])$ and $A \in \cF^{P^{t,x,\mu}}_s$,
  Lemma \ref{lem_X_alpha} implies that
    \beas
   \big( \Pi^{x,\mu}_{t,s} \big)^{-1} \big( \cE \times A \big)
  = \big\{ (r,\o) \in [t,s] \times \O^t: \big(r, \cX^{t,x,\mu}(\o) \big) \in \cE \times A  \big\}
  =  \cE \times  \big(\cX^{t,x,\mu} \big)^{-1} (A) \in \sB([t,s]) \otimes \ol{\cF}^t_s  .
  \eeas
  Hence, the rectangular measurable set
  $  \cE \times A \in \L_s \dfnn \big\{ \cD \subset [t,s] \times \O^t : \big( \Pi^{x,\mu}_{t,s} \big)^{-1}  ( \cD  )
      \in \sB([t,s]) \otimes \ol{\cF}^t_s  \big\}   $, which is clearly a $\si-$field of $[t,s] \times \O^t$.
  It follows that $ \sB([t,s]) \otimes \cF^{P^{t,x,\mu}}_s \subset \L_s$.
  For any $\cM \in \sB(\hM)$, since $ \cZ^{-1}(\cM) \dfnn \big\{ (r,\o) \in [t,s] \times \O^t: \cZ_r (\o) \in \cM  \big\}
 \in  \sB([t,s]) \otimes \cF^{P^{t,x,\mu}}_s \subset \L_s$,
   \beas
 &&        \big\{ (r,\o) \in [t,s] \times \O^t: \cZ_r (\cX^{t,x,\mu}(\o)) \in \cM  \big\}
=  \big\{ (r,\o) \in [t,s] \times \O^t:   (r,\cX^{t,x,\mu}(\o)) \in \cZ^{-1} (\cM)  \big\}   \\
 && \hspace{2cm} =  \big( \Pi^{x,\mu}_{t,s} \big)^{-1}  ( \cZ^{-1}(\cM) )
 \in \sB([t,s]) \otimes \ol{\cF}^t_s .
  \eeas
   Hence, $\cZ (\cX^{t,x,\mu}) $ is  $\ol{\bF}^t  -$progressively measurable.

\ss If  $\cZ  \in \hH^{p, loc}_{\bF^{P^{t,x,\mu}}} \big([t,T], \hE, P^{t,x,\mu}  \big) $ for some $p  \in [1,\infty)$,
let  $A \dfnn  \big\{\o \in \O^t:      \int_t^T   |\cZ_s (\o)|^p ds < \infty       \big\}$.
    Similar to  \eqref{eq:xd511},  we see from \eqref{eq:d265}    that
   $  1 = P^{t,x,\mu} (A)    = P^t_0
   \big( \big\{ \o \in \O^t: \cX^{t,x,\mu} (\o) \in A  \big\} \big) $, i.e.,
 $\cZ (\cX^{t,x,\mu})$ has $P^t_0-$a.s. $p-$integrable paths. Thus $ \cZ (\cX^{t,x,\mu}) \in \hH^{p, loc}_{\ol{\bF}^t} \big([t,T], \hE  \big)$.

 \ss  Moreover,    if $\cZ  \in \hH^{p, \wh{p}}_{\bF^{P^{t,x,\mu}}}
  \big([t,T], \hE, P^{t,x,\mu}  \big) $ for some $p,\wh{p} \in [1,\infty)$, we can deduce  from   \eqref{eq:d265}   that
    \beas
      E_t \Bigg[ \bigg( \int_t^T    \big| \cZ_s \big( \cX^{t,x,\mu} \big)  \big|^p ds \bigg)^{\wh{p}/p} \Bigg]
 = E_{P^{t,x,\mu}} \Bigg[ \bigg( \int_t^T \big| \cZ_s \big|^p ds \bigg)^{\wh{p}/p} \Bigg] < \infty   .
    \eeas
 Therefore, $\cZ (\cX^{t,x,\mu})  \in \hH^{p, \wh{p}}_{\ol{\bF}^t} \big([t,T], \hE  \big)$.   \qed

\begin{lemm}  \label{lemm_mart_part_fI}
 Let $(t,x) \in [0,T] \times \hR^d$.   For any  $P \in \cQ^{t,x}_S $, the $\bF^t-$adapted continuous process
  $M^t_s  \dfnn B^t_s - \int_t^s b(r, B^{t,x}_r) dr$, $s \in [t,T]$
 is a continuous martingale with respect to  $(\bF^t, P  )$. Consequently, 
   \bea  \label{eq:mart_part_fI}
 W^P_s =     \int^{P}_{[t,s]}      \big(\hat{q}^t_r\big)^{-1  }    d  M^t_r  , \q    s \in [t,T ]   .
   \eea
\end{lemm}

      \ss \no {\bf Proof:} We fix   $ \mu \neg \in \neg  \cU^t $ and
 let $\wt{\cN}_\mu \dfnn \{ X^{t,x,\mu}_r \ne \wt{X}^{t,x,\mu}_r, \exists \; r \in [t,T] \} \in \sN^{P^t_0}$.
    Given  $t  \neg \le \neg  s  \neg < \neg  r  \neg \le \neg  T$,  since $     \int_t^s  \mu_r  \, dB^t_r  $, $s \in [t,T]$
    is a martingale with respect to $( \ol{\bF}^t, P^t_0)$,
         for any finite subset $\{t_1 < \cds < t_m \}$    of $\hQ \cap [t,s]$
         and any $\{ (x_i, \l_i) \}^m_{j=1} \subset \hQ^d \times \hQ_+ $,
         we can deduce from    \eqref{FSDE3} and the $\ol{\bF}^t-$adaptedness of $X^{t,x,\mu}$ that
        \beas
    && \hspace{-0.7cm}     \int_{\o' \in \underset{i=1}{\overset{m}{\cap}}  (B^t_{t_i}  )^{-1}   ( O_{\l_i}(x_i) ) }
 \big( M^t_r  (\o') \neg  -  \neg   M^t_s (\o') \big) \, dP^{t,x,\mu} (\o')  \\
         && \hspace{-0.3cm}  = \neg \int_{ \o \in (\cX^{t,x,\mu} )^{-1} \neg \big( \underset{i=1}{\overset{m}{\cap}}  (B^t_{t_i}  )^{-1}   ( O_{\l_i}(x_i) )  \big)}
         \big( B^t_r (\cX^{t,x,\mu}(\o)) \neg   -  \neg   B^t_s(\cX^{t,x,\mu}(\o))
 \neg - \neg   \hb{$\int_s^r \neg b \big(r' \neg ,(x \neg + \neg B^t_{r'}(\cX^{t,x,\mu}(\o)) \big) dr'$}  \big)  dP^t_0 (\o)   \\
 && \hspace{-0.3cm}  = \neg  \int_{ \o \in \big(\cN^{\,t,x}_\mu \cup \wt{\cN}_\mu \big)^c \cap \, (X^{t,x,\mu} )^{-1} \neg \big( \underset{i=1}{\overset{m}{\cap}}  (B^t_{t_i}  )^{-1}   ( O_{\l_i}(x_i) )  \big)}
         \big(  X^{t,x,\mu}_r(\o) \neg   -  \neg    X^{t,x,\mu}_s (\o)
  \neg - \neg   \hb{$\int_s^r \neg b(r' \neg ,  X^{t,x,\mu}_{r'}(\o)  dr'$}
     \big)  dP^t_0 (\o)   \\
   && \hspace{-0.3cm}  = \neg  \int_{ \o \in    \underset{i=1}{\overset{m}{\cap}} (X^{t,x,\mu}_{t_i})^{-1}  (  O_{\l_i}(x_i) ) }
 \big(\hb{$\int_s^r \mu_{r'} (\o) \, d B^t_{r'}(\o) $} \big) \,  dP^t_0 (\o) = 0.
    \eeas
This shows
   $  \cC^t_s \subset \L^t_s  \dfnn \big\{A  \neg \in  \neg  \cF^t_s:   \int_A    ( M^t_r    \neg   -   \neg  M^t_s  ) \, dP^{t,x,\mu} =0 \big\}$. As $ \L^t_s $   is clearly a Dynkin system,
  we see from Lemma \ref{lem_countable_generate1} and Dynkin system theorem  that $ \cF^t_s =  \si  \big(\cC^t_s \big) \subset \L^t_s $,
      which implies that $  E_{P^{t,x,\mu}} [ M^t_r | \cF^t_s ] = M^t_s $, $P^{t,x,\mu}-$a.s.  Hence, $M^t$ is a continuous martingale with respect to  $(\bF^t, P^{t,x,\mu} )$. By    Lemma \ref{lem_F_version} (1),
 $M^t$ is also a continuous martingale with respect to  $\big(\bF^{P^{t,x,\mu}}, P^{t,x,\mu} \big)$. As
 \beas
  \fI^{P^{t,x,\mu}}_s =     \int^{P^{t,x,\mu}}_{[t,s]}      \big(\hat{q}^t_r\big)^{-1  }    d  B^t_r
 =  \int_t^s  \big(\hat{q}^t_r\big)^{-1  }  b(r, x + B^t_r) dr
 + \int^{P^{t,x,\mu}}_{[t,s]}  \big(\hat{q}^t_r\big)^{-1  } d  M^t_r, \q s \in [t,T],
 \eeas
 we see that the process $\Big\{\int^{P^{t,x,\mu} }_{[t,s]}  \big(\hat{q}^t_r\big)^{-1  } d  M^t_r \Big\}_{s \in [t,T]}$
 is exactly the martingale part $W^{P^{t,x,\mu}}$ of $\fI^{P^{t,x,\mu}}$.
  \qed

     \ss \no {\bf Proof of Proposition \ref{prop_PS-in-PW}:}
  Fix $(t,x) \in [0,T] \times \hR^d$ and  $ \mu \neg \in \neg  \cU^t $.
  As  $M^t_s = B^t_s - \int_t^s b(r, x+B^t_r) dr$, $s \in [t,T]$
 is a continuous martingale with respect to  $(\bF^t, P  )$ by Lemma \ref{lemm_mart_part_fI}, we see that
   $B^t$ is a continuous semi-martingale with respect to  $(\bF^t, P^{t,x,\mu} )$.  Thus, $P^{t,x,\mu} \in \cQ^t $.

 \ss    Given $i,j \in \{1,\cds,d\}$,  it follows from \eqref{eq:d265}   that
 \beas
   0& \dneg  \tneg =& \dneg  \tneg P^{t,x,\mu} \neg - \neg \lmt{n \to \infty} \,
  \underset{ s \in [t,T] }{\sup} \bigg| \lan B^{t,i}, B^{t,j} \ran^{P^{t,x,\mu}}_s
   \neg - \neg  \sum^{\lfloor 2^n \neg  s\rfloor}_{k=1 }  
   \Big(B^{t,i}_{ \frac{k}{2^n}  }  \neg - \neg  B^{t,i}_{ \frac{k-1}{2^n}} \Big)
    \Big(B^{t,j}_{ \frac{k}{2^n}  } \neg  - \neg  B^{t,j}_{ \frac{k-1}{2^n}}   \Big)  \bigg|    \nonumber   \\
   & \dneg  \tneg = & \dneg  \tneg P^t_0 - \lmt{n \to \infty} \,
   \underset{ s \in [t,T] }{\sup} \bigg| \lan B^{t,i}, B^{t,j} \ran^{P^{t,x,\mu}}_s \big( \cX^{t,x,\mu}   \big)
    \neg - \neg   \sum^{\lfloor 2^n \neg  s\rfloor}_{k=1 }  
   \Big(\cX^{t,x,\mu,i}_{ \frac{k}{2^n}  }  \neg -  \neg  \cX^{t,x,\mu,i}_{ \frac{k-1}{2^n}} \Big)\Big(\cX^{t,x,\mu,j}_{ \frac{k}{2^n}  }
   \neg -   \neg   \cX^{t,x,\mu,j}_{ \frac{k-1}{2^n}}   \Big)  \bigg|     \nonumber  \\
 & \dneg  \tneg = & \dneg  \tneg P^t_0 - \lmt{n \to \infty} \,
   \underset{ s \in [t,T] }{\sup} \bigg| \lan B^{t,i}, B^{t,j} \ran^{P^{t,x,\mu}}_s \big( \cX^{t,x,\mu}   \big)
    \neg - \neg   \sum^{\lfloor 2^n \neg  s\rfloor}_{k=1 }  
   \Big(X^{t,x,\mu,i}_{ \frac{k}{2^n}  }  \neg -  \neg  X^{t,x,\mu,i}_{ \frac{k-1}{2^n}} \Big)\Big(X^{t,x,\mu,j}_{ \frac{k}{2^n}  }
   \neg -   \neg   X^{t,x,\mu,j}_{ \frac{k-1}{2^n}}   \Big)  \bigg|  \,  ,   
 \eeas
 which   implies that $P^t_0-$a.s.,
 \bea  \label{eq:xc021}
   \lan B^{t,i}, B^{t,j} \ran^{P^{t,x,\mu}}_s \big( \cX^{t,x,\mu}   \big)
    =   \lan X^{t,x,\mu,i}, X^{t,x,\mu,j} \ran^{P^t_0}_s
    =   \int_t^s \sum^d_{\ell=1} \mu^{i\ell}_r \mu^{j\ell}_r  \, dr  , \q   \fa s \in [t,T] .
 \eea
 Since   $  \mu_r \in \hU_0 = \hS^{>0}_d$, $dr \times dP^t_0-$a.s. (so is $\mu^2_r = \mu_r \cd \mu_r $),
   we can deduce from  \eqref{eq:xc021}, \eqref{eq:d265} and \eqref{eq:a075}    that
  \bea
  1   & \dneg  \tneg =& \dneg  \tneg  P^t_0 \bigg\{ \o \in \O^t:  s \to \lan B^t, B^t \ran^{P^{t,x,\mu}}_s \big( \cX^{t,x,\mu} (\o) \big) \hb{ is absolutely continuous and }  \nonumber \\
   &  \dneg  \tneg & \qq \qq
   \lmt{m \to \infty}   m  \Big(  \lan B^t  \ran^{P^{t,x,\mu}}_s - \lan B^t  \ran^{P^{t,x,\mu}}_{(s-1\neg /m)^{\neg {}^+} }  \Big)
      \big( \cX^{t,x,\mu} (\o) \big)  =  \mu^2_s   (\o)  \in   \hS^{>0}_d  \hb{ for a.e. }  s \in [t,T] \bigg\}
   \qq   \label{eq:xc031a}  \\
      & \dneg  \tneg =& \dneg  \tneg  P^{t,x,\mu} \bigg\{ \o' \in \O^t:  s \to \lan B^t, B^t \ran^{P^{t,x,\mu}}_s  (  \o' )   \hb{ is absolutely continuous and }  \nonumber \\
   &  \dneg  \tneg & \qq \qq
  \hat{a}^t_s  ( \o' ) = \lmt{m \to \infty} m \Big( \lan B^t  \ran^{P^{t,x,\mu}}_s - \lan B^t  \ran^{P^{t,x,\mu}}_{(s-1\neg /m)^{\neg {}^+} }
  \Big)  (  \o')  \in  \hS^{>0}_d  \hb{ for a.e. }  s \in [t,T] \bigg\} .    \label{eq:xc031b}
  \eea
   Therefore, $ P^{t,x,\mu} \in \cQ^t_W$.  \qed

   \begin{lemm} \label{lem_2equivalence}
   Let  $(t,x) \in [0,T] \times \hR^d$ and $\mu \in \cU^t$. It holds $P^t_0-$a.s. that
    \bea
   \hat{q}^t_s \big( \cX^{t,x,\mu} \big)      =  \mu_s     \hb{\;  for a.e.  }s \in [t,T] .  \label{eq:d251a}
   \eea
  And for any $\cZ \in \hH^{p,loc}_{\bF^{t,x,\mu}}   ([t,T], \hR^d, P^{t,x,\mu})$ with $p \in[1,\infty)$, it holds $P^t_0-$a.s. that
 \bea
      \bigg( \int^{P^{t,x,\mu}}_{[t,s]} \cZ_r  d W^{P^{t,x,\mu}}_r \bigg) \big(\cX^{t,x,\mu}\big)
    = \int_t^s \cZ_r (\cX^{t,x,\mu} )  d B^t_r     ,    \q   \fa    s    \in   [t,T]   ,  \label{eq:d251b}
        \eea
        where   $\cZ (\cX^{t,x,\mu}) \in \hH^{p, loc}_{\ol{\bF}^t} \big( [t,T], \hR^d  \big)$ by Lemma \ref{lem_Z_transform}.
   \end{lemm}

       \ss \no {\bf Proof:}
      We can deduce from Lemma \ref{lem_hat_a_half},  \eqref{eq:xc031b} and \eqref{eq:d265}   that
  \beas
 1& \tneg =& \tneg  P^{t,x,\mu} \left\{ \o' \in \O^t:  \big( \hat{q}^{t }_s \big)^2 (\o') = \hat{a}^{t }_s (\o')
 =   \lmt{m \to \infty} m \Big( \lan B^t  \ran^{P^{t,x,\mu}}_s - \lan B^t  \ran^{P^{t,x,\mu}}_{(s-1\neg /m)^{\neg {}^+} }
  \Big)  (  \o')     \hb{ for a.e. }  s \in [t,T] \right\}\\
 & \tneg =& \tneg  P^t_0 \left\{ \o \in \O^t: \big( \hat{q}^{t }_s \big)^2 \big(\cX^{t,x,\mu}(\o)\big)
 =  \lmt{m \to \infty} m \Big( \lan B^t  \ran^{P^{t,x,\mu}}_s - \lan B^t  \ran^{P^{t,x,\mu}}_{(s-1\neg /m)^{\neg {}^+} }
  \Big) \big(\cX^{t,x,\mu}(\o)\big)
   \hb{ for a.e. }  s \in [t,T] \right\}   .
  \eeas
  This together with \eqref{eq:xc031a} yields that
 \bea  \label{eq:xd021}
 1=  P^t_0 \left\{ \o \in \O^t: \big( \hat{q}^{t }_s \big)^2 \big(\cX^{t,x,\mu}(\o)\big)
 =   \mu^2_s(\o ) \hb{ and } \mu_s(\o )  \in   \hS^{> 0}_d
   \hb{ for a.e. }  s \in [t,T] \right\} ,
 \eea
 where we used the fact that $\mu_s \in \hU_0$, $ ds \times d P^t_0 -$a.s.
 Since for each   $\G \in \hS^{> 0}_d $ there exists a unique   element in $ \hS^{> 0}_d $
 such that $\vs^2 =\vs \cd \vs = \G $ (see e.g. Theorem VI.9 of \cite{Reed_Simon_1972}),
  \eqref{eq:xd021} leads to    \eqref{eq:d251a}.

 Now, let $\cZ \in \hH^{p,loc}_{\bF^{t,x,\mu}}   ([t,T], \hR^d, P^{t,x,\mu})$ for some  $p \in[1,\infty)$.
 We have seen from the proof  of Proposition \ref{prop_PS-in-PW} that
 the   process    $M^t_s   =   B^t_s -    \int_t^s b(r, x    +   B^t_r) dr$, $s    \in    [t,T]$
 is a continuous martingale with respect to  $(\bF^t, P^{t,x,\mu} )$. It is thus a   continuous martingale with respect to
 $ \big(\bF^{P^{t,x,\mu}}, P^{t,x,\mu} \big)$ by Lemma \ref{lem_F_version} (1).
 Then we know   (see e.g. Problem 3.2.27 of \cite{Kara_Shr_BMSC}) that
 there exists a sequence of $ \hR^d -$valued, $\bF^{P^{t,x,\mu}}-$simple processes
   $  \Big\{\F^n_s \neg = \neg  \sum^{\ell_n}_{i=1} \xi^n_i \b1_{ \big\{s \in (t^n_i, t^n_{i+1}] \big\} } ,
  \,  s  \neg \in  \neg  [t,T] \Big\}_{n \in \hN}$ \big(where $t  \neg = \neg  t^n_1  \neg <  \neg  \cds
 \neg <  \neg  t^n_{\ell_n+1} \neg = \neg T$
   and $\xi^n_i   \neg \in  \neg  \cF^{P^{t,x,\mu}}_{ t^n_i}$  for $i \neg = \neg 1,  \neg \cds  \neg  , \ell_n$\big) such that
    \bea
  && {P^{t,x,\mu}} \neg - \neg  \lmt{n \to \infty} \int_t^{T} \Big( \big(\F^n_s\big)^T  - \big( \hat{q}^t_s \big)^{-1} \cZ^T_s \Big)
 \,  d \lan M^t  \ran^{P^{t,x,\mu}}_s \, \Big(  \F^n_s - \cZ_s  \big( \hat{q}^t_s \big)^{-1  } \Big)  =0 , \label{eq:d211a} \\
         \hb{and}  &&
             {P^{t,x,\mu}}  \neg - \neg   \lmt{n \to \infty} \, \underset{s \in [t,T]}{\sup} \left|  \,  \sum^{\ell_n}_{i=1} \xi^n_i  \big(  M^t_{s \land t^n_{i+1}} -  M^t_{s \land t^n_i} \big)
            - \int^{P^{t,x,\mu}}_{[t,s]}   \cZ_r  \big(\hat{q}^t_r\big)^{-1  }    d M^t_r \right| =0   . \qq \qq  \label{eq:d211b}
    \eea
 Define a martingale with respect to $ \big( \ol{\bF}^t, P^t_0 \big)$:  $ \U^t_s \dfnn \int_t^s \mu_r dB^t_r$, $s \in [t,T]$.
  We can deduce from \eqref{eq:d211a}, \eqref{eq:xc073},  \eqref{eq:d265} and   \eqref{eq:d251a}  that
    \bea
 0 &\tneg =&\tneg   P^{t,x,\mu} \neg -\neg \lmt{n \to \infty} \int_t^{T}     \Big( \big( \F^n_s \big)^T -   \big( \hat{q}^t_s
  \big)^{-1  } \cZ^T_s \Big)   \,   d \lan B^t  \ran^{P^{t,x,\mu}}_s \, \Big( \F^n_s - \cZ_s  \big( \hat{q}^t_s \big)^{-1  } \Big)
      \nonumber  \\
     &\tneg =&\tneg    P^t_0 \neg -\neg \lmt{n \to \infty} \int_t^{T}      \Big( \big( \F^n_s (\cX^{t,x,\mu}) \big)^T
    - \big( ( \hat{q}^t_s )^{-1  }   \cZ^T_s \big) (\cX^{t,x,\mu})     \Big)    \, \hat{a}^t_{s} (\cX^{t,x,\mu})\,
    \Big( \F^n_s (\cX^{t,x,\mu}) -
 \big( \cZ_s ( \hat{q}^t_{s} )^{-1 } \big) (\cX^{t,x,\mu})    \Big)    d s   \nonumber  \\
 &\tneg =&\tneg    P^t_0\neg -\neg \lmt{n \to \infty} \int_t^{T}  \Big( \big( \F^n_s (\cX^{t,x,\mu}) \big)^T- \mu^{-1}_{s}
 \cZ^T_s  (\cX^{t,x,\mu}) \Big)
       \,  d  \lan \U^t  \ran^{P^t_0}_{s}    \,  \Big( \F^n_s (\cX^{t,x,\mu}) -\cZ_s (\cX^{t,x,\mu})  \mu^{-1}_{s}    \Big)  \; .  \label{eq:d281}
            \eea
   Also,     \eqref{eq:d211b}, \eqref{eq:mart_part_fI}  and  \eqref{eq:d265} yield  that
      \bea
  0 & \tneg  \dneg =& \tneg  \dneg   P^{t,x,\mu}\neg -\neg \lmt{n \to \infty} \, \underset{s \in [t,T]}{\sup}  \,
 \bigg| \, \sum^{\ell_n}_{i=1} \xi^n_i  \big(  M^t_{s \land t^n_{i+1}} \neg  - \neg   M^t_{s \land t^n_i} \big)
 \neg  -  \neg \int^{P^{t,x,\mu}}_{[t,s]}   \cZ_r d W^{P^{t,x,\mu}}_r \bigg|  \nonumber \\
   & \tneg  \dneg =&  \tneg  \dneg  P^t_0\neg -\neg \lmt{n \to \infty} \,  \underset{s \in [t,T]}{\sup}  \,
  \bigg|  \, \sum^{\ell_n}_{i=1} \xi^n_i (\cX^{t,x,\mu}) \bigg(  \cX^{t,x,\mu}_{s \land t^n_{i+1}}
           \neg   -  \neg  \cX^{t,x,\mu}_{s \land t^n_i}
 \neg - \neg \int_{s \land t^n_i}^{s \land t^n_{i+1}} b \big(r, x \neg + \neg \cX^{t,x,\mu}_r \big) dr \bigg)
  \neg  -  \neg  \bigg(\int^{P^{t,x,\mu}}_{[t,s]} \neg  \cZ_r d W^{P^{t,x,\mu}}_r \bigg) (\cX^{t,x,\mu}) \bigg|   \nonumber \\
            & \tneg  \dneg =&  \tneg  \dneg   P^t_0\neg -\neg \lmt{n \to \infty} \,  \underset{s \in [t,T]}{\sup}  \,
  \bigg|  \, \sum^{\ell_n}_{i=1} \xi^n_i (\cX^{t,x,\mu}) \bigg(  X^{t,x,\mu}_{s \land t^n_{i+1}}
            \neg  -  \neg  X^{t,x,\mu}_{s \land t^n_i}
  -  \int_{s \land t^n_i}^{s \land t^n_{i+1}} b \big(r, X^{t,x,\mu}_r \big) dr  \bigg)  \neg   - \bigg(\int^{P^{t,x,\mu}}_{[t,s]}   \cZ_r d W^{P^{t,x,\mu}}_r \bigg) (\cX^{t,x,\mu}) \bigg| \nonumber \\
  & \tneg  \dneg =&  \tneg  \dneg   P^t_0\neg -\neg \lmt{n \to \infty} \,  \underset{s \in [t,T]}{\sup}  \,
  \bigg|  \, \sum^{\ell_n}_{i=1} \xi^n_i (\cX^{t,x,\mu}) \big(  \U^t_{s \land t^n_{i+1}}
            \neg  -  \neg  \U^t_{s \land t^n_i}
     \big)  \neg   - \bigg(\int^{P^{t,x,\mu}}_{[t,s]}   \cZ_r d W^{P^{t,x,\mu}}_r \bigg) (\cX^{t,x,\mu}) \bigg| .   \label{eq:d287}
    \eea

     For $i = 1,\cds, \ell_n$, since  $  \xi^n_i  \in  \cF^{P^{t,x,\mu}}_{ t^n_i}$, we see  from Lemma  \ref{lem_X_alpha} that
     $ \xi^n_i (\cX^{t,x,\mu}) \in \ol{\cF}^t_{ t^n_i}$. Thus, $  \F^n (\cX^{t,x,\mu})$ is also an  $\ol{\bF}^t -$simple process.
            Proposition 3.2.26 of \cite{Kara_Shr_BMSC} and \eqref{eq:d281} then imply  that
      \beas
         P^t_0\neg -\neg \lmt{n \to \infty} \underset{s \in [t,T]}{\sup} \left|  \,
               \sum^{\ell_n}_{i=1} \xi^n_i (\cX^{t,x,\mu}) \big(  X^{t,x,\mu}_{ s \land t^n_{i+1}}
            -  X^{t,x,\mu}_{s \land t^n_i } \big)  - \int_t^s \cZ_r (\cX^{t,x,\mu}) \mu^{-1 }_r d\,  \U^t_r  \right| =0   ,
      \eeas
     which     together with \eqref{eq:d287}   shows that $P^t_0-$a.s.
       \beas
    \hspace{1.6cm}     \bigg(\int^{P^{t,x,\mu}}_{[t,s]}   \cZ_r d W^{P^{t,x,\mu}}_r \bigg) (\cX^{t,x,\mu})  
  = \int_t^s \cZ_r (\cX^{t,x,\mu}) \mu^{-1 }_r d\, \U^t_r =  \int_t^s \cZ_r (\cX^{t,x,\mu}) d B^t_r , \q \fa  s \in [t,T] .
    \hspace{1.5cm}  \hb{\qed}
       \eeas

 \ss \no {\bf Proof of Proposition \ref{prop_filtration_coincide}:}
  Let $(t,x)  \in  [0,T] \times \hR^d$ and $\mu    \in    \cU^t$.
  For any $s    \in    [t,T]$,  since $W^{P^{t,x,\mu}}_s  \neg \in  \neg  \cF^{P^{t,x,\mu}}_s$, one can easily deduce that
   \beas
    \cG^{P^{t,x,\mu}}_s  = \si \Big( \si \big(W^{P^{t,x,\mu}}_r, r \in [t,s]  \big) \cup \sN^{P^{t,x,\mu}} \Big)  \subset \cF^{P^{t,x,\mu}}_s  .
    \eeas
         As to  the inverse inclusion,     \eqref{eq:d251b} implies that
         $    B^t_\cd     =   W^{P^{t,x,\mu}}_\cd (\cX^{t,x,\mu})  $
         holds except on a $P^t_0-$null set $  \cN_0    $.
   Given   $r \neg \in \neg  [t,s]$ and $\cE  \neg \in  \neg   \sB(\hR^d)$,  one has
 \beas
  \big(B^t_r\big)^{-1} (\cE) = \{ \o \in \O^t: B^t_r (\o) \in \cE \}
 =  \big\{ \o \in \O^t: W^{P^{t,x,\mu}}_r \big( \cX^{t,x,\mu}  (\o) \big)  \in \cE \big\}  \, \D \,   \cN_\cE
 \eeas
with  $\cN_\cE   \dfnn    \big( \{ \o \in \O^t: B^t_r (\o) \in \cE \}  \, \D \,   \{ \o \in \O^t: W^{P^{t,x,\mu}}_r \big( \cX^{t,x,\mu}  (\o) \big)  \in \cE \} \big)    \subset \cN_0   $. To wit,
      \beas
    \big(B^t_r\big)^{-1} (\cE)  \in   \L^t_s  \dfnn \Big\{ A \subset \O^t:
     A =  \Big( ( \cX^{t,x,\mu})^{-1}\big(\wt{A}\big) \Big)  \,\D\, \cN
     \,    \hb{ for some }   \wt{A} \in  \cG^{P^{t,x,\mu}}_s  \hb{ and }    \cN \in     \sN^{P^t_0} \Big\} .
       \eeas
   Clearly, $\Big\{(\cX^{t,x,\mu})^{-1}\big(\wt{A}\big): \wt{A} \in  \cG^{P^{t,x,\mu}}_s  \Big\} $ is a $\si-$field of $\O^t$.
      Similar to Problem 2.7.3 of  \cite{Kara_Shr_BMSC}, one can show  that $ \L^t_s $ forms a $\si-$field of $\O^t$.
   It follows that      $ \ol{\cF}^t_{\neg s}  \subset \L^t_s $.
            As $ \sN^{P^t_0} \subset \L^t_s $,  we further have   $ \ol{\cF}^t_{\neg s}  \subset \L^t_s $.

      \ms Now, for any   $A \in \cF^{P^{t,x,\mu}}_s $, we know from Lemma \ref{lem_X_alpha} that
  $  (\cX^{t,x,\mu})^{-1}  (   A  )  \in    \ol{\cF}^t_{\neg s}      \subset \L^t_s $.
     Hence, there exists $\wt{A} \in  \cG^{P^{t,x,\mu}}_s \subset \cF^{P^{t,x,\mu}}_s  $ and $ \cN \in     \sN^{P^t_0}$ such that
   $   (\cX^{t,x,\mu})^{-1}  ( A  )    =   (\cX^{t,x,\mu})^{-1}  \big(  \wt{A} \, \big) \,\D \, \cN $,
          which leads to  that
       \beas
   \cN=     (\cX^{t,x,\mu})^{-1} ( A  )  \,\D\,  (\cX^{t,x,\mu})^{-1} \big(  \wt{A} \,  \big)
  =  (\cX^{t,x,\mu})^{-1} \big(  A \,\D\, \wt{A}  \, \big).
                 \eeas
 Then    \eqref{eq:d265} shows that
 $     P^{t,x,\mu}  \big(  A \,\D\, \wt{A}  \, \big) 
 = P^t_0 \big(\cN\big) =0 $, namely, $A \,\D\, \wt{A}  \in \sN^{P^{t,x,\mu}}$.
 It follows that $A = \wt{A} \, \D \, \big( A \,\D\, \wt{A} \big) \in \cG^{P^{t,x,\mu}}_s$,
 thus $ \cF^{P^{t,x,\mu}}_s = \cG^{P^{t,x,\mu}}_s$.
  \if{0}
\ss If $s < T$,  since the $P-$augmented   filtration $\big\{ \cG^{P^{t,x,\mu}}_s \big\}_{s \in [t,T]} $    generated
 by the $P-$Brownian motion $W^{P^{t,x,\mu}}$
 satisfies the usual condition, we see that    $\cF^t_{s+}
 \subset \cF^{P^{t,x,\mu}}_{s+} =  \cF^{P^{t,x,\mu}}_s$. If $s=T$, we clearly have $ \cF^t_{T+}= \cF^t_T \subset \cF^{P^{t,x,\mu}}_T$.
 \fi
 \qed

    \ss \no {\bf Proof of Proposition \ref{prop_DRBSDE_transform}:}  Let
 us simply denote $\Big(\cY^{t,x, {P^{t,x,\mu}}}   \big(T,  h ( B^{t,x}_T  ) \big)  ,
 \cZ^{t,x, {P^{t,x,\mu}}}    \big(T,  h ( B^{t,x}_T  ) \big)  ,
 \ul{\cK}^{t,x, {P^{t,x,\mu}}}    \big(T,  h ( B^{t,x}_T  ) \big)  , \\
   \ol{\cK}^{t,x, {P^{t,x,\mu}}}    \big(T,  h ( B^{t,x}_T  ) \big)  \Big) $  by     $(\cY, \cZ, \ul{\cK},\ol{\cK})  $.
  Lemma \ref{lem_Y_transform} and Lemma \ref{lem_Z_transform} shows that
  \beas
 (\sY, \sZ, \ul{\sK},\ol{\sK}) \dfnn \big(\cY (\cX^{t,x,\mu}), \cZ (\cX^{t,x,\mu}), \ul{\cK}
  (\cX^{t,x,\mu}),\ol{\cK}(\cX^{t,x,\mu}) \big) \in \hG^q_{\ol{\bF}^t}  \big([t,T]  \big).
 \eeas
 And we can deduce from \eqref{eq:d265}, Lemma \ref{lem_2equivalence}  that
 $(\sY, \sZ, \ul{\sK},\ol{\sK})$ satisfies
 \beas
    \left\{\ba{l}
 \dis   \sY_s= h  \big(  \wt{X}^{t,x,\mu}_T \big) \neg + \neg  \int_s^T  \neg
 f  \big(r, \wt{X}^{t,x,\mu}_r,  \sY_r, \sZ_r, \mu_r \big)  \, dr  \neg+ \neg  \ul{\sK}_{\,T} \neg-\neg \ul{\sK}_{\,s}
  \neg- \neg \big( \ol{\sK}_T \neg - \neg \ol{\sK}_s \big)
   \neg-\neg  \int_s^T \neg \sZ_r  \, d B^t_r  , \q    s \in [t,T] ,   \vspace{1mm} \q \\
   \dis \ul{L}^{t,x,\mu}_{  s}   \le \sY_s
   \le  \ol{L}^{t,x,\mu}_{  s}  , \q    s \in [t,T] ,
   \q \hb{and} \q  \int_t^T \neg  \big(  \sY_s   -  \ul{L}^{t,x,\mu}_{  s}  \big)     d \ul{\sK}_{\,s}
  = \int_t^T \neg  \big(  \ol{L}^{t,x,\mu}_{  s} - \sY_s  \big)     d \ol{\sK}_s   = 0
     \ea \right.
    \eeas
 on  the probability space $ \big(\O^t, \ol{\cF}^t_T, P^t_0 \big)$.
 Since this doubly Reflected BSDE admits a unique solution in $\hG^q_{\ol{\bF}^t}  \big([t,T]  \big)$
 according to Section \ref{sec:zs_drgame}, we have
  \beas
 \hspace{0.6cm} (\sY, \sZ, \ul{\sK},    \ol{\sK})  =  \Big( Y^{t,x,\mu} \big(T, h  \big(  \wt{X}^{t,x,\mu}_T \big) \big),
 Z^{t,x,\mu} \big(T, h  \big(  \wt{X}^{t,x,\mu}_T \big) \big), \ul{K}^{t,x,\mu} \big(T, h  \big(  \wt{X}^{t,x,\mu}_T \big) \big)
, \ol{K}^{t,x,\mu} \big(T, h  \big(  \wt{X}^{t,x,\mu}_T \big) \big)\Big) .  \hspace{0.5cm}  \hb{\qed}
 \eeas

\medskip

\bibliographystyle{siam}
\bibliography{SDGVU}

\begin{thebibliography}{10}

\bibitem{Atar_Budhiraja_2010}
{\sc R.~Atar and A.~Budhiraja}, {\em A stochastic differential game for the
  inhomogeneous {$\infty$}-{L}aplace equation}, Ann. Probab., 38 (2010),
  pp.~498--531.

\bibitem{BH11}
{\sc E.~Bayraktar and Y.~Huang}, {\em On the multi-dimensional controller and
  stopper games},  (2011).
\newblock Available at \url{http://arxiv.org/abs/1009.0932}.

\bibitem{OS_CRM}
{\sc E.~Bayraktar, I.~Karatzas, and S.~Yao}, {\em Optimal stopping for dynamic
  convex risk measures}, to appear in the Illinois Journal of Mathematics
  (special issue on honor of Don. Burkholder),  (2012), pp.~1--45.
\newblock Available at \url{http://lanl.arxiv.org/abs/0909.4948}.

\bibitem{BS12}
{\sc E.~Bayraktar and M.~S\^{i}rbu}, {\em Stochastic perron's method and
  verification without smoothness using viscosity comparison: obstacle problems
  and {D}ynkin games},  (2011).
\newblock Available at \url{http://arxiv.org/abs/1112.4904}.

\bibitem{SDGVU}
{\sc E.~Bayraktar and S.~Yao}, {\em ({A} longer version of) {O}n zero-sum
  stochastic differential games},  (2011).
\newblock Available at \url{http://arxiv.org/abs/1112.5744}.

\bibitem{OSNE1}
{\sc E.~Bayraktar and S.~Yao}, {\em Optimal stopping for non-linear
  expectations---{P}art {I}}, Stochastic Process. Appl., 121 (2011),
  pp.~185--211.

\bibitem{OSNE2}
\leavevmode\vrule height 2pt depth -1.6pt width 23pt, {\em Optimal stopping for
  non-linear expectations---{P}art {II}}, Stochastic Process. Appl., 121
  (2011), pp.~212--264.

\bibitem{Browne_2000}
{\sc S.~Browne}, {\em Stochastic differential portfolio games}, J. Appl.
  Probab., 37 (2000), pp.~126--147.

\bibitem{Buckdahn_Li_1}
{\sc R.~Buckdahn and J.~Li}, {\em Stochastic differential games and viscosity
  solutions of {H}amilton-{J}acobi-{B}ellman-{I}saacs equations}, SIAM J.
  Control Optim., 47 (2008), pp.~444--475.

\bibitem{Buckdahn_Li_3}
\leavevmode\vrule height 2pt depth -1.6pt width 23pt, {\em Probabilistic
  interpretation for systems of {I}saacs equations with two reflecting
  barriers}, NoDEA Nonlinear Differential Equations Appl., 16 (2009),
  pp.~381--420.

\bibitem{Buckdahn_Li_2}
\leavevmode\vrule height 2pt depth -1.6pt width 23pt, {\em Stochastic
  differential games with reflection and related obstacle problems for {I}saacs
  equations}, Acta Math. Appl. Sin. Engl. Ser., 27 (2011), pp.~647--678.

\bibitem{CSTV_2007}
{\sc P.~Cheridito, H.~M. Soner, N.~Touzi, and N.~Victoir}, {\em Second-order
  backward stochastic differential equations and fully nonlinear parabolic
  {PDE}s}, Comm. Pure Appl. Math., 60 (2007), pp.~1081--1110.

\bibitem{Cvitanic_Karatzas_1996}
{\sc J.~Cvitani{\'c} and I.~Karatzas}, {\em Backward stochastic differential
  equations with reflection and {D}ynkin games}, Ann. Probab., 24 (1996),
  pp.~2024--2056.

\bibitem{Dugundji_1966}
{\sc J.~Dugundji}, {\em Topology}, Allyn and Bacon Inc., Boston, Mass., 1966.

\bibitem{EHW_2011}
{\sc B.~El~Asri, S.~Hamad\`ene, and H.~Wang}, {\em ${L}^p$-solutions for doubly
  reflected backward stochastic differential equations}, to appear in
  Stochastic Analysis and Applications,  (2011).

\bibitem{EH-03}
{\sc N.~El-Karoui and S.~Hamad{\`e}ne}, {\em B{SDE}s and risk-sensitive
  control, zero-sum and nonzero-sum game problems of stochastic functional
  differential equations}, Stochastic Process. Appl., 107 (2003), pp.~145--169.

\bibitem{Fleming_2011}
{\sc W.~H. Fleming and D.~Hern\'andez~Hern\'andez}, {\em On the value of
  stochastic differential games}, Communications on Stochastic Analysis, 5
  (2011), pp.~341--351.

\bibitem{Fleming_1989}
{\sc W.~H. Fleming and P.~E. Souganidis}, {\em On the existence of value
  functions of two-player, zero-sum stochastic differential games}, Indiana
  Univ. Math. J., 38 (1989), pp.~293--314.

\bibitem{Hamadene_Hassani_2005}
{\sc S.~Hamad\`ene and M.~Hassani}, {\em {BSDE}s with two reflecting barriers:
  the general result.}, Probab. Theory Relat. Fields, 132 (2005), pp.~237--264.

\bibitem{Hamadene_Lepeltier_1995a}
{\sc S.~Hamad{\`e}ne and J.~P. Lepeltier}, {\em Backward equations, stochastic
  control and zero-sum stochastic differential games}, Stochastics Stochastics
  Rep., 54 (1995), pp.~221--231.

\bibitem{Hamadene_Lepeltier_1995b}
\leavevmode\vrule height 2pt depth -1.6pt width 23pt, {\em Zero-sum stochastic
  differential games and backward equations}, Systems Control Lett., 24 (1995),
  pp.~259--263.

\bibitem{Hamadene_Lepeltier_2000}
{\sc S.~Hamad{\`e}ne and J.-P. Lepeltier}, {\em Reflected {BSDE}s and mixed
  game problem}, Stochastic Process. Appl., 85 (2000), pp.~177--188.

\bibitem{Hamadene_Lepeltier_Wu_1999}
{\sc S.~Hamad{\`e}ne, J.-P. Lepeltier, and Z.~Wu}, {\em Infinite horizon
  reflected backward stochastic differential equations and applications in
  mixed control and game problems}, Probab. Math. Statist., 19 (1999),
  pp.~211--234.

\bibitem{Hamadene_Popier_2008}
{\sc S.~Hamad\`ene and A.~Popier}, {\em ${L}^p$-solutions for reflected
  backward stochastic differential equations}, to appear in Stochastics and
  Dynamics.

\bibitem{Hamadene_Rotenstein_Zlinescu_2009}
{\sc S.~Hamad{\`e}ne, E.~Rotenstein, and A.~Z{\u{a}}linescu}, {\em A
  generalized mixed zero-sum stochastic differential game and double barrier
  reflected {BSDE}s with quadratic growth coefficient}, An. \c Stiin\c t. Univ.
  Al. I. Cuza Ia\c si. Mat. (N.S.), 55 (2009), pp.~419--444.

\bibitem{Karandikar_1995}
{\sc R.~L. Karandikar}, {\em On pathwise stochastic integration}, Stochastic
  Process. Appl., 57 (1995), pp.~11--18.

\bibitem{Kara_Shr_BMSC}
{\sc I.~Karatzas and S.~E. Shreve}, {\em Brownian motion and stochastic
  calculus}, vol.~113 of Graduate Texts in Mathematics, Springer-Verlag, New
  York, second~ed., 1991.

\bibitem{MR1872738}
{\sc I.~Karatzas and W.~D. Sudderth}, {\em The controller-and-stopper game for
  a linear diffusion}, Ann. Probab., 29 (2001), pp.~1111--1127.

\bibitem{MR2435857}
{\sc I.~Karatzas and I.-M. Zamfirescu}, {\em Martingale approach to stochastic
  differential games of control and stopping}, Ann. Probab., 36 (2008),
  pp.~1495--1527.

\bibitem{Krylov_CDP}
{\sc N.~V. Krylov}, {\em Controlled diffusion processes}, vol.~14 of Stochastic
  Modelling and Applied Probability, Springer-Verlag, Berlin, 2009.
\newblock Translated from the 1977 Russian original by A. B. Aries, Reprint of
  the 1980 edition.

\bibitem{LLPU_1994}
{\sc A.~J. Lazarus, D.~E. Loeb, J.~G. Propp, and D.~Ullman}, {\em Richman
  games}, in Games of no chance ({B}erkeley, {CA}, 1994), vol.~29 of Math. Sci.
  Res. Inst. Publ., Cambridge Univ. Press, Cambridge, 1996, pp.~439--449.

\bibitem{Matoussi_Possamai_Zhou_2012}
{\sc A.~Matoussi, D.~Possamai, and C.~Zhou}, {\em Second Order Reflected
  Backward Stochastic Differential Equations}, 2012.
\newblock Available at \url{http://arxiv.org/abs/1201.0746}.

\bibitem{Nutz_2011}
{\sc M.~Nutz}, {\em A quasi-sure approach to the control of non-{M}arkovian
  stochastic differential equations}, tech. rep., ETH Zurich, 2011.
\newblock Available at \url{http://arxiv.org/abs/1106.3273}.

\bibitem{Peng_G_2007b}
{\sc S.~Peng}, {\em {G}-{B}rownian motion and dynamic risk measure under
  volatility uncertainty}, 2007.
\newblock Available at \url{http://lanl.arxiv.org/abs/0711.2834}.

\bibitem{Peng_G_2007a}
\leavevmode\vrule height 2pt depth -1.6pt width 23pt, {\em {$G$}-expectation,
  {$G$}-{B}rownian motion and related stochastic calculus of {I}t\^o type}, in
  Stochastic analysis and applications, vol.~2 of Abel Symp., Springer, Berlin,
  2007, pp.~541--567.

\bibitem{Peng_XuXm_2010}
{\sc S.~Peng and X.~Xu}, {\em B{SDE}s with random default time and related
  zero-sum stochastic differential games}, C. R. Math. Acad. Sci. Paris, 348
  (2010), pp.~193--198.

\bibitem{PSSW_2009}
{\sc Y.~Peres, O.~Schramm, S.~Sheffield, and D.~B. Wilson}, {\em Tug-of-war and
  the infinity {L}aplacian}, J. Amer. Math. Soc., 22 (2009), pp.~167--210.

\bibitem{Reed_Simon_1972}
{\sc M.~Reed and B.~Simon}, {\em Methods of modern mathematical physics. {I}.
  {F}unctional analysis}, Academic Press, New York, 1972.

\bibitem{revuz_yor}
{\sc D.~Revuz and M.~Yor}, {\em Continuous martingales and {B}rownian motion},
  vol.~293 of Grundlehren der Mathematischen Wissenschaften [Fundamental
  Principles of Mathematical Sciences], Springer-Verlag, Berlin, third~ed.,
  1999.

\bibitem{Royden_real}
{\sc H.~L. Royden}, {\em Real analysis}, Macmillan Publishing Company, New
  York, third~ed., 1988.

\bibitem{STZ_2011b}
{\sc H.~Soner, N.~Touzi, and J.~Zhang}, {\em Dual formulation of second order
  target problems}, preprint,  (2011).
\newblock Available at \url{http://arxiv.org/abs/1003.6050}.

\bibitem{STZ_2011d}
\leavevmode\vrule height 2pt depth -1.6pt width 23pt, {\em Martingale
  representation theorem for the {$G$}-expectation}, Stochastic Process. Appl.,
  121 (2011), pp.~265--287.

\bibitem{STZ_2011a}
\leavevmode\vrule height 2pt depth -1.6pt width 23pt, {\em Quasi-sure
  stochastic analysis through aggregation}, to appear in Electronic Journal of
  Probability,  (2011).
\newblock Available at \url{http://arxiv.org/abs/1003.4431}.

\bibitem{STZ_2011c}
\leavevmode\vrule height 2pt depth -1.6pt width 23pt, {\em Wellposedness of
  second order {B}ackward {SDE}s}, to appear in Probability Theory and Related
  Fields,  (2011).
\newblock Available at \url{http://arxiv.org/abs/1003.6053}.

\bibitem{Stroock_Varadhan}
{\sc D.~W. Stroock and S.~R.~S. Varadhan}, {\em Multidimensional diffusion
  processes}, Classics in Mathematics, Springer-Verlag, Berlin, 2006.
\newblock Reprint of the 1997 edition.

\end{thebibliography}

\end{document}